%% file: munk-def-HAL.tex
\documentclass{amsbook}

\usepackage{amsmath}
\usepackage{amsthm}
\usepackage{amssymb}
\usepackage{mathtools}

\usepackage{graphicx}

     \usepackage[active]{srcltx}
\usepackage[frenchb,english]{babel}
\usepackage[utf8]{inputenc}
\usepackage{yfonts}
\usepackage{enumerate}
\usepackage[shortlabels]{enumitem}
\numberwithin{equation}{chapter}
\numberwithin{figure}{chapter}

\usepackage{xcolor}

\newtheorem{Thm}{Theorem}
\newtheorem{Cor}{Corollary}
\newtheorem{Rmk}{Remark}[section]
\newtheorem{Prop}[Rmk]{Proposition}
\newtheorem{Lem}[Rmk]{Lemma}
\newtheorem{Def}[Rmk]{Definition}

\newtheorem{Lemma}{Lemma}

\def\bC {\mathbf{C}}

\def\N {\mathbf{N}}
\def\bR {\mathbf{R}}

\def\Z {\mathbf{Z}}

\def\R {\mathbf{R}}

\newcommand{\mean}[1]{\left\langle #1\right\rangle}

\def\cD {\mathcal{D}}

\def\cF {\mathcal{F}}
\def\ccM {\mathcal{M}}

\def\cT {\mathcal{T}}

\def\eps {{\epsilon}}

\def\Om {{\Omega}}

\def\indc {{\bf 1}}

\def\d {{\partial}}

\newcommand{\Div}{\operatorname{div}}

\newcommand{\curl}{\operatorname{curl}}

\newcommand{\sgn}{\operatorname{sign}}

\newcommand{\Supp}{\operatorname{Supp}}

\newcommand{\supp}{\operatorname{Supp}}
\newcommand{\psib}{\bar \psi}
\newcommand{\ga}{\textgoth{a}}
\newcommand{\gb}{\textgoth{b}}
\newcommand{\ba}{\begin{aligned}}
\newcommand{\ea}{\end{aligned}}

\newcommand{\T}{\mathbb T}
\newcommand{\be}{\begin{equation}}
\newcommand{\ee}{\end{equation}}

\newcommand{\psbl}{\psi^{BL}}
\newcommand{\psibl}{\psi^{BL}}

\newcommand{\psapp}{\psi^{app}}
\newcommand{\chib}{\chi_0}

\newcommand{\psil}{\psi^{\mathrm{lift}}}

\newcommand{\psisig}{\psi^\Sigma}

\newcommand{\dtl}{\delta\tau^{\mathrm{lift}}}
\newcommand{\viscosite}{\mathfrak{E}}
\renewcommand{\nu}{\mathfrak{E}}

\def\psiapp {{\psi_{app}}}
\def\dpsi {{\delta \psi}}
\def\dtau {{\delta \tau}}

\let\ds=\displaystyle

\date{\today}

\begin{document}

\frontmatter

\title[Degenerate boundary layers ]
{Mathematical Study\linebreak[1] of Degenerate Boundary Layers:\linebreak[1]
	A Large Scale Ocean Circulation Problem}
\author{Anne-Laure Dalibard}
\address{UPMC Univ Paris 06, UMR 7598 Laboratoire Jacques-Louis Lions,
	Paris, F-75005 France --- and ---
	CNRS, UMR 7598 LJLL, Paris, F-75005 France}
\email{dalibard@ann.jussieu.fr}
\author{Laure Saint-Raymond}
\address{UMPA- UMR 5669, \'Ecole normale sup\'erieure de Lyon, 46 all\'ee d'Italie,
	69364 Lyon Cedex 07, France}
\email{laure.saint-raymond@ens-lyon.fr}

\setcounter{page}{3}

\begin{abstract} This paper is concerned with a complete asymptotic
	analysis as $\viscosite \to 0$ of the stationary Munk equation $\partial_x\psi-\viscosite \Delta^2 \psi=
	\tau$ in a domain $\Omega\subset \mathbf{R}^2$, supplemented with boundary
	conditions for $\psi $ and $\partial_n \psi$. This equation is a simple
	model for the circulation of currents in closed basins, the variables
	$x$ and $y$ being respectively the longitude and the latitude. A crude
	analysis shows that as $\viscosite \to 0$, the weak limit of $\psi$ satisfies
	the so-called Sverdrup transport equation inside the domain, namely
	$\partial_x \psi^0=\tau$, while boundary layers appear in the vicinity of
	the boundary.
	
	These boundary layers, which are the main center of interest of the
	present paper, exhibit several types of peculiar behaviour. First, the
	size of the boundary layer on the western and eastern boundary, which
	had already been computed by several authors, becomes formally very
	large as one approaches northern and southern portions of the boudary,
	i.e. pieces of the boundary on which the normal is vertical. This
	phenomenon is known as geostrophic degeneracy. In order to avoid such
	singular behaviour, previous studies imposed restrictive assumptions
	on the domain $\Omega$ and on the forcing term $\tau$. Here, we prove
	that a superposition of two boundary layers occurs in the vicinity of
	such points: the classical western or eastern boundary layers, and
	some northern or southern boundary layers, whose mathematical
	derivation is completely new. The size of northern/southern boundary
	layers is much larger than the one of western boundary layers
	($\viscosite^{1/4}$ vs. $\viscosite^{1/3}$). We explain in detail how the superposition
	takes place, depending on the geometry of the boundary.
	
	Moreover, when the domain $\Omega$ is not connex in the $x$ direction,
	$\psi^0$ is not continuous in $\Omega$, and singular layers appear in
	order to correct its discontinuities. These singular layers are
	concentrated in the vicinity of horizontal lines, and therefore
	penetrate the interior of the domain $\Omega$. Hence we exhibit some kind
	of boundary layer separation. However, we emphasize that we remain
	able to prove a convergence theorem, so that the singular layers
	somehow remain stable, in spite of the separation.
	
	Eventually, the effect of boundary layers is non-local in several
	aspects. On the first hand, for algebraic reasons, the boundary layer
	equation is radically different on the west and east parts of the
	boundary. As a consequence, the Sverdrup equation is endowed with a
	Dirichlet condition on the East boundary, and no condition on the West
	boundary. Therefore western and eastern boundary layers have in fact
	an influence on the whole domain $\Omega$, and not only near the
	boundary. On the second hand, the northern and southern boundary layer
	profiles obey a propagation equation, where the space variable $x$
	plays the role of time, and are therefore not local.
\end{abstract}


\maketitle

\tableofcontents

\mainmatter
\input{intro}

\input{Multiscale}

\input{approximate}

\input{convergence}

\input{conclusion}

\section*{Acknowledgments} We acknowledge the support of the ANR
(grant ANR-08-BLAN-0301-01 and grant ANR Dyficolti
ANR-13-BS01-0003-01). This project has received funding from the European Research Council (ERC) under the European Union's Horizon 2020 research and innovation program Grant agreement No 637653, project BLOC ``Mathematical Study of Boundary Layers in Oceanic Motion''.
We thank Franck Sueur and David G\'erard-Varet
for pointing out to us the literature on convection diffusion
equations (articles \cite{grasman,jung-temam}), and Didier Bresch for
his insight on the question of islands. We also thank the referees for
their valuable comments.


\appendix
\input{app-islands}

\input{appendixB}
\input{appendixC}
\input{appendixD}

\input{notations}
\backmatter

\end{document}

%% file: intro.tex
\chapter{Introduction}

The present paper is mainly concerned with mathematical methods  investigating singular behaviours on the boundary of a bounded domain $\Omega$, when the size of the boundary layer
depends strongly on its localization, and more precisely when it becomes degenerate on some part of the boundary.

Such a situation has been depicted in many former works, but never really dealt with insofar as additional assumptions were often made to guarantee that  boundary terms vanish in the vicinity of the singularity. This is the case for instance of  our study \cite{DSR} of the $\beta$-plane model for rotating fluids in a thin layer where we suppose that the wind forcing vanishes at the equator. The same type of assumption was used in the paper  \cite {Ro2} of F. Rousset which investigates the behaviour of Ekman-Hartmann boundary layers on the sphere. This holds also true for the work of Desjardins and Grenier \cite{DG} on Munk and Stommel layers where it is assumed that the Ekman pumping (which is directly related to the wind forcing) is zero in the vicinity of the Northern and Southern coasts. The difficulty was pointed out by D. G\'erard-Varet and T. Paul in \cite{GVP}.

We intend here to get rid of this non physical assumption on the forcing, and to obtain a mathematical description of the singular boundary layers. More generally, we would like to understand how to capture the effects of the geometry in such problems of singular perturbations on domains with {possibly characteristic} boundaries.

\section{Munk boundary layers}

The equation we will consider, the so-called Munk equation, can be written as follows when the domain $\Om$ is simply connected:
\begin{equation}
\label{M}
\begin{aligned}
 \d_x \psi -\viscosite \Delta^2 \psi =  \tau \hbox{ in } \Omega,\\
\psi_{|\d \Omega} = 0,\quad (n\cdot \nabla \psi)_{|\d \Omega} =0.
\end{aligned}
\end{equation}
This model comes from large-scale oceanography and is expected to provide a good approximation of the stream function $\psi$ of oceanic currents, assuming that
\begin{itemize}
\item  the motion of the (incompressible) fluid is purely two-dimensional $\Omega \subset \R^2$ (shallow-water approximation),
\item the wind forcing is integrated as a source term $\tau \in W^{4,\infty}(\bar \Omega)$ (Ekman pumping),\label{tau}
\item  the Coriolis parameter depends linearly on latitude ($\beta$-plane approximation),
\item nonlinear effects are negligible,
\item the Ekman pumping at the bottom is negligible.
\end{itemize}
We refer to Chapter 5 for a derivation of the equation together with a discussion of the physical approximations involved.
When $\Om$ is not simply connected, the boundary conditions are slightly different. We refer to paragraph \ref{ssec:islands} for more details.

\subsection{State of the art}$ $
\label{formal-asymptotic}
Desjardins and Grenier studied in \cite{DG} a time-dependent and nonlinear version of \eqref{M}. Their result implies  in particular the following  statement for the linear Munk  equation:

\begin{Thm}\label{thm:DG}
Let  $\Omega $ be a smooth  domain, defined by
\be\label{de:om-DG}
\Om:=\left\{(x,y)\in \R^2,\ x_W(y)<x<x_E(y),\ y_{min}<y< y_{max}\right\},
\ee
where $x_E, x_W \in \mathcal C^2(y_{min}, y_{max})$.

Assume  that the wind forcing $\tau \in  H^s(\Omega)$   vanishes identically in the vicinity of the North and of the South
$$\exists\lambda>0,\quad \tau (x,y) \equiv 0 \hbox{ if } y\leq y_{min}+\lambda \hbox{ or } y\geq y_{max}-\lambda\,.\leqno (A0)$$

  Denote by $u_\viscosite =\nabla^\perp \psi_\viscosite $  any solution to the vorticity formulation of the 2D Stokes-Coriolis system
\be\label{Munk-evol}
\begin{aligned}
\viscosite \d_t \Delta \psi  +  \d_x \psi  -\viscosite \Delta^2 \psi = \tau ,\\
\psi_{|\d \Omega} = 0,\quad (n\cdot \nabla \psi)_{|\d \Omega} =0.
\end{aligned}
\ee
Then,  for all $N\in \mathbb N$, if $s$ is sufficiently large, there exists $$u_{app}^N=\sum_{i=0}^N \viscosite^{i/3} (u^{int}_i + u^{BL}_i)$$ such that
$$
\|u_\viscosite-u_{app}^N\|_{L^\infty((0,T), L^2(\Om))}\leq C \viscosite^{(N-1)/6}.
$$
The main order term $u^{int}_0=\nabla^\perp \psi^0$ is the weak limit of $u_\viscosite$ in $L^2$ as $\viscosite\to 0$ and satisfies the Sverdrup relation
\begin{equation}\label{sverdrup}
 \d_x  \psi ^0=\tau, \qquad  \psi ^0_{|\Gamma_E}=0,
\end{equation}
where  $\Gamma_E$ is  the East boundary
$$\Gamma_E = \{ (x,y) \in \d\Omega,x =x_E(y)\}\,.$$
The boundary layers $u_i^{BL}$ are located in a band of width $\viscosite^{1/3}$ in the vicinity of {  the East and West components of} $\d\Om$.
\end{Thm}

The approximate solution $u_{app}$ is computed starting from  asymptotic expansions in terms of the  parameter $\viscosite$.
\begin{itemize}
\item At main order  in the interior of the domain, we get the Sverdrup relation
$$ \d_x\psi^0 =\tau\,.$$
We can then prescribe only one boundary condition, either on the Eastern coast or on the Western coast.

\item Since $u^{int}_0 $ does not vanish on the boundary, we then introduce boundary layer corrections, which are given by the balance between
$ \d_x \psi^{BL}$ and $\viscosite \Delta^2 \psi^{BL}$.
Because the space of admissible (localized) boundary corrections is of dimension 2 on the Western boundary, and of dimension 1 on the Eastern boundary,
we can recover all the boundary conditions except $\psi^0_{|\Gamma_E} =0$, which is chosen as the boundary condition for the Sverdrup equation.

\end{itemize}
The proof relies then on a standard energy method. The assumption on $\tau$ guarantees that $\psi^0$ is smooth, and that the coefficients arising in the definition of the boundary layer terms are uniformly bounded, and  small enough to be absorbed in the viscosity term. Desjardins and Grenier actually study a nonlinear version of \eqref{Munk-evol}, for which a smallness condition on $\tau$ is required  to deal with the nonlinear convection term.
%
%

\subsection{Boundary layer degeneracies} $ $

Our main goal in this paper is to get rid of the assumption $(A0)$ on the forcing $\tau$, and to treat more general domains $\Om$ (e.g. domains with islands, or which are not convex in the $x$  direction). { This problem related to singular boundary layers has already been studied (at the formal level) by De Ruijter \cite{DR} for some simplified geometries involving typically rectangles.

The main difficulty to get rid of assumption $(A0)$ }  is that the incompatibility between the formal limit equation and the boundary conditions
$$\psi_{|\d \Omega} = 0,\quad (n\cdot \nabla \psi)_{|\d \Omega} =0$$
is not of the same nature depending on the part of the boundary to be considered:
 at the North and at the South, the transport term is tangential to the boundary  and therefore is not expected to be singular even in boundary layers.

Introducing boundary layer terms to restore the boundary conditions, we will find typically
$$ \d_x \psi_{E,W} -\viscosite \d_x^4 \psi_{E,W} =0\,,$$
on the lateral boundaries, and
$$ \d_x \psi_{N,S} -\viscosite \d_y^4 \psi_{N,S}=0\,,$$
on the horizontal boundaries.
This implies in particular that western and eastern boundary layers should be of size $\viscosite^{1/3}$ while northern and southern boundary layers should be much larger, of size $\viscosite^{1/4}$. Of course there are superposition zones, which have to be described rather precisely if we want to get an accurate approximation, { and prove a rigorous convergence result.}

This is actually  the first step towards the understanding of more complex geometries. Additional difficulties will be discussed in the next section.

\subsection{Stability of the stationary Munk equation}

Another important difference with \cite{DG} comes from the fact that (\ref{M}) is a stationary equation, so that classical energy methods (\cite{grenier, schochet}) are irrelevant.

Denote by $\psiapp$ the solution to the following approximate equation
\begin{equation}
\label{Mapp}
\begin{aligned}
 \d_x \psiapp -\viscosite \Delta^2 \psiapp =\tau +\dtau \hbox{ on } \Omega,\\
\psi_{app|\d \Omega} = 0,\quad (n\cdot \nabla \psi_{app})_{|\d \Omega} =0.
\end{aligned}
\end{equation}
In particular, denoting $\dpsi =\psi-\psiapp$, we have
\begin{equation}
\label{Mstab}
\begin{aligned}
 \d_x \dpsi -\viscosite \Delta^2 \dpsi =\dtau \hbox{ on } \Omega,\\
\dpsi_{|\d \Omega} = 0,\quad (n\cdot \nabla \dpsi)_{|\d \Omega} =0.
\end{aligned}
\end{equation}

\bigskip
{ The stability estimates used in this paper rely on the use of}  {\bf  weighted spaces} as suggested by  Bresch and Colin \cite {BC}. From the identities
$$-\int_\Omega( \d_x \dpsi ) e^x \dpsi =\frac12 \int_\Omega (\dpsi)^2 e^x,$$
$$\int_\Omega(\Delta ^2 \dpsi ) e^x \dpsi =\int_\Omega(\Delta  \dpsi )^2 e^x+2 \int_\Omega e^x \Delta \dpsi\; \d_x \dpsi +\int_\Omega e^x \Delta \dpsi\; \dpsi,$$
together with Cauchy-Schwarz inequality, we indeed deduce that
$$ \frac{1 - \viscosite}2  \int_\Omega \dpsi^2 e^x+\frac \viscosite2 \int_\Omega(\Delta  \dpsi )^2 e^x\leq -\int_\Omega \delta\tau \dpsi e^x -2\viscosite \int_\Omega \Delta \dpsi \d_x \dpsi e^x \,.$$
Assume that $\delta\tau$ can be decomposed into
$$
\delta\tau=\delta \tau_1 + \delta \tau_2,
$$
with $\delta\tau_1\in L^2(\Om)$, $\delta \tau_2\in H^{-2}(\Om)$. Then, using the Poincar\'e inequality,
$$
\left|\int_\Omega \delta\tau \dpsi e^x\right|\leq C(\|\delta\tau_1\|_{L^2}\|\dpsi\|_{L^2} + \|\delta \tau_2\|_{H^{-2}}\|\Delta\dpsi\|_{L^2}),
$$
while
$$
\left|\int_\Omega e^x  \Delta \dpsi\; \d_x \dpsi \right|\leq C \|\Delta \dpsi\|_{L^2}\|\d_x \dpsi\|_{L^2}\leq C  \|\Delta \dpsi\|_{L^2}^{3/2}\|\dpsi\|_{L^2}^{1/2}.
$$

It comes finally, assuming that $\viscosite\ll 1$,
\begin{equation}
 \label{weight-est}
\|  \dpsi\|_{L^2}^2 +\viscosite \| \Delta  \dpsi\|_{L^2}^2 \leq C\left({\| \dtau_1 \|_{L^2}^2}+\frac{\| \dtau_2 \|_{H^{-2}}^2}{\viscosite}\right).
   \end{equation}

\bigskip
We then define some relevant notion of approximate solution.

\begin{Def}\label{def:app-general}
A function $\psiapp \in H^2(\Omega)$ is an approximate solution to (\ref{M}) if it satisfies the approximate equation (\ref{Mapp}) for some $\dtau \in H^{-2}(\Omega)$ such that $\dtau=\dtau_1+\dtau_2$ with
\begin{equation}
\label{app-def}
\lim_{\viscosite \to 0}  \frac{{\| \dtau_1 \|_{L^2}^2}+\viscosite^{-1} {\| \dtau_2 \|_{H^{-2} }^2}}{  \|  \psiapp\|_{L^2}^2 +\viscosite \| \Delta \psiapp \|_{L^2}^2 } =0\,.
\end{equation}
\end{Def}
Plugging \eqref{app-def} in  \eqref{weight-est}  we obtain
$$
  \| \psi- \psiapp\|_{L^2}^2 +\viscosite \| \Delta(\psi- \psiapp )\|_{L^2}^2 =o (  \|  \psiapp\|_{L^2}^2 +\viscosite \| \Delta \psiapp \|_{L^2}^2 )$$
 meaning that $\psiapp\sim \psi$.

By such a method, we exhibit the {\it dominating phenomena in terms of their contribution to the energy balance}. In particular, for the oceanic motion, we expect  to justify the crucial role of boundary currents as they account for a macroscopic part of the energy.

{
In view of  the energy estimate \eqref{weight-est} and of the expected sizes of boundary layers (see \eqref{est:L2} and \eqref{est:H2} below), the idea is to prove that $\psiapp$ satisfies equation \eqref{Mapp} with an error term $\dtau=\dtau_1+\dtau_2$ such that
\be\label{hyp:dtau}
\|\dtau_1\|_{L^2(\Om)}=o(\viscosite^{1/8}),\quad \|\dtau_2\|_{H^{-2}(\Om)}=o(\viscosite^{5/8}).
\ee
Notice that the energy estimate \eqref{weight-est} allows us to capture both types of boundary layers and the interior term, as we will explain in Theorem \ref{thm}.

\begin{Def}\label{remainder-def}
In the rest of the paper, we  say that $\dtau=\dtau_1+\dtau_2$ is an admissible error term if it satisfies \eqref{hyp:dtau}.
\label{def:admissible}
\end{Def}

Our main result, Theorem \ref{thm}, is the construction of an explicit approximate solution with an admissible error term in the sense of Definition \ref{remainder-def}.
}
%
%
%
%
%
%
%
%
%
%
%
%
%
%

 \section{Geometrical preliminaries}
\label{sec:geom}

The goal of this section is to state the precise assumptions regarding the domain,
under which we are able to derive a nice approximation for the Munk problem.

As we will see in the course of the proof, error estimates depend strongly on the geometry of the domain, and especially
on the flatness of the boundary near horizontal parts.

We will therefore consider
\begin{itemize}
\item the generic case, when the slope of the tangent vector vanishes polynomially;

\item  a special case of flat boundary, for which the exponential decay is prescribed;

{
\item the case of a corner between the East/West boundary and the horizontal part  (as considered by De Ruijter in \cite{DR}).}
\end{itemize}

The specific assumptions in the case of corners will be considered separately in subsection \ref{corner-par}, in an attempt to keep the presentation as simple as possible. Therefore the regularity and geometry assumptions detailed in subsections \ref{ssec:regularity}, \ref{ssec:singularity-lines}, and \ref{ssec:islands} apply to the first two cases above.

\subsection{Regularity and flatness assumptions}$ $
\label{ssec:regularity}
First of all we assume, without loss of generality,  that  $\Omega$ is bounded and connected\footnote{If $\Om$ is not connected, since \eqref{M} is a local equation, we can perform our study on every connected component, and we therefore obtain a result on the whole domain $\Om$.}. We further require that
$$ \hbox{ the boundary is a compact $C^4$ manifold,  described by a finite number of charts.}\leqno (H1)$$
In particular, the perimeter $L$  is finite, and  the boundary of $\Om$ can be parametrized by the arc-length $s$. If $\d\Om$ is connected, we will often write  $\d\Om=\{(x(s),y(s))\,/\, 0\leq s \leq L\}$; if $\d\Om$ has several connected components, such a representation holds on every connected component of $\d\Om$. We further choose the orientation of the arc-length in such a way that the local frame $(t,n)$ with $n$ the exterior normal is direct.
\label{localcoord}

Furthermore, by the local inversion theorem, there exists $\delta>0$ such that
any point $(x,y)\in \Omega$ at a distance $z$ smaller than $\delta$ from $\d \Omega$ has  a unique projection on the boundary. In other words, on such a tube, the arc-length $s$ and  the distance  to the boundary  $z$ form a nice system of coordinates.

We thus introduce some truncation function $\chib\in C^\infty_c([0,+\infty), [0,1])$ such that\label{chib}
\begin{equation}
\label{chi0-def}
 \hbox{Supp} \chib \subset [0,\delta) \hbox{ and }  \chib\equiv 1 \hbox{ on } [0,\delta/2]\,.
 \end{equation}

\bigskip
As mentioned in the previous heuristic study, we expect the boundary layers to be singular on the horizontal parts of the boundary. In order to understand the connection between  both types of boundary layers and to get a nice approximation of the solution to the Munk problem, we therefore need precise information on the profile of the boundary near horizontal parts.

We assume that the horizontal part $\Gamma_N \cup \Gamma_S$ of the boundary consists in a finite number of intervals (possibly reduced to points where the tangent to the boundary is horizontal).
Denote by $\theta(s)$\label{theta} the oriented angle  between the horizontal vector $e_x$ and the exterior normal $n$. By definition, we set\label{GammaNS}
$$
\begin{aligned}
\Gamma_N:=\{s\in \d\Om,\   \cos \theta(s) =0\hbox{ and } \sin \theta(s) =1\},\\
 \Gamma_S:=\{s\in\d\Om,\   \cos \theta(s) =0\hbox{ and } \sin \theta(s) =-1\}.\\
\end{aligned}$$
\begin{figure}
 \begin{center}
 \includegraphics[height = 8cm]{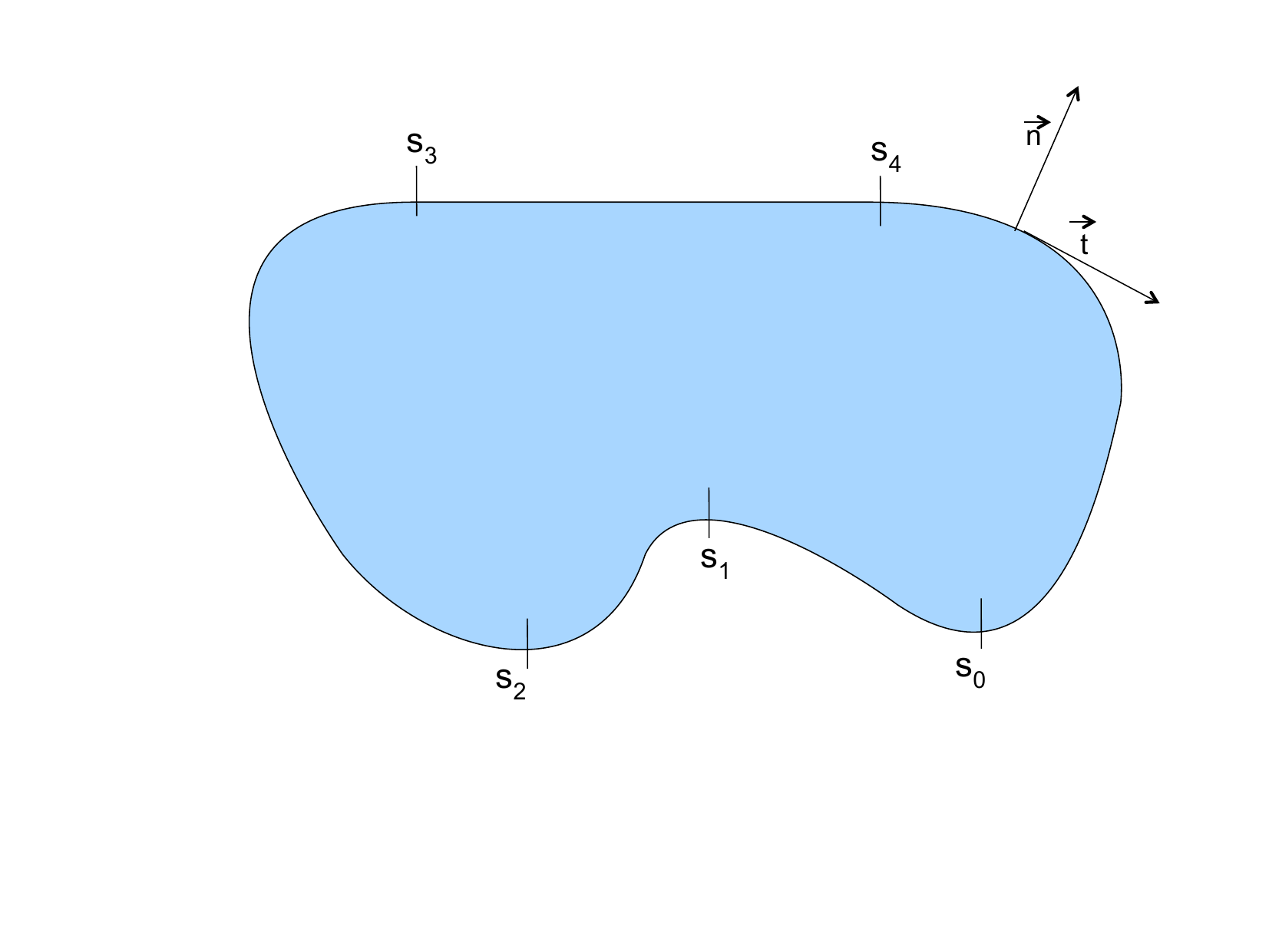}
\end{center}
\caption{Positions of the abscissa $s_i$}\label{fig:abscisses}
\end{figure}

We also introduce the following notation (see figure \ref{fig:abscisses}): let $s_1<s_2<\cdots<s_k$\label{si} such that
$$
        \cos(\theta(s_i))=0\quad 1\leq i\leq k
        $$
and for all $i\in\{1,\cdots, k\}$, either $\cos(\theta(s))$ is identically zero  or $\cos(\theta(s))$ does not vanish between $s_i$ and $s_{i+1}$. Throughout the article, we use the conventions $s_{k+1}=s_1$, $s_0=s_k$. We further define some partition of unity $(\rho_i)_{1\leq i \leq k}$
\begin{equation}
\label{rho-def}
\begin{aligned}
\rho_i \in C^\infty_c(\d\Om,[0,1]), \quad \sum_{i=1}^k \rho_i(s) =1,\\
\supp \rho_i\subset ]s_{i-1}, s_{i+1}[.
\end{aligned}
\end{equation}

We denote by $\Gamma_E$, $\Gamma_W$\label{GammaEW} the East and West boundaries of the domain:
$$
\begin{aligned}
\Gamma_E:=\{s\in \d\Om,\quad \cos \theta(s)>0\},\\
\Gamma_W:=\{s\in \d\Om,\quad \cos \theta(s)<0\}.
\end{aligned}
$$
Eventually, let
\be\label{I+}
I_+:=\{i\in \{1,\cdots, k\},\ s_i\in \d\Gamma_E\}.
\ee

\bigskip
The profile assumption states then as follows: for any $s_i$, for $\sigma =\pm$ such that $\cos \theta(s) \neq 0$ on $[s_i, s_{i\sigma1}]$,
\begin{enumerate}[(i)]
\item either    there exists $n\geq 1$ and $C\neq 0$ (depending on $i$ and $\sigma$) such that as $s\to s_i$, $s\in (s_i, s_{i\sigma1}) $,
        $$
        \begin{aligned}
                \cos \theta(s)\sim \frac{C}{n!}(s-s_i)^n,\\
                \text{and }\theta^{(l)}(s) \sim \frac{C}{(n-l)!}(s-s_i)^{n-l}\text{ for }1\leq l\leq \inf(3,n);
        \end{aligned} \leqno (H2(i))
        $$

        \item or there exists $C\neq 0$ and $\alpha>0$ (also depending on $i$ and $\sigma$) such that  as $s\to s_i$, $s\in (s_i, s_{i\sigma1}) $,
        $$
        \begin{aligned}
                \cos \theta(s)\sim C e^{-\alpha/|s-s_i|},\\
                \text{and } \theta'(s)\sim\frac{C\alpha}{(s-s_i)^2}e^{-\alpha/|s-s_i|}.
        \end{aligned}\leqno (H2(ii))
        $$
\end{enumerate}
The first situation corresponds to the generic case when the cancellation is of finite order. The second one is an example of infinite order cancellation: in that case, which is important since it is the archetype of $C^\infty$ boundary with flat parts, we prescribe the exponential decay because there is no general formula for error estimates. Notice that we do not require the behaviour of $\theta$ to be the same on both sides of $s_i$, provided the function $\theta$ belongs to $\mathcal C^3(\d \Om)$, so that $(x(s), y(s))\in \mathcal C^4(\d\Om)$.

This profile assumption will essentially guarantee that, up to a small truncation, we will be able to lift boundary conditions either by East/West boundary layers, or by North/South boundary layers at any point of the boundary. This is therefore  the main point to get rid of assumption $(A0)$.

\subsection{Singularity lines}$ $
\label{ssec:singularity-lines}
\begin{figure}[h]
 \includegraphics[width=0.8\textwidth]{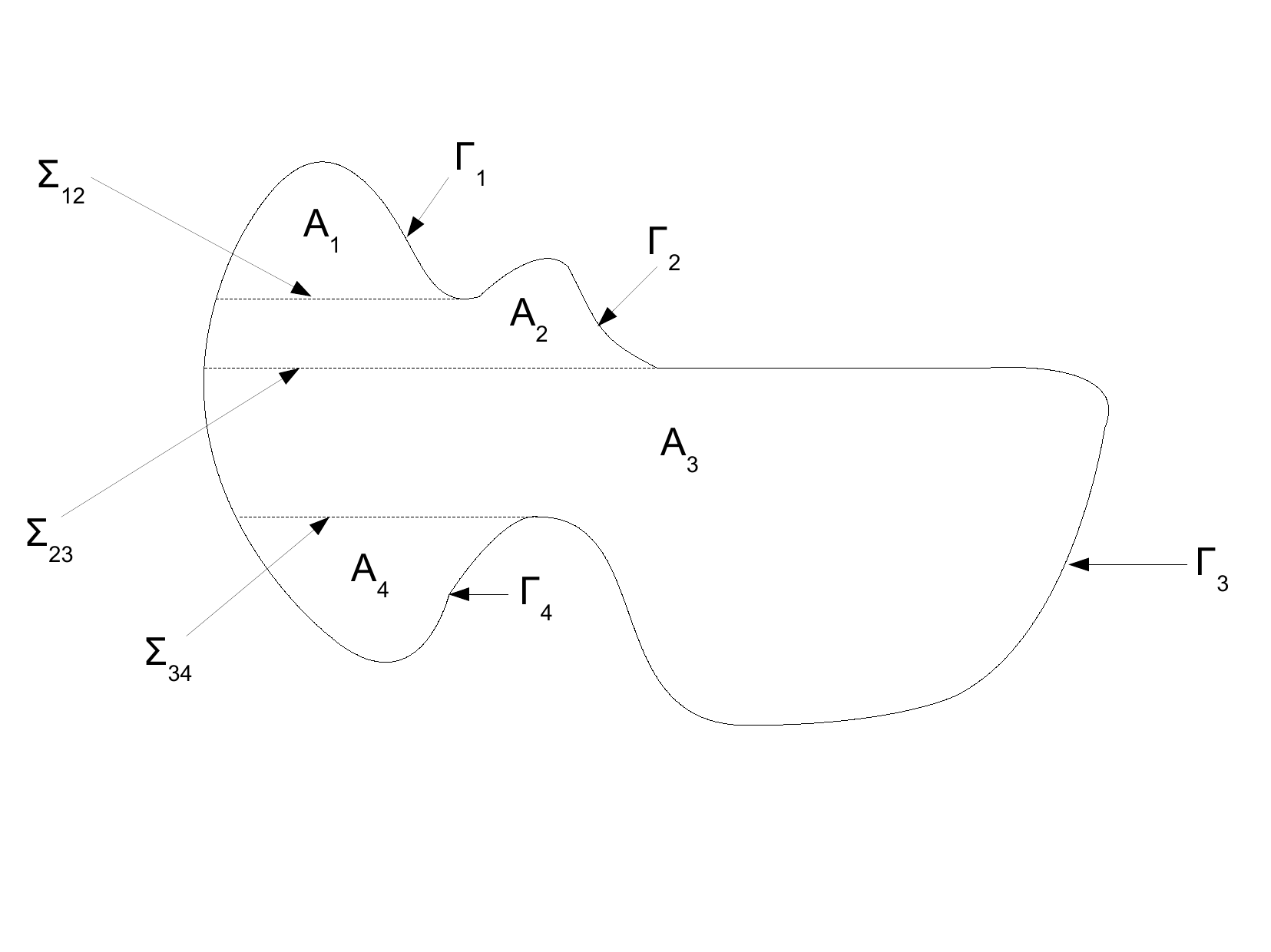}
\caption{Discontinuity lines of a general domain $\Om$}\label{fig:complexe_general}
\end{figure}
In the case when the domain $\Omega$ is not convex in the $x$ direction, we will see that the asymptotic picture is much more complex, especially because the solution $\psi^0$ to the Sverdrup equation (\ref{sverdrup}) is discontinuous as soon as $\bar \Gamma_E$ has more than two connected components (notice that if $\cos \theta$ has an isolated point of cancellation in the interior of $\bar\Gamma_E$, so that $\Gamma_E$ has two connected components but $\bar \Gamma_E$ is connected, then there is no discontinuity in $\psi^0$.)

We therefore introduce the lines $\Sigma_{ij}$, across which the main order term will be discontinuous, which will give rise to boundary layer singularities: we set
$$
\bar \Gamma_E:=\bigcup_{j=1}^M \Gamma_j,
$$
where $\Gamma_1, \cdots, \Gamma_M$ are the closed connected components of $\bar \Gamma_E$. For $j\in \{1,\cdots, M\}$, we set\label{A-def}
$$
A_j:=\left\{(x,y)\in \Om,\ \exists x'\in \bR,\ (x',y)\in \Gamma_j\text{ and }(tx+(1-t)x',y)\in \Om\ \forall t\in ]0,1[\right\}.
$$
We have clearly
$$
\Om=\bigcup_{j=1}^M A_j.
$$
We also define (see Figure \ref{fig:complexe_general}, and also Figures 3.2 and 3.3 in \cite{DR})
\be\label{def:Sigma}
\Sigma_{ij}:=\bar A_i\cap \bar A_j\quad\text{for }i\neq j, \quad \Sigma:=\bigcup_{i,j} \Sigma_{ij}.
\ee
It can be easily checked that every set $\Sigma_{ij}$ is either empty or a horizontal line with ordinate $y_{ij}$ such that there exist $x_i, x_j\in \bR$ with  $(x_i, y_{ij})\in   \Gamma_i$,  $(x_j, y_{ij})\in\Gamma_j$, and either $(x_i, y_{ij})\in \d\Gamma_i$ or $(x_j, y_{ij})\in\d\Gamma_j$.

Eventually, we parametrize every set $\Gamma_j$ by a graph $x_E^j$, namely
$$
\Gamma_j=\{(x_E^j(y), y),\ y_{min}^j\leq y \leq y_{max}^j\}.
$$

We will therefore need to build singular correctors which are not localized in the vicinity of the boundary. This construction is rather technical, and for the sake of simplicity, we will use two  additional assumptions, which are quite general:

\begin{figure}[h]
	 \includegraphics[height=6cm]{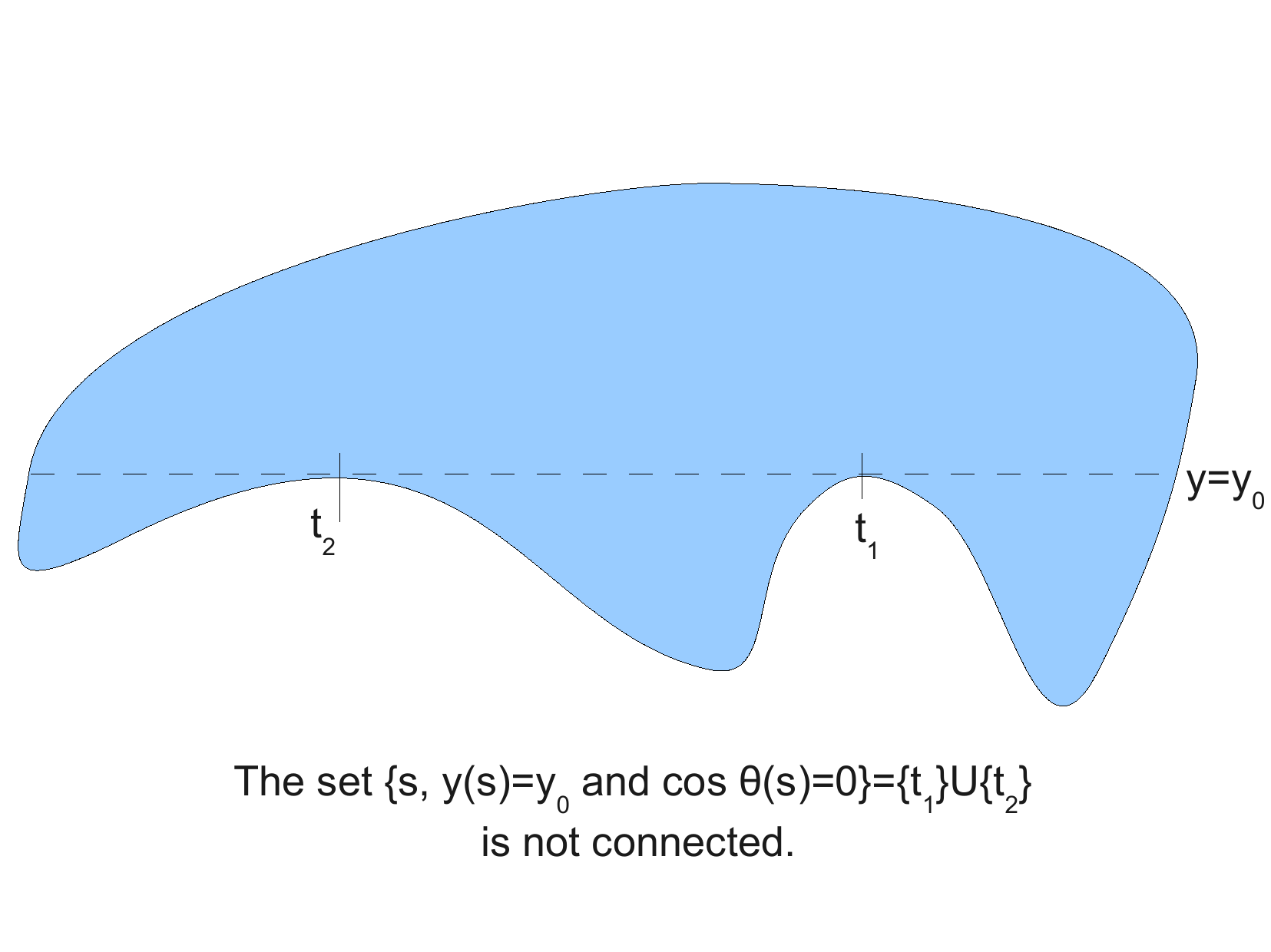}
\caption{Example of a domain when assumption (H3) is not satisfied}\label{fig:H3}
\end{figure}

(H3) $\forall y_0\in \R$, the set
$$
\{s\in \d\Om,\ y(s)=y_0,\ \cos \theta(s)=0\}
$$
is a connected set (see Figure \ref{fig:H3}).

(H4) Let $s_j\in \{s_1, \cdots, s_k\}$ be a boundary point such that
\begin{itemize}
 \item $s_j\in\d\Gamma_E\cap \d\Gamma_W$;
\item and $\Om$ is not convex in $x$ in a neighbourhood of $(x(s_j), y(s_j))$.
\end{itemize}
Then $\cos \theta(s)=O(|s-s_j|^4)$  for $s\to s_j$.
\smallskip

Assumption (H4) will be discussed in Remark \ref{rem:hyp-h4}.

\subsection{Domains with islands} $ $
\label{ssec:islands}

When the domain $\Om$ is not simply connected, the boundary conditions on $\d\Om$ are slightly different. This case has been studied in particular in \cite{BGGRB}, where the authors investigate the weak limit of the Munk equation in a domain with islands. Let $\Om_1, \cdots \Om_K$\label{omega-i} be $\mathcal C^4$ simply connected domains of $\R^2$, such that
\begin{itemize}
 \item $\Om_i\Subset \Om_1$ for $i\geq 2$;
\item $\Om_i\cap \Om_j=\emptyset$ for $2\leq i<j\leq K$;
\item $\Om:=\Om_1\setminus \cup_{i\geq 2} \Om_i$ satisfies (H1)-(H4).
\end{itemize}
Let $C_i=\d\Om_i$ for $i\geq 1$. Notice that of course, the presence of islands gives rise to discontinuity lines $\Sigma$ as described in the preceding paragraph.

Then the Munk equation can be written as
\be
\label{Munk-islands}\begin{aligned}
\d_x \psi-\viscosite \Delta^2 \psi=\tau\text{ in }\Om,\\
\d_n \psi=0\text{ on } \d\Om,\\
\psi_{|C_1}=0,\quad \psi_{|C_i}=c_i\text{ for }i\geq 2.
\end{aligned}
\ee
The constants $c_i$ are different from zero in general: indeed, the condition $u\cdot n_{|\d\Om}=0$, where $u= \nabla^\bot \psi$ is the current velocity, becomes $\psi=\text{constant}$ on every connected component of $\d\Om$. However, the constants are not required to be all equal. In fact, the values of $c_2, \cdots, c_K$ are dictated by compatibility conditions, namely
\be\label{compatibility}
\forall j\geq 2,\quad \viscosite \int_{C_j} \d_n \Delta \psi - \int_{C_j}\cT^\bot\cdot n=0,
\ee
where $\tau=\curl\cT$.

We explain in Appendix A where condition \eqref{compatibility} comes from. The constants $c_2, \cdots, c_K$ are then uniquely determined, as shows the following
\begin{Lem}
$\bullet$ Following \cite{BGGRB}, define $\psi_1,\cdots \psi_K$ by
\begin{equation}
\label{psi1-eq}
\begin{aligned}
\d_x \psi_1- \viscosite \Delta^2 \psi_1= \tau,\\
\psi_{1|\d\Om }=0,\quad\d_n \psi_{1|\d\Om }=0,
\end{aligned}
\end{equation}
and for $i\geq 2$,
\begin{equation}
\label{psii-eq}
\begin{aligned}
\d_x \psi_i- \viscosite \Delta^2 \psi_i= 0,\\
\d_n \psi_{i|\d\Om }=0,\\
\psi_{i|C_j}=\delta_{ij}\quad\forall j\in \{1,\cdots, K\}.
\end{aligned}
\end{equation}
Then  $\psi$ satisfies equation \eqref{Munk-islands} for some constants $c_1, \cdots c_K$ if and only if
\be\label{dec:psi}
\psi= \psi_1 + \sum_{i\geq 2} c_i \psi_i.
\ee
$\bullet$ Define the matrix $M^\viscosite$ and the vector $D^\viscosite$ by
$$
\begin{aligned}
M^\viscosite:=\left(\viscosite \int_{C_i}\d_n \Delta \psi_j\right)_{2\leq i,j\leq K},\\
D^\viscosite:=\left(-\viscosite \int_{C_j}\d_n \Delta \psi_1 + \int_{C_j} \cT^\bot \cdot n\right)_{2\leq j\leq K}.
\end{aligned}
$$
Then the following facts hold:
\begin{itemize}
\item $M^\viscosite$ is invertible;
\item $\psi$ is a solution of \eqref{Munk-islands} satisfying the compatibility condition \eqref{compatibility} if and only if $\psi $ is given by \eqref{dec:psi} with $M^\viscosite c=D^\viscosite$.
\end{itemize}

\label{lem:def-M,D}
\end{Lem}
The proof of Lemma \ref{lem:def-M,D} is postponed to Appendix A.

$\bullet$ The formula \eqref{dec:psi} shows that it is sufficient to understand the asymptotic behaviour of the functions $\psi_i$. Indeed, the coefficients $c_i$ (which depend on $\viscosite$) are obtained as the solutions of a linear system involving the functions $\psi_i$. It is proved in  \cite{BGGRB} that for any domain $V\subset \Omega$ such that $\bar V \cap \Sigma=\emptyset$, where $\Sigma$ is defined by \eqref{def:Sigma}
$$\begin{aligned}
\psi_1\rightharpoonup \psi^0\quad \text{in }L^2(\Om),\\
\psi_1 \to \psi^0\quad \text{in }L^2(V),
  \end{aligned}
$$
where $\psi^0$ is the solution of the Sverdrup equation \eqref{sverdrup}, and
$$\begin{aligned}
\psi_i\rightharpoonup \mathbf 1_{B_i}\quad \text{in }L^2(\Om),\\
\psi_i \to  \mathbf 1_{B_i}\quad \text{in }L^2(V),
  \end{aligned}
$$
where $B_i$ is defined by\label{B-def}
\begin{multline*}
B_i:=\left\{(x,y)\in \Om,\ \exists x'>x,\ (x',y)\in C_i\cap \bar \Gamma_E \right.\\\left.\text{ and }(tx+(1-t)x',y)\in \Om\ \forall t\in (0,1)\right\}.
\end{multline*}
\begin{figure}
	 \includegraphics[width=\textwidth]{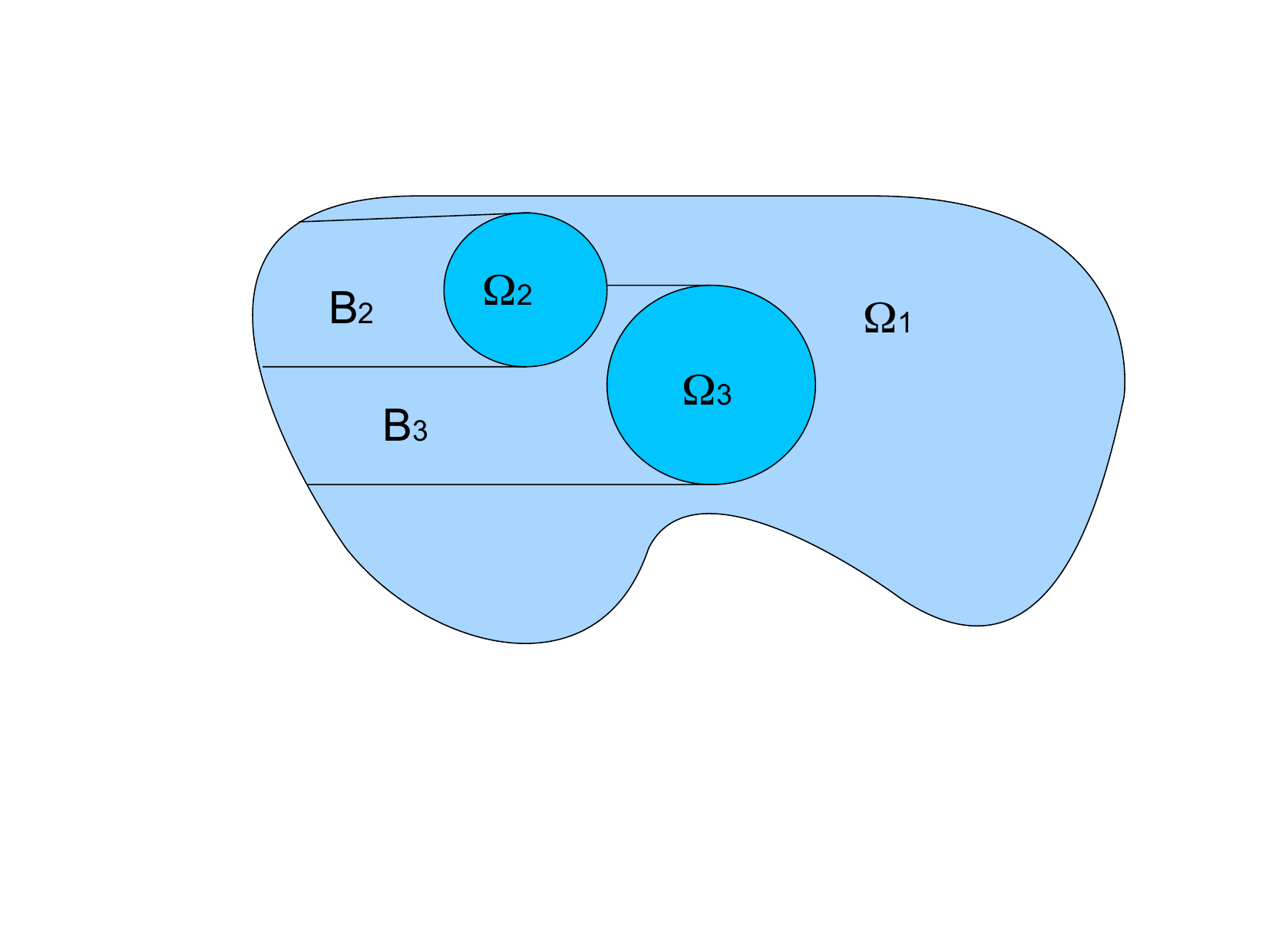}
\caption{Definition of the domains $B_i$}\label{fig:islands}
\end{figure}

This result is in fact sufficient to compute the asymptotic limit of the coefficients $c_i$, which converge towards some constants $\bar c_i$. We refer to section \ref{sec:islands} for more details. We will go one step further in the present paper, since we are able to compute an asymptotic development for the functions $\psi_i$, and therefore give a rate of convergence for the coefficients $c_i$ and the function $\psi$.

{{
\subsection{Periodic domains and domains with corners}\label{corner-par}$ $

The case when the connection between the horizontal part of the boundary and the East or West part is a corner, which is typically the case of rectangles considered by De Ruijter, is actually  easier to deal with, because there is no superposition zone for the boundary layers.

This corresponds to have a parametrization of the boundary by a function which is piecewise $\mathcal C^4$, with some jumps for the angle $\theta(s)$
$$ \cos \theta(s) \equiv 0 \hbox{ on } \Gamma_N \cup \Gamma_S, \qquad \inf_{\Gamma_E \cup \Gamma_W}|\cos \theta |>0 \,.$$
We need only to suppose that angles in West corners are obtuse.

Our arguments also allow us to investigate domains of the type $\Om= \T\times (y_-, y_+)$,  where $\T=\R/\Z$. The interest for such domains  stems from the analysis of circumpolar currents: indeed, realistic ocean basins consist of a part of a spherical surface, and for some latitudes $y\in ]y_c^-, y_c^+[$ there may be no continent. The latter part of the fluid domain is called the ``circumpolar'' component. More precisely, De Ruijter considers domains of the form
\be\label{def:circumpolar}
\Om=\Omega_{circ}\cup \Sigma_c^+ \cup \Sigma_c^- \cup \Omega_c^+\cup \Omega_c^-,\ee
where (see Figure \ref{fig:circumpolar})
$$
\begin{aligned}
\Omega_{circ} =  \T\times ]y_c^-, y_c^+[,\\
\Om_c^+ =\{(x,y)\in \T\times \R, y_c^+<y<\gamma_c^+(x) \},\\
\Om_c^- = \{(x,y)\in \T\times \R , y_c^->y>\gamma_c^-(x) \},
\\
\Sigma_c^\pm=\T^*\times \{y_c^\pm\},
\end{aligned}
$$
where $\gamma_c^\pm$ are smooth periodic functions such that $\gamma_c^\pm(0)=y_c^\pm$, and $\gamma_c^+(x)>y_c^+$ for all $ x\in (0,1)$, $\gamma_c^-(x)<y_c^-$ for all $x\in(0,1)$.

\begin{figure}
 \begin{center}
 \includegraphics[height = 8cm]{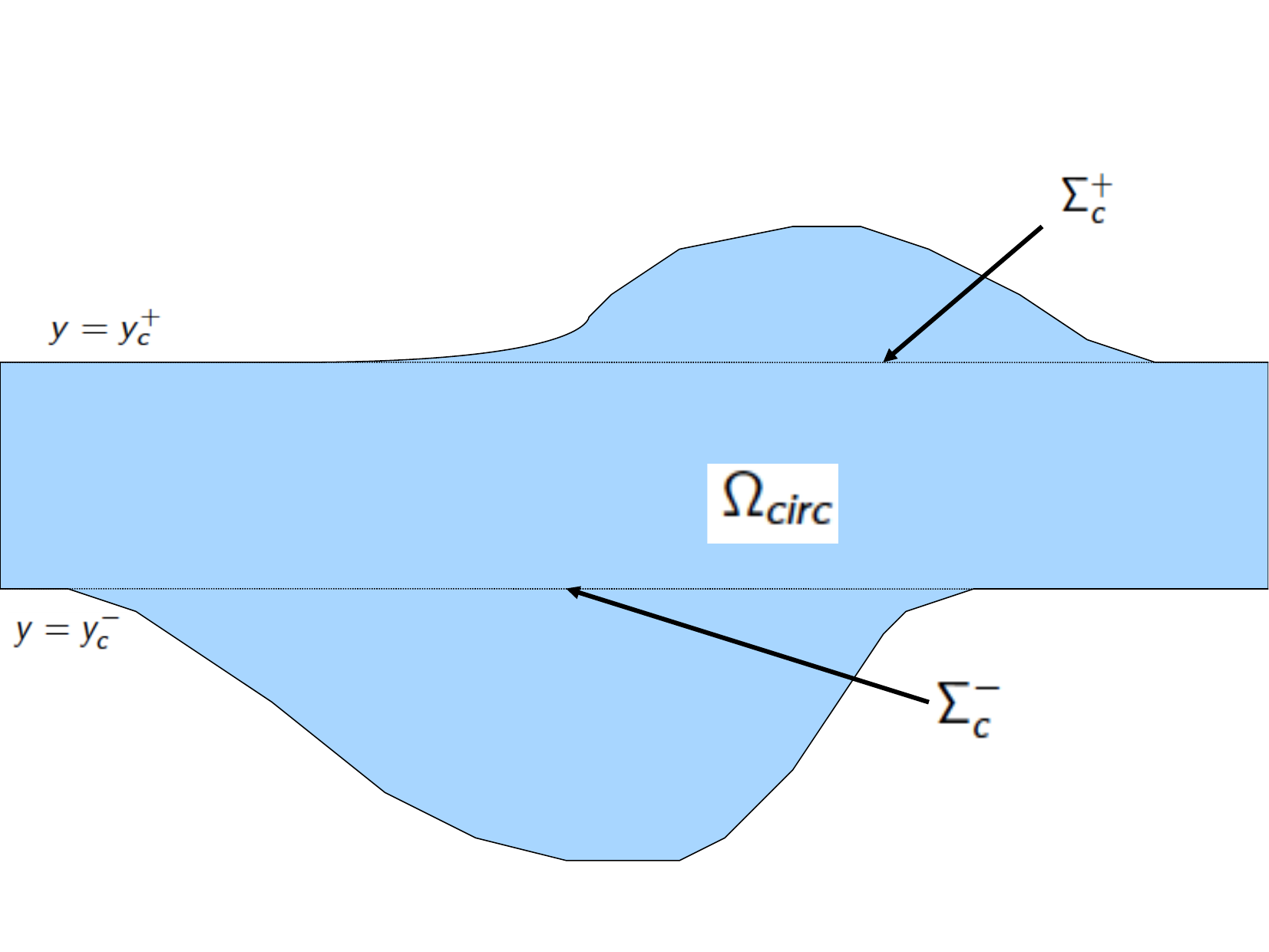}
\end{center}
\caption{Circumpolar domain}\label{fig:circumpolar}
\end{figure}

It is very likely that the analysis of the case \eqref{def:circumpolar} is in fact a combination of the arguments for ``standard'' smooth domains in $\R^2$, which are the main concern of this article, and periodic domains of the type $\T\times (y_-,y_+)$. However, because of strong singularities near the junction points $(0,y_c^\pm)$, the construction and the energy estimates become very technical, without seemingly exhibiting any new mathematical behaviour or ideas. Therefore we will focus on the periodic case and explain, within this simplified geometry, why circumpolar currents appear.

}}

\bigskip\noindent
\section{Main approximation results}

\subsection{General case}
We first describe our result in a domain $\Om\subset \R^2$ satisfying (H1)-(H4), and then explain how our result can be extended to other types of domains.

We will  prove the existence of  approximate solutions in the form
\be\label{def:psiapp}\psiapp =  \psi^{int}+\psi_{E,W}+\psi_{N,S}+ \psi^\Sigma\ee
where

\begin{itemize}
\item $\psi^{int}$ is a regularization of $\psi^0$ (or of $\psi^0+\sum_{i=2}^K \bar c_i \mathbf 1_{B_i}$ when the domain $\Om$ has islands), and $\psi^0$ is  the solution to the Sverdrup equation \eqref{sverdrup}. Let us emphasize that this regularization includes boundary layer correctors  located in horizontal bands of width $\viscosite^{1/4}$ in the vicinity of every singular line $\Sigma_{ij}$.
Furthermore, $\|\psi^{int}-\psi^0\|_{L^2}=o(\|\psi^{int}\|_{L^2})$;

\item $\psi_{E,W}$ groups together  eastern and western boundary terms, which decay on a distance of order $\viscosite^{1/3}$;
\item   $\psi_{N,S}$  is the  contribution of southern and northern boundary layers terms, which decay on a distance of order $\viscosite^{1/4}$;

\item and $\psi^\Sigma$ is an additional boundary layer term, located in horizontal bands of width $\viscosite^{1/4}$ in the vicinity of every singular line $\Sigma_{ij}$.

\end{itemize}

Notice that $\psi^{int}$, $\psi_{E,W}$, $\psi_{N,S}$ and $\psi^\Sigma$ do not have the same sizes in $L^2$ and $H^2$: typically, if the boundary condition to be lifted is of order 1, the size in $L^2$ of a boundary layer type term is $\lambda^{-1/2}$, where $\lambda^{-1}$ is the size of the boundary layer, while its size in $H^2$ is $\lambda^{3/2}$.
As a consequence, we roughly expect that
\be\label{est:L2}
\underbrace{\|\psi^{int}\|_{L^2}}_{\sim 1}\gg \underbrace{\|\psi_{N,S}\|_{L^2} \sim\|\psi^\Sigma\|_{L^2}}_{\sim \viscosite^{1/8}} \gg  \underbrace{ \|\psi_{E,W}\|_{L^2}}_{\sim \viscosite^{1/6}},
\ee
while
\be\label{est:H2}
\underbrace{\|\psi^\Sigma\|_{H^2}\sim \|\psi_{N,S}\|_{H^2} \sim\|\psi^{int}\|_{H^2}}_{\sim \viscosite^{-3/8}} \ll  \underbrace{ \|\psi_{E,W}\|_{H^2}}_{\sim \viscosite^{-1/2}},
\ee

Our result is the following;
\begin{Thm}\label{thm}
Assume that the domain $\Omega\subset \R^2$ satisfies assumptions (H1)-(H4). Consider a non trivial forcing $\tau=\curl \cT\in W^{4,\infty}(\Om)$, and let $\psi$  be the solution of the Munk equation \eqref{M} (or \eqref{Munk-islands}-\eqref{compatibility} when $\Om$ has islands).

 Then there exists a function $\psiapp$ of the form
$$\psiapp = \psi^{int}+\psi_{E,W}+\psi_{N,S}+ \psi^\Sigma$$
satisfying the approximate equation (\ref{Mapp}) with an admissible remainder in the sense of (\ref{hyp:dtau}).

More specifically, $\psi^{app}$
 is a good approximation of $\psi$, in the following sense: let $V\subset\Om$ be a non-empty open set. Then, for a generic forcing term $\tau$, the following properties hold:
\begin{itemize}
\item Approximation of the interior term: if $V \Subset \Om$ is such that $\bar V \cap \Sigma=\emptyset$, then
$$
\|\psi -\psi^0\|_{L^2(V)}=o(\|\psi^0\|_{L^2(V)});
$$
\item Approximation of the $\Sigma$ boundary layers: if $V \Subset \Om$ is such that $V\cap \Sigma_{ij}\neq \emptyset$ for some $i\neq j\in \{1,\cdots, M\}$, then
$$
\begin{aligned}
\| \psi- (\psi^{int} + \psi^{\Sigma})\|_{L^2(V)}= o( \| \psi^{\Sigma}\|_{L^2(V)}) ,\\
\| \psi- (\psi^{int} + \psi^{\Sigma})\|_{H^2(V)}= o( \|\psi^{int} + \psi^{\Sigma}\|_{H^2(V)}).
\end{aligned}
$$
\item Approximation of the North and South boundary layers: if $\bar V \cap \mathring{\Gamma}_{N,S}\neq \emptyset$ and $\bar V \cap \Gamma_{E,W}=\emptyset$, $\bar V \cap \Sigma=\emptyset$, then
$$
\begin{aligned}
\| \psi- (\psi^0 + \psi_{N,S})\|_{L^2(V)}= o( \| \psi_{N,S}\|_{L^2(V)}) ,\\
\| \psi -\psi_{N,S}\|_{H^2(V)}= o( \|\psi_{N,S}\|_{H^2(V)}).
\end{aligned}
$$
\item Approximation of the West boundary layer: if $\bar V \cap \Gamma_W\neq \emptyset$ and $\bar V \cap \Gamma_{N,S}=\emptyset$, $\bar V \cap \Sigma=\emptyset$, $\bar V \cap \Gamma_E=\emptyset$,
$$
\|\psi-\psi_W\|_{H^2(V)}= o(\|\psi_W\|_{H^2(V)}).
$$

\end{itemize}

\end{Thm}
In the above Theorem, the term ``generic forcing'' is necessarily unprecise at this stage. It merely ensures that the terms constructed in the approximate solution are not identically zero. We will give a more precise assumption in the next chapter (see Definition \ref{def:tau-generic}). The interior of $\Gamma_{N,S}$, namely $ \mathring{\Gamma}_{N,S}$, is to be understood through the induced topology on $\d\Om$. If the forcing is generic, some rough estimates on the sizes of the different terms are given in the table page \pageref{tableaux-tailles}.

{

\begin{Rmk}\label{rem:east-corr}
Note that the present result does not say anything about the validity of East boundary correctors. Indeed, their amplitude is typically $O(\viscosite^{1/3})$ on the zones where $\cos \theta$ is bounded away from zero, and therefore they are too small to be captured by the energy estimate. We will comment more on this point in section \ref{sec:macro-corrector} (see Remark \ref{rmk:east-corr2}). Note however that in non-degenerate settings, a solution can be built at any order (see \cite{DG}), and therefore the energy estimates may capture  the east corrector  in this case. Moreover, the east boundary layer equation is somewhat indirectly justified by the fact that the interior term vanishes on the east boundary.
\end{Rmk}

The construction of an approximate solution relies on a local asymptotic expansion of the form
$$
\psi\simeq \psi^0 + \psbl(s,\lambda(s) z),
$$
where $\psi^0$ is the interior term, which solves the Sverdrup equation \eqref{sverdrup}, $s$ is the arc-length, $z$ is the distance to the boundary, and $\lambda(s)\gg 1$ is the inverse of the boundary layer size. The boundary layer term $\psbl(s,Z)$ is also assumed to vanish as $Z\to \infty$. Plugging this expansion into equation \eqref{M}, we find an equation for $\psbl$, in which both $\psbl$ and $\lambda$ are unknown. The idea is then to choose $\lambda$ in a clever way, depending on the zone of the boundary under consideration, and then to solve the corresponding equation on $\psbl$. As a consequence, different equations for the boundary layer are obtained on different zones of the boundary, and a matching between these zones must be performed: this is the difficult part of the mathematical analysis, which is absent from \cite{DR}. Eventually, once the approximate solution is defined, we check that the corresponding error terms are all admissible in the sense of Definition \ref{def:admissible}.

\begin{Rmk}[About the geometric assumptions]
The geometric assumptions are mostly used in the zones of transition between  the different types of boundary layers.  Note that they may not all be necessary, but we found no systematic way of dealing with all possible cases simultaneously.
\begin{itemize}
\item $(H1)$ can be weakened since we are able to deal with corners (corners involving only East or West boundaries are handled by a simple truncation while corners with horizontal boundaries will be discussed in  chapter \ref{chap:construction});
\item $(H2)$ is used to have explicit estimates for the transition zone between East/West boundary layers and North/South boundary layers. In particular, the rate of cancellation of $\cos \theta$ determines:
\begin{itemize}
\item The size of the truncation of the forcing $\tau$ in singular zones, see \eqref{def:deltax1}, \eqref{def:deltax2};
\item The domain of validity of the east and west boundary layers, see \eqref{hyp:valid-E/W-2}, \eqref{equiv_si+_alg}, \eqref{equiv_si+_exp};
\item The zone on which energy is injected in the north and south boundary layers, see \eqref{Wsigmaipm}
\end{itemize}

\item $(H3)$ allows to avoid the connection between some $\Sigma$ layer and a horizontal boundary. It is not completely clear whether or not we would be able to handle such a singular transition;
\item $(H4)$ is a technical assumption (comparable to $(H2)$) to get good estimates on the $\Sigma$ layer, and in particular, on the impact of the $\Sigma$ layer on the west boundary. It is discussed in detail in Remark \ref{rem:hyp-h4}.
\end{itemize}
\end{Rmk}

}

\begin{Rmk}[Comparison with previous results]
$\bullet$ One of the main features of our construction lies in the precise description of the connection between boundary layers. In particular, we prove that the sizes and profiles of the North/South and East/West boundary layers are unrelated; the transition between both types of boundary layers occurs through their amplitude only. Notice that this result is rather unexpected: indeed, most works on boundary layers (see for instance \cite{GVP}) assume that an asymptotic expansion of the form
$$
\psi\sim \psi^{int}+\psbl(s,\lambda(s) z)
$$
holds, where the size of the boundary layer, $\lambda^{-1}$, is defined on the whole boundary and is continuous. Here, we exhibit a different type of asymptotic expansion, which shows that a superposition of two types of boundary layers occurs on the transition zone. Also, on this transition zone, the ratio between the sizes of the two boundary layers is very large, so that the asymptotic expansion above cannot hold.

$\bullet$ The analysis we present in this paper can be extended without difficulty to very general anisotropic degenerate elliptic equations, in particular to the convection-diffusion equation
$$\d_x \psi -\viscosite \Delta \psi = \tau$$
which has been studied for special domains by Eckhaus and de Jager in  \cite{eckhaus-jager}, followed by Grasman in \cite{grasman}, and more recently by Jung and Temam \cite{jung-temam}. Note that in the case of the convection-diffusion equation,   the maximum principle (which does not hold anymore in our case) can be used to prove convergence in $L^\infty$. In the paper \cite{eckhaus-jager}, the authors exhibit parabolic boundary layers on the North and South boundaries, but only treat the case when the domain $\Om$ is a rectangle, which turns out to be easier, as explained in the next chapter.

Note that the tools we develop in the present paper allow to consider more general geometries, for which
\begin{itemize}[label=-]
\item there is a continuous transition between lateral  and horizontal boundaries,
\item there are singular interfaces $\Sigma_{ij}$, and even islands.
\end{itemize}

\end{Rmk}

{\subsection{Periodic and rectangle cases}
 \label{ssec:rectangle-case}

 Our arguments also allow to consider  the cases when the domain $\Omega$ is a rectangle $(x_-, x_+)\times(y_-, y_+)$ or an  $x$-periodic domain  $\T \times (y_-, y_+)$, which are in fact much less involved than the case of smooth domains.

In the case of a rectangle, the approximate solution is given by
\be\label{psiapp-rect}
\psiapp=\psi^0_t + \psi_{E,W} + \psi_{N,S},
\ee
exactly as in the case of a smooth domain. It turns out that the interaction between east/west boundary layers and north/south boundary layers in the corners is rather simple.

 In the periodic case, the definition of approximate solution has to be changed a little bit because the circumpolar current $\psi^{circ}$ (see definition below) is generically very large in all norms, so that comparing the size of the error with $\psi^{circ}$ does not give any precise information on the asymptotic expansion of the solution $\psi$. Moreover, the weight $\exp(x)$ (or, more generally, any increasing function of $x$ whose derivatives are bounded from below) cannot be used in the energy estimates, since only periodic weights are allowed. Therefore  the energy estimate \eqref{weight-est} for equation \eqref{Mstab} becomes
\be\label{est-per}
\viscosite \|\Delta \delta \psi\|_{L^2(\T \times (y_-, y_+))}\leq \|\delta \tau\|_{H^{-2}(\T \times (y_-, y_+))}.
\ee
We therefore replace Definition \ref{def:app-general} by the following:
\begin{Def}\label{def:app-per}
A function $\psiapp \in H^2(\T \times (y_-, y_+))$ is an approximate solution to (\ref{M}) if it satisfies the approximate equation (\ref{Mapp}) for some $\dtau \in H^{-2}(\T \times (y_-, y_+))$ with
\begin{equation}
\label{app-def-circ}
\lim_{\viscosite \to 0}  \frac{ {\| \dtau \|_{H^{-2} }}}{  \viscosite \| \Delta (\psiapp-\psi^{circ} ) \|_{L^2} } =0\,.
\end{equation}
\end{Def}

In this case, the approximate solution is defined by
\be\label{psiapp-per}\psiapp =\psi^{circ} + \psi^0_{per} + \psbl_{N,S} \ee
where
\begin{itemize}
\item  $\psi^{circ}=\psi^{circ}(y)$  is the circumpolar current, due to the average forcing by the wind, namely
$$
\begin{aligned}
 -\viscosite \d_y^4 \psi^{circ} = \mean{\tau} (y) , \quad y\in (y_-, y_+),\\
  \psi^{circ}_{|y=y_\pm} = 0,\\
  \d_y \psi^{circ}_{|y=y_c^\pm} = 0,
  \end{aligned}
  $$
where $\mean{f}:=\int_\T f$ for any periodic function $f$.

\item $\psi^0_{per}$ is the classical Sverdrup current in $\Omega_{circ}$, defined by
$$
\d_x \psi^0_{per}(x,y)= \tau(x,y)- \mean {\tau}(y),\quad x\in \T, \  y\in (y_-, y_+),
$$
with $\mean{\psi^0_{per}}(y)=0$.

\item $\psbl_{N,S}$ are periodic North and South boundary layers, whose definition differ slightly from usual North and South boundary layers.

\end{itemize}

We then have the following results:
\begin{Prop}
\begin{enumerate}
\item Assume that $\Om=(x_-, x_+)\times (y_-, y_+)$ and that $\tau \in W^{4,\infty}(\Om)$. Then  $\psiapp$ defined by \eqref{psiapp-rect} satisfies \eqref{Mapp} with an admissible remainder in the sense of \eqref{hyp:dtau}, and is therefore an approximate solution in the sense of Definition \ref{def:app-general}.

\item Assume that $\Om=\T\times (y_-, y_+)$ and that $\tau\in H^s(\Om)$ for $s$ sufficiently large. Then  $\psiapp$ defined by \eqref{psiapp-per} satisfies \eqref{Mapp} with a remainder $\dtau$ satisfying
$$
\|\dtau\|_{H^{-2}}=O(\viscosite),
$$
so that for generic forcing $\tau$,
$$
\|\psi - (\psi^{circ} + \psbl_{N,S})\|_{H^2}= o(\| \psbl_{N,S}\|_{H^2}).
$$
In particular, $\psiapp$ is an approximate solution in the sense of Definition \ref{def:app-per}.

\end{enumerate}
\label{cor:per-rect}
\end{Prop}

\begin{Rmk}
Notice that because of the lack of an $L^2$ estimate, the Sverdrup part of the solution in the periodic case is not captured by the energy estimate. However, if $\tau$ is sufficiently smooth, the construction can be iterated and an approximate solution can be built at any order, so that the existence of every term in the expansion can be justified. Indeed, there is no degeneracy in the problem, and therefore no singularity. The construction of North and South boundary layers only costs a finite, quantifiable number of derivatives on $\tau$.
\end{Rmk}

\subsection{Outline of the paper}

Since the proof of Theorem \ref{thm} is very technical, we have chosen to separate as much as possible the construction of $\psiapp$ and the proof of convergence. Therefore the organization of the paper is the following.

 In Chapter \ref{chap:multiscale}, we expose the main lines of the construction of the boundary layer type terms, namely $\psi_{E,W}$, $\psi_{N,S}$ and $\psi^\Sigma$, without going into the technicalities. An important point is that, while the East and West boundary layers $\psi_{E,W}$ are defined by some local operator, $\psi_{N,S}$ and $\psi^\Sigma$ are obtained as the solutions of some parabolic equations, which accounts for the terminology of parabolic boundary layers used by De Ruijter in his book \cite{DR}: ``when the southern boundary of the continent A coincides with a characteristic the first approximation in the free shear layer is not only dependent on the matching conditions but also on the initial condition at the rim of the parabolic boundary layer. In this way the information about the processes along the southern coast of A is reflected in the interior of the basin. "

 In Chapter \ref{chap:construction}, we give all the necessary details for the construction insisting on the connection between the different types of boundary layers, and we estimate the sizes of the four terms defining $\psiapp$. Even though this part can seem essentially technical, there are two important features of the construction to be noted. The first one is that the connection between boundary layers has to be understood as a superposition: amplitudes are matched, but not profiles. The second ingredient is that the order for the construction of the different correctors is prescribed and this has something to do with the disymmetry between East and West: the construction is essentially Westwards, as are the transport by the Sverdrup equation and the diffusion in the parabolic layers. We also explain how rectangular or periodic domains can be handled.

 Eventually, in Chapter \ref{chap:convergence}, we prove the estimate on the error term \eqref{hyp:dtau}, which entails that $\psiapp$ is an approximate solution.
}
Moreover, as a very large number of notations are introduced throughout the paper, an index of notations is available after the Appendix. We also included a table summarizing the sizes of the different parameters and terms.

%% file: Multiscale.tex
\chapter{Multiscale analysis}

\label{chap:multiscale}

Searching as usual an approximate solution to the Munk equation \eqref{M}
$$\d_x \psi -\viscosite \Delta^2 \psi =\tau$$ in the form
$$
\begin{aligned}
\psiapp = \psi^0 +\psi^{BL}\quad \text{ if } \Om \text{ is simply connected (no islands)},\\
\text{and } \psiapp = \psi^0 + \sum_{i=2}^K \bar c_i \mathbf1_{B_i}+\psi^{BL}\quad \text{ if } \Om \text{ has islands,}
\end{aligned}
$$
with $\psi^0$ satisfying the Sverdrup relation
$$\d_x \psi^0 =\tau \text{ in }\Om ,\quad \psi^0_{|\Gamma_E}=0,$$
we see that both $\psi^0$ and $\psi^{BL}$ present  singularities near the North and South boundaries of the domain, i.e. as $\cos \theta$ vanishes, as well as on the interfaces $\Sigma_{ij}$.

\begin{itemize}
 \item First,  the main term of the approximate solution, $\psi^0$, is singular near all smooth ``North-East'' and ``South-East corners''. More precisely, combining the integral definition of $\psi^0$ together with equivalents for the coordinates of boundary points  given in Appendix B, we see that  the $y$ derivatives of $\psi^0$ explode near such corners.

\item Moreover, as we explained in paragraph \ref{formal-asymptotic}, the size of the boundary layer becomes much larger as $\cos \theta\to 0$, going from $\viscosite^{1/3}$ to $\viscosite^{1/4}$. This also creates strong singularities in the boundary layer terms, which are completely independent from the singularity described above. In fact, there is a small zone in which both boundary layers coexist and are related to one another through their amplitude.

\item Finally, for complex domains, i.e. for domains where the closure of the East boundary $\bar \Gamma_E$  is not connected, the solution of the interior problem $\psi^0$ has discontinuities across the horizontal lines $\Sigma_{ij}$.  These discontinuities gives rise to a \textbf{``boundary layer singularity''}, which  is apparented to the North and South boundary terms, the  size of which is therefore  $\viscosite^{1/4}$.

%

\end{itemize}

\begin{figure}[h]
 \includegraphics[height=8cm]{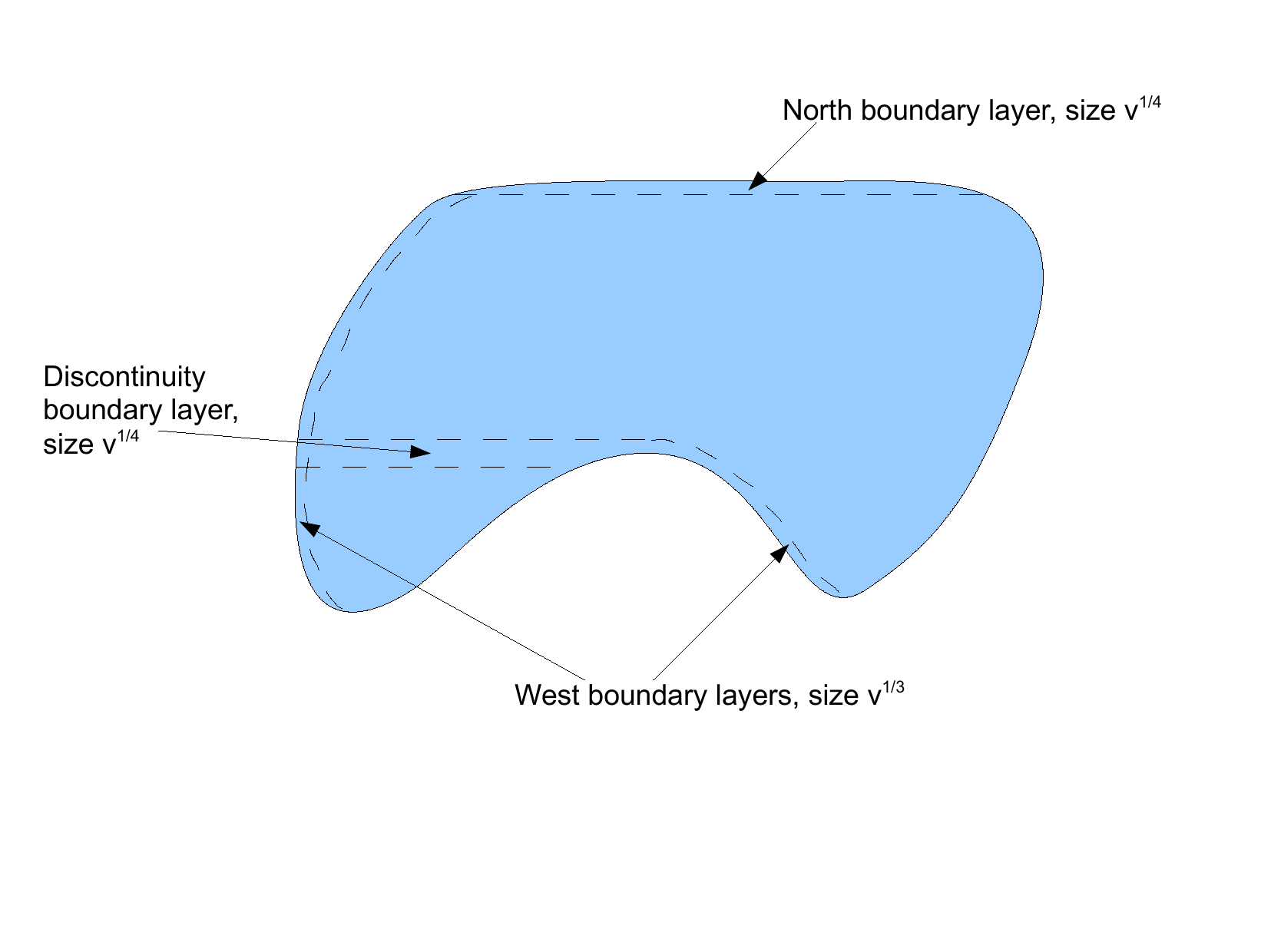}
\caption{Overall view of the boundary layers}
\label{fig:recap}
\end{figure}

Because of these three types of singularities, the construction of the solution is quite technical.
We therefore start with a brief description of all kinds of boundary type terms (see Figure \ref{fig:recap}).

\section{ Local coordinates and  the boundary layer equation}

As usual in linear singular perturbation problems, we build boundary layer correctors as solutions to the homogeneous linear equation
\begin{equation}
\label{HM}
 \d_x \psi -\viscosite \Delta^2 \psi =0\,,
 \end{equation}
localized in the vicinity of the boundary (therefore depending in a singular way of the distance $z$ to the boundary)\label{capZ}
$$\begin{aligned}
\psi \equiv \psi^{BL} (s, \lambda z) \hbox{ with } \d_z \psi^{BL} =\lambda  \d_Z \psi^{BL} \gg \d_s \psi^{BL}\\
 \psi^{BL} \to 0 \hbox{ as } Z\to \infty,
 \end{aligned}$$
 and lifting boundary conditions
 $$\psi ^{BL} (s,0)=-\psi^{int} (s,0)  , \quad \d_Z \psi^{BL} (s,0)=- \d_n \psi^{int} (s,0).$$
 Note that the parameter $\lambda\gg1$ is expected to measure the inverse size of the boundary layer, it can therefore depend on $s$.

It is then natural to rewrite the homogeneous Munk equation in terms of the local coordinates $(z,s)$.
We have by definition of $\theta$, $s$ and $z$,
\be\label{jacobien}\begin{pmatrix}\ds { \d x \over \d z} & \ds{\d x\over \d s}\\ & \\
\ds{ \d y \over \d z} & \ds{\d y\over \d s}\end{pmatrix} = \begin{pmatrix} -\cos \theta & (1+z\theta') \sin \theta\\
- \sin \theta & -(1+z\theta') \cos \theta\end{pmatrix}\ee
so that
$$\begin{pmatrix}\ds{ \d z \over \d x} & \ds{\d z\over \d y}\\&\\
\ds{ \d s \over \d x} & \ds{\d s\over \d y}\end{pmatrix} = \begin{pmatrix} -\cos \theta & - \sin \theta\\
(1+ z \theta')^{-1} \sin \theta & -(1+z\theta')^{-1} \cos \theta\end{pmatrix}.$$
We therefore deduce that the jacobian of the change of variables $(x,y)\to (z,s)$ is equal to $(1+z\theta')^{-1}$. Furthermore,
$$\d_x \psi^{BL} =- \cos \theta \lambda  \d_Z \psi^{BL} +{\sin\theta \over 1+ \lambda^{-1} Z \theta'}\left[ \d_s \psi^{BL}+ {\lambda' \over \lambda} Z\d_Z \psi^{BL} \right].$$
On most part of the boundary the second term is negligible compared to the first one. Nevertheless we see immediately that on horizontal parts the first term is zero, so that we have to keep the second one. Since $\lambda\gg 1$, we approximate the jacobian term $(1 + \lambda^{-1} Z \theta')^{-1}$ by $1$.

Similar computations allow to express the bilaplacian in terms of the local coordinates (which involves more or less twenty terms). However in the boundary layers, we expect that the leading order term is the fourth derivative with respect to $z$
$$\Delta^2 \psi^{BL} =\lambda^4 \d_Z^4 \psi^{BL} + O(\lambda^3).$$
We will thus consider only this term and check \textit{a posteriori} in Chapter \ref{chap:convergence} that the contribution of other terms is indeed negligible.
Note  that, as we want to estimate the $H^{-2}$ norm of the remainder, we will only need to express the laplacian (rather than the bilaplacian) in local coordinates.

We will therefore define boundary layer correctors as (approximate) solutions to the equation
\begin{equation}
\label{BL}
-\lambda \cos \theta \d_Z \psi^{BL} +\sin \theta \Big(\d_s \psi^{BL} +{\lambda' \over \lambda} Z \d_Z \psi^{BL}\Big) -\viscosite \lambda^4 \d_Z^4 \psi^{BL} =0.
\end{equation}
This equation remaining still complicated, we will actually consider two regimes depending on the precise  localization on the boundary.
For each one of these regimes, we will neglect one of the first two terms in (\ref{BL}) (i.e. part of the $\d_x$ derivative), so that the size of $\lambda$ as well as the profile of $\psi^{BL}$ will be different. Of course we will need to check \textit{a posteriori} that the term which has been neglected can be dealt with as a remainder in the energy estimate.

\section{East and West boundary layers}
\label{sec:EW}
In this section, we construct the boundary layers on the lateral sides of the domain. We retrieve rigorously the result announced in the introduction, namely the intensification of Western boundary currents and the dissymetry between the East and West coasts.

\subsection{The scaled equation}$ $

Along the East and West coasts of the domain, on intervals on which $\cos \theta$ remains bounded away from zero, it can be expected that the main terms in the boundary layer equation \eqref{BL} are $\lambda \cos \theta\d_Zf(s,\lambda z)$ and $\viscosite \lambda^4\d_Z^4f(s,\lambda z)$. Hence we take $\lambda$ such that
\be\label{def:lambda_EW}
\lambda_{E,W}= \left(\frac{|\cos \theta|}{\viscosite}\right)^{1/3},
\ee
and $f$ such that
\be\label{eq:BL-EW}
\d_Z^4 f = -\sgn(\cos(\theta))\d_Z f.
\ee
We recall that $Z$ is the rescaled boundary layer variable ($Z=\lambda z$), so that $Z\in (0, \infty)$.

We look for solutions of the above equation which decay as $Z\to \infty$. Consequently, the dimension of the vector space of solutions of the simplified boundary layer equation \eqref{eq:BL-EW} depends on the sign of $\cos(\theta)$:
\begin{itemize}
        \item If $\cos \theta>0$ (East coast),  decaying solutions of equation \eqref{eq:BL-EW} are of the form
$$
f(s,Z)= A(s) \exp(-Z)
$$

        \item If $\cos \theta<0$ (West coast), decaying solutions of equation \eqref{eq:BL-EW} are of the form
$$
f(s,Z)=A_+(s) \exp(-e^{i\pi/3} Z) + A_- (s) \exp(-e^{-i\pi/3} Z).
$$

\end{itemize}
Notice that we retrieve the dissymetry between the East and West coasts: indeed, only one boundary condition can be lifted on the East boundary, whereas two boundary  conditions (namely, the traces of $\psi^0$ and $\d_n\psi^0$) can be lifted on the West boundary. As a consequence, $\psi^0$ must vanish at first order on the East coast, so that the role of $\psi^{BL}$ on $\Gamma_E$ is merely to correct the trace of $\d_n\psi^0$.

Note that, in order that the trace  and the normal derivative of $\psi^0+\psi^{BL}$ are exactly zero on $\Gamma_E$, we will actually need an additional corrector, which is built in the next chapter.

In first approximation, the boundary layer terms on the East and West coasts, denoted respectively by $\psi_E$ and $\psi_W$, are thus defined by\label{psiEW}
$$
\begin{aligned}
        \psi_E(s,Z)=A(s) \exp(-Z),\\
        \psi_W(s,Z)=\sum_\pm A_\pm(s) \exp(-e^{\pm i \pi/3} Z),
\end{aligned}
$$
where  the coefficients $A, A_+, A_-$ ensure that the trace and the normal derivative of $\psi^0 (s,z)+ \psbl (s,\lambda_{E,W}z) $  vanish at main order on the East and West coasts. This leads to
\begin{eqnarray}
 A(s)&=&\lambda_{E}^{-1} \d_n \psi^0_{|\d \Om}(s),\label{def:AE}\\
\begin{pmatrix}
        A_+(s)\\A_-(s)
        \end{pmatrix} &= &
        \begin{pmatrix}
                1&1\\\lambda_{W} e^{i\pi/3}&\lambda_{W}e^{-i\pi/3}
        \end{pmatrix}^{-1}
        \begin{pmatrix}
                - \psi^0_{|\d \Om}(s)\\ \d_n\psi_{|\d \Om}^0(s)
        \end{pmatrix}\nonumber\\&=&{1\over \sqrt{3}}
        \begin{pmatrix}
                 e^{i\pi/6} &-i\lambda_W^{-1}\\ e^{-i\pi/6} &i \lambda_{W}^{-1}
        \end{pmatrix}
        \begin{pmatrix}
                - \psi^0_{|\d \Om}(s)\\ \d_n\psi_{|\d \Om}^0(s)
        \end{pmatrix}\label{def:AW}
\end{eqnarray}
Let us emphasize that the precise value of $A$ on East boundaries is in fact irrelevant in energy estimates, since equation \eqref{def:AE} implies that $A=O(\viscosite^{1/3})$ on zones where $\cos \theta$ does not vanish. Therefore the East boundary layer itself is not captured by energy estimates. But its incidence on the interior term, through the fact that $\psi^0_{|\Gamma_E}=0$, is clearly seen in the $L^2$ estimate.

\subsection{ Domain of validity} $ $

Simplifying \eqref{BL} into \eqref{eq:BL-EW}, we have neglected the terms $\lambda'\lambda^{-1} Z\d_Z f$  and $\d_s f$.
As $f$ is an exponential profile, the term corresponding to $\lambda'\lambda^{-1} Z\d_Z f$ is smaller than $\viscosite \lambda^4 \d_Z^4f$ and {$\lambda\cos(\theta)\d_Z  f$} as long as
$$
\left|\frac{\lambda'}{\lambda}\right|\ll |\viscosite\lambda^4|=\frac{|\cos\theta|^{4/3}}{\viscosite^{1/3}},
$$
which leads to
\be\label{hyp:valid-E/W-2}
        \viscosite^{1/3}|\theta'|\; |\cos \theta|^{-7/3}\ll 1.
\ee

As for the term $\d_s f$, {using (\ref{def:AE}) and (\ref{def:AW}), we obtain that,}  if $\psi^0_{|\d\Om}$ and $\d_n\psi^0_{|\d\Om}$ are smooth with respect to $s$, the corresponding error terms can be neglected as long as
$$
\lambda|\cos(\theta)|\gg1, \text{ i.e. }|\cos(\theta)|\gg \viscosite^{1/4}.
$$
This last condition is  less stringent than \eqref{hyp:valid-E/W-2}. Hence we only keep \eqref{hyp:valid-E/W-2} in order to determine the interval of validity of the construction.

\bigskip
The intervals on which the East and West boundary layers are defined follow from the validity condition \eqref{hyp:valid-E/W-2}. More precisely, if $\cos \theta\neq0$ on $(s_i, s_{i+1})$ (West or East coasts), we set\label{sipm}
$$
\begin{aligned}
        s_i^+:=\sup\left\{s\in\left(s_i, \frac{s_i+s_{i+1}}{2}\right), \viscosite^{1/3}|\theta'|\; |\cos \theta|^{-7/3}\geq 1\right\},\\
        s_{i+1}^-=\inf\left\{s\in\left( \frac{s_i+s_{i+1}}{2}, s_{i+1}\right), \viscosite^{1/3}|\theta'|\; |\cos \theta|^{-7/3} \geq 1\right\}.
\end{aligned}
$$

Notice that $s_i^+, s_{i+1}^-$ depend on $\viscosite$ and are well-defined as long as $\theta\in\mathcal C^1([0,L])$. By definition,
$$
\viscosite^{1/3}|\theta'|\; |\cos \theta|^{-7/3}<1 \quad \forall s\in(s_i^+, s_{i+1}^-),$$
and it is easily proved that $\cos(\theta(s_i^\pm))$ vanish as $\viscosite\to 0$, so that
$$
\lim_{\viscosite\to 0} s_i^+=s_i,\quad\lim_{\viscosite\to 0} s_{i+1}^-=s_{i+1}.
$$

\begin{Rmk}
Notice that \eqref{hyp:valid-E/W-2} is not satisfied in the vicinity of $s_i^+, s_{i+1}^-$. Moreover, the derivatives of $\psi^0_{|\d\Om}$ and $\d_n \psi^0_{|\d\Om}$ with respect to $s$ may become very large as $s$ approaches $s_i$ and $s_{i+1}$. Chapter \ref{chap:convergence} is devoted to the estimation of the corresponding error terms and to the proof of their admissibility.
\end{Rmk}

\section{North and South boundary layers}
\label{sec:NS}

In this section, we construct the boundary layer terms near the intervals where $\cos \theta$ vanishes. Notice that these intervals may in fact be reduced to single points.
{ On the horizontal parts of the boundary, we get a parabolic equation of order 4 with constant coefficients, as obtained  by De Ruijter in \cite{DR}. However to account for the part of the boundary which is ``almost horizontal" we will need to consider a more complicated equation with variable coefficients. }
%
%
%
%

\medskip

\subsection{The scaled equation}$ $
We consider in this section an interval $(s_i, s_{i+1})$
on which $\cos \theta$ is identically zero and $\sin \theta\equiv -1$ (South boundary), so that equation \eqref{BL} becomes
\be\label{eq:North1}
\d_s f + \frac{\lambda'}{\lambda}Z \d_Z f +\viscosite \lambda^4\d_Z^4 f=0.
\ee
As in section \ref{sec:EW}, we have to choose $\lambda$ and $f$ so that the above equation is satisfied. The simplest choice is to take
\be
\lambda_{N,S}:=\viscosite^{-1/4}\label{def:lambdaNS}
\ee
so that $\lambda'=0$, and the boundary layer equation becomes a diffusion-like equation, with the arc-length $s$ playing the role of the time variable.  Note that such degenerate parabolic boundary layers have been exhibited in \cite{eckhaus-jager} for instance.

This raises several questions:
\begin{itemize}
 \item What is the direction of propagation of ``time'' in \eqref{eq:North1}? In other words, what is the initial data for \eqref{eq:North1}?

\item Is equation \eqref{eq:North1} well-posed?

\item What is the domain of definition (in $s$) of the South boundary layer term?

\item How are the South, East and West boundary layers connected?
\end{itemize}

We will prove in this section that equations of the type
$$\begin{aligned}
\d_s f +\d_Z^4f=0, \quad s\in (0,T), \ Z>0,\\
f_{|s=0}=f_{{in}}\\
f_{|Z=0}=f_0,\quad\d_Z f_{|Z=0}=f_1,
  \end{aligned}
$$
with $f_{\text{in}}\in L^2(0, \infty)$, $f_0, f_1\in W^{1,\infty}(0,T)$,  are well-posed, and we will give some energy estimates on the solutions of such equations. In the present context, this means that on South boundary layers, equation \eqref{eq:North1} is a forward equation (in $s$), while on North boundary layers it becomes a backward equation. This is consistent with the definition of the interior term $\psi^0$: in all cases, the boundary condition in $s$ is prescribed on the East end of the interval.

\begin{figure}[h]
 \includegraphics[width=\textwidth]{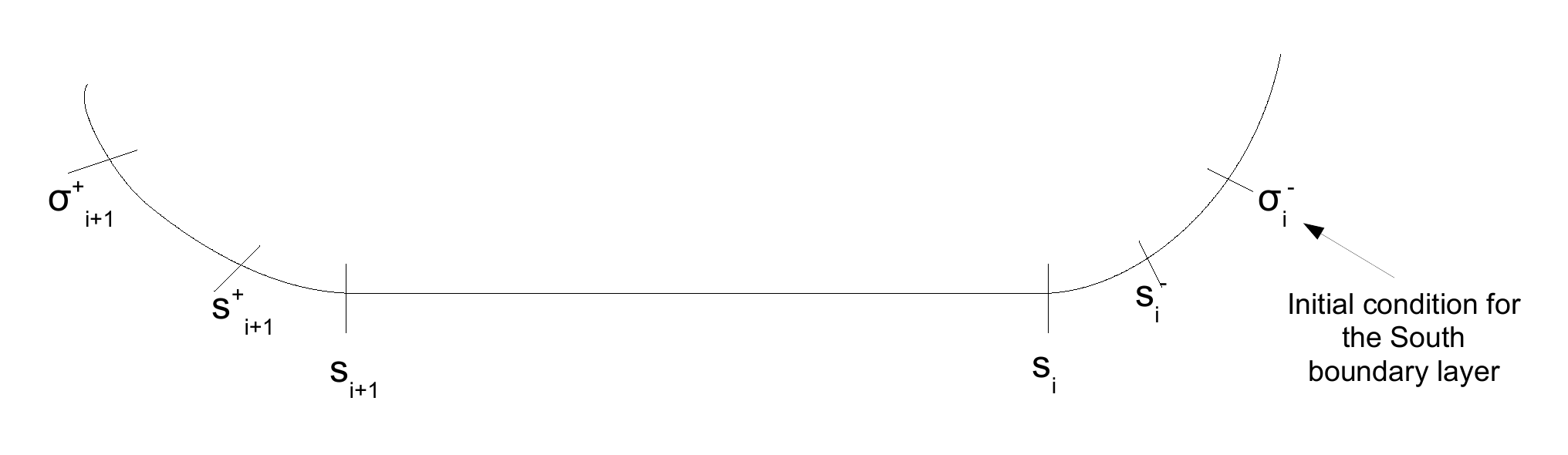}
\caption{The arc-lengths $\sigma_i^-, s_i^-, s_i, s_{i+1}, s_{i+1}^+, \sigma_{i+1}^+$ on a South boundary}
\label{fig:bord-sud}
\end{figure}

Let us also recall that the domain of validity of the West and East boundary layer terms does not reach the zone where $\cos \theta=0$ (see \eqref{hyp:valid-E/W-2}). As a consequence, the South boundary layer term must be defined for $s\in (\sigma_i^-, s_i)$ and $s\in (s_{i+1}, \sigma_{i+1}^+)$ where the arc-lengths $\sigma_i^-, \sigma_{i+1}^+$ will be defined later on and should satisfy (see Figure \ref{fig:bord-sud})
$$\sigma_i^-<s_i^-,\quad \sigma_{i+1}^+>\sigma_{i+1}^+.$$

%
%
This  requires to slightly modify the boundary layer equation. From now on, we define the South boundary layer term $\psi_S$ as the solution of\label{psiNS}
\be
\begin{aligned}
 \d_s \psi_S - \viscosite^{-1/4}\frac{\cos \theta}{\sin\theta}\d_Z \psi_S -\frac{1}{\sin \theta}\d_Z^4 \psi_S=0,\ s>\sigma_{i}^-,\ Z>0,\\
\psi_{S|s=\sigma_{i}^-}=0,\\
\psi_{S|Z=0}=\Psi_0,\quad \d_Z\psi_{S|Z=0}=\Psi_1,
\end{aligned}
\label{eq:North}
\ee
with $\Psi_0$, $\Psi_1$ to be defined later on (see \eqref{def:Psi01}). Notice that we choose $\lambda_{N,S}= \viscosite^{-1/4}$ even in the zones where $\cos \theta \neq 0$, i.e. outside the interval $(s_i, s_{i+1})$.

\subsection{Study of the boundary layer equation \eqref{eq:North}:}$ $

In this paragraph, we give a well-posedness result for equation
\be\label{eq:North2}
\begin{aligned}
\d_s f +\gamma(s)\d_Z f + \mu(s)\d_Z^4f=0,\quad s\in(0,T),\ Z>0, \\
f_{|s=0}=f_{in},\\
f_{|Z=0}=f_0,\quad\d_Zf_{|Z=0}=f_1\,.
\end{aligned}
\ee

We first note that, up to lifting the boundary conditions, we are brought back to the study of the same parabolic equation with homogeneous boundary conditions and with a source term: indeed $f$ is a solution of \eqref{eq:North2} if and only if
$$
g(s,Z):=f(s,Z)-f_0(s)(Z+1)\exp(-Z)-f_1(s) Z \exp(-Z)
$$
is a solution of
$$\begin{aligned}
\d_s g + \gamma(s)\d_Z g + \mu(s)\d_Z^4g= S(s,Z),\quad s\in(0,T),\ Z>0, \\
g_{|s=0}=g_{in},\\
g_{|Z=0}=0,\quad\d_Zg_{|Z=0}=0,
\end{aligned}
$$
where
$$
g_{in}:=f_{in}-f_0(0)(Z+1)\exp(-Z)-f_1(0) Z \exp(-Z)
$$
and
\begin{eqnarray*}
 S(s,Z)&:= &-f_0'(s)(Z+1)\exp(-Z)- f_1'(s)Z\exp(-Z)\\
&&+ \gamma(s) f_0(s)Z\exp(-Z)-\gamma(s)f_1(s)(1-Z)\exp(-Z)\\
&&-\mu(s) f_0(s)(Z-3)\exp(-Z) - \mu(s) f_1(s)(Z-4)\exp(-Z).
\end{eqnarray*}

Our precise result is then the following:

\begin{Lem}\label{Cauchy-NS}
Let $T>0$, and let $\mu\in L^{\infty}(0,T)$, $\gamma\in L^{\infty}(0,T)$, $f_0, f_1\in \mathcal C^1([0,T])$, $g_{in}\in L^2(0,\infty)$, and $S\in L^1((0,T), L^2(\R_+))$.

Assume that there exists $\mu_0>0$ such that
$$
\mu(s)\geq \mu_0\quad \text{for a.e. }s\in(0,T).
$$

 Then the   equation
\be\label{eq:Northg}
\begin{aligned}
\d_s g +\gamma(s)\d_Z g +\mu(s)\d_Z^4g= S(s,Z),\quad s\in(0,T),\ Z>0, \\
g_{|s=0}=g_{in},\\
g_{|Z=0}=0,\quad\d_Zg_{|Z=0}=0,
\end{aligned}
\ee
has a unique solution $g\in \mathcal C([0,T], L^2(\R_+))\cap L^2((0,T), H^2(\R_+))$ which satisfies the energy estimate
\begin{eqnarray*}
  &&\frac12\| g(s)\|_{L^2(\R^+)}^2 + \int_0^s \mu(s')\| \d_Z^2 g(s')\|_{L^2(\R^+)} ^2ds' \\&\leq& \| g_{in}\|_{L^2(\R^+)}^2 +\left(\int _0^s \| S(s')\|_{L^2(\R^+)}ds'\right)^2\,.
\end{eqnarray*}

\end{Lem}
\begin{proof}
Let
$$
V:=\{v\in H^2(0,\infty),\ v_{|Z=0}=0,\ \d_Z v_{|Z=0}=0\}.
$$
For $u,v\in V$, $s\in(0,T)$, define the quadratic form
$$
a(s,u,v)=\int_0^\infty \gamma(s) \d_Z u(Z) \; v(Z)\:dZ + \mu(s)\int_0^\infty\d_Z^2 u(Z)\: \d_Z^2 v(Z)\:dZ.
$$
Then for all $u,v\in V$, the map
$$
s\mapsto a(s,u,v)
$$
is measurable, and for almost every $s\in (0,T)$
\begin{eqnarray*}
|a(s,u,v)|&\leq &\Big(\|\gamma\|_{L^\infty(0,T)}+ \|\mu\|_{L^\infty(0,T)}\Big)\|u\|_{H^2(\R_+)}\|v\|_{H^2(\R_+)}\\
a(s,u,u)&\geq& \mu_0 \|\d_Z^2 u\|_{L^2(\R_+)}^2  \\
&\geq &  \mu_0( \|u\|_{H^2(\R_+)}^2 -  \|\d_Zu\|_{L^2(\R_+)}^2 - \|u\|_{L^2(\R_+)}^2 ) .
\end{eqnarray*}
Notice that if $u\in V$,
\begin{eqnarray*}
\| \d_Z u\|_{L^2}^2 &=& \int_0^\infty (\d_Z u)^2 = -\int_0^\infty u \d_Z^2 u\\
&\leq & \frac{\|u\|_{L^2}^2}{2}+  \frac{\|\d_Z^2u\|_{L^2}^2}{2}.
\end{eqnarray*}
We infer eventually that
$$
a(s,u,u)\geq \frac{\mu_0}{2}\|u\|_{H^2(\R_+)}^2 -\frac{3\mu_0}{2}\|u\|_{L^2(\R_+)}^2.
$$
Using Theorem 10.9 by J.L. Lions in \cite{Brezis}, we infer that there exists a unique solution $g\in L^2((0,T), V)\cap \mathcal C([0,T], L^2(\R_+))$ of equation \eqref{eq:Northg} such that $\d_s g \in L^2((0,T), H^{-2}(\R_+))$.

\bigskip
Multiplying equation \eqref{eq:Northg} by $g$ and integrating on $\R_+$, we infer
$$
\frac{1}{2}\frac{d}{ds}\int_0^\infty |g(s)|^2+\mu(s) \int_0^\infty |\d_Z^2 g(s)|^2=\int_0^\infty  S(s,Z)g(s,Z)\:dZ.
$$
Therefore
\begin{eqnarray*}
&&\frac{d}{ds}\left[ \frac 12 \| g(s)\|_{L^2}^2 + \int_0^s \mu(s')\|\d_Z^2 g(s')\|_{L^2}^2ds'\right]\\
&\leq & \| S(s)\|_{L^2} \| g(s)\|_{L^2}\\
&\leq&\sqrt{2} \| S(s)\|_{L^2}\left[ \frac 12 \| g(s)\|_{L^2}^2 + \int_0^s \mu(s')\|\d_Z^2 g(s')\|_{L^2}^2ds'\right]^{1/2}.
\end{eqnarray*}
Integrating with respect to $s$ leads then to the desired inequality.
\end{proof}

We deduce easily that \eqref{eq:North} has a unique solution on any interval of the form $(\sigma^-_i,\sigma_{i+1}^+)$ such that $\sin \theta$ does not vanish on  $[\sigma^-_i,\sigma_{i+1}^+]$, and that it satisfies
\begin{eqnarray*}
\|\psi_S(s)\|_{L^2(\R_+)}&\leq&C (|\Psi_0(s)| + |\Psi_1(s)|+ |\Psi_0(\sigma_i^-)| + |\Psi_1(\sigma_i^-)|) \\
&&+C\int_{\sigma_i^-}^s\left(|\Psi_0'(s')| + |\Psi_1'(s')| \right)\:ds'\\
&&+C\int_{\sigma_i^-}^s(\mu(s')+|\gamma(s')| ) (|\Psi_0(s')| + |\Psi_1(s')|)\:ds',
\end{eqnarray*}
with $\mu(s)=-1/\sin \theta(s)$, $\gamma(s)=-\viscosite^{-1/4} \cos \theta(s)/\sin \theta(s)$.

\subsection{Boundary conditions for $s\in (s_i, s_{i+1})$}
 In order to satisfy  the boundary conditions
$$
\psi_{|\d \Om}=0,\quad \d_n\psi_{|\d \Om}=0,
$$
the South boundary layer term constructed above must be such that for $s\in (s_i, s_{i+1})$,
$$
\begin{aligned}
 \psi_S(s,Z=0)=-\psi^0_{|\d\Om}(s),\\
\viscosite^{-1/4}\d_Z  \psi_S(s,Z=0) = -\d_n \psi^0_{|\d\Om}(s).
\end{aligned}
$$
In equation \eqref{eq:North}, we therefore take, for $s\in (s_i, s_{i+1})$,
\begin{align}\label{def:Psi01}
\Psi_0(s):= -\psi^0_{|\d\Om}(s),\\
\Psi_1(s):=\viscosite^{1/4} \d_y  \psi^0_{|\d\Om}(s).
\end{align}
There remains to define $\Psi_0$ and $\Psi_1$ for $s\in (\sigma_i^-, s_i]\cup [s_{i+1}, \sigma_{i+1}^+)$.

\subsection{Connection with East and West boundary layers}$ $

$\bullet$ In the most simple cases, the North or South boundaries are connected on the one hand to the East boundary, and on the other hand to the West boundary. This corresponds to the situation when

\noindent
-  the North boundary $ y=y_i$ is a local maximum of the ordinate;

\noindent
-  the South boundary $ y=y_i$ is a local minimum of the ordinate.

Without loss of generality, we assume that the corresponding piece of $\d \Om$ is a ``South'' boundary, meaning that $\cos\theta<0$ in a neighbourhood on the right of $s_{i+1}$,  $\cos\theta>0$ in a neighbourhood on the right of $s_{i}$, and $\sin \theta=-1$ for $s\in [s_i, s_{i+1}]$.

The connection with the East boundary is fairly simple: we will introduce some  truncation $\chi_\viscosite$ of $\tau$ close to the East corner (parametrized by $s_{i}$). Therefore the solution to the transport equation $\d_x\psi^0_t=\tau \chi_\viscosite$ is identically zero in a vicinity of $(x_{i}, y_{i})$, as well as the East boundary layer (which lifts the trace of $\d_n \psi^0_t$), so that we merely require $\psi_{S|s=\sigma_{i}^-}=0$.

Concerning the connection with the West boundary layer, the situation is not as straightforward. In order to satisfy  the boundary conditions
$$
\psi_{|\d \Om}=0,\quad \d_n\psi_{|\d \Om}=0,
$$
the West and South boundary layer terms must be such that
$$
\begin{aligned}
\psi_W(s,Z=0) + \psi_S(s,Z=0)=-\psi^0_{t|\d\Om}(s),\\
\lambda_{W} \d_Z  \psi_W(s,Z=0)  + \lambda_{S} \d_Z  \psi_S(s,Z=0) = -\d_n \psi^0_{t|\d\Om}(s)
\end{aligned}
$$
in a vicinity of $s=s_{i+1}$. As a consequence, we take
$$
\begin{aligned}
 \psi_W(s,Z=0)= - \varphi_{i+1}(s) \psi^0_{t|\d\Om}(s),\\
 \d_Z  \psi_W(s,Z=0) =-\lambda_{W}^{-1}\varphi_{i+1}(s) \d_n\psi^0_{t|\d\Om}(s),
\end{aligned}
$$
and
$$
\begin{aligned}
 \psi_S(s,Z=0)= -(1-\varphi_{i+1}(s)) \psi^0_{t|\d\Om}(s),\\
 \d_Z  \psi_S(s,Z=0) =-\lambda_{S}^{-1}(1-\varphi_{i+1}(s) )\d_n\psi^0_{t|\d\Om}(s),
\end{aligned}
$$
where $\varphi_{i+1}$ is a truncation function which we will define precisely in the next chapter. We emphasize that with this definition, the South and West boundary layers are related via their amplitude only: in particular, the sizes of the boundary layers are not related.

\bigskip
{$\bullet$ The case when the boundary has non smooth corners between the horizontal part and the meridional boundaries (see paragraph \ref{corner-par})
$$ \cos \theta(s) \equiv 0 \hbox{ on } \Gamma_N \cup \Gamma_S, \qquad \inf_{\Gamma_E \cup \Gamma_W} |\cos \theta| >0   ,$$
typically the case of rectangles studied by De Ruijter \cite{DR},  could seem more singular at first sight but is actually easier to deal with.

Near East corners, there is no need to introduce the truncation $\chi_\viscosite$ since the solution to the Sverdrup equation is smooth.

Now, on the West boundary, if we lift the function $\psi^0+ \psi_E +\psi_{N,S}$, we obtain local boundary terms which vanish identically when $s\to s_i$. Hence there is  no error due to the trace of West boundary terms on the horizontal part of the boundary. Notice that:
\begin{itemize}
\item If the angle between a horizontal and a western boundary is exactly $\pi/2$, our construction works without any adaptation;
\item If the angle is obtuse, the North/South boundary layer needs to be extended beyond the corner, in the spirit of the connection between the discontinuity boundary layers and western boundary layers (see Remark \ref{rmk:sigma-west}).

\item If the angle is acute, the correct rescaled boundary layer variable on the western boundary is $(x-x_W(y))/\viscosite^{1/3}$ instead of $z\viscosite^{-1/3}$, so as not to pollute the trace on the North/South boundary. Notice that this is the boundary layer variable used by Desjardins and Grenier in \cite{DG}.

\end{itemize}
We leave the technical details to the reader, and treat in complete detail  the case of rectangles in the present article.

}

\bigskip
$\bullet$ In more complex cases, the line $y=y_i$ intersects the interior of $\Omega$, and we get a singularity which is not only localized in the vicinity of $\d \Omega$. Techniques of boundary layers allow however to understand the qualitative behaviour of the solution in the vicinity of the line of singularity. This heuristic approach is presented in the next paragraph.

{\section{Discontinuity zones}

\label{sec:disc}

Discontinuity zones occur when the East boundary $\bar \Gamma_E$ is not connected.}
At leading order, we expect the solution to (\ref{M}) to be approximated by the solution $\psi_t^0$ to the transport equation with suitable truncations $\chi_\viscosite$ near East corners
$$
\begin{aligned}
\d_x \psi^0_t=\tau \chi_\viscosite,\\
\psi^0_{t|\Gamma_E}=0.
\end{aligned}
$$
If the domain $\Om$ has islands, we add to $\psi^0_t$ the quantity $\sum_{j=2}^K \bar c_j \mathbf 1_{B_j}$. We will give more details on this case in the next section.
Note that, because of the truncation $\chi_\viscosite$, the main order approximation $\psi^ 0_t$ now depends (weakly) on $\viscosite$.

Therefore on every set $A_i$, $\psi^0_t$ takes the form
$$
\psi^0_t(x,y)=-\int_x^{x_E^i(y)}\tau \chi_\viscosite(x',y)\:dx'.
$$
Unlike the main order term in Theorem \ref{thm:DG},
\textit{$\psi^0_t$ does not belong to $H^2(\Om)$ in general}. It is indeed obvious that $\psi^0_t$ and $\d_y \psi^0_t$ may be discontinuous across every (nonempty) line $\Sigma_{ij}$. More precisely, the jump of $\psi^0_t$ across a given line $\Sigma_{ij}$ takes the following form (up to an inversion of the indexes $i$ and $j$)
$$
[\psi^0_t]_{\Sigma_{ij}}=- \int_{x_i}^{x_j}(\tau \chi_\viscosite)(x',y_{ij})\:dx'.
$$
The existence of such discontinuities in the main interior term is a serious impediment to energy methods. Indeed, inequality \eqref{weight-est}, for instance,  requires the approximate solution to be at least in $H^2(\Om)$. In fact, we will prove  that this discontinuity gives rise to a \textbf{``boundary layer singularity''}: the corresponding corrector is a boundary layer term located in the vicinity of $\Sigma_{ij}$. Since the normal vector to $\Sigma_{ij}$ is parallel to $e_y$, this boundary layer term is apparented to the North and South boundary terms constructed in Section \ref{sec:NS}, and therefore the  size of the boundary layer is $\viscosite^{1/4}$.

We now define what we have called in Theorem \ref{thm} a ``generic forcing $\tau$'':
\begin{Def}
We say that the forcing $\tau$ is \textbf{generic} if all the following conditions are satisfied:
\begin{itemize}
 \item For all $i\neq j$ such that  $\Sigma _{ij}\neq \emptyset$, there exists  $c_\Sigma>0$ such that
$$
\hbox{ for all $\viscosite$ sufficiently small, }\left|\left[\psi^0_t + \sum_{l=2}^K \bar c_l \mathbf 1_{B_l}\right]_{\Sigma_{ij}}\right|\geq c_\Sigma;
$$
\item If $\mathring{\Gamma}_{N,S}\neq \emptyset$ (for the induced topology on $\d\Om$), then for any $V\subset \Gamma_{N,S}$, there exists $c_V>0$ such that
$$
\hbox{ for all $\viscosite$ sufficiently small, }\|\psi^0_{t|\d\Om}\|_{L^2(V)}\geq c_V;
$$
\item For any nonempty $V\subset \Gamma_{W}$, there exists $c_V>0$ such that
$$
\hbox{ for all $\viscosite$ sufficiently small, }\left\|\left(\psi^0_t + \sum_{l=2}^K \bar c_l \mathbf 1_{B_l}\right)_{|\d\Om}\right\|_{L^2(V)}\geq c_V.
$$
\end{itemize}

\label{def:tau-generic}
\end{Def}

\smallskip

$\bullet$ From now on, for the sake of simplicity, we restrict the presentation to the case when $M=2$, i.e. \textbf{$\bar \Gamma_E$ has two connected components}, and $\psi^0_t$ has exactly one line of discontinuity. Of course, our construction can be immediately generalized to the case when there are more than two connected components, and therefore several lines of discontinuity: we merely add up the local correctors constructed  in the vicinity of every $\Sigma_{ij}$. In particular, when $\Om$ has islands, there are always at least two discontinuity lines. We will sketch the necessary adaptations  in the next paragraph and leave the details to the reader.

\begin{figure}
\includegraphics[width=\textwidth]{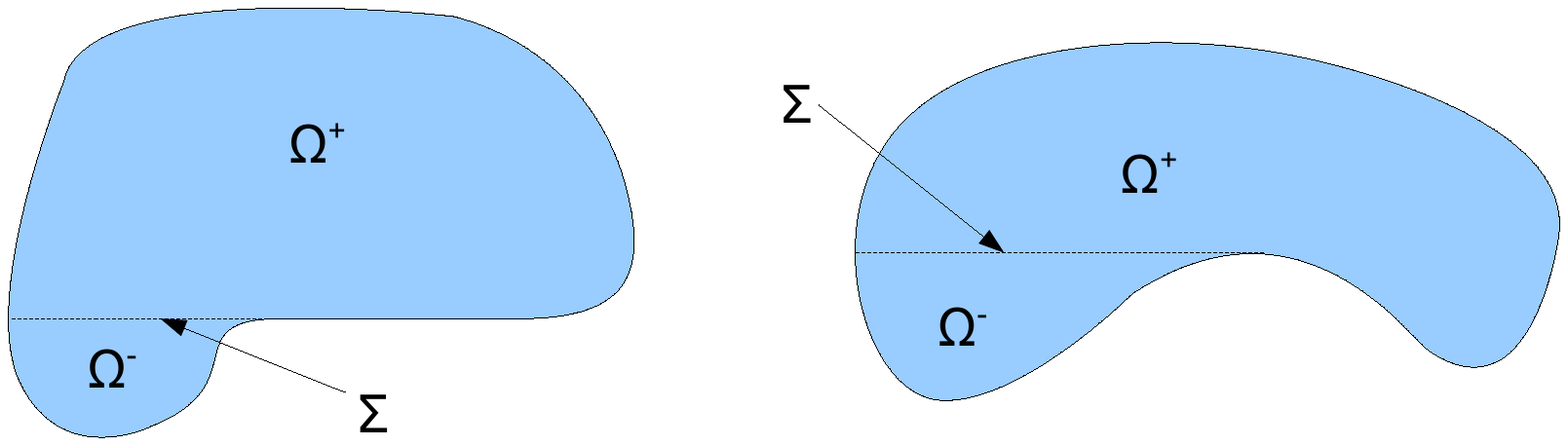}
\caption{Two possible configurations for the domain $\Om$}\label{fig:ex-complexe}
\end{figure}

 Thus we henceforth assume that $\Om$ has the following form (see Figure \ref{fig:ex-complexe}):
$$
\Om=\Om^+ \cup \Om^- \cup \Sigma,
$$\label{Sigma}
where
\begin{itemize}
 \item $\Om^\pm$ are non empty, open and convex in $x$;
\item $\Om^+= \{(x,y)\in \Om,\ y>y_1\text{ or } x\geq x_1\}$;
\item $\Om^-= \{(x,y)\in \Om,\ y<y_1\text{ and } x< x_1\}$;
\item $\Sigma=\{(x,y_1)\in \Om,\  x<x_1\}$
\end{itemize}
and $(x_1, y_1)=(x(s_1), y(s_1))$ for some $s_1$ such that $\cos \theta(s_1)=0$. Without any loss of generality, we further assume that $\sin \theta(s_1)=-1$.
This type of domain can correspond to two different configurations (see Figure \ref{fig:ex-complexe}).

We will use the following notations
$$\chi^- = \indc_{x<x_1} \indc_{y<y_1},\quad \chi^+ =1- \chi^-$$
which somehow stand for $\indc_{\Omega^\pm}$ : their role is to avoid any artificial singularity on $\d \Omega$. We also set 
$\Gamma_E^\pm= \d \Om^\pm\cap \Gamma_E$. We parametrize each set $\Gamma_E^\pm$ by a graph $(x_E^\pm(y), y)$.
\label{x-E-pm}\label{chipm}

As explained above, the function $\psi^0_t$ and its $y$ derivative are discontinuous across $\Sigma$.
 More precisely, for $x<x_1$,
\begin{eqnarray*}
[\psi^0_t]_{|\Sigma}(x,y_1)&=&\psi^0_t(x, y_1^+ ) - \psi^0_t(x, y_1^- )\\
&=&-\int_x^{x_E^+(y_1)}(\tau\chi_\viscosite)(x',y_1)\:dx' + \int_x^{x_E^-(y_1)}(\tau\chi_\viscosite)(x',y_1)\:dx' \\
&=&-\int_{x_E^-(y_1)}^{x_E^+(y_1)}(\tau\chi_\viscosite)(x',y_1)\:dx' .
\end{eqnarray*}
Notice in particular that the jump is constant along $\Sigma$.
In a similar way, since $(\tau \chi_\viscosite)(x_E^-(y), y)$ vanishes in a neighbourhood on the left of $y_1$,
\begin{eqnarray*}
[\d_y\psi^0_t]_{|\Sigma}(x,y_1)&=&- (x_E^+)'(y_1) (\tau \chi_\viscosite)(x_E^+(y_1), y_1)\\
&&-\int_{x_E^-(y_1)}^{x_E^+(y_1)}\d_y(\tau\chi_\viscosite)(x',y_1)\:dx' .
\end{eqnarray*}
\medskip

We now construct the boundary layer type correctors $\psil, \psi^\Sigma$.
The role of $\psil$ is
 to counterbalance the jump of $\psi^0_t$, and of its normal derivative $\d_y \psi^0_t$, across $\Sigma$:
 $$ [ \psil]_{\Sigma} = -[\psi^0_y]_{\Sigma},\quad [\d_y \psil]_{\Sigma}=-[\d_y \psi_t^0]_{\Sigma}\,.$$
When we lift these boundary conditions, we  introduce a source term in the equation. This source term is then handled by a boundary layer type term $\psi^\Sigma$, which has no discontinuity across $\Sigma$: $\psi^\Sigma\in H^2(\Om)$.

 At first sight, if we consider the whole singular  corrector $\psil+ \psi^\Sigma$, this problem could seem underdetermined. Indeed  there are two jump conditions to be satisfied, and possibly two boundary layers (on $\Om^+$ and on $\Om^-$), each of which is the solution of an equation of the type \eqref{eq:North}, and therefore having each two entries ($\Psi_0$ and $\Psi_1$).

  However, when looking at the error terms, we see that
 the traces of $\d_y^2 \psil _{|\Sigma}$ and $\d_y^3 \psil_{|\Sigma}$ appear in the energy estimates (see the proof of Lemma \ref{lem:eq-psil} below, and in particular the derivation of equation \eqref{eq:der23}). Therefore we have further  to request that
$$[\d_y^2 \psil ]_{\Sigma}=[\d_y^3 \psil ]_{\Sigma}=0\,.$$
 Because of this constraint, the energy of the boundary layer discontinuity is distributed on both sides of $\Sigma$. The boundary term $\psil+ \psi^\Sigma$ is unequivocally defined, since roughly speaking, there are four jump conditions and four entries for the boundary layer terms.


\subsection{Lifting the discontinuity}
\begin{figure}
 \includegraphics[width=\textwidth]{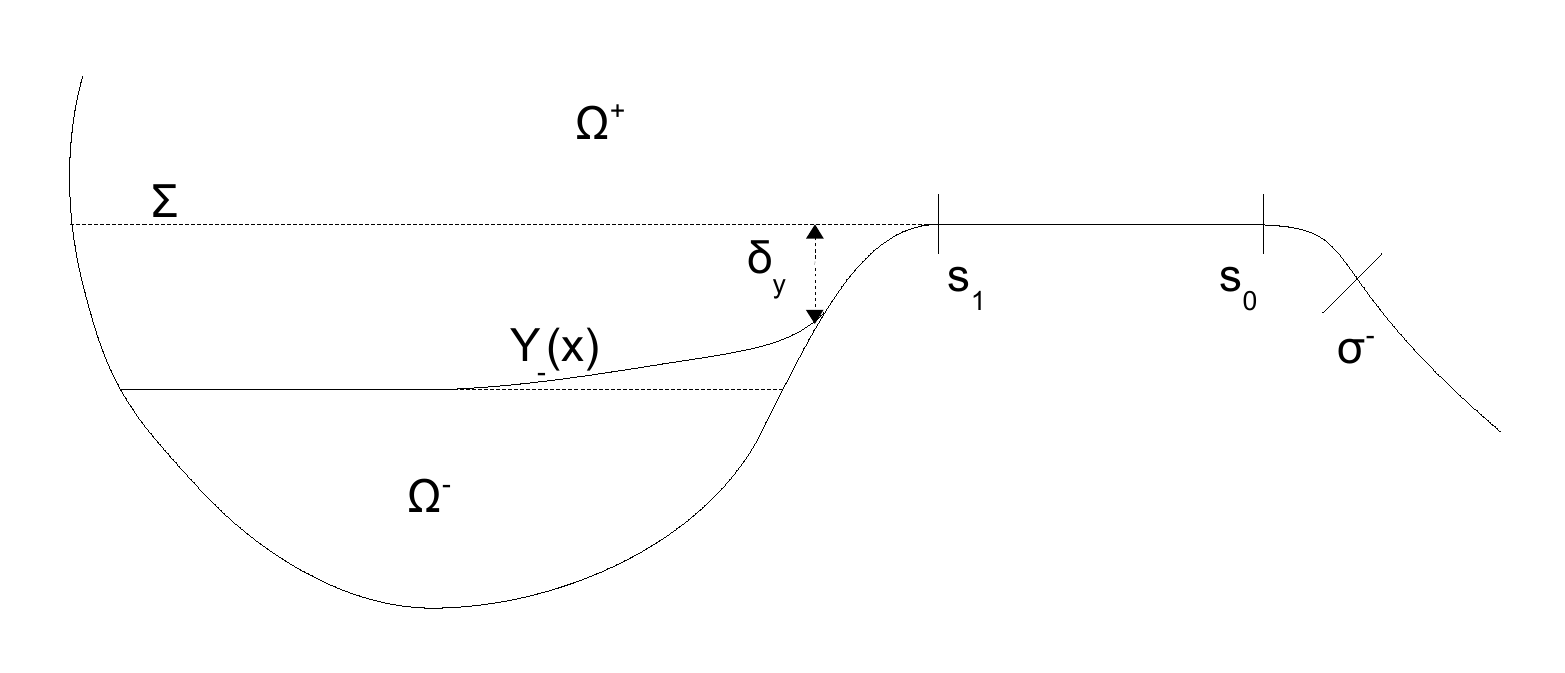}
\caption{The interior singular layer}
\label{fig:sigma-layer-ch2}
\end{figure}

$ $Let us start by introducing some notation. Consider the closed set $\{s\in \d\Om,\ y(s)=y_1\}$. We denote by $I_1$ its connected component containing $s_1$. Then $I_1$ is a closed interval (possibly reduced to a single point). Without loss of generality, we assume that $s_1=\sup I_1$, which means (recall that $\sin \theta(s_1)=-1$) that $s_1$ is the West end of $I_1$. Hence we have either $I_1=\{s_1\}$, or $I_1=[s_0, s_1]$ with $s_0<s_1$.
As in equation \eqref{eq:North}, we introduce a point $\sigma_{in}$ which will be the initial point for $\psil$ and $\psi^{\Sigma}$, and which we will define more precisely in the next chapter (see \eqref{def:sigma-} and Figure \ref{fig:sigma-layer-ch2}).

In order to keep the construction as simple as possible, we use that the cancellation of $\cos \theta$ near $\inf I_1$ is strong (assumption (H4) in section \ref{sec:geom}), which implies in particular  that the boundary of $\Om^+$ near the junction point of arc-length $\inf I_1 $ has $\mathcal C^4$ regularity.


The role of $\psil$ is to lift both the jumps of $\psi^0_t$ and $\d_y\psi^0_t$ across $\Sigma$, and the traces of $\psi^0_t$ and $\d_n \psi^0_t$ on the portion of $\d\Om$ between $\sigma_{in} $ and $s_1$, so that $\psi^0_t +\psil$ belongs to $H^2(\Omega)$.

We have indeed the following

\begin{Prop}\label{H2-prop}
Let $\psi$ be  a function on $\Omega$ such that
$$
\psi_{|\Omega^+} \in H^2(\Omega^+), \quad \psi_{|\Omega^-} \in H^2(\Omega^-),\
\hbox{ and } [\psi]_{\Sigma} =[\d_n \psi]_\Sigma =0\,.
$$
Then $\psi$ belongs to $H^2(\Omega)$.
\end{Prop}

\begin{proof}
The squared $L^2$ norm being additive, the only point to be checked is that
\be\label{dec-DeltaPsi}\Delta \psi = \Delta(\psi\chi^+) + \Delta (\psi\chi^-)\ee
in the sense of distributions, which  is obtained by a simple duality argument.

For any  $\phi\in C^\infty_c(\Omega)$, because of the jump conditions on $\d_y \psi$ and $\psi$, we have
$$
\begin{aligned}
\langle \phi, \Delta \psi\rangle_{\mathcal D, \mathcal D'}&=\int_\Omega \Delta \phi ( \psi \chi^+ +\psi \chi^-)\\
&= -\int_\Sigma \d_y \phi  \psi_{|\Sigma^+}- \int_{\Omega^+}  \nabla \phi \cdot \nabla (\psi \chi^+) \\&+\int _\Sigma \d_y \phi \psi_{|\Sigma^-} - \int_{\Omega^-}  \nabla \phi \cdot \nabla (\psi \chi^-) \\
&= \int_{\Omega^+}   \phi\Delta  (\psi \chi^+) + \int_\Sigma \phi \d_y \psi_{|\Sigma^+}+\int_{\Omega^-}   \phi\Delta  (\psi \chi^-) - \int_\Sigma \phi \d_y \psi_{|\Sigma^-}\\
&=\int _\Omega \phi (\Delta  (\psi \chi^+)+\Delta  (\psi \chi^-))
\end{aligned}
$$
from which we deduce that $\Delta \psi \in L^2(\Omega)$, and $\Delta \psi= \Delta(\psi\chi^+) + \Delta (\psi\chi^-)$.
\end{proof}

\bigskip

We therefore seek $\psil$ in the form
\be\label{def:psil}
\psil(x,y)=\chi^+ \chi\left(\frac{y-y_1}{\viscosite^{1/4}}\right)\left[a(x) + b(x) (y-y_1)\right],
\ee
where the truncation $\chi\in \mathcal C^\infty_0(\R)$ is such that
\be\label{chi}
\supp\chi\subset[-1,1],\quad
\chi(y)=1\text{ for }|y|\leq 1/2,
\ee
and satisfying the conditions
$$\begin{aligned}
&(i) \quad \psil_{|\d\Om^+\cap \d\Om}=-(1-\varphi)\psi^0_{t|\d\Om},\text{ for }s\in (\sigma_{in}, s_1),\\
& \qquad \d_n \psil_{|\d\Om^+\cap \d\Om}=-(1-\varphi)\d_n \psi^0_{t|\d\Om}\text{ for }s\in (\sigma_{in}, s_1),\\
&(ii) \quad \psil_{|\Sigma}=   -[\psi^0_t]_\Sigma , \ \d_y \psil_{|\Sigma}=-[\d_y\psi^0_t]_\Sigma.
  \end{aligned}
$$
As in section \ref{sec:NS}, $\varphi$ is a partition function in $s$, whose role is to ensure a smooth transition between the East/West  boundary layer on $[\sigma_{in},\inf I_1]$ and the discontinuity boundary layer.

A straightforward computation then provides the following formulas for the coefficients $a$ and $b$
\begin{equation}
\label{def:ab1}
\begin{aligned}
b(x)=&(1-\varphi(s))\left[-\sin \theta(s) \d_n \psi^0_{t|\d\Om}(s) +\cos \theta(s)\d_s\psi^0_{t|\d\Om}(s)\right]\\& - \cos \theta(s) \varphi' \psi^0_{t|\d\Om}(s),\\
 a(x)=&-(1-\varphi(s) )\psi^0_{t|\d\Om}(s) - b(x(s))(y(s)-y_1),
\end{aligned}
\end{equation}
for $x_1<x=x(s)<x(\sigma_{in})$,
and by
\begin{equation}
\label{def:ab2}
\begin{array}{l}
a(x)= -[\psi^0_t]_\Sigma,\\
b(x)=- [\d_y\psi^0_t]_\Sigma
\end{array}
\quad\text{for }x\leq x_1.
\end{equation}
For $x\geq x(\sigma_{in})$, we merely take $a(x)=b(x)=0$.

\bigskip
Of course, however, $\psil$ is not a solution of \eqref{M}, even in the approximate sense. More precisely, we have the following proposition :
\begin{Prop}
With the previous definition and notations,
$$(\d_x- \viscosite \Delta^2) (\psi^0_t + \psil)= \tau+\dtl\left(x,\frac{y-y_1}{\viscosite^{1/4}}\right)  + r_{\mathrm{lift}}^1 + r_{\mathrm{lift}}^2,
$$
where
$$
\| r_{\mathrm{lift}}^1\|_{L^2(\Om)} =o(\viscosite^{1/8}), \quad \| r_{\mathrm{lift}}^2\|_{H^{-2}(\Om)} =o(\viscosite^{5/8}),
$$
and
\begin{eqnarray*}
 \dtl(x,Y)&=&\chi^+ (x, y_1 + \viscosite^{1/4}Y)\chi\left(Y\right)\left[a'(x) + \viscosite^{1/4}b'(x) Y\right]\\
&&- \chi^+ (x, y_1 + \viscosite^{1/4}Y)\chi^{(4)}(Y) \left[a(x) + \viscosite^{1/4} b(x) Y\right]\\
&&-4\chi^+ (x, y_1 + \viscosite^{1/4}Y)\viscosite^{1/4} \chi^{(3)}(Y)b(x).
\end{eqnarray*}

\label{lem:eq-psil}
\end{Prop}

\begin{Def}\label{psiint-def}
In the rest of the paper, we set
\be \label{def:psiint}
\psi^{int}:= \psi^0_t + \psil;
\ee
which is consistent with the statement of Theorem \ref{thm}, in which $\psi^{int}$ was described as an $H^2$ regularization of $\psi^0$.
\end{Def}

\begin{proof}[Sketch of proof]
By Proposition \ref{H2-prop}, $\psi^0_t + \psil$ belongs to $H^2(\Omega)$ and, using \eqref{dec-DeltaPsi},
\begin{eqnarray*}
 \d_x(\psi^0_t + \psil)&=& \tau\chi_\viscosite+\chi^+\chi\left(\frac{y-y_1}{\viscosite^{1/4}}\right)(a'(x)+ b'(x)(y-y_1)),\\
\Delta (\psi^0_t + \psil)&=& \chi^+(\Delta(\psi^0_t+\psil))+\chi^-( \Delta\psi^0_t).
\end{eqnarray*}

Moreover, as we have the following trace identities on $\d\Omega^+$
$$(\Delta \psil)_{|\Sigma} =0, \quad (\d_y \Delta \psil)_{|\Sigma} =0,$$
using again \eqref{dec-DeltaPsi}, we get
\begin{eqnarray}
\nonumber \viscosite \Delta(\chi^+ \Delta\psil)&=&\viscosite\chi^+ \chi\left(\frac{y-y_1}{\viscosite^{1/4}}\right)(a^{(4)}(x)+ b^{(4)}(x)(y-y_1))\\
\nonumber&&+ 4 \viscosite^{3/4}\chi^+ \chi'\left(\frac{y-y_1}{\viscosite^{1/4}}\right)b''(x)\\
\label{eq:der23}&&+ 2 \viscosite^{1/2}\chi^+\chi''\left(\frac{y-y_1}{\viscosite^{1/4}}\right)(a'' (x)+ b''(x)(y-y_1)) \\&&
\nonumber+\chi^+\chi^{(4)}\left( \frac{y-y_1}{\viscosite^{1/4}}\right)\left[a(x) +  b(x) (y-y_1)\right]\\
\nonumber&&+4\viscosite^{1/4} \chi^+ \chi^{(3)}\left(x, \frac{y-y_1}{\viscosite^{1/4}}\right)b(x).
\end{eqnarray}
We thus define
$$
\begin{aligned}
r_{\mathrm{lift}}^1&=  4 \viscosite^{3/4}\chi^+ \chi'\left(\frac{y-y_1}{\viscosite^{1/4}}\right)b''(x)+ 2 \viscosite^{1/2}\chi^+\chi''\left(\frac{y-y_1}{\viscosite^{1/4}}\right)(a'' (x)+ b''(x)(y-y_1))\\
r_{\mathrm{lift}}^2&=\tau (\chi_\viscosite - 1)+\viscosite\chi^+ \chi\left(\frac{y-y_1}{\viscosite^{1/4}}\right)(a^{(4)}(x)+ b^{(4)}(x)(y-y_1))\\
&- \viscosite \Delta^2(\chi^+ \psi^0_t) -  \viscosite \Delta^2(\chi^- \psi^0_t).
\end{aligned}
$$
so that
\begin{eqnarray*}
&&(\d_x- \viscosite \Delta^2) (\psi^0_t + \psil)\\
&=& \tau+\dtl\left(x,\frac{y-y_1}{\viscosite^{1/4}}\right)  + r_{\mathrm{lift}}^1 + r_{\mathrm{lift}}^2,
\end{eqnarray*}

\begin{itemize}
\item The remainder terms will be dealt with in the last chapter, using some technical estimates on $a$, $b$ and their derivatives (see Lemmas \ref{lem:trunc} and \ref{lem:est-A-B} and paragraph \ref{ssec:proof-eq-psil}).
\item
The source term $\dtl$ is not an admissible error term in the sense of Definition \ref{remainder-def}. Therefore it has to be corrected by yet another boundary layer type term.
\end{itemize}
\end{proof}

\subsection{The interior singular layer}$ $
The lifting term $\psil$ has introduced a remainder  $\dtl$ which we now treat as a source term. By analogy with North and South boundary layers, we therefore define $\psisig$ as the solution of the following equation
\be\label{eq:psisig}
\begin{aligned}
\d_x \psisig - \d_Y^4 \psisig=-\dtl,\quad x<x(\sigma_{in}), \quad Y>Y_-(x)\\
\psisig_{|x=x(\sigma_{in})}=0,\quad
\psisig_{|Y=Y_-(x)}= \d_n  \psisig_{|Y=Y_-(x)}=0,
\end{aligned}
\ee
where the function $Y_-$  parametrizes the lower boundary of the interior singular layer (see Figure \ref{fig:sigma-layer-ch2})
\begin{itemize}
\item which has typical size $\delta_y$, where $\delta_y$ is a parameter such that $\delta_y\gg \viscosite^{1/4}$,
\item which is of course included in $\Omega$.
\end{itemize}
Note that we use here cartesian coordinates $x$ and $Y = (y-y_1)/\viscosite^{1/4}$, since the discontinuity line $\Sigma$ is essentially horizontal (which is not the case for North and South boundary layers which are extended on macroscopic parts of the adjacent East and West boundaries).
The precise definitions of $Y_-$ and $\sigma_{in}$ will be given in the next chapter (see page \pageref{Y-}).

\section{The case of islands}\label{sec:islands}$ $

When $\Om$ is not simply connected, the solution of \eqref{Munk-islands}-\eqref{compatibility} is given by (see Lemma \ref{lem:def-M,D})
$$
\psi= \psi_1 + \sum_{i=2}^K c_i \psi_i,
$$
where  $\psi_i$ is the solution of a Munk-like equation,
and the constants $c_i$ satisfy a linear system whose coefficients depend on the $\psi_i$.

Therefore, as explained in the introduction  (paragraph \ref{ssec:islands}) the main issue is to compute the asymptotic behaviour of every function $\psi_i$. Following the steps described in the preceding sections (and which will be developed in the next chapter), we are able to construct a function $\psi_i^{app}$ for every $1\leq i \leq K$. Of course, there are a few minor changes, due to the non homogeneous boundary condition, but we leave those to the reader, as they do not bear any additional difficulty.

Let us admit for the time being that Theorem \ref{thm} holds for the functions $\psi_i$, namely that we are able to construct $\psi_i^{app}$ such that
\begin{equation}
\label{thm-consequence}
\|\psi_i-\psi_i^{app}\|_{L^2}=o(\viscosite^{1/8}),\quad \|\psi_i-\psi_i^{app}\|_{H^2}=o(\viscosite^{-3/8}).
\end{equation}
Hence, in order to prove Theorem \ref{thm} for the function $\psi$, it suffices to define
$$\psiapp =\psi_1^{app} + \sum _{i=2}^K c_i^{app} \psi_i^{app}$$
where the constants $c_i^{app}$ are such that
\be\label{ciapp}
|c_i-c_i^{app}|=o(\viscosite^{1/8}), \quad c_i^{app}=O(1).
\ee
According to \eqref{est:L2}, \eqref{est:H2}, we have $\|\psi_i^{app}\|_{L^2}= O(1)$, $\|\psi_i^{app}\|_{H^2}= O(\viscosite^{-1/2})$, so that \eqref{thm-consequence} and \eqref{ciapp} imply
$$
\|\psi-\psiapp\|_{L^2}=o(\viscosite^{1/8}),\quad \|\psi-\psiapp\|_{H^2}=o(\viscosite^{-3/8}).
$$

We now turn to the proof of \eqref{ciapp}, which relies on the following Lemma:
\begin{Lem} For all $i\in \{2,\cdots, K\}$, let $g_i\in \mathcal C^\infty(\bar \Om)$ such that $g_i\equiv 1$ in a neighbourhood of $C_i$, and $\supp g_i\cap C_j=\emptyset$ for $i\neq j$.

Let us define the matrix
$$
M^{app}:=\left(-\int_{\Om} \psi_j^{app} \d_x g_i\right)_{2\leq i,j\leq K}
$$
and the vector
$$
D^{app}:=\left(\int_{\Om} (\tau g_j + \psi_1^{app}\d_x g_j) + \int_{C_j} \cT^\bot\cdot n\right)_{2\leq j\leq K}.
$$
Then the following facts hold (recall that $M^\viscosite, D^\viscosite$ are defined in Lemma \ref{lem:def-M,D}):
\begin{itemize}
 \item $|M^\viscosite-M^{app}|=o(\viscosite^{1/8})$, $|D^\viscosite-D^{app}|=o(\viscosite^{1/8})$;
\item $M^\viscosite$, $M^{app}$ are invertible matrices, with $(M^\viscosite)^{-1}=O(1), (M^{app})^{-1}=O(1)$;
\item $D^{app}=O(1)$.

\end{itemize}

\label{lem:capp}

\end{Lem}

Before addressing the proof of the lemma, let us explain why \eqref{ciapp} follows: by definition,
$$
M^\viscosite \begin{pmatrix} c_2\\\vdots\\ c_K\end{pmatrix}= M^\viscosite c= D^\viscosite,
$$
so that
$$
M^\viscosite (c-c^{app})= (M^{app}-M^\viscosite)c^{app}+ D^\viscosite - D^{app},$$
where $c^{app} $ is defined by
$$
M^{app} c^{app}= D^{app}.
$$
Using the last two items of Lemma \ref{lem:capp}, we infer that $c^{app}=O(1)$. Therefore
$$
M^\viscosite (c-c^{app})=o(\viscosite^{1/8}).
$$
Since $(M^\viscosite)^{-1}=O(1)$,  \eqref{ciapp} holds.

\begin{proof}[Proof of Lemma \ref{lem:capp}]
$\bullet$ The first step is to use the equation satisfied by every $\psi_i$ in order to express the coefficients of $M^\viscosite$, $D^\viscosite$ as  integrals on $\Om$. More precisely, we have, since $\nabla g_i=0$ on a neighbourhood of $\d\Om$,
\begin{eqnarray*}
 \viscosite \int_{C_i}\d_n \Delta \psi_j&=& \viscosite \int_{\d\Om} g_i \d_n \Delta \psi_j\\
&=& \viscosite \int_{\Om}g_i \Delta^2 \psi_j + \viscosite \int_{\Om } \nabla \Delta \psi_j \cdot \nabla g_i\\
&=& \int_{\Om} g_i \d_x \psi_j -\viscosite \int_{\Om} \psi_j \Delta^2 g_i\\
&=& - \int_{\Om} \psi_j \d_x g_i + \delta_{ij}\int_{C_i} e_x\cdot n-\viscosite \int_{\Om} \psi_j \Delta^2 g_i.
\end{eqnarray*}
Since $C_i$ is a closed curve, we obtain eventually
$$
 \viscosite \int_{C_i}\d_n \Delta \psi_j=- \int_{\Om} \psi_j \d_x g_i -\viscosite \int_{\Om} \psi_j \Delta^2 g_i.
$$
In a similar way,
$$
-\viscosite \int_{C_j}\d_n \Delta \psi_1 = \int_{\Om}( \tau g_j + \psi_1\d_x g_j) +  \viscosite \int_{\Om} \psi_1 \Delta^2 g_i.
$$
The energy estimate \eqref{weight-est} shows that $\|\psi_j\|_{L^2}=O(1)$ for $1\leq j\leq K$, so that
\be\label{MDnu}\begin{aligned}
M^\viscosite= \left(- \int_{\Om} \psi_j \d_x g_i\right)_{2\leq i,j\leq K} + O(\viscosite),\\
D^\viscosite= \left(   \int_{\Om}( \tau g_j + \psi_1\d_x g_j)+ \int_{C_j} \cT^\bot\cdot n\right)_{2\leq j\leq K} + O(\viscosite).
  \end{aligned}
\ee
Replacing every function $\psi_j$ by $\psi_j^{app}$ and using the error estimate (\ref{thm-consequence}), we obtain the first point of the lemma.

$\bullet$ To prove the second point of the Lemma, we will use identity \eqref{in:psia} in  Appendix A, from which we will deduce the following coercivity inequality: there exists a constant $C>0$, independent of $\viscosite$, such that
\be \label{coercivity}
\forall a=(a_2,\cdots, a_K)\in \R^{K-1},\quad \sum_{i,j} M^\viscosite_{ij} a_ia_j \leq -C |a|^2.
\ee
Of course this entails easily that $(M^\viscosite)^{-1}$ is bounded. Using the first item of the Lemma, we infer that $M^{app}$ satisfies the same type of coercivity property, and therefore $(M^{app})^{-1}$ is bounded as well.

Let us now prove \eqref{coercivity}.
Let $a\in \R^{K-1}$ be a fixed vector. Let $V\subset \Om$ be an open set such that $V$ is a neighbourhood of $\Gamma_W$, namely
\begin{itemize}
 \item $\d V\cap \d \Om\subset \Gamma_W$ (notice that this implies that $\bar V$ does not intersect $\Gamma_{N,S}$, nor $\Gamma_E$);

\item $\bar V \cap \Sigma=\emptyset$ ($V$ does not meet the discontinuity boundary layer) ;

\item $V\cap B_i\neq \emptyset$, and $\d V\cap \bar B_i\cap \Gamma_W\neq \emptyset$ for all $2\leq i\leq K$ ($V$ intersects  the West boundary of every set $B_i$).

\end{itemize}

Then we may write, with $\psi^{app}_a= \sum a_i \psi^{app}_i$,
\begin{eqnarray*}
 \viscosite \int_{\Om}|\Delta \psi_a|^2&=&\viscosite \int_{\Om} |\Delta \psi_a^{app}|^2 + o(\viscosite^{1/8}|a|^2)\\
&\geq & \viscosite \int_{V} |\Delta \psi_a^{app}|^2 + o(\viscosite^{1/8}|a|^2)\\
&\geq & \viscosite\sum_{i=2}^K  \int_{V\cap B_i} |\Delta \psi_a^{app}|^2 + o(\viscosite^{1/8}|a|^2).
\end{eqnarray*}
Notice that $ \viscosite \int_{\Om}|\Delta \psi_a|^2=O(1)$ (thanks to energy estimates), while $\psi_a^{app}$ is proportional to $a_i$ on $V\cap B_i$. These are the two main ingredients for the proof of the coercivity inequality. More precisely, on every set $V\cap B_i$, by definition\footnote{Notice that on $B_i$, the main order of $\psi_i^{app}$ is $\mathbf 1_{B_i}$, whose normal derivative on $\d \Om$ is zero (away from the end point of $\Sigma$). As a consequence, several corrector terms, whose construction is necessary in the general case, are identically zero here.},
$$
\psi_a^{app}=a_i \left(1+ \chi_0(z)\sum_{\pm} \frac{e^{\pm i\pi/6}}{\sqrt{3}}\exp(-\lambda_W e^{\pm i \pi/3} z)\right).
$$
Consequently, at the main order,
\begin{eqnarray*}
 \viscosite  \int_{V\cap B_i} |\Delta \psi_a^{app}|^2&\simeq& \viscosite a_i^2 \int_{V\cap B_i}\lambda_W^4\left| \chi_0(z)\sum_{\pm} \frac{e^{\pm 5 i\pi/6}}{\sqrt{3}}\exp(-\lambda_W e^{\pm i \pi/3} z)\right|^2ds\:dz\\
&\geq & C \viscosite a_i^2 \int_{\bar V\cap \bar B_i \cap \Gamma_W} \lambda_W^3(s)\:ds\\
&\geq & C a_i^2.
\end{eqnarray*}
Gathering all the terms, we obtain \eqref{coercivity}.

\smallskip

$\bullet$ The fact that $D^{app}=O(1)$ is obvious and follows from \eqref{weight-est}, which yields $\|\psi_1\|_{L^2}=O(1)$.
\end{proof}

For the sake of completeness, we also give the system satisfied by $\bar c=\lim_{\viscosite \to 0} c= \lim_{\viscosite \to 0} c^{app}$:

\begin{Cor}
Let $\bar c$ be the solution of
$$
\bar M \bar c= \bar D,
$$
where
$$
\begin{aligned}
\bar M:=\left(-\int_{\d B_j\cap C_i} e_x\cdot n\right)_{2\leq i,j\leq K},\\
\bar D:= \left(\int_{\Om} (\tau g_j + \psi^0\d_x g_j) + \int_{C_j} \cT^\bot\cdot n\right)_{2\leq j\leq K}.
\end{aligned}
$$
Then $\bar c=\lim_{\viscosite \to 0} c_i$.

\end{Cor}
\begin{proof}
We pass to the limit in the formulas \eqref{MDnu}.
 It is proved in \cite{BGGRB}, and it can be easily deduced from our construction, that
$$
\psi_j\rightharpoonup \mathbf 1_{B_j} \quad\text{in }L^2(\Om)\text{ for }2\leq j\leq K,
$$
and therefore
$$
- \int_{\Om} \psi_j \d_x g_i\to - \int_{B_j} \d_x g_i= - \int_{\d B_j\cap C_i} e_x\cdot n.
$$
Hence the identities \eqref{MDnu} imply easily that
$$
M^\viscosite \to \bar M, \quad D^\viscosite \to \bar D.
$$
Passing to the limit in inequality \eqref{coercivity}, we infer that $\bar M$ is also coercive, and therefore invertible. Thus $\bar c$ is well-defined. The corollary follows.

\end{proof}

{\section{North and South periodic boundary layers}
\label{sec:BL-per}
When the domain $\Om$ is $\Om:=\T\times (y_-, y_+)$, boundary layers occur in the vicinity of $y_-$, $y_+$. Because of the periodicity, they can be easily computed by using Fourier series in the $x$ variable. Therefore this kind of construction is covered by article \cite{GVP}; we include the computations here for the sake of completeness.

By definition of the Sverdrup term in \eqref{psiapp-per}, the boundary conditions to be lifted on the north and south boundary have zero average. Let us focus for instance on the south boundary; setting $Z=(y-y_-)/\viscosite^{1/4}$, the equation satisfied by the boundary layer term is
$$
\begin{aligned}
\d_x \psbl_S - \d_Z^4 \psbl_S=0\quad\text{in } \T\times \R_+,\\
\psbl_{S|Z=0}=f_0,\ \d_x\psbl_{S|Z=0}=f_1,
\end{aligned}
$$
with $\mean{f_0}= \mean{f_1}=0$. Looking for $\psbl_S$ in the form
$$
\psbl_S(x,Z)=\sum_{k\in 2\pi \Z\setminus\{0\}} \hat \psi_k(Z)\exp(ik x),
$$
we infer that $\hat \psi_k$ is given by
$$
\hat \psi_k(Z)= \sum_\pm A_k^\pm \exp(-|k|^{1/4}\lambda_k^\pm Z),
$$
where
$$
\begin{aligned}
\lambda_k^+= e^{i\pi/8}, \ \lambda_k^-=e^{-3i\pi/8}\quad \text{if }k>0,\\
\lambda_k^+= e^{3i\pi/8}, \ \lambda_k^-=e^{-i\pi/8}\quad \text{if }k<0,
\end{aligned}
$$
and
$$
\begin{aligned}
A_k^+ + A_k^-= \hat f_0(k),\\
\lambda_k^+ A_k^+ + \lambda_k^- A_k^- = -|k|^{-1/4} \hat f_1(k).
\end{aligned}
$$
Hence
$$
|A_k^\pm| \leq C (|\hat f_0(k)| + |k|^{-1/4}|\hat f_1(k)|  ).
$$
It follows easily that
\be\label{est:psblNS-per}
\begin{aligned}
\|\psbl_S\|_{L^2(\T, H^2(\R_+))}\leq C \left(\sum_{k\neq 0} |k|^{3/4} (|\hat f_0(k)| + |k|^{-1/4}|\hat f_1(k)|  )^2\right)^{1/2},\\
\|\d_x^l \psbl_S\|_{L^2(\T, H^2(\R_+))}\leq C \left(\sum_{k\neq 0} |k|^{2l+3/4} (|\hat f_0(k)| + |k|^{-1/4}|\hat f_1(k)|  )^2\right)^{1/2}.
\end{aligned}
\ee

}

%% file: approximate.tex
\chapter{Construction of the approximate solution}
\label{chap:construction}

The approximate solutions to the Munk equation (\ref{M}) are obtained by gathering together the different elementary pieces described in the previous chapter, namely
\begin{itemize}
\item the interior term which is essentially the solution to the transport equation (\ref{sverdrup}), regularized in the vicinity of ``East corners'' and of the interface $\Sigma$;
\item East and West boundary layer terms, which lift locally the boundary conditions but become singular in the vicinity of the points $s_i$;
\item North and South boundary layer terms, which lift the boundary conditions on the horizontal parts of the boundary but in a non local way;
\item singular layer terms, which make up for the source terms introduced by the regularization of the discontinuity at $\Sigma$.
\end{itemize}
We have now to understand the interplay between those different elementary pieces.

Of course, the equation (\ref{M}) being linear, the errors induced by all these terms are simply added, so that the control on the remainders in the approximate equation will be rather simple to obtain. The point to be stressed is that we need  each elementary term to be smooth enough (namely $H^2$), with suitable controls on the corresponding derivatives.

\bigskip

More precisely, we will consider here only one basic problem among  (\ref{psi1-eq})(\ref{psii-eq}), say the case (\ref{psi1-eq}) with homogeneous boundary conditions $\psi_{|\Omega} =\d_n \psi_{|\Omega} =0$ and forcing $\tau$. We will further assume, without loss of generality, that we have the simple geometry described in the previous chapter
$$
\Om=\Om^+ \cup \Om^- \cup \Sigma,
$$
where
\begin{itemize}
 \item $\Om^\pm$ are non empty, open and convex in $x$;
\item $\Om^+= \{(x,y)\in \Om,\ y>y_1\text{ or } x\geq x_1\}$;
\item $\Om^-= \{(x,y)\in \Om,\ y<y_1\text{ and } x< x_1\}$;
\item $\Sigma=\{(x,y_1)\in \Om,\  x<x_1\}$
\end{itemize}
and $(x_1, y_1)=(x(s_1), y(s_1))$ for some $s_1$ such that $\cos \theta(s_1)=0$ and $\sin \theta(s_1)=-1$.
We emphasize that the only simplification here regards notations, and that more complex domains satisfying assumptions $(H1)-(H4)$ are handled exactly in the same way.

Periodic and rectangle cases will  be dealt with separately in Section 3.6.
\bigskip

We will actually focus  on the following technical difficulties
\begin{itemize}
\item the precise regularization  process for the interior term in $\Omega^\pm$  (section \ref{sec:interior});
\item the fact that the East  boundary layer does not lift simultaneously both boundary conditions (correcting the normal derivative introduces indeed an error on the trace). This implies that one has to introduce an additional corrector defined on a macroscopic domain in the vicinity of the East boundary but far from corners $(x_i,y_i)$ where it would be singular (section \ref{sec:macro-corrector});
\item the connection with North/South boundary layers (even if the corresponding horizontal part of the domain is reduced to one point $(x_i,y_i)$): we indeed expect these boundary layers to lift  the boundary conditions both on the horizontal parts and on the East/West boundaries close to the corners  $(x_i,y_i)$  (section  \ref{sec:NSboundary}).
Note that, since these North/South boundary layers are defined by a non local equation, we need to check that they do not carry any more energy beyond  some point, so that they can be truncated and considered  as local contributions;
\item the truncation of the surface layer near its West end, in order that its trace on the West boundary is not too singular (see \ref{sec:Sigmaboundary});
\item the precise definition of the West boundary layers, which have to lift the traces of all previous contributions (see \ref{sec:Wboundary});
\end{itemize}
We will also derive estimates on the terms constructed at every step.

Note that, {\bf for boundary layers, the construction must be performed in a precise order, starting with  East boundaries, then defining North/South (and surface) boundary layers, and finally lifting the West boundary conditions}. This dissymmetry between East and West boundaries is similar to what happens at the macroscopic level for the interior term.

\section{The interior term}
\label{sec:interior}

As in the previous chapter, we define
\be\label{def:psi0}
\psi^0(x,y)=-\int_x^{x_E(y)}\tau(x',y)\:dx',
\ee\label{psi0}
where $x_E(y)$  is the abscissa (or longitude) of the projection of $(x,y)$ on  $\Gamma_E$ alongside $e_x$.
We have seen that $\psi^0$ does not belong to $H^2(\Om)$ in general: indeed,

\begin{itemize}
\item
 The function $x_E$ has singularities near the points $y_i= y(s_i), i\in I_+$;
 \item Since $x_E$ takes different values on $\Om^+$ and $\Om^-$, the function $\psi^0_t$ and its $y$ derivative are discontinuous across $\Sigma$.
 \end{itemize}
  Therefore, $\psi^0$ cannot be used as such in the definition of the approximate solution $\psapp$ of \eqref{M}.

  \bigskip

  Let us first consider separately the domains $\Omega^\pm$. We remove the singularities  by truncating the function $\tau$ near the end points of $\Gamma_E$ (with abscissa $s_i, i\in I_+$).
\begin{figure}\label{fig:troncature}
 \includegraphics[height=8cm]{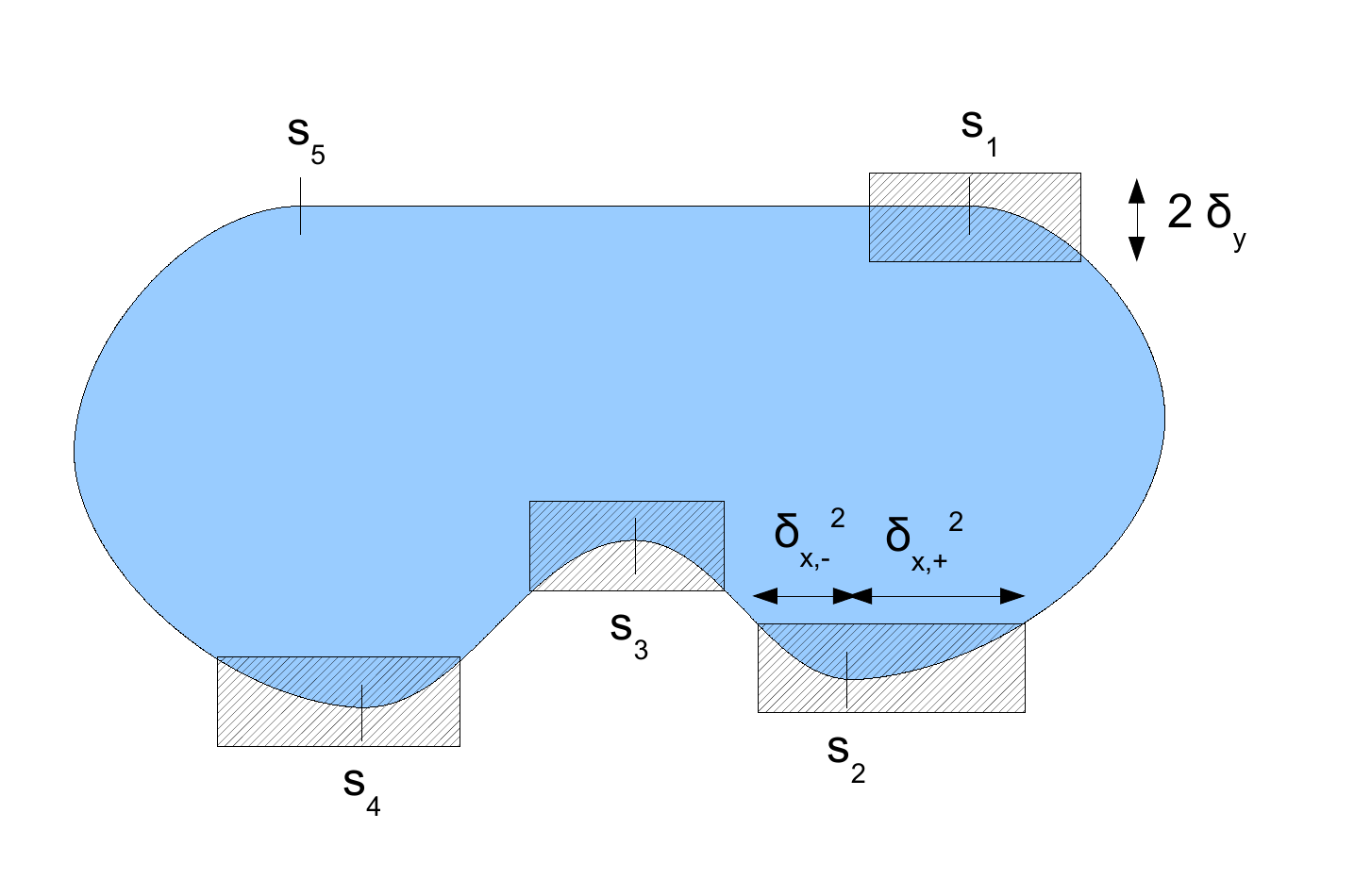}
\caption{The function $\tau$ is truncated in the hatched zones. In the picture, $I_+=\{1,2,3,4\}$.}
\end{figure}

 The size of the truncation depends on the rate of cancellation of $\cos \theta$ near $s_i$. Since the rate of cancellation may be different on the left and on the right of $s_i$, we seek for truncation functions with different behaviours on the left and on the right of every cancellation point.
More precisely, let
\be\label{def:chi}
\chi_\viscosite(x,y):=\prod_{i\in I_+} \chi_i (x-x_i, y-y_i)
\ee
where $(x_i,y_i)$ are the coordinates of the point of $\d\Om$ with arc length $s_i$. Each function $\chi_i$ takes the following form:
$$
\chi_i(x,y):=\left\{ \begin{array}{lc}
                                \bar \chi\left(\frac{x}{\delta_{x,+}^{i}}, \frac{y}{\delta_y}\right)&\text{ if }x>0,\\
                                \bar \chi\left(\frac{x}{\delta_{x,-}^{i}}, \frac{y}{\delta_y}\right)&\text{ if }x\leq 0.
                     \end{array}\right.
$$

$\bullet$
The function $\bar \chi$ is defined by
$$
1-\bar \chi(x,y):=(1-\bar \chi_1(x))(1-\bar\chi_2(y))\quad (x,y)\in \R^2,
$$
with $\bar \chi_1, \bar\chi_2\in \mathcal C^\infty(\R)$, and
$$
\begin{aligned}
\bar\chi_1(x)=1\text{ if } |x|\geq  1+\frac{1}{|\ln \delta_y|},\quad \bar\chi_1(x)=0\text{ if } |x|\leq  1,\\
\bar \chi_2(y)=1 \text{ if } |y|\geq  1,\quad \bar\chi_2(y)=0\text{ if } |y|\leq  1/2 \hbox{ which does not depend on } \viscosite.
\end{aligned}
$$
\label{chii}
 It can be easily checked that all the $x$ derivatives of $\bar \chi$ vanish for $x=0$, so that $\chi$ belongs to $\mathcal C^\infty_0(\R^2)$.
 Furthermore, it  satisfies obviously
$$
\begin{aligned}
0\leq \bar \chi\leq 1,\\
\bar \chi(x,y)=0\text{ if }|y|\leq \frac{1}{2}\text{ and }|x|\leq  1,\\
\bar \chi (x,y)=1\text{ if }|y|\geq 1\text{ or }|x|\geq  1+\frac{1}{|\ln \delta_y|}.
\end{aligned}
$$

$\bullet$
Let us now give the definition of the truncation rates. We take\label{deltay}
\be\label{def:deltay}
 \delta_y:=\viscosite^{1/4} |\ln \viscosite  |^{1/5} (\ln |\ln \viscosite |)^{-\beta},
\ee
for some arbitrary exponent $\beta>0$.   Note that a rate of this type is mandatory in the case of an exponential cancellation (assumption (H2ii)). In the case of an algebraic cancellation (assumption (H2i)), there is more flexibility, but the above choice still works.

The definition of $\delta_x$ is a little more involved.  We define $\delta_x$ so that if we start from the point $(x_E(y_i), y_i)$ (i.e. from the end point of $\Gamma_E$) and perform a shift of size $\delta_x$ in the $x$ direction, $-\sin \theta(s_i)\delta_y$ in the $y$ direction, the end point still belongs to $\Gamma_E$. This leads to the following definitions:
\begin{itemize}
\item If $x_E(y)$ is defined on both sides of $y_i$, i.e. if $s_i$ belongs to the interior of $\bar \Gamma_E$, then we set
\be\label{def:deltax1}
\begin{aligned}
\delta_{x, +}^i:=|x_E(y_i -\sin \theta(s_i). \delta_y)- x_i|,\\
\delta_{x, -}^i:=|x_E(y_i+ \sin \theta(s_i).  \delta_y)- x_i| \,.
\end{aligned}
\ee
\item If $x_E(y)$ is defined only on one side of $y_i$, say for $y>y_i$, then we take
\be\label{def:deltax2}
\delta_{x, +}^i= \delta_{x, -}^i= | x_E(y_i+ \delta_y)-x_i|.
\ee
\label{deltax}
\end{itemize}

In particular,  if $(x_E(y), y)\in \supp  \chi_\viscosite $, then
$$
\text{either }\inf_{i\in I_+}\frac{|y-y_i|}{\delta_y}\geq \frac{1}{2} \text{ or } \inf_{i\in I_+}\sup_\pm\frac{(x_E(y)-x_i)_\pm}{\delta_{x,\pm}^i}\geq 1.
$$
According to \eqref{def:deltax1}-\eqref{def:deltax2}, the above condition becomes
$$
\text{either }\inf_{i\in I_+}\frac{|y-y_i|}{\delta_y}\geq \frac{1}{2}\text{ or }\inf_{i\in I_+}\sup_\pm\frac{(x_E(y)-x_i)_\pm}{|x_E(y_i\mp\sin \theta(s_i) \delta_y)-x_i)|}\geq 1.
$$
Since $x_E$ is locally a monotonous function (near $y=y_i$), this amounts to $$\inf_{i\in I_+}|y-y_i|\geq \delta_y/2.$$
We infer that
\be\label{est:xE-suppchi}
(x_E(y), y)\in \supp \chi_\viscosite  \ \Rightarrow\ \forall i\in I_+,\ |y-y_i|\geq \frac{\delta_y}{2}.
\ee

\bigskip

We then set
\be\label{def:psi0t}
\psi^0_t(x,y):=-\int_x^{x_E(y)}\chi_\viscosite (x',y)\tau(x',y)\:dx'
\ee
\label{psi0t}
With this definition, $\psi^0_t\in H^2(\Om^\pm)$, but $\psi^0_t \notin H^2(\Om)$: indeed $\psi^0_t$ and $\d_y \psi^0_t$ are still discontinuous across $\Sigma$. The term $\psil$ constructed in paragraph \ref{ssec:psil} of the present chapter will lift this discontinuity.

 Moreover, by definition of $\psi^0_t$, we have
$$
\d_x \psi^0_t - \viscosite \Delta^2 \psi^0_t= \tau + \delta\tau \quad\text{on }\Om^+ \cup \Om^-
$$
with
\be
\label{dtauint}
 \delta \tau = \tau (\chi_\viscosite -1)  - \viscosite \Delta^2 \psi^0_t.
\ee

        \begin{Prop}
                Assume that $\Om$ satisfies assumptions {(H1), (H2)}. Let $\beta>0$ be arbitrary, and let $\delta_y$, $\delta_x$ as in \eqref{def:deltay}, \eqref{def:deltax1}-\eqref{def:deltax2}.
Then
$$
\begin{aligned}
\| \tau (\chi_\viscosite -1)\|_{H^{-2}(\Om)}= o (\viscosite ^{5/8}),\\
\| \Delta \psi^0_t\|_{L^2(\Om^\pm)}= o(\viscosite ^{-3/8}).
\end{aligned}
$$

\label{lem:trunc}

        \end{Prop}
This proposition will be proved in Chapter \ref{chap:convergence}.

\begin{Rmk}
 Notice that the rate of cancellation in $y$, namely $\delta_y$, is the same for every function $\chi_i$. It is far from obvious \textit{a priori }that such a choice can lead to a suitable truncation. However, in order not to further burden the notations, we have chosen to anticipate on this result, which follows from the proof of Lemma \ref{lem:trunc} below.

The $x$ derivative of the function $\bar \chi$ is   unbounded (it is of order $|\ln  \delta_y|$ in $L^{\infty}$). This choice is mandatory in the case of an exponential cancellation of $\cos \theta$ around $s_i$ (assumption (H2ii)). If the cancellation around $s_i$ is algebraic (assumption (H2i)), the function $\bar \chi$ can be any smooth function vanishing near zero and such that $1-\bar \chi$ has compact support.
\end{Rmk}

\bigskip
We conclude this paragraph by giving some estimates on the trace of $\psi^0_t$ and $\d_n \psi^0_t$. These estimates will be useful when we construct the boundary layer terms lifting these boundary conditions. For the sake of readability, we introduce the following majorizing functions\label{cM-def}
$$\begin{aligned}
 \ccM(y) =  \prod_{i\in I_+} {1\over \max ( |(y-y_i)\ln |y-y_i|, \delta _y |\ln \delta_y|) } \\
\ccM^*(y) = \ccM(y) \prod_{i\in I_+}\mathbf{1}_{|y-y_i|\geq \delta_y/2 }\,.
 \,.
\end{aligned}$$

\begin{Lem}\label{lem:est-psi0}
Let $\psi^0_t$ be defined by
$$
\psi^0_t(x,y)=-\int_x^{x_E(y)}(\tau\chi_\viscosite )(x',y)\:dx'.
$$

\medskip
\noindent
$\bullet$
\textbf{Trace estimates on $\d\Om$:}
There exists a constant $C$, depending only on $\Om$ and $\tau$, such that for all $s\in  \d\Om$,
\begin{eqnarray*}
\|\psi^0_{t|\d\Om}\|_{L^\infty(\d\Om)}&\leq & C,\\
|\d_n \psi^0_{t|\d\Om}(s)|&\leq&  C\ccM(y)
\end{eqnarray*}

Moreover, for all $s\in \d\Om$ such that $(x(s), y(s))\notin \Sigma$,
$$  |\d_s \psi^0_{t|\d\Om}(s)| + \delta_y  |\d_s \d_n\psi^0_{t|\d\Om}(s)|\leq  C|\cos \theta(s)|\ccM(y(s))
$$
 and
\begin{eqnarray*}
|\d_s^2 \psi^0_{t|\d\Om}(s)| + \delta_y|\d_s^2 \d_n\psi^0_{t|\d\Om}(s)|
&\leq & C\left(\frac{1}{\delta_y} + \frac{|\cos \theta|^2}{\delta_y^2 |\ln \delta_y|}\right),\\
|\d_s^3 \psi^0_{t|\d\Om}(s)| + \delta_y|\d_s^3 \d_n\psi^0_{t|\d\Om}(s)|
&\leq & \frac{C}{\delta_y}\left(1 + \frac{|\cos \theta|}{\delta_y} + \frac{|\cos \theta|^3}{\delta_y^2}  \right).
\end{eqnarray*}
Eventually, on the East coast, we have the following more precise estimates:
\begin{eqnarray*}
|\d_s^2 \d_n\psi^0_{t|\Gamma_E}(s)|
&\leq &  C\frac{|\theta'(s)|}{\delta_y|\cos \theta(s)|}\prod_{i\in I_+} \mathbf{1}_{|y(s)-y_i|\geq \delta_y/2},\\
|\d_s^3 \d_n\psi^0_{t|\Gamma_E}(s)|
&\leq & C\frac{|\theta'(s)|^2}{\delta_y|\cos \theta(s)|^2}\prod_{i\in I_+} \mathbf{1}_{|y(s)-y_i|\geq \delta_y/2}.
\end{eqnarray*}

%
%

\medskip
\noindent
$\bullet$
\textbf{Jump estimate on $\Sigma$:}
The jumps $[\psi^0_t]_\Sigma$, $[\d_y \psi^0_t ]_\Sigma$ are constant along $\Sigma$. Moreover, there exists a constant $C$ such that
$$
\left|[\psi^0_t]_\Sigma\right|\leq C,\quad \left|[\d_y \psi^0_t ]_\Sigma\right|\leq  \frac{C}{\delta_y |\ln \delta_y|}.
$$

\end{Lem}

\begin{proof}
The $L^\infty$ bound on $\psi^0_{t|\d\Om}$ is obvious: we merely observe that
$$
\|\psi^0_{t|\d\Om}\|_{L^\infty(\d\Om)}\leq \|\tau\|_{L^\infty(\Om)}\sup_{(x,y),(x',y)\in \Om}|x-x'|.
$$

As for the bound on $\d_n\psi^0_t$, we have $\d_n=\cos \theta \d_x+\sin \theta \d_y$, so that
\begin{eqnarray}
\d_n\psi^0_{t|\d\Om}(s)&=& - \sin \theta(s) x_E'(y(s)) (\tau \chi_\viscosite )(x_E(y(s)), y(s))\\
&&-\sin \theta(s) \int_{x(s)}^{x_E(y(s))}\d_y( \tau\chi_\viscosite )(x,y(s))\:dx\\
&&+ \cos \theta (s) (\tau \chi_\viscosite )(x(s), y(s)).
\end{eqnarray}

The formulas in Appendix B together with (\ref{est:xE-suppchi}) imply that  for all $i\in I_+$, in the vicinity of $y_i$
\be\label{est:x_E'}
|x_E'(y)|\leq \frac{C}{|y-y_i| |\ln |y-y_i||^2}
\ee
so that
$$
(x_E(y), y)\in \supp \chi_\viscosite  \ \Rightarrow\
|x_E'(y)|\leq C\ccM(y) \leq {C\over \delta_y|\ln\delta_y|^2}.
$$

Recalling eventually that for all $c>1$, for $\viscosite $ sufficiently small,
$$\begin{aligned}
 \|\d_y \chi_\viscosite \|_\infty= O(\delta_y^{-1}),\\
\supp \d_y \chi_\viscosite \subset \bigcup_{i\in I_+}\left[x_i-c\delta_{x,-}^i, x_i+c\delta_{x,+}^i\right]\times\left[ y_i - \delta_y, y_i+\delta_y\right],
\end{aligned}
$$
with $\delta_{x,\pm}^i= O(|\ln \delta_y|^{-1}),$ this leads  to the estimate on $\d_n \psi^0_{t|\d\Om}$. We have indeed
\begin{equation}\label{est:nabla_chi}
\begin{aligned}
\Big| \int_{x(s)}^{x_E(y(s))}\!\!\!\d_y( \tau\chi_\viscosite )(x,y(s))\:dx
\Big| \leq &C\sum_{i \in I_+}  \| \d _y \chi_\viscosite \|_\infty  \indc_{
  |y-y_i| \leq \delta_y}  \int_{x(s)}^{x_E(y(s))} \!\!\!\indc_{ |x-x_i| \leq c\delta_x} dx \\
& + \sum_{i \in I_+} \int_{x(s)}^{x_E(y(s))} \!\!\!\chi_\viscosite  | \d_y \tau | (x,y(s)) dx\\
&\leq C\left(1+ {\delta_x \over \delta_y}\indc_{ |y-y_i| \leq \delta_y} \right)\,.
\end{aligned}
\end{equation}

The estimate on $\d_s \psi^0_t$ is obtained by differentiating the identity
$$
\psi^0_{t|\d\Om}(s)= - \int_{x(s)}^{x_E(y(s))} (\tau\chi_\viscosite )(x,y(s))\:dx
$$
with respect to $s$. Using \eqref{jacobien} and therefore the relation
$$
\frac{d}{ds}(x(s), y(s))= (\sin \theta(s), -\cos \theta(s)),
$$
 we infer that
\begin{eqnarray*}
\d_s  \psi^0_{t|\d\Om}(s)&=& \cos \theta(s) x_E'(y(s)) (\tau\chi_\viscosite )(x_E(y(s)), y(s))\\
&&+ \sin \theta(s)  (\tau\chi_\viscosite )(x(s), y(s))\\
&&-\cos \theta(s)\int_{x(s)}^{x_E(y(s))} \d_y(\tau\chi_\viscosite )(x,y(s))\:dx.
\end{eqnarray*}
 The estimates \eqref{est:x_E'} and \eqref{est:nabla_chi} above yield the desired inequality.

The other trace estimates are derived in a similar fashion. Notice that is is hard to derive a sharp global estimate, since the angle $\theta(s)$ for $s\in \Gamma_W$ is in general different from the angle $\theta(s')$, where $s'\in \Gamma_E$ is such that $y(s)=y(s')$. Therefore, there is no simplification for terms of the type
$
\cos^2(\theta(s)) x_E''(y(s))
$ when $s\in \Gamma_W$. On the East coast, however, we can use the formulas in Appendix B, from which we deduce that
$$
\cos^2(\theta(s)) x_E''(y(s))= \frac{\theta'(s)}{\cos \theta(s)}\quad \forall s\in \Gamma_E.
$$
Notice also that $\psi^0_t$ vanishes on $\Gamma_E$ by definition, that $\d_n\psi^0_{t|\Gamma_E}$ is supported in $\{|y(s)-y_i|\geq \delta_y/2\}$, and that all terms of the type
$$
\int_{x(s)}^{x_E(y(s))} \d_y^k (\tau \chi_\viscosite )(x,y(s))\:dx
$$
are zero for $s\in \Gamma_E$.
The upper-bounds for $\d_s^2\d_n \psi^0_{t|\Gamma_E}$, $\d_s^3\d_n \psi^0_{t|\Gamma_E}$ follow.

As for the jump estimates, we recall that with the notations of Chapter \ref{chap:multiscale}, Section \ref{sec:disc},
$$
[\psi^0_t]_\Sigma=-\int_{x_1}^{x_E^+(y_1)}(\tau \chi_\viscosite )(x,y_1)\:dx,
$$
and therefore the jump of $\psi^0_t$ is of order one. Notice that the assumptions of Chapter \ref{chap:multiscale}, Section \ref{sec:disc} imply that $x_1=x_E^-(y_1)$.
The jump of the $y$ derivative is given by
\begin{eqnarray*}
[\d_y \psi^0_t]_\Sigma&=&-{x_E^+}'(y_1)(\tau \chi_\viscosite )(x_E^+(y_1), y_1),\\
&&+ {x_E^-}'(y_1)(\tau \chi_\viscosite )(x_E^-(y_1), y_1)\\
&&-\int_{x_1}^{x_E^+(y_1)}\d_y(\tau \chi_\viscosite )(x,y_1)\:dx.
\end{eqnarray*}
The points $(x_E^\pm(y_1), y_1)$ do not depend on $\viscosite $. Therefore, either $(x_E^\pm(y_1), y_1)$ is an East corner ($1\in I_+$), and then $(\tau \chi_\viscosite )(x_E^\pm(y_1), y_1)=0$ for all $\viscosite >0$, or $(x_E^\pm(y_1), y_1)$ is not a corner, and then ${x_E^\pm}'(y_1)$ is bounded by a constant independent of $\viscosite $. Hence the first two terms of the right-hand side above are bounded as $\viscosite \to 0$. Using inequality \eqref{est:nabla_chi}, we infer  that
$$
\left|\int_{x_1}^{x_E^+(y_1)}\d_y(\tau \chi_\viscosite )(x,y_1)\:dx\right|\leq \frac{C}{\delta_y |\ln \delta_y|}.
$$

\end{proof}

%
%
%
%
%

At this stage, we have built on each subdomain $\Omega^\pm$ an interior term $ \psi^0_t $
\begin{itemize}
\item which vanishes on $\Gamma_E$,
\item which belongs to $H^2(\Omega^\pm)$ (typically with derivatives of order $O(\delta_y^{-1})$ in the vicinity of East corners),
\item and which approximately satisfies the transport equation $\d_x \psi_t^0 = \tau$.
\end{itemize}
Nevertheless
 $\psi^0_t$ does not fit the boundary conditions, and therefore boundary layer correctors must be defined. Following the direction of propagation of the equation at main order, we start with the East boundary layers.


\section{Lifting the East boundary conditions}
\label{sec:macro-corrector}

In this section, we focus on an interval $(s_i, s_{i+1})\subset \Gamma_E$.
We  recall that the domain $[s_i^+, s_{i+1}^-]\subset [s_i, s_{i+1}] $ of validity of the East boundary layer is given by
$$
\begin{aligned}
        s_i^+:=\sup\left\{s\in\left(s_i, \frac{s_i+s_{i+1}}{2}\right), \viscosite ^{1/3}|\theta'|\; |\cos \theta|^{-7/3} \geq 1\right\},\\
        s_{i+1}^-=\inf\left\{s\in\left( \frac{s_i+s_{i+1}}{2}, s_{i+1}\right), \viscosite ^{1/3}|\theta'|\; |\cos \theta|^{-7/3} \geq 1\right\}.
\end{aligned}
$$
Easy computations based on the explicit rate of cancellation of $\cos \theta$ near $s_i$ provide then the following Lemma, which shows in particular that, because of the truncation, the trace of $\d_n \psi^0_t$ is zero outside  the domain of validity of the East boundary layer.

\begin{Lem}\label{sigma-s-est}
Under assumption (H2i),
$$s_i^+-s_i \sim C_i \viscosite ^\frac{1}{4n+3}\,.$$
Under assumption (H2ii),
$$
s_i^+-s_i= -\frac{4\alpha}{\ln\viscosite }\left(1 - 6 \frac{\ln |\ln \viscosite |}{\ln \viscosite } + O\left(\frac{1}{\ln \viscosite }\right)\right). $$
In particular,
$$|y(s_i^+)-y_i|\ll \delta_y$$
and
\begin{equation}
\label{lE-est}
\lambda_E(s_i^+) \geq C\viscosite ^{-1/4}\,.
\end{equation}
\end{Lem}

\begin{proof}
Lemma \ref{sigma-s-est} is obtained by straightforward computations~:
\begin{itemize}
\item
If  (H2i) is satisfied, then for $s>s_i$
$$
\viscosite ^{1/3}|\theta'(s)|\; |\cos \theta(s)|^{-7/3} \sim C' \viscosite ^{1/3} |s-s_i|^{-\frac{4n+3}{3}}
$$
so that there exists a positive constant $C_i$ such that
\be\label{equiv_si+_alg}
s_i^+-s_i \sim C_i \viscosite ^\frac{1}{4n+3}.
\ee
\item  If  (H2ii) is satisfied, then for $s>s_i$
$$
\viscosite ^{1/3}|\theta'(s)|\; |\cos \theta(s)|^{-7/3} \sim C' \viscosite ^{1/3} \frac{\exp\left(\frac{4\alpha}{3(s-s_i)}\right)}{(s-s_i)^2}
$$
so that
$$
-2 \ln (s_i^+-s_i) + \frac{4\alpha}{3(s_i^+-s_i)}=- \ln C' - \frac{\ln \viscosite }{3} +o(1).
$$
Computing an asymptotic development of $s_i^+-s_i$ leads to
\be\label{equiv_si+_exp}
s_i^+-s_i= -\frac{4\alpha}{\ln\viscosite }\left(1 - 6 \frac{\ln |\ln \viscosite |}{\ln \viscosite } + O\left(\frac{1}{\ln \viscosite }\right)\right).
\ee

\end{itemize}

Using the formulas in Appendix B, we infer in particular that in case (H2i),
$$
y(s_i^+)-y_i\sim C \viscosite ^{\frac{n+1}{4n+3}},
$$
and $(n+1)/(4n+3)>1/4$, so that $|y(s_i^+)-y_i|\ll \delta_y$. In a similar way, if (H2ii) is satisfied,
\begin{eqnarray*}
|y(s_i^+)-y_i|&\sim& C |\ln \viscosite |^{-2}\exp\left(\alpha\frac{\ln \viscosite }{4\alpha}\left(1+6 \frac{\ln|\ln \viscosite |}{\ln \viscosite }+ O (|\ln \viscosite |^{-1})\right)\right)\\
&\sim & C |\ln \viscosite |^{-2} \exp\left(\frac{\ln \viscosite }{4}+ \frac{3}{2}\ln |\ln \viscosite | + O(1)\right)\\
&\sim & C |\ln \viscosite |^{-1/2}\viscosite ^{1/4}\ll \delta_y.
\end{eqnarray*}

Therefore, in all cases, we have $|y(s_i^+)-y_i|\ll \delta_y$, so that the extremities of the domain of validity   of the East boundary layer are in the truncation zone.
\end{proof}

\subsection{Traces of the East boundary layers}$ $

By definition of $\psi^0_t$, we have
$$(\psi_t^0)_{|\Gamma_E} = 0.$$
We thus define the East boundary layer to lift the trace of $\d_n \psi_t^0$, i.e.
$$
\psi_E(s,Z)=-\lambda_E^{-1}(s) \d_n \psi^0_{t|\Gamma_E}(s) \exp(-Z)
$$
where
$$
\lambda_E(s)= \left(\frac{\cos \theta}{\viscosite }\right)^{1/3}.
$$

But then $\psi_t^0 + \psibl_E(s,\lambda_E z)$ no longer satisfies the zero trace condition on $\Gamma_E$.
More precisely, since $\d_n \psi^0_{t|\Gamma_E}$ vanishes for $|y(s)-y_i|\leq \delta_y/2$, we have the following trace estimate:
\begin{Lem}\label{lem:psiE-trace}
The trace of $\psi_E$ satisfies the following bound
\be \label{est:psiE}
|\psi_{E|Z=0}(s)|\leq C\frac{\viscosite ^{1/3}}{|\cos \theta(s)|^{1/3}} \ccM^*(y(s)).
\ee
Moreover, its  derivatives with respect to $s$ satisfy
\begin{eqnarray*}
 |\d_s\psi_{E|Z=0}(s)|&\leq &C\viscosite ^{1/3}\delta_y^{-1}|\cos \theta|^{2/3}\ccM^*(y(s)),\\
  |\d_s^2\psi_{E|Z=0}(s)|&\leq &C\viscosite ^{1/3}\delta_y^{-1} \frac{|\theta'(s)|}{|\cos \theta(s)|^{1/3}}\ccM^*(y(s)),\\
    |\d_s^3\psi_{E|Z=0}(s)|&\leq &C\viscosite ^{1/3}\delta_y^{-1} \frac{|\theta'(s)|^2}{|\cos \theta(s)|^{4/3}} \ccM^*(y(s)).
\end{eqnarray*}

\end{Lem}
These estimates are a straightforward consequence of Lemma \ref{lem:est-psi0} and of the formula defining $\psi_E$.

Hence the remaining trace on $\Gamma_E$ is non-zero, and must be corrected. Note in particular that the bounds above are too singular in the vicinity of $s_i, i\in I_+$ in order that we can lift the remaining trace by a simple macroscopic corrector. We lift the trace in two different ways depending on the value of $s$:
\begin{itemize}
\item If $s$ is far from any point $s_i, i\in I_+$, we lift the trace thanks to a macroscopic corrector $\psi^{corr}_\Sigma $, which we construct in the next paragraph;

\item If $s$ is in a neighbourhood of size one of any $s_i$, we lift the trace of $\psi_E$ thanks to the North and South boundary layer terms; we will explain the latter construction in the next section.

\end{itemize}

\begin{Rmk}\label{rmk:east-corr2}
Since the East corrector is not captured by energy estimates, as explained in Remark \ref{rem:east-corr}, it would be tempting to construct a boundary layer corrector which is not a solution of the East boundary layer equation, but which lifts the trace of the normal derivative of $\psi^0_t$ without perturbing the zero order trace on the boundary, namely a corrector of the type
$$
A(s) z \exp(-z/\viscosite ^p),
$$
for some $p>0$ and for an adequate choice of $A(s)$. However, because of the strong singularities of $\d_n\psi^0_{t|\d\Om}(s)$ near $s_i$ and $s_{i+1}$, it can be proved that no choice of the parameter $p$ leads to admissible error terms. In other words, close to the singularity zones, the corrector should be an approximate solution of the equation. This justifies the need for an elaborate construction, even though the corresponding terms are negligible in the final energy. This is in fact classical in multi-scale problems: quite often, it is necessary to construct high  order correctors, whose sole purpose is to ensure that the remainder terms are admissible, but which are not seen by the total energy.

\end{Rmk}

\subsection{Definition of the East corrector}$ $

The idea is therefore to split each East boundary component $[s_i, s_{i+1}]$ in three subdomains (independent of $\viscosite $)
\begin{itemize}
\item  one which is far from the singularities
$$[\sigma_i^+,\sigma_{i+1}^-] \subset ]s_i, s_{i+1}[$$
so that we have uniformly small bounds on the trace  of $\psibl_E$ on $[\sigma_i^+,\sigma_{i-1}^-]$;

\item  two which are (macroscopic) neighbourhoods of the singularities and such that
\begin{equation}
\label{Esigmaipm}
\forall s \in [s_i,\sigma_i^+] \cup [\sigma_{i+1}^-, s_{i+1}],\quad |\sin\theta(s)| \geq \frac12\,,
\end{equation}
so that we can extend the North or South boundary layers on these parts of the boundary.

\end{itemize}

We therefore define suitable truncation functions  in $s$, namely $\varphi_i^+, \varphi_{i+1}^- \in C^\infty(\R, [0,1])$
\begin{equation}
\label{Evarphii}
\begin{aligned}
\varphi_i^+ (s)= 1 \hbox{ if  }s \geq \sigma_i^+ ,\quad \varphi_i^+(s) = 0  \hbox{ if } s \leq \frac12( s_i +\sigma_i^+) ,\\
\varphi_{i+1}^-(s) = 1 \hbox{ if  }s \leq \sigma_{i+1}^- ,\quad \varphi_{i+1}^-(s)  = 0  \hbox{ if } s \geq \frac12( s_{i+1} +\sigma_{i+1}^-) .
\end{aligned}
\end{equation}

\bigskip
The East corrector is then expected to lift  $ \varphi_i^+\varphi_{i+1}^- \psi_{E|Z=0}$, which is the part of the trace which is known to be small.
More precisely, we set, for $s\in [s_{i-1}, s_i]$,
\be \label{def:psi1}
\psi^{corr}_E(s,z)= - \psi_{E|Z=0}(s)\varphi_i^+(s)\varphi_{i-1}^- (s)\chi_0(z).
\ee

Using Lemmas \ref{sigma-s-est} and  \ref{lem:psiE-trace}, we infer that
\be\label{est:psi1}
\|\psi^{corr}_E\|_{H^m(\Om)}=O(\viscosite ^{1/3})\quad\forall m\leq 3.
\ee
Moreover, $\psi^{corr}_{E|\d \Omega\setminus \Gamma_E}=0$, and $\d_n \psi^{corr}_{E|\d \Omega}=0$ by definition of $\chi_0$. Notice  also that
$$
\begin{aligned}
(\psi^0_t+ \psi_E(s,\lambda_E z)+ \psi^{corr}_E)_{|[\sigma_i^+, \sigma_{i-1}^-]}=0,\\
\d_n (\psi^0_t+ \psi_E(s,\lambda_E z)+ \psi^{corr}_E)_{|\Gamma_E}=0.
\end{aligned}
$$

Therefore, at this point, we have restored both boundary conditions on the subdomain $[\sigma_i^+,\sigma_{i-1}^-] $ of $ [s_i, s_{i+1}]$, but not the condition on the trace  in the neighborhoods of the singularities. This is handled by the North and South boundary layers, which we now address.

In the next section, we set
\be\label{def:psib}
\psib:=\psi^0_t + \chi_0 ( \psi^{corr}_E+ \psi_E(s, \lambda_E(s) z)),
\ee
where $\chi_0$ is the macroscopic truncation defined in (\ref{chi0-def}), so that on the boundary
$$
\begin{aligned}
\psib_{|\d\Om\setminus \Gamma_E}= \psi^0_{t|\d\Om\setminus \Gamma_E},\\
\d_n\psib_{|\d\Om\setminus \Gamma_E}= \d_n\psi^0_{t|\d\Om\setminus \Gamma_E},\\
\psib_{|(s_i, s_{i+1})}= (1-\varphi_i^+ \varphi_{i+1}^-)\psi_{E|Z=0}(s),\\
\d_n \psib_{|\Gamma_E}=0,
\end{aligned}
$$
where $(s_i, s_{i+1})$ is a connected component of $\Gamma_E$.

\section{North and South boundary layers}
 \label{sec:NSboundary}

 In order to lift the boundary conditions both on the horizontal parts and on the East/West boundaries close to the points  $(x_i,y_i)$ for $1\leq i \leq k$, we then define North and South boundary layers.
 Without loss of generality, we focus on the case of South boundaries (the case of North boundaries can be deduced by a simple symmetry).

 Denote by $s_i$ and $s_{i-1}$ the curvilinear abscissa of the endpoints of the horizontal part to be considered, with $s_i=s_{i-1}$ if it is just an isolated point of the boundary with horizontal tangent.

 There are several point to be discussed in the present section:
\begin{itemize}
\item In the first paragraph, we give a precise definition of the interval on which energy is injected in the South boundary layer (i.e. we define the support of the functions $\Psi_0$ and $\Psi_1$ appearing in \eqref{eq:North});
\item The second paragraph is devoted to the derivation of regularity and moment estimates on $\psi_S$;
\item Eventually, we explain how we truncate $\psi_S$ beyond the support of $\Psi_0$ and $\Psi_1$.
\end{itemize}

 \subsection{Definition of the initial boundary value problem}$ $

\begin{figure}[h]
\includegraphics[height=8cm]{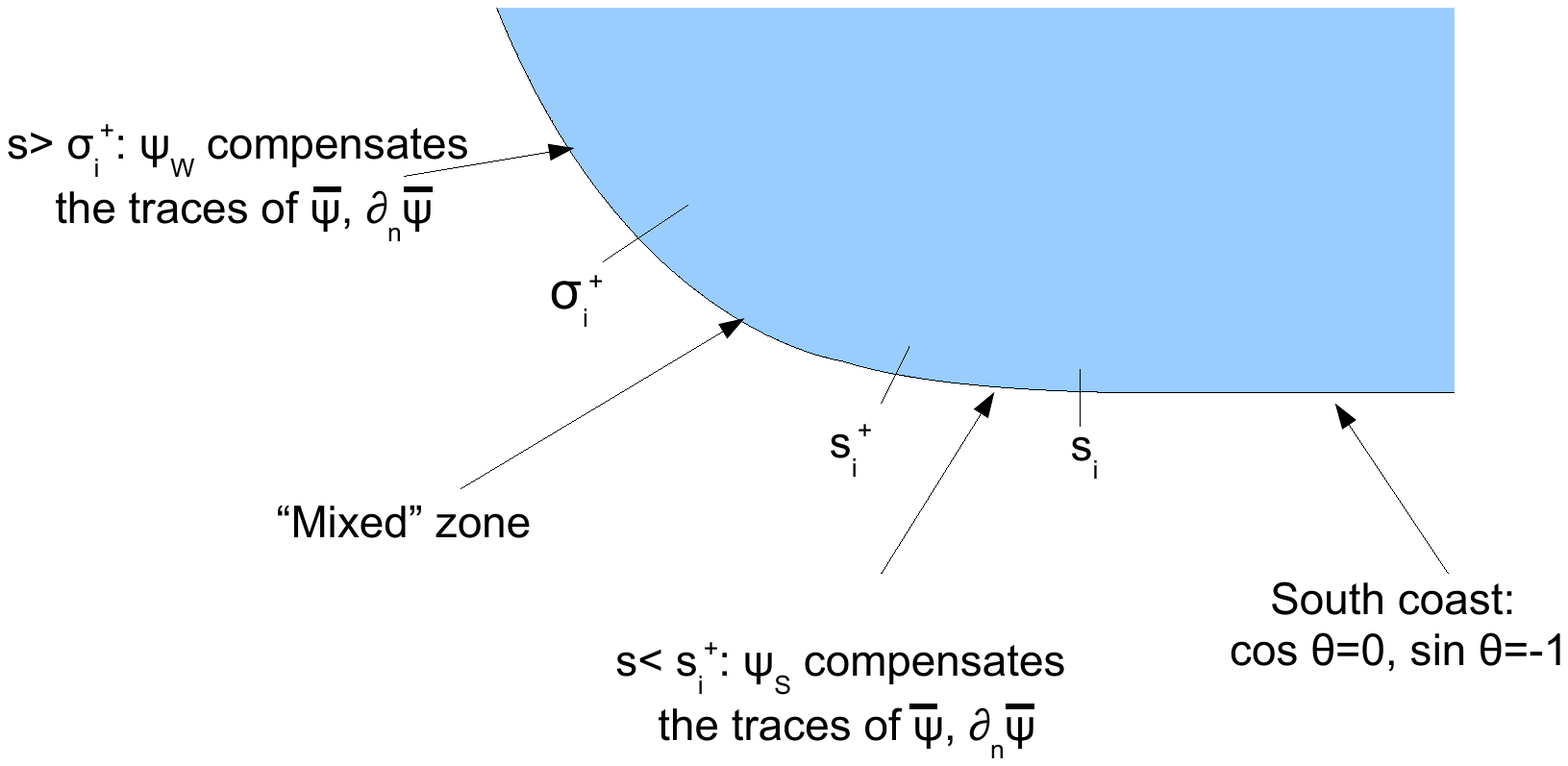}
\caption{The interplay between the North/South and East/West boundary layers}\label{fig:couches_NS}
\end{figure}

 We first define  the extremal points of the {\bf interval on which some energy is injected in the boundary layer}~:
 \begin{itemize}
 \item
 If $\cos \theta(s) >0$ for $s>s_i$ (East boundary), we need to lift the East boundary condition on a macroscopic neighborhood of $s_i$, namely $[s_i,\sigma_i^+]$ where $\sigma_i^+$ is defined by (\ref{Esigmaipm}).

\noindent Similarly,  if $\cos \theta(s) >0$ for $s<s_{i-1}$, we define
 $\sigma_{i-1}^- $ by (\ref{Esigmaipm}). Truncation functions $\varphi_i^+, \varphi_{i-1}^-$ are then defined by \eqref{Evarphii}.

  \item
 If $\cos \theta(s) <0$ for $s>s_i$ (West boundary), we only need to lift the West  boundary condition on a small neighbourhood of $s_i$, precisely when the West boundary layer becomes too singular (Recall that  the trace of $\psi^{0}_t$ is not  zero in general in the vicinity of $s_i$ on West boundaries). We therefore denote by $\sigma_i^+$ the arc-length beyond which no energy is injected in the boundary layer (see Figure \ref{fig:couches_NS}), and we expect that $\sigma_i^+>s_i^+>s_i$, and $|\sigma_i^+-s_i|\ll 1$.

 Note that, in order that the transition between the two types of boundary layers is not too singular, we have to choose $|\sigma_i^+-s_i^+|$ as large as possible. But in order that the transport coefficient $b$ in (\ref{eq:North}) remains bounded in $L^1$, we need that $\cos \theta(s)$ is small for $s \leq \sigma_i^+$.   We thus choose $\sigma_i^+$ satisfying  the condition
\begin{equation}
\label{Wsigmaipm}       \int_{s_i}^{\sigma_i^+}|\cos \theta(s)|\:ds=\viscosite ^{1/4}.
\end{equation}
We further define $\varphi_{i}^+ \in C^\infty(\R, [0,1])$ such that
  \begin{equation}
 \label{Wvarphii}
\varphi_i^+(s) = 0  \hbox{ if } s \leq s_i^+ \text{ and }\Supp (1-\varphi_i^+ )\subset (-\infty, \sigma_i^+ [.
\end{equation}

 Similarly,  if $\cos \theta(s) <0$ for $s<s_{i-1}$, we define $\sigma_{i-1}^-$ by
 $$\int_{\sigma_{i-1}^-}^{s_{i-1}}|\cos \theta(s)|\:ds=\viscosite ^{1/4},       $$
 and $\varphi_{i-1}^- \in C^\infty(\R, [0,1])$ such that
$$\Supp (1-\varphi_{i-1}^-)\subset ]\sigma_{i-1}^-, \infty) \text{ and } \varphi_{i-1}^-(s)  = 0  \hbox{ if } s \geq s_{i-1}^- .
$$

%
%

We will take
$$
\varphi_i^+ (s)= \Phi\left(s_i^+ + \frac{s-s_i^+}{\sigma_i^+ - s_i^+}\right) $$
where $\Phi \in \mathcal C^\infty(\R)$ is such that $\Phi(s)=0$ for $s\leq 0$, $\Phi(s)=1$ for $s\geq \frac12$, $0\leq \Phi\leq 1$ and all the derivatives of $\Phi$ are bounded.

\end{itemize}

\begin{Rmk}
There are two cases of South boundary layer terms which will not be considered in this section, since they both give rise to interface boundary layer terms:
\begin{itemize}
\item when $\cos \theta(s) <0$ for $s<s_{i-1}$ (West boundary on the right) and $\cos \theta(s) >0$ for $s>s_{i}$ (East boundary on the left);

\item when $\cos \theta(s)>0 $ for $s<s_{i-1}$ and for $s>s_i$ (East boundary on both sides), with $s_{i-1}<s_i$;
\end{itemize}

 We will indeed consider the corresponding South boundary  together with   the singular interface $\Sigma$ in the next section.
\end{Rmk}

As in Lemma \ref{sigma-s-est}, we can compute asymptotic developments for $\sigma_i^+$ and $\sigma_{i-1}^-$.

\begin{Lem}\label{sigma-s-est-bis}  The transition functions $\varphi_i^\pm$ have controlled variations:
\begin{itemize}
\item If $\cos \theta>0$ for $s\in (s_i, s_{i+1})$ (East coast), then there exists a constant $c>0$ such that
$$
|s_i-\sigma_i^+|\geq C,
$$
and therefore
$$
\|\d_s^k \varphi_i^+\|_\infty=O(1)\quad\forall k.
$$

\item If $\cos \theta<0$ for $s\in (s_i, s_{i+1})$ (West coast), then
$$
\sigma_i^+-s_i\sim C_i' \viscosite ^\frac{1}{4n+4}
$$
if (H2i) is satisfied in a neighbourhood on the right of $s_i$, and
$$
 \sigma_i^+- s_i=-\frac{4\alpha}{\ln \viscosite }\left(1-\frac{8\ln |\ln \viscosite |}{\ln \viscosite }+ O\left(|\ln \viscosite |^{-1}\right)\right)
$$
if (H2ii) is satisfied in a neighbourhood on the right of $s_i$.

In this case,
$$\begin{aligned}
\|\d_s^k \varphi_i^+\|_\infty \leq C \viscosite ^{-k/7},\\
\|\d_s^k \varphi_i^+\|_1\leq C\viscosite ^{-(k-1)_+/7}.
\end{aligned}
$$

\end{itemize}

\end{Lem}
\begin{proof}
The proof follows from the definition of $\sigma_i^+$ and $\varphi_i^+$ on East and West coasts.

On the West coast, if $\cos \theta$ vanishes algebraically near $s_i$, then $$s_i^+-s_i\sim \viscosite ^{\frac{1}{4n+3}}\hbox{ and } \sigma_i^+-s_i\sim \viscosite ^{\frac{1}{4n+4}}.$$
  On the other hand, if $\cos \theta$ vanishes exponentially near $s_i$, then $$s_i^+-s_i \sim \sigma_i^+-s_i\sim -4 \alpha/\ln \viscosite  \hbox{ and }\sigma_i^+-s_i^+=O(\ln (|\ln \viscosite |)/(\ln \viscosite )^2)\ll |s_i^+-s_i|.$$
Therefore  we always have$|s_i^+-s_i|\ll |\sigma_i^+-\sigma_i|$.

Furthermore, by construction, we have
\be\label{est:ds-varphi_i}
\begin{aligned}
\|\d_s^k\varphi_i^\pm\|_\infty=O(|s_i^\pm-\sigma_i^\pm|^{-k} + |s_i^\pm-s_i|^{-k})\text{ if }(s_i, s_{i\pm 1})\subset \Gamma_W,\\
\|\d_s^k\varphi_i^\pm\|_\infty=O(|s_i^\pm-\sigma_i^\pm|^{-k})\text{ if }(s_i, s_{i\pm 1})\subset \Gamma_E.
\end{aligned}
\ee
so that
$$
\|\d_s^k\varphi_i^+\|_\infty= O(\viscosite ^{-k/(4n+3)})\quad \text{(resp. } \|\d_s^k\varphi_i^+\|_\infty= O\left( \left(\frac{|\ln \viscosite |^2}{\ln |\ln \viscosite |}\right) )^k \right) \,.
$$
To derive $L^1$ estimates, we merely multiply the $L^\infty$ estimate
by the size of \linebreak $\Supp \d_s^k \varphi_i^+$.
\end{proof}

\bigskip
Let us then define the {\bf initial boundary value problem on $s \geq \sigma_{i-1}^-$, $Z\geq 0$}.

Denote
\begin{equation}
\label{gamma-i}
\gamma_{i-1,i} = (1-\varphi_i^+) (1-\varphi_{i-1}^-)\,.
\end{equation}

 The South boundary layer is therefore described by
\be\label{eq:South}
\begin{aligned}
\d_s \psi_S + \gamma (s)\d_Z \psi_S + \mu(s)\d_Z^4\psi_S=0,\quad s\geq \sigma_{i-1}^-,\ Z>0, \\
\psi_{S|Z=0}=\Psi_0,\quad\d_Z\psi_{S|Z=0}=\Psi_1
\end{aligned}
\ee
 where $\gamma(s)=-\viscosite ^{-1/4} \cos \theta(s)/\sin \theta(s)$, $\mu(s)=-1/\sin \theta(s)$, with boundary condition
\begin{equation}
\label{S-bc}
\begin{aligned}
 \Psi_0(s):= -\gamma_{i-1,i}\psib_{|\d \Omega}(s),\\
\Psi_1(s): =-\viscosite ^{1/4}\gamma_{i-1,i}\d_n\psib_{|\d \Omega}(s),
\end{aligned}
\end{equation}
and zero initial data  prescribed at $s = \sigma_{i-1}^-$. We recall that $\psib $ is defined by \eqref{def:psib}.

Note that, when the two types of boundary layers meet (that is on the support of $(\gamma_{i-1,i})'$), the width of the North/South boundary layer is much larger than the one of the East/West boundary layer. Indeed, the width of the North/South boundary layer is always $\viscosite ^{1/4}$, while that of the East/West boundary layer is at most, using Lemma \ref{sigma-s-est} and hypothesis (H2),
$$
\lambda_{E,W}^{-1}(s_i^+)= \left(\frac{\viscosite }{\cos \theta(s_i^+)}\right)^{1/3}\leq C\left\{
\begin{array}{ll}
\viscosite ^{\frac{n+1}{4n+3}}&\text{ in case (H2i)},\\
\viscosite ^{1/4}|\ln \viscosite |^{-3/2}&\text{ in case (H2ii)}.
\end{array}
\right.
$$
Therefore, it seems more accurate to talk about a {superposition} of the boundary layers, rather than a connection.

 \subsection{Estimates for $\psi_{N,S}$}$ $

 \begin{Lem}[Trace estimates]\label{trace-lem}
The functions $\Psi_0, \Psi_1$ defined by
(\ref{S-bc})
for $s\in (\sigma_{i-1}^-, s_{i+1})$ are such that
$$
\begin{aligned}
\|\d_s^k \Psi_0\|_{L^\infty} + \|\d_s^k \Psi_1\|_{L^\infty}\leq C\viscosite ^{-k/7}\quad \text{for }k\in \{0,\cdots, 3\},\\
\|\d_s \Psi_0\|_{L^1} + \|\d_s \Psi_1\|_{L^1}\leq C,\\
\|\d_s^2 \Psi_0\|_{L^1} + \|\d_s^2 \Psi_1\|_{L^1}\leq C\delta_y^{-1},\quad \|\d_s^3 \Psi_0\|_{L^1} + \|\d_s^3 \Psi_1\|_{L^1}\leq C\delta_y^{-1}\viscosite ^{-1/7}.
\end{aligned}
$$
Moreover,
$$
\begin{aligned}
\|\gamma \Psi_0\|_{L^1}+  \|\gamma \Psi_1\|_{L^1}=O(1),\\
\|\d_s(\gamma \Psi_0)\|_{L^1}+  \|\d_s(\gamma \Psi_1)\|_{L^1}=o(\viscosite ^{-1/6}),\\
\|\d_s^2(\gamma \Psi_0)\|_{L^1}+  \|\d_s^2(\gamma \Psi_1)\|_{L^1}=O(\viscosite ^{-2/7}).
\end{aligned}
$$

\end{Lem}

\begin{proof}
We recall that $\psib=\psi^0_t+ \chi_0( \psi^{corr}_E+ \psi_E(s, \lambda_E(s)z))$. Notice that thanks to assumption (H3), $\psi^0_t$ is continuous (and even $\mathcal C^2$) on the interval $(\sigma_{i-1}^-, \sigma_i^+)$. The estimates are slightly different depending on whether $(s_{i-2}, s_{i-1})$ and $(s_i, s_{i+1})$ are portions of $\Gamma_E$ or $\Gamma_W$. We focus for instance on the portion $(\sigma_{i-1}^-, s_{i-1})$, keeping in mind that the portion $(s_{i}, \sigma_i^+)$ is analogous.

\noindent
{\bf
Connection with West boundaries.}

\noindent
 {If $(\sigma_{i-1}^-, s_{i-1})\subset \Gamma_W$,} then by definition (\ref{Wsigmaipm}) of $\sigma_{i-1}^-$,
\be\label{in:ys}
\forall s\in (\sigma_{i-1}^-, s_{i-1}),\quad |y(s)-y_i|\leq \viscosite ^{1/4}\ll \delta_y.
\ee
Moreover, $\bar \psi = \psi^0_t$ on $(\sigma_{i-1}^-, s_i)$ in this case. Therefore Lemma \ref{trace-lem} is  a  consequence of the trace estimates on $\psi_t^0$ stated in Lemma \ref{lem:est-psi0}.

First, the $L^\infty$ bounds together with the definition of $\varphi_{i-1}^-$ imply that $\|\Psi_0\|_{L^\infty}=O(1)$.
Moreover, for $s\in (\sigma_{i-1}^-, s_i)$, using Lemma \ref{lem:est-psi0} and inequality \eqref{in:ys}, we have
$$
|\d_s  \psi^0_{t|\d\Om}(s)|\leq C \left( 1+ \frac{|\cos \theta(s)|}{\delta_y |\ln \delta_y|}\right),
$$
so that, using the definition \eqref{Esigmaipm} of $\sigma_i^\pm$ together with the definition \eqref{def:deltay} of $\delta_y$,
\begin{eqnarray*}
\|  \d_s \psi^0_{t|\d\Om}\|_{L^1 (\sigma_{i-1}^-, s_i)}&\leq &C\Big(1+ \viscosite ^{-1/4}|\ln \viscosite |^{-6/5}(\ln|\ln \viscosite |)^{\beta}\int_{\sigma_{i-1}^-}^{\sigma_i^+}|\cos \theta(s)|\:ds\Big)\\
&\leq & C(1+ |\ln \viscosite |^{-6/5}(\ln|\ln \viscosite |)^{\beta}).
\end{eqnarray*}
We infer that
$$
\|\Psi_0'\|_{L^1(\sigma_{i-1}^-, s_i)}\leq \| \psi^0_{t|\d\Om}\|_\infty \|\d_s \gamma_{i-1,i}\|_{L^1( \sigma_{i-1}^-, s_i)} + \|  \d_s  \psi^0_{t|\d\Om}\|_{L^1 (\sigma_{i-1}^-, s_i)}\leq C.
$$
In a similar fashion, Lemma \ref{sigma-s-est-bis} and assumption (H2) yield
\begin{eqnarray*}
\|\Psi_1\|_{L^\infty(\sigma_{i-1}^-, s_i)}&\leq & \viscosite ^{1/4} \| \d_n   \psi^0_{t|\d\Om}\|_{L^\infty(\sigma_{i-1}^-, s_i)}\\
&\leq & C \viscosite ^{1/4} \left(1+ \frac{1}{\delta_y |\ln \delta_y|}\right)\leq C,
\end{eqnarray*}
and eventually
\begin{eqnarray*}
  \|\Psi_1'\|_{L^1}&\leq & \viscosite ^{1/4}  \|\d_s \gamma_{i-1,i}\|_{L^1( \sigma_{i-1}^-, s_i)}  \| \d_n  \psi^0_{t|\d\Om}\|_{L^\infty(\sigma_{i-1}^-, s_i)}\\
&&+ C \viscosite ^{1/4}   \| \d_s \d_n  \psi^0_{t|\d\Om}\|_{L^1(\sigma_{i-1}^-, s_i)}\\
&\leq & C \left(1 + \viscosite ^{1/4} \int_{\sigma_{i-1}^-}^{s_i}\left(\frac{1}{\delta_y}+ \frac{|\cos \theta|}{\delta_y^2|\ln \delta_y|}\right)\right)\\
&\leq & C  \left(1 +\frac{\viscosite ^{1/4}}{\delta_y} + \frac{\viscosite ^{1/2}}{\delta_y^2|\ln \delta_y|}\right)\leq C.
\end{eqnarray*}

The higher order estimates are obtained in a similar way. We use the estimates of Lemma \ref{sigma-s-est-bis}, which lead to
$$
\| \d_s^k\gamma_{i-1,i}\|_\infty = O(\viscosite ^{-k/7}),\quad\| \d_s^k\gamma_{i-1,i}\|_1= O(\viscosite ^{-(k-1)_+/7});
$$
 In a similar fashion,
$$
\|\cos^k \theta\|_{L^\infty(\sigma_{i-1}^-, s_i)}= O(\viscosite ^{k/8}), \quad \|\cos^k \theta\|_{L^1(\sigma_{i-1}^-,s_i)}= O(\viscosite ^{(k+1)/8}).
$$
Lemma \ref{lem:est-psi0}  implies that for $s\in( \sigma_{i-1}^-, s_i)$
\begin{eqnarray*}
| \Psi_0''(s)| +  | \Psi_1''(s)| &\leq & C \left(\delta_y^{-1} \left(1+ \frac{\cos^2 \theta(s)}{\delta_y}\right) + \viscosite ^{-1/7}  \left(1+ \frac{|\cos \theta(s)|}{\delta_y}\right) \right)\\
&&+ C \left(|\d_s^2 \varphi_i^+| +  |\d_s^2 \varphi_{i-1}^-| \right),\\
| \d_s^3\Psi_0(s)| +  | \d_s^3\Psi_1(s)| &\leq & C \delta_y^{-1}\left(1 + \frac{|\cos \theta|}{\delta_y} + \frac{|\cos \theta|^3}{\delta_y^2}\right)\\
&&+C \viscosite ^{-1/7} \delta_y^{-1} \left(1+ \frac{\cos^2 \theta(s)}{\delta_y}\right) \\
&&+ C \viscosite ^{-2/7}  \left(1+ \frac{|\cos \theta(s)|}{\delta_y}\right) + C \left(|\d_s^3 \varphi_i^+| +  |\d_s^3 \varphi_{i-1}^-| \right),
\end{eqnarray*}
from which we easily infer the estimates of the Lemma.

The estimates on $\gamma\Psi_0$ and $\gamma \Psi_1$ are a consequence of the following estimates on $\gamma$:
$$
\begin{aligned}
\gamma\equiv 0 \text{ on } (s_{i-1}, s_i),\\
\| \gamma\|_{L^1(\sigma_{i-1}^-, s_{i-1})}=1
\end{aligned}
$$
by definition of $\sigma_{i-1}^-$ on West coasts, and
$$\begin{aligned}
\|\gamma \|_{L^\infty(\sigma_{i-1}^-, s_{i-1})}, \|\d_s \gamma\|_{L^1(\sigma_{i-1}^-, s_{i-1})}= O(\viscosite ^{-1/8}),\\
\| \d_s^2 \gamma \|_{L^1(\sigma_{i-1}^-, s_{i-1})}= O(\viscosite ^{-1/6}).
\end{aligned}
$$
Indeed
$$\|\d_s^k \gamma\|_{L^1(\sigma_{i-1}^-, s_{i-1})} \leq C (1+\viscosite ^{-1/4}\| \d_s^k \theta\|_{L^1(\sigma_{i-1}^-, s_{i-1})}) .$$
Using assumption (H2) together with the definition of $\sigma_{i-1}^-$, we infer the desired result; notice that the most singular case corresponds to $n=1$ in (H2i) for the estimates on $\gamma$ and $\d_s \gamma$, and to $n=2$ for the estimate on $\d_s^2 \gamma$.

\smallskip

\medskip
\noindent
{\bf Connection with East boundaries}

\noindent
 {If $(\sigma_{i-1}^-, s_{i-1})\subset \Gamma_E$,} the estimates are different in several regards:
\begin{itemize}
\item The function $\gamma_{i-1, i}$ has bounded derivatives;

\item Because of the truncation $\chi_\viscosite $, the traces of $\psi_t^0$ and $\d_n \psi_t^0$ are identically zero on the vicinity of $s_{i-1}$;

\item The trace of $\psi^0_t$ is zero on $(\sigma_{i-1}^-, s_{i-1})$.

\item The normal derivative of $\psib=\psi^0_t+ \psi_E + \psi_E^{corr}$ is identically zero on $(\sigma_{i-1}^-, s_{i-1})$ by definition of $\psi_E$, so that $\Psi_1=0$ on $(\sigma_{i-1}^-, s_{i-1})$;

\end{itemize}
Therefore it suffices to prove that $\gamma_{i-1,i}\psi_{E|Z=0}$ satisfies the desired estimates on the interval $(\sigma_{i-1}^-, s_{i-1})$. Notice that the estimates on $(s_{i-1}, s_i)$ can be treated with the same arguments as in the first case.

We use the estimates of Lemma \ref{lem:psiE-trace} and  assumption (H2), together with the formulas of Appendix B, in order to compute $y(s)$, $\cos \theta$ and $\theta'$ in terms of $s$. We infer that on $(\sigma_{i-1}^-, s_{i-1})$, $\Psi_0$ satisfies the following bounds:
\begin{itemize}
\item If (H2i) is satisfied in a neighbourhood on the left of $s_{i-1}$, then for $k=1,2,3$, $s\in (\sigma_{i-1}^-, s_{i-1})$,
\begin{eqnarray*}
|\Psi_0|&\leq & C \viscosite ^{1/3}\left(|s-s_{i-1}|^{-\frac{4n+3}{3} }|\ln (s_{i-1}-s)|^{-1} \mathbf{1}_{C \delta_y^{1/(n+1)}\leq |s-s_{i-1}|\leq 1/2} +1\right),\\
|\d_s^k \Psi_0|&\leq &  \viscosite ^{1/3} \delta_y^{-1}\left(|s-s_{i-1}|^{-\frac{n}{3}-k }|\ln(s_{i-1}-s)|^{-1} \mathbf{1}_{C \delta_y^{1/(n+1)}\leq |s-s_{i-1}|\leq 1/2} +1\right).
\end{eqnarray*}
We infer in particular that
$$
\|\Psi_0\|_{L^\infty(\sigma_{i-1}^-, s_{i-1})}\leq C \viscosite ^{1/3}\delta_y^{-\frac{4n+3}{3(n+1)}}|\ln \delta_y|^{-1}.
$$
Since $\frac{4n+3}{3(n+1)}<4/3$ for all $n\geq 1$, we infer that
$$
\|\Psi_0\|_{L^\infty(\sigma_{i-1}^-, s_{i-1})}\leq C \viscosite ^{1/3}\delta_y^{-4/3}|\ln \viscosite |^{-1} \leq C |\ln \viscosite |^{-1}.
$$
As for the other estimates, we have, for $k=1,2,3$,
$$
\begin{aligned}
\| \d_s^k \Psi_0\|_{L^\infty(\sigma_{i-1}^-, s_{i-1})}\leq C \viscosite ^{1/3}\delta_y^{-1} \delta_y^{-\frac{n+3k}{3(n+1)}}|\ln \viscosite |^{-1},\\
\| \d_s^k \Psi_0\|_{L^1(\sigma_{i-1}^-, s_{i-1})}\leq C \viscosite ^{1/3}\delta_y^{-1} \delta_y^{-\frac{n+3(k-1)}{3(n+1)}}|\ln \viscosite |^{-1}.
\end{aligned}
$$
It can be checked that the most singular estimates on the $L^\infty$ norm, and on the $L^1$ norm as soon as $k\geq 2$, correspond to $n=1$. We then obtain the following (non optimal) upper bounds
$$
\begin{aligned}
\|\d_s \Psi^0\|_{L^1(\sigma_{i-1}^-, s_{i-1})}\leq C |\ln \viscosite |^{-1},\\
\|\d_s \Psi^0\|_{L^\infty(\sigma_{i-1}^-, s_{i-1})},\ \|\d_s^2 \Psi^0\|_{L^1(\sigma_{i-1}^-, s_{i-1})}\leq C \viscosite ^{-1/12}|\ln \viscosite |^{-1},\\
\|\d_s^2 \Psi^0\|_{L^\infty(\sigma_{i-1}^-, s_{i-1})},\ \|\d_s^3 \Psi^0\|_{L^1(\sigma_{i-1}^-, s_{i-1})}\leq C \viscosite ^{-5/24}|\ln \viscosite |^{-1},\\
\|\d_s^3 \Psi^0\|_{L^\infty(\sigma_{i-1}^-, s_{i-1})}\leq C \viscosite ^{-1/3}|\ln \viscosite |^{-1}.
\end{aligned}
$$

\item If (H2ii) is satisfied in a neighbourhood of the left of $s_{i-1}$, then for $s\in (\sigma_{i-1}^-, s_{i-1})$ and $k=1,2,3$,
$$\begin{aligned}
|\Psi_0(s)|\leq C \viscosite ^{1/3}|s-s_{i-1}|^{-1} \exp\left(\frac{4\alpha}{3|s-s_{i-1}|}\right)\mathbf{1}_{|y(s)-y_{i-1}|\geq \delta_y/2},\\
|\d_s^k \Psi_0(s)|\leq C \viscosite ^{1/3} \delta_y^{-1} |s-s_{i-1}|^{-1-2k}\exp\left(\frac{\alpha}{3|s-s_{i-1}|}\right)\mathbf{1}_{|y(s)-y_{i-1}|\geq \delta_y/2}.
\end{aligned}
$$
Moreover, in case (H2ii), we recall that
$$
y(s)-y_{i-1}\sim C(s-s_{i-1})^2 \exp\left(-\frac{\alpha}{|s-s_{i-1}|}\right),
$$
so that for $k=0,1,2$
\begin{eqnarray*}
&&\|\d_s^k \Psi^0\|_{L^\infty(\sigma_{i-1}^-, s_{i-1})},\ \|\d_s^{k+1} \Psi^0\|_{L^1(\sigma_{i-1}^-, s_{i-1})}\\
&\leq& C \viscosite ^{1/3}\delta_y^{-4/3}|\ln \viscosite |^{2k-\frac{5}{3}}\leq C |\ln \viscosite |^{2k-\frac{5}{3}}
\end{eqnarray*}
and similarly
$$
\|\d_s^3 \Psi^0\|_{L^\infty(\sigma_{i-1}^-, s_{i-1})}\leq C|\ln \viscosite |^{13/3}.
$$
Therefore we retrieve the desired on the norm of $\Psi_0$. Concerning the estimates on $\gamma \Psi_0$, we use the following inequalities on $(\sigma_{i-1}^-, s_{i-1})$
\begin{eqnarray*}
|\gamma \Psi_0(s)|&\leq &C \viscosite ^{1/12}|\cos \theta|^{2/3}\ccM(y(s)),\\
|\d_s(\gamma \Psi_0(s))|&\leq & C \viscosite ^{1/12}\Big (\delta_y^{-1} |\cos \theta|^{5/3}+ |\theta'(s)||\cos \theta|^{-1/3}\Big)\ccM(y(s)) \\
|\d_s^2(\gamma \Psi_0(s))|&\leq & C \viscosite ^{1/12}\Big (\delta_y^{-1} |\cos \theta|^{2/3}|\theta'(s)|+ \frac{|\theta''(s)|}{|\cos \theta|^{-1/3}}\Big)\ccM(y(s)).
\end{eqnarray*}
\end{itemize}
Once again, we distinguish between  (H2i) and (H2ii), and we obtain in the worst situation
$$
\begin{aligned}
\|\gamma  \Psi_0\|_{L^1}\leq C |\ln \viscosite |^{-1},\\
\|\d_s(\gamma \Psi_0)\|_{L^1}\leq C\viscosite ^{1/12}\delta_y^{-1}\leq C \viscosite ^{-1/6}|\ln \viscosite |^{-1/5}|\ln |\ln \viscosite ||^{\beta},\\
\|\d_s^2(\gamma  \Psi_0)\|_{L^1}\leq C \viscosite ^{1/12} \delta_y^{-7/6}\leq C \viscosite ^{-5/24}.
\end{aligned}
$$

\end{proof}

\bigskip
By Lemma \ref{Cauchy-NS}, we then have the well-posedness of  equation \eqref{eq:South}, as well as suitable a  priori estimates for the solution $\psi_S$
\begin{equation}
\label{psiS-est1}
\begin{aligned}
\|\psi_S(s)\|_{L^2(\R_+)}&\leq C (|\Psi_0(s)| + |\Psi_1(s)|)\\
&+C\int_{\sigma_{i-1}^-}^s\left(|\Psi_0'(s')| + |\Psi_1'(s')| \right)\:ds'\\
&+C\int_{\sigma_{i-1}^-}^s(\mu(s')+|\gamma(s')| ) (|\Psi_0(s')| + |\Psi_1(s')|)\:ds'.
\end{aligned}
\end{equation}

In order to show that the South boundary layer is stable, we will need additional estimates on the boundary layer terms, namely regularity and moment estimates.

\begin{Lem}[Regularity and moment estimates] \label{reg-lem}
Let $\psi_S$ be the solution to (\ref{eq:South}) with boundary conditions (\ref{S-bc}), defined on
$$I_S:=\{s \geq \sigma_{i-1}^-\,/\, \forall s' \in [\sigma_{i-1}^-, s], \, |\sin \theta(s)| \geq \frac{1}{4}\}.$$
Then
\begin{equation}
\label{moment}
\begin{aligned}
\| (1+Z ^k) \psi_S (s,Z)\|_{L^\infty (I_S,L^2_Z)}  +\| (1+Z^k) \d_Z^2 \psi_S \|_{L^2(I_S\times \R_+)}=O(1),\\
\| (1+Z^k) \d_s  \psi_S\|_{L^\infty (I_S,L^2_Z)} +\| (1+Z^k) \d_Z^2\d_s \psi_S \|_{L^2(I_S\times \R_+)}=O(\delta_y^{-1}),\\
\| (1+Z^k) \d_s^2  \psi_S\|_{L^\infty (I_S,L^2_Z)} +\| (1+Z^k) \d_Z^2\d_s^2 \psi_S \|_{L^2(I_S\times \R_+)}=O(\viscosite ^{-1/2}|\ln \viscosite |^{-1}).
\end{aligned}
\end{equation}

\end{Lem}

\begin{Rmk}
The condition $|\sin \theta|\geq 1/4$ in the definition of $I_S$ is somewhat arbitrary: we just want to prevent $\sin \theta$ from touching zero, but obviously any fixed positive number will do.
\end{Rmk}

\begin{proof}
Such a priori estimates are obtained by combining precised energy estimates for the parabolic equation \eqref{eq:phi-sigma}, together with
  high order estimates on the traces.

\medskip
$\bullet$
Define as in the previous chapter
$$
\begin{aligned}
 S(s,Z)&:= -\Psi_0'(s)(Z+1)\exp(-Z)- \Psi_1'(s)Z\exp(-Z)\\
&+ \gamma(s) \Psi_0(s)Z\exp(-Z)-\gamma(s)\Psi_1(s)(1-Z)\exp(-Z)\\
&-\mu(s) \Psi_0(s)(Z-3)\exp(-Z) - \mu(s) \Psi_1(s)(Z-4)\exp(-Z).
\end{aligned}
$$
From the trace estimates in Lemma \ref{trace-lem}, we deduce that, for all $k\geq 0$,
 $$
\begin{aligned}
  \|Z^k S\|_{L^1}\leq C,\\
\|Z^k \d_s S\|_{L^1} \leq C \delta_y^{-1},\quad \|Z^k \d_s^2 S\|_{L^1} \leq C \delta_y^{-1}\viscosite ^{-1/7}.
\end{aligned}
$$
Notice also that by definition \eqref{Evarphii}, \eqref{Wvarphii} of $\varphi_{i}^\pm$, we have $\supp (1-\varphi_{i-1}^-)\subset ]\sigma_{i-1}^-,\infty[$, so that $\Psi_0$ and $\Psi_1$ are identically zero on a neighbourhood of $\sigma_{i-1}^-$. Therefore $\d_s^k\psi_{S|s=\sigma_{i-1}^-}=0$.

\medskip
$\bullet$
In order to derive estimates on $\d^k_s \psi_S $ for $k\geq 1$, we lift the boundary conditions to get
\be\label{eq:Southg}
\begin{aligned}
\d_s g + \gamma(s)\d_Z g + \mu(s)\d_Z^4g= S(s,Z),\quad s\in(0,T),\ Z>0, \\
g_{|s=0}=0,\\
g_{|Z=0}=0,\quad\d_Zg_{|Z=0}=0,
\end{aligned}
\ee
 then differentiate the equation and  proceed by induction on $k$. We start with
$$
\d_s(\d_s^k g) + \d_s^k (\gamma(s) \d_Z g) +\d_s^k (\mu (s) \d_Z^4 g) =\d_s^k S\,.$$
Integrating against $\d_s^k g$, we get
$$
\begin{aligned}
&\frac12 \| \d_s^k g(s) \|_{L^2_Z}  ^2 + \int_{\sigma_{i-1}^-}^s \mu(s') \| \d_Z^2 \d_s^k g(s') \|_{L^2_Z}^2 ds' \\
& \leq  \int_{\sigma_{i-1}^-}^s\|  \d_s^k S(s')\| _{L^2_Z}  \| \d_s^k g(s') \|_{L^2_Z} ds'\\
&+ \sum_{j=0}^{k-1} C^j_k \int_{\sigma_{i-1}^-}^s  | \d_s^{k-j} \gamma(s') | \| \d_Z\d_s^j g(s') \|_{L^2_Z} \| \d_s^k g(s') \|_{L^2_Z} ds' \\
&+ \sum_{j=0}^{k-1} C^j_k \int_{\sigma_{i-1}^-}^s  \| \d_Z^2 \d_s^k g(s') \| _{L^2_Z} | \d_s^{k-j}\mu(s') | \| \d_Z^2\d_s^j g(s') \|_{L^2_Z} ds' \,,
\end{aligned}
$$
and therefore, since $\mu$ is bounded from below,
\begin{eqnarray*}
& \| \d_s^k g\|_{L^\infty((\sigma_{i-1}^-, s),L^2_Z)}^2 +  \| \d_s^k \d_Z^2 g\|_{L^2((\sigma_{i-1}^-, s),L^2_Z)}^2\\
&\leq  C \|\d_s^k S \|_{L^1((\sigma_{i-1}^-, s),L^2_Z)} \| \d_s^k g\|_{L^\infty((\sigma_{i-1}^-, s),L^2_Z)}\\
&+C \sum_{j=0}^{k-1}\|\d_s^{k-j} \gamma \|_{L^{4/3}(\sigma_{i-1}^-, s)} \|\d_s^jg\|_{L^\infty_s (L^2_Z)}^{1/2}  \| \d_s^j \d_Z^2 g\|_{L^2_s (L^2_Z)}^{1/2} \| \d_s^k g\|_{L^\infty_s (L^2_Z)}\\
&+C  \sum_{j=0}^{k-1}\|\d_s^{k-j} \mu\|_{L^{\infty}(\sigma_{i-1}^-, s)} \| \d_s^j \d_Z^2 g\|_{L^2_s (L^2_Z)} \| \d_s^k \d_Z^2 g\|_{L^2_s (L^2_Z)}.
\end{eqnarray*}
Eventually, using the Cauchy-Schwarz inequality, we obtain for $k\leq 2$
\begin{eqnarray}
\label{in:dskg} & \| \d_s^k g\|_{L^\infty((\sigma_{i-1}^-, s),L^2_Z)}^2 +  \| \d_s^k \d_Z^2 g\|_{L^2((\sigma_{i-1}^-, s),L^2_Z)}^2\\
&\leq  C \|\d_s^k S \|_{L^1((\sigma_{i-1}^-, s),L^2_Z)}^2\nonumber \\
&+ C \sum_{j=0}^{k-1}\|\d_s^{k-j} \gamma \|_{L^{4/3}(\sigma_{i-1}^-, s)}^2 \|\d_s^jg\|_{L^\infty((\sigma_{i-1}^-, s),L^2_Z)}  \| \d_s^j \d_Z^2 g\|_{L^2((\sigma_{i-1}^-, s)\times\R_+)}\nonumber\\
&+ C \sum_{j=0}^{k-1}\|\d_s^{k-j} \mu\|^2_{L^{\infty}(\sigma_{i-1}^-,s)} \| \d_s^j \d_Z^2 g\|_{L^2((\sigma_{i-1}^-, s)\times\R_+)}^2.\nonumber
\end{eqnarray}

Next recall that $\mu$ is bounded, and that all its derivatives are of order $O(1)$
$$ \mu \geq \mu_0>0,\qquad \| \d_s ^j \mu \| _\infty \leq C.$$
In order to bound the term involving $\gamma$, we have to distinguish between the cases when $(s_{i-2}, s_{i-1})\subset \Gamma_E$ or $\Gamma_W$:

\begin{itemize}
\item If $(s_{i-2}, s_{i-1})\subset \Gamma_W$, then $\gamma= \viscosite ^{-1/4} \cos \theta \mu$ satisfies
$$
\begin{aligned}
\|\d_s\gamma\|_{L^{4/3}(\sigma_{i-1}^-, \sigma_i^+)}=O(\viscosite ^{- \frac{1}{4}+ \frac{3}{4}\frac{1}{8}})= O(\viscosite ^{-5/32}),\\
\| \d_s^2 \gamma \|_{L^{4/3}(\sigma_{i-1}^-, \sigma_i^+)}=O(\viscosite ^{- \frac{1}{4}+ \frac{3}{4}\frac{1}{12}})= O(\viscosite ^{-3/16}).
\end{aligned}
$$

Notice that in this case, we also have $(s_i, s_{i+1})\subset \Gamma_W$ since we have excluded for the moment the case leading to a discontinuity line $\Sigma$.

Inequality \eqref{in:dskg} therefore yields, for $k=0$,
$$
\| g\|_{L^\infty_s (L^2_Z)} + \| \d_Z^2 g \|_{L^2_{s,Z}} \leq  C \| S\|_{L^1_s(L^2_Z)}\leq C,$$
then, for $k=1$,
$$
\| \d_s g\|_{L^\infty_s (L^2_Z)} + \| \d_Z^2\d_s g \|_{L^2_{s,Z}} \leq C \delta_y^{-1},
$$
and, for $k=2$,
$$
\| \d_s^2 g\|_{L^\infty_s (L^2_Z)} + \| \d_Z^2\d_s^2 g \|_{L^2_{s,Z}}
\leq C \delta_y^{-1} \viscosite ^{-3/16}=o(\viscosite ^{-7/16}).
$$

\item If $(s_{i-2}, s_{i-1})\subset \Gamma_E$, then for $s\in (\sigma_{i-1}^-, s_{i-1})$
$$
\|\d_s^{j} \gamma \|_{L^{4/3}(\sigma_{i-1}^-, s)}=O(\viscosite ^{-1/4})
$$
for all $j$ and for all $s>\sigma_{i-1}^-$. In this case, we refer to the proof of Lemma \ref{trace-lem} and we state that
$$
\begin{aligned}
\|S\|_{L^1((\sigma_{i-1}^-, s_{i-1}), L^2_Z)}= O(|\ln\viscosite |^{-1}),\\
\|\d_s S\|_{L^1((\sigma_{i-1}^-, s_{i-1}), L^2_Z)}= o(\viscosite ^{-1/6}),\\
\|\d_s^2 S\|_{L^1((\sigma_{i-1}^-, s_{i-1}), L^2_Z)}= o(\viscosite ^{-5/24}).
\end{aligned}
$$

Inequality \eqref{in:dskg} therefore yields, for $k=0$,
$$
\| g\|_{L^\infty ((\sigma_{i-1}^-, s_{i-1}),L^2_Z)} + \| \d_Z^2 g \|_{L^2 ((\sigma_{i-1}^-, s_{i-1}),L^2_Z)} \leq  C |\ln \viscosite |^{-1},$$
then, for $k=1$,
$$
\| \d_s g\|_{L^\infty ((\sigma_{i-1}^-, s_{i-1}),L^2_Z)} + \| \d_Z^2\d_s g \|_{L^2 ((\sigma_{i-1}^-, s_{i-1}),L^2_Z)} \leq C \viscosite ^{-1/4} |\ln \viscosite |^{-1},
$$
and, for $k=2$,
$$
\| \d_s^2 g\|_{L^\infty ((\sigma_{i-1}^-, s_{i-1}),L^2_Z)} + \| \d_Z^2\d_s^2 g \|_{L^2 ((\sigma_{i-1}^-, s_{i-1}),L^2_Z)}
\leq C \viscosite ^{-1/2} |\ln \viscosite |^{-1}.
$$
There remains to derive bounds for $s>s_{i-1}$. There are essentially two cases:
\begin{itemize}
\item either $(s_i, s_{i+1})\subset \Gamma_W$, and in this case the estimates on $(s_{i-1}, \sigma_i^+)$ are similar to the ones of the case $(s_{i-2}, s_{i-1})\subset \Gamma_W$. The only change lies in the initial data: we merely replace
$$
\d_s^k g_{|s=\sigma_{i-1}^-}=0
$$
by
$$
\|\d_s^k g_{|s=s_{i-1}}\|_{L^2_Z}= O(\viscosite ^{-k/4}|\ln \viscosite |^{-1})\quad \text{for }k=0,1,2.
$$
We infer that
$$\begin{aligned}
\| g \|_{L^\infty_s( L^2_Z)} +\|g \|_{L^2_{s,Z}} = O(1),\\
\| \d_s g \|_{L^\infty_s( L^2_Z)} +\|\d_s g \|_{L^2_{s,Z}} = O(\delta_y^{-1}),\\
\| \d_s^2 g \|_{L^\infty_s( L^2_Z)} +\|\d_s^2 g \|_{L^2_{s,Z}} = O(\viscosite ^{-1/2}|\ln \viscosite |^{-1}).
\end{aligned}
$$

\item or $(s_i, s_{i+1})\subset \Gamma_E$. Notice that in this case, $s_i=s_{i-1}$ (i.e. there is no flat horizontal boundary), otherwise the corresponding piece of flat boundary would be the beginning of a $\Sigma$ layer, which will be treated in the next section. Whence we can treat the interval $(s_i, \sigma_i^+)$ exactly as $(\sigma_{i-1}^-, s_{i-1})$, and we obtain
$$\begin{aligned}
\| g \|_{L^\infty_s( L^2_Z)} +\|g \|_{L^2_{s,Z}} = O(|\ln \viscosite |^{-1}),\\
\| \d_s^k g \|_{L^\infty_s( L^2_Z)} +\|\d_s^k g \|_{L^2_{s,Z}} = O(\viscosite ^{-k/4}|\ln \viscosite |^{-1})\text{ for }k=1,2.
\end{aligned}
$$

\end{itemize}

\end{itemize}

$\bullet$
We get moment estimates in a similar fashion
$$
\d_s(Z^k g) + Z^k (\gamma(s) \d_Z g) +Z^k (\mu (s) \d_Z^4 g) =Z^k S\,.$$
Integrating against $Z^k g$, we get
$$
\begin{aligned}
\frac12 \| Z^k g(s) \|_{L^2_Z}  ^2 &+ \int_{\sigma_{i-1}^-}^s \mu(s') \| Z^k \d_Z^2  g(s') \|_{L^2_Z}^2 ds' \\
& \leq  \int_{\sigma_{i-1}^-}^s\|  Z^k S(s')\| _{L^2_Z} \| \| Z^k g(s') \|_{L^2_Z} ds'\\
&+ 2k \int_{\sigma_{i-1}^-}^s  |\gamma(s') | \| Z^k g(s') \|_{L^2_Z} \| Z^{k-1} g(s')  \|_{L^2_Z} ds' \\
&+ 2k \int_{\sigma_{i-1}^-}^s  | \mu(s') | \| Z^k \d_Z^2  g(s') \| _{L^2_Z}  \| Z^{k-1} \d_Z g(s') \|_{L^2_Z} ds' \\
&+ k(k-1)\int_{\sigma_{i-1}^-}^s  | \mu(s') | \| Z^k \d_Z^2  g(s') \| _{L^2_Z}  \| Z^{k-2}  g(s') \|_{L^2_Z} ds'  \,,
\end{aligned}
$$
A simple recursion then leads to the expected uniform bound for the moments of $\psi_S$ and $\d_{Z}^2 \psi_S$.

$\bullet$ Combining both arguments we obtain the moment estimates for higher order derivatives.

\end{proof}

 \subsection{Extinction and truncation} $ $
 \label{ssec:extinction}

 We  emphasize that $\psi_{S|s=\sigma_i^+}$ is not small in general: there is still some energy in the South boundary layer, even though the boundary conditions are absorbed by the West boundary layer only for $s\geq \sigma_i^+$. We deal with this somewhat unexpected phenomenon by propagating the boundary layer term $\psi_S$ beyond $\sigma_i^+$. We prove that when $\cos \theta $ becomes bounded away from zero, $\psi_S$ vanishes. As before, there are two cases, depending on whether $(s_i, s_{i+1})$ is a part of $\Gamma_W$ or $\Gamma_E$:

 $ \bullet$ First case: $(s_i, s_{i+1})\subset \Gamma_E$: looking closely at the proof of Lemma \ref{reg-lem}, we see that in this case,
\be\label{est:psi-S-Gamma-E}
 \|\psi_S\|_{L^\infty((s_i, \bar s), L^2_Z)}= O(|\ln \viscosite |^{-1})
 \ee
for all fixed  $\bar s>\sigma_{i}^+$ (independent of $\viscosite $) such that $|\sin \theta|\geq 1/2$ on $(s_{i}, \bar s)$. Therefore $\psi_S$ is small on the West boundary layer, on intervals of size one.

$\bullet $ Second case: $(s_i, s_{i+1})\subset \Gamma_W$: we use the following result:

\begin{Lem}[Extinction of the South boundary layer]
Let $\psi_S$ be the solution to (\ref{eq:South}) with boundary conditions (\ref{S-bc}), defined for $s\geq s_i$ as long as $\sin \theta (s)<0$.
Then
$$
\int_{s_i}^s  \viscosite ^{-1/4}|\cos \theta(s')|\|\psi_S(s')\|_{L^2}^2 \:ds'\leq C.
$$

\label{lem:extinction}
\end{Lem}

Before proving Lemma \ref{lem:extinction}, we now explain how we truncate the South boundary layer on $(s_i, s_{i+1})$. We introduce a function $\tilde \gamma_{i-1,i}\in \mathcal C^\infty(\d\Om)$ by
\be\label{def:tilde-gamma-i}\begin{aligned}
\supp \tilde\gamma_{i-1,i}\subset \{s \in  ]s_{i-2}, s_{i+1}[ \,/\, \sin \theta(s) <0\},\\
\tilde\gamma_{i-1,i}\equiv 1 \text{ on a neighbourhood of } \Supp(1-\varphi_i^+)\cap \Supp (1-\varphi_{i-1}^-).
  \end{aligned}
\ee
Notice than we can always assume that $\tilde \gamma_{i-1,i}$ has bounded derivatives.

The South boundary layer  on $(s_{i-2}, s_{i+1})$ is then defined by $$\tilde \gamma_{i-1,i}(s) \psi_S(s, z \viscosite ^{-1/4}).$$
Using either Lemma \ref{lem:extinction} or estimation \eqref{est:psi-S-Gamma-E}, we infer that $\|\psi_S\|_{L^2_{s,Z}}= O(|\ln \viscosite |^{-1})$ on the support of $\tilde \gamma_{i-1,i}'$. We will prove in the next chapter that this is enough to ensure that the error terms generated by $\tilde \gamma_{i-1,i}$ are admissible in the sense of Definition \ref{def:admissible}.

\begin{proof}[Proof of Lemma \ref{lem:extinction}]

By  Lemma \ref{reg-lem}, for every $s\in (s_i, s_{i+1})$ such that $\sin \theta(s) <0$,
$$
\|(1+Z) \psi_S\|_{L^2_Z}=O(1),
$$
and
$$
 \int_{\sigma_{i-1}^-}^{s} \left(|\Psi_0'(s)| + |\Psi_1'(s)| + (\mu(s)+|\gamma (s)| ) (|\Psi_0(s)| + |\Psi_1(s)|)\right)\:ds= O(1).
$$
Moreover, for $s\geq \sigma_i^+$, $\Psi_0(s)=\Psi_1(s)=0$. Integrating \eqref{eq:South} against $Z\psi_S$, we infer that there exists a constant $C>0$, which does not depend on $\viscosite $, such that for all $s\geq s_i$ such that $\sin \theta(s) <0$
$$\int_0^\infty Z\psi_S^2(s)\:dZ + \int_{s_i}^s  \viscosite ^{-1/4}|(\tan \theta(s'))^{-1}\|\psi_S(s')\|_{L^2}^2 \:ds'\leq C.
$$
which gives the desired estimate.
\end{proof}


\section{The interface layer}
 \label{sec:Sigmaboundary}

  Let us now focus on the interface $\Sigma = \d \Omega^+ \cap \d \Omega^-$. We recall that, by construction, $\psi_t^0$ is discontinuous across $\Sigma$ so that
  \begin{itemize}
  \item we need to lift this discontinuity to get a smooth ($H^2$) approximate solution,
  \item we then have to define a corrector (in the form of a singular layer term) to remove the energy created by the lifting term.
  \end{itemize}

  \subsection{The lifting term $\psil$}$ $

\label{ssec:psil}
We first construct  the term $\psil$ which lifts the discontinuity across $\Sigma$. The construction has already been explained in Chapter \ref{chap:multiscale}.  Since $\sigma_i^\pm$ is now defined precisely by \eqref{Wsigmaipm} for West boundaries and \eqref{Esigmaipm} for East boundaries, we can define $\sigma_{in}$ in the following way
\be\label{def:sigma-}\sigma_{in}:=
\left\{
\begin{array}{ll}
\sigma_1^- & \text{ if } I_1 =\{s_1\},\\
\sigma_0^- &\text{ if } I_1 =[s_0, s_1]\text{ and }(\sigma_0^-, s_0)\subset \Gamma_W,\\
s_0&\text{ if } I_1 =[s_0, s_1]\text{ and }(\sigma_0^-, s_0)\subset \Gamma_E.
\end{array}
\right.
\ee
Notice that if $I_1=\{s_1\}$, then we always have $(\sigma_1^-,s_1 )\subset \Gamma_W$ (an isolated point of cancellation of $\cos \theta$ on the East boundary does not create a discontinuity line). The reason why we choose to take $s_0$ instead of $\sigma_0^-$ as an initial point in the last line is essentially technical:  we have chosen to write the equation on the discontinuity type boundary layer in cartesian coordinates, rather than curvilinear ones (see equation \eqref{eq:psisig}). This is legitimate as long as $\cos \theta\ll 1$, so that the normal vector coincides almost with $-e_y$, while the tangent vector coincides with $-e_x$. However, by definition of $\sigma_i^\pm$ on East boundaries, $\cos \theta$ takes ``large" values on $(\sigma_0^-, s_0)$, and therefore the present construction cannot be used. We by-pass this problem by considering separately the East boundary for $s<s_0$ and the discontinuity boundary layer (which will be justified in Remark \ref{asymptotic-rmk}).

We thus define a truncation $\varphi$ by
$$
\varphi:=\left\{
\begin{array}{ll}
\varphi_1^- & \text{ if } I_1 =\{s_1\},\\
\varphi_0^- &\text{ if } I_1 =[s_0, s_1]\text{ and }(\sigma_0^-, s_0)\subset \Gamma_W,\\
0&\text{ if } I_1 =[s_0, s_1]\text{ and }(\sigma_0^-, s_0)\subset \Gamma_E.
\end{array}
\right.
$$
Notice that in the last case, $\psi^0_t$ vanishes identically in a neighbourhood of $(x_0,y_0)$, and thus there is no need for a truncation in $s$.

We then define as in  \eqref{def:psil}
$$
\psil(x,y)=\chi^+ \chi\left(\frac{y-y_1}{\viscosite ^{1/4}}\right)\left[a(x) + b(x) (y-y_1)\right],
$$\label{psil}
where $\chi$ is a truncation function as in \eqref{chi} and $a$ and $b$ are defined by
\eqref{def:ab1} and \eqref{def:ab2}.

Eventually, in order to simplify the present analysis, we  truncate the function $\tau$ in the vicinity of $\inf I_1$ and $\sup I_1$, even in the case when these points do not belong to $I_+$. Note that this does not change anything to the estimates of Lemma \ref{lem:trunc}.

 \subsection{The interior singular layer $\psi^\Sigma$}

The lifting term $\psil$ has introduced a remainder  $\dtl$  in the approximate Munk equation (see Lemma \ref{lem:eq-psil}), which we have to treat as a source term.
As suggested in the previous chapter, we therefore define $\psisig$ as the solution of  equation  \eqref{eq:psisig}
$$
\begin{aligned}
\d_x \psisig - \d_Y^4 \psisig=-\dtl,\quad x<x(\sigma_{in}), \quad Y>Y_-(x)\\
\psisig_{|x=x(\sigma_{in})}=0,\quad
\psisig_{Y=Y_-(x)}= \d_n  \psisig_{Y=Y_-(x)}=0.
\end{aligned}
$$
The function $Y_-$ \label{Y-}   is any $C^4$ function satisfying the following conditions:
$$
\begin{aligned}
	Y_-(x)=\frac{y(x)-y_1}{\nu^{1/4}}\quad \text{for } x \in (x(\varsigma_1^+), x(\sigma^-)),\\
Y_-(x)\equiv\bar Y_->\delta_y\nu^{-1/4}\quad \text{for }x<x(t_1^+),\\
\forall k\leq 4, \forall p\geq 1,\ \exists C>0,\  \quad\| Y_-\|_{W^{k,p}(\bar x, x(\sigma^-))}\leq C \| Y_-\|_{W^{k,p}(x(\varsigma_1^+), x(\sigma^-))},
\end{aligned}
$$
where
\begin{equation}
\label{t1+}
 \varsigma_1^+= \inf\{s>s_1,\quad |y(s) -y(s_1) |=\delta_y\},\quad t_1^+ = \inf\{s>s_1,\quad |y(s) -y(s_1) |=2\delta_y\},
\end{equation}
 $y(x)$ is the ordinate of the point of $\d\Om$ with abscissa $x$ in the vicinity of $(x_1, y_1)$ and of $(x(s_0),y(s_0))$, and  $\bar x$ is an arbitrary fixed point such that the vertical line $x=\bar x$ does not intersect $\bar \Om$.

%
%

\begin{figure}
 \includegraphics[width=\textwidth]{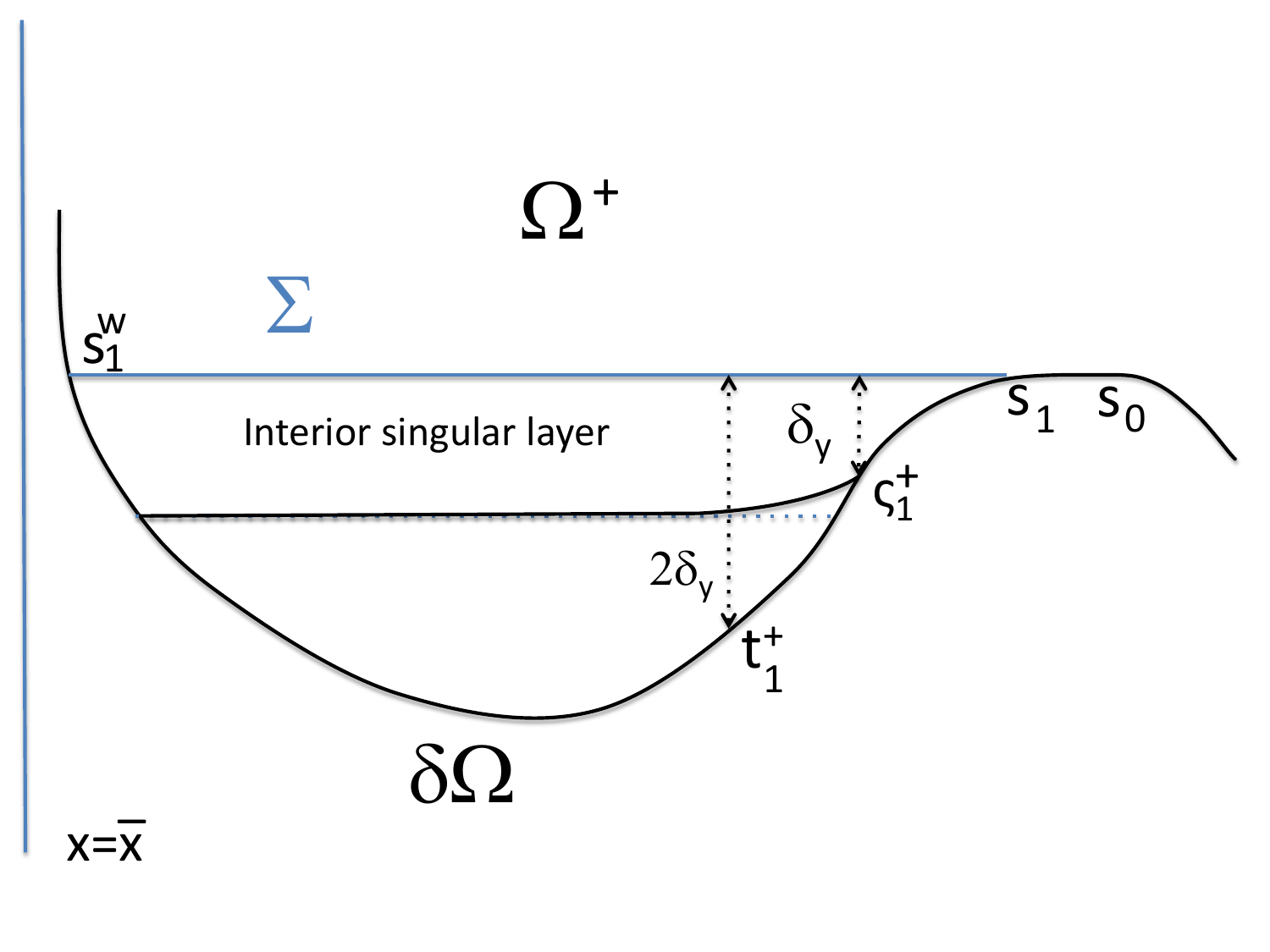}
\caption{The interior singular layer}
\end{figure}

Then the solution $\psisig$ of equation \eqref{eq:psisig} can be expressed as
$$
\psisig(x,Y)=\phi^\Sigma(x,Y-Y_-(x)),
$$
where $\phi^\Sigma$ is the solution of the backward equation
\be\label{eq:phi-sigma}
\left\{\begin{array}{l}
\d_x \phi^\Sigma -Y_-'(x)\d_Z \phi^\Sigma - \d_Z^4\phi^\Sigma=\\\quad= - \dtl(x, Y_-(x) + Z), \quad x<x(\sigma_{in}),\quad Z>0,\\
\phi^\Sigma_{|x=x(\sigma_{in})}=0,\quad
\phi^\Sigma_{|Z=0}=\d_Z\phi^\Sigma_{|Z=0}=0.
\end{array}\right.
\ee
which is precisely an equation of the type \eqref{eq:Northg}. In the equation above, we still denote the rescaled boundary layer coordinate by $Z$, in order to be consistent with \eqref{eq:Northg}. However, in the present context, $Z$ is just a scaled cartesian coordinate.

Notice that $\psisig$ and $\phi^\Sigma$ depend in fact on $\viscosite $ through the function $Y_-$ and the coefficients $a$ and $b$. The asymptotic value of $\psisig$ as $\viscosite $ vanishes will be given later on.

This indeed requires a refined version of Lemma \ref{Cauchy-NS} giving high order \textit{a priori }estimates on the solutions to the singular layer equation (\ref{eq:psisig}), as well as Lemma \ref{lem:est-A-B} in Appendix C  controlling the sizes of $a$, $b$ and their derivatives.

%


\begin{Lem}
Let $\phi^\Sigma$ be the solution to (\ref{eq:phi-sigma}), and $\psi^\Sigma (x,Y)=\phi^\Sigma(x,Y-Y_-(x))$.

Denote by $(x_1^W, y_1)= (x(s^W_1),y(s^W_1))$ the projection of $(x_1,y_1)$ on  the West boundary. \label{projW-def}

Then, for all $\bar x <x_1^W$,
$$
\begin{aligned}
\|\psi^\Sigma\|_{ L^{\infty}((\bar x, x(\sigma_{in})), L^2(Y_-,\infty))}=O(1),\\
\|\d_x\psi^\Sigma \|_{ L^{4}((\bar x, x(\sigma_{in})), L^2(Y_-,\infty))}=O(\viscosite ^{-1/19}),\\
\|\d_x^2\psi^\Sigma \|_{ L^{2}((\bar x, x(\sigma_{in})), L^2(Y_-,\infty))}=O(\viscosite ^{-2/19}).
\end{aligned}
$$
Furthermore, for all integers $k,l$ with $k+l\leq 3$,
$$
\| \d_x^k \d_Z^l \psi^\Sigma\|_{L^\infty(\bar x, x(t_1^+))} \leq C \viscosite ^{-\frac{1+2l + 8k}{8\times 19}},
$$
where $t_1^+$ is defined by \eqref{t1+}.
\label{lem:est-psi-sigma}
\end{Lem}

\begin{Rmk}\label{rmk:sigma-west}
The function $\psi^\Sigma$ is defined beyond the West end of $\{ (x,y) \in \Omega) \,/\, y = \bar Y_-\}$: indeed, for some values of $s$ in a neighbourhood of $s^W_1$, $x(s)$ may be smaller than $x_1^W$ (notice that this is the case as soon as $\sin \theta(s^W_1)\neq 0$, i.e. as soon as the tangent to $\d\Om$ at $s=s^W_1$ is not vertical). Thus the trace $\psi^\Sigma_{|\Gamma_W}$ involves values of $\psi^\Sigma$ at points $x<x^W_1$.
\end{Rmk}

\begin{proof}
Such a priori estimates are obtained by combining precised energy estimates for the parabolic equation \eqref{eq:phi-sigma}, together with
a refined description of the geometry $Y_-$, and high order estimates on the source terms.

\medskip \noindent
\textbf{$\bullet$ Estimates on $\phi^\Sigma$.}

 Denoting by $S(x,Z)= -\dtl( x, Y_-(x) +Z)$ the right-hand side of \eqref{eq:phi-sigma}, we recall that
\begin{eqnarray*}
 S(x,Z)&=&- \left(1- \mathbf{1}_{x<x_1,Z + Y_-(x)<0}\right) \chi(Z + Y_-(x))\left[a'(x) + \viscosite ^{1/4}b'(x) (Z+ Y_-(x))\right]\\
&+& \left(1- \mathbf{1}_{x<x_1,Z + Y_-(x)<0}\right) \chi^{(4)} (Z + Y_-(x))\left[a(x) + \viscosite ^{1/4}b(x) (Z+ Y_-(x))\right]\\
&+& 4 \left(1- \mathbf{1}_{x<x_1,Z + Y_-(x)<0}\right) \chi^{(3)} (Z + Y_-(x))\viscosite ^{1/4} b(x).
\end{eqnarray*}
Notice that the indicator function does not create any singularity on the domain of integration $x\in (\bar x, x(\sigma_{in}))\times (0,\infty)$: indeed, $\chi^{(3)}$ and $\chi^{(4)}$ are identically zero in a neighbourhood of zero, so that $\mathbf 1_{Y<0} \chi^{(k)}(Y)$ is in fact a $\mathcal C^\infty$ function for $k=3,4$.
Moreover,  for $x<x_1$, $a$ and $b$ are constant, so that
$$
 \mathbf{1}_{x<x_1} \mathbf 1_{Z + Y_-(x)<0}\, \chi(Z + Y_-(x))\left(a'(x) + \viscosite ^{1/4}b'(x) (Z+ Y_-(x))\right)\equiv 0.
$$
Eventually, for $x=x_1$, $Y_-=0$, so that $\mathbf 1_{Z>0}\mathbf 1_{Z + Y_-(x_1)<0}=0$.

It is easily checked that
$$
\|S(x,\cdot)\|_{L^2(\R_+)}\leq C(|a|+ |a'|+ \viscosite ^{1/4} |b|+ \viscosite ^{1/4}|b'|),
$$
so that, using Lemma \ref{lem:est-A-B} in  Appendix C,
$$
\|S\|_{L^1((\bar x, x(\sigma_{in})), L^2(\R_+))}\leq C.
$$
According to Lemma \ref{Cauchy-NS}, we then have
$$
\|\phi^\Sigma\|_{L^{\infty}(\bar x, x(\sigma_{in})), L^2(\bR_+))} + \|\d_Z^2\phi^\Sigma\|_{L^2((\bar x, x(\sigma_{in}))\times\bR_+)}=O(1).
$$

In order to derive estimates on $\d^k_x \phi^\Sigma $ for $k\geq 1$, we differentiate equation \eqref{eq:phi-sigma} and we proceed by induction on $k$. Notice that because of the various truncations and of the choice of $\sigma_{in}$, $\phi^\Sigma $ is identically zero in a neighbourhood of $x= x(\sigma_{in})$, so that the initial data is always zero, together with the boundary conditions. Moreover, as explained before, the indicator function
in the right-hand side does not raise any singularity. Using the same kind of estimates as in Lemma \ref{reg-lem}, we infer that
\begin{eqnarray}\label{est:uk}
 &&\| \d_x^k\phi^\Sigma\|_{L^\infty_x(L^2_Z)}+  \|\d_Z^2\d_x^k \phi^\Sigma \|_{ L^2_{x,Z}}\\
\nonumber&\leq &  C\|\d_x^k S\|_{L^1((\bar x, x(\sigma_{in}), L^2(\R_+))} \\
\nonumber&&+C \sum_{l=0}^{k-1}\|Y_-^{(k-l)}\|_{L^{4/3}(\bar x, x(\sigma_{in}))}\left(\| \d_x^l\phi^\Sigma\|_{L^\infty_x(L^2_Z)}+  \|\d_Z^2\d_x^l \phi^\Sigma \|_{L^2_{x,Z}}\right).
\end{eqnarray}

\smallskip
\noindent
$\bullet$ \textbf{ Estimates on $\d_x^k S$.}

As in the proof of Lemma \ref{sigma-s-est-bis}, the most singular estimates for $Y_-^{(k)}$ are obtained when $\cos \theta$ vanishes algebraically near $s=s_1$, with the lowest possible exponent $n$.
Moreover, by assumption (H4), the exponent $n$ above satisfies $n\geq 4$. Hence we have
\begin{equation}
\label{est:Y-}
\begin{aligned}
\|Y_-^{(k)}\|_\infty\leq C |\ln \viscosite |^{1/5}\viscosite ^{-k/20},\\
\|Y_-^{(k)}\|_{L^{4/3}}\leq C  |\ln \viscosite |^{1/5}\viscosite ^{(\frac{3}{4}-k)\frac{1}{20}}.
\end{aligned}
\end{equation}

\medskip
We now use the estimates of Lemma \ref{lem:est-A-B} in Appendix C together with \eqref{est:Y-}. Denoting 
$$
A_k:=\|a^{(k)}\|_{L^1(x_{min}, x_{max})} + \viscosite ^{1/4} \|b^{(k)}\|_{L^1(x_{min}, x_{max})}=O(\viscosite ^\frac{1-k}{19}).
$$
for all $k\in \{1,\cdots, 4\},$ we have
\begin{eqnarray*}
\|\d_x S\|_{L^1L^2}&\leq&C\left((A_1+A_2) + \|Y_-'\|_\infty(A_0+A_1)\right),\\
\|\d_x^2 S\|_{L^1L^2}&\leq&C(A_2+A_3) + \|Y_-'\|_\infty(A_1+A_2) \\&&+C (\|Y_-''\|_\infty+\|Y_-'\|_\infty^2)(A_0+A_1) ,\\
  \|\d_x^3 S\|_{L^1L^2}&\leq&C(A_3+A_4) + \|Y_-'\|_\infty(A_2+A_3) \\&&+C (\|Y_-''\|_\infty+\|Y_-'\|_\infty^2)(A_1+A_2)\\
&&+ C(\|Y_-^{(3)}\|_\infty+ \|Y_-'\|_\infty\|Y_-''\|_\infty+ \|Y_-'\|_\infty^3)(A_0+ A_1),
\end{eqnarray*}
so that
\begin{equation}
\label{est:S}
\|\d_x^k S\|_{L^1_x(L^2_Z)}\leq C \viscosite ^{-k/19}.
\end{equation}

\smallskip
\noindent
 \textbf{$\bullet\  W^{k, \infty}$ estimates on $\psi^\Sigma$:}

$\bullet$ Gathering \eqref{est:uk}, \eqref{est:Y-} and \eqref{est:S}, we prove easily by induction that for $k\in \{0,\cdots 3\}$,
$$
\| \d_x^k\phi^\Sigma\|_{L^\infty((\bar x, x(\sigma_0^-), L^2(\R_+))}+  \|\d_Z^2\d_x^k \phi^\Sigma \|_{ L^2((\bar x, x(\sigma_0^-), L^2(\R_+))}=O(\viscosite ^{-k/19}).
$$
The first estimates on $\psisig$ follow from the identities
\begin{eqnarray*}
\d_x \psisig(x,Y)&=& \d_x \phi^\Sigma(x, Y-Y_-)-Y_-'\d_Z \phi^\Sigma(x, Y-Y_-),\\
\d_x^2    \psisig(x,Y)&=&  \d_x^2 \phi^\Sigma(x, Y-Y_-)-2 Y_-'\d_x\d_Z \phi^\Sigma(x, Y-Y_-)\\&&-Y_-''\d_Z \phi^\Sigma(x, Y-Y_-)-{Y_-'}^2\d_Z^2\phi^\Sigma(x, Y-Y_-).
\end{eqnarray*}

Then, notice that for $x<x(t_1^+)$, $Y_-\equiv\bar Y_-$, so that
$$
\| \psi^{\Sigma}\|_{W^{k,\infty}((\bar x , x(t_1^+))\times (\bar Y_-, \infty))}= \| \phi^{\Sigma} \|_{W^{k,\infty}((\bar x , x(t_1^+))\times (0, \infty))}.
$$
We then use Sobolev inequalities in order to bound the $L^\infty$ norms. More precisely, we use repeatedly the following Lemma :
\begin{Lem}There exists a constant $C>0$ such that for any $f\in H^4(\R_+)$ satisfying $f_{|Z=0}= \d_Z f_{|Z=0}=0$, we have $f\in W^{3,\infty}(\R_+)$ and for $k=0, 1, 2, 3$,
$$
\|\d_Z^kf \|_{L^\infty(\R_+)}\leq C \|f\|_{L^2(\R_+)}^\frac{7-2k}{8} \|\d_Z^4f\|_{L^2(\R_+)}^\frac{1+2k}{8}.
$$
\label{lem:sobolev}

\end{Lem}
Lemma \ref{lem:sobolev} is proved in Appendix D.
Notice that for $x<x(t_1^+)$, $\phi^\Sigma$ satisfies
$$
\begin{aligned}
\d_x \phi^\Sigma- \d_Z^4 \phi^\Sigma= & (\mathbf{1}_{Z+ \bar Y^-<0}-1 )\chi^{(4)}(Z+ \bar Y^-) (\bar a + \viscosite ^{1/4}\bar b (Z+ \bar Y^-))\\
& +4 \left(1- \mathbf{1}_{x<x_1,Z + Y_-(x)<0}\right) \chi^{(3)} (Z + Y_-(x))\viscosite ^{1/4} b(x)
,
\end{aligned}
$$
so that for $k\geq 2$,
$$
\d_x^k \phi^\Sigma= \d_Z^4 \d_x^{k-1} \phi^\Sigma,\quad x<x(t_1^+), Z>\bar Y^-.
$$
Lemma \ref{lem:sobolev} and equation \eqref{eq:phi-sigma} imply that for $k\geq 1$,
\begin{eqnarray*}
\|\d_x^k\phi^{\Sigma}\|_{L^\infty ((\bar x,x(t_1^+))\times (0, \infty))}&\leq & \|\d_x^k\phi^{\Sigma} \|_{L^\infty_x(L^2_Z)}^{7/8} \|\d_Z^4\d_x^k\phi^{\Sigma}\|_{L^\infty_x(L^2_Z)}^{1/8}\\
&\leq &  \|\d_x^k\phi^{\Sigma} \|_{L^\infty_x(L^2_Z)}^{7/8} \|\d_x^{k+1}\phi^{\Sigma}\|_{L^\infty_x(L^2_Z)}^{1/8}\\
&\leq & C \viscosite ^{-\frac{7k}{8\times 19}}\viscosite ^{-\frac{k+1}{8\times 19}}\\
&\leq & C\viscosite ^{-\frac{k}{19}-\frac{1}{8\times 19}},
\end{eqnarray*}
while
\begin{eqnarray*}
\|\d_Z^l \phi^{\Sigma}\|_{L^\infty ((\bar x,x(t_1^+))\times (0, \infty))}&\leq &C \|\phi^{\Sigma}\|_{L^\infty_x(L^2_Z)}^{(7-2l)/8}\|\d_Z^4 \phi^{\Sigma}\|_{L^\infty_x(L^2_Z)}^{(1+2l)/8}\\
&\leq & C\left(\|\d_x \phi^{\Sigma} \|_{L^\infty_x(L^2_Z)} + \|S\|_{L^\infty_x(L^2_Z)}\right)^{(1+2l)/8}\\
&\leq & C \viscosite ^{-\frac{1+2l}{8\times 19}}.
\end{eqnarray*}
Similar estimates can be derived for terms of the type $\d_x^k \d_Z^l \phi^{\Sigma}$, with $k+l\leq 3$. This concludes the proof of Lemma \ref{lem:est-psi-sigma}.

\end{proof}

\medskip
\begin{Rmk}\label{asymptotic-rmk}
$\bullet$
Because of the function $Y_-$ and the coefficients $a$ and $b$, the profile $\psi^\Sigma$ depends on $\viscosite $ in general, which is not entirely satisfactory for an approximate solution. However, the estimates of Lemma \ref{lem:est-psi-sigma}
show that $\psi^\Sigma $ is uniformly bounded in $L^\infty(L^2)\cap L^2(H^2)$.

Let $\underbar{x}\in (\bar x , x_1)$ arbitrary. Then $x(t_1^+)>\underbar{x}$ for $\viscosite  $ small enough, so that $Y_-(x)=\bar Y_-$ for $x\in (\bar x, \underbar{x})$. Extending $\psi^\Sigma$ by zero below $\bar Y_-$, we infer that $\psi^\Sigma$ is uniformly bounded in $L^2( (\bar x, \underbar{x}), H^2(\R))$. Up to the extraction of a sub-sequence, $\psi^\Sigma$ converges  towards a function $\bar \psi^\Sigma$ weakly in $L^2( (\bar x, \underbar{x}), H^2(\R))$. Passing to the limit in the equation satisfied by $\psi^\Sigma$, we obtain
$$
\d_x \bar \psi^\Sigma- \d_Y^4\bar \psi^\Sigma=\bar S(Y),\text{ in } \mathcal D'( (\bar x, \underbar{x})\times \R),
$$
where
$$\begin{aligned}
\bar S(Y)=-\bar a \chi^{(4)}(Y),
\text{with }\bar a=-\int_{x_E^-(y_1)}^{x_E^+(y_1)} \tau(x,y_1)\:dx.
  \end{aligned}
$$

$\bullet$ Notice that with this construction, there remains a non zero trace on the East part of the boundary $(s_1, s_2)$, and also on $(s_{k-1}, s_0)$ if $(s_{k-1}, s_0)\subset \Gamma_E$ (since $\sigma_{in}= s_0$  according to \eqref{def:sigma-}, and therefore the trace of $\psi_E$ on $(\sigma_0^-, s_0)$ is not lifted by $\psil$).
These traces correspond to $\psi_{E|Z=0} (1- \varphi_1^+ \varphi_2^-)$ and $\psi_{E|Z=0} (1- \varphi_0^- \varphi_{k-1}^+)$. 

The right way to proceed is to lift first $\psi_{E|Z=0} (1- \varphi_0^- \varphi_{k-1}^+)$ thanks to  a South type boundary layer $\psi_S$, exactly as explained in the previous section, then to construct the 
 $\psil$ and $\psi^\Sigma$ corresponding to $\psi^{int}+\psi_S$, and finally to construct another South boundary layer to lift  $\psi_{E|Z=0} (1- \varphi_1^+ \varphi_2^-)$. The principle as usual is to restore boundary conditions starting from the East. This will introduce minor corrections in the approximate solution, and will not change the a priori estimates on the different terms, nor the connection with the West boundary.

\end{Rmk}

 \subsection{ Connection with the West boundary}
 \label{ssec:sigma-west}
 The question is then to describe the West junction, that is the connection with the West boundary layer. The traces of the discontinuity terms $\psil+ \psi^\Sigma\chi\left(\frac{y-y_1}{\delta_y}\right) $ on the West boundary create yet another singularity in the West boundary layer term: indeed, the coefficient  $A(s)$ in \eqref{def:AW} now lifts the traces of $\psi^{int}+ \psi^{\Sigma}$, and thus changes abruptly near the point  $s_1^W\in \Gamma_W$. We recall that $\psi^{int}= \psi^0_t + \psil$.

At this stage we encounter another difficulty: the function $\d_s A$ can be written as
\begin{eqnarray*}
\d_s A_\pm(s)&=&\pm \left|\frac{\viscosite}{\cos \theta}\right|^{1/3}\cos \theta(s) \sin \theta(s) \d_y^2 \psi^{int}(x(s),y(s))\\
&&+\text{differentiable terms.}
\end{eqnarray*}
The function $\d_y^2 (\psi^0_{t}+\psil)$ is discontinuous across $\Sigma$: more precisely,
$$
[\d_y^2 (\psi^0_{t}+\psil)]_{\Sigma}= [\d_y^2 \psi^0_{t}]_{\Sigma},
$$
and the jump of $\d_y^2 \psi^0_{t}$ across $\Sigma$ is constant along $\Sigma$ and of order $\delta_y^{-2}$. As a consequence, $A'$ is not continuous, which is problematic, because the error estimates for the West boundary layer involve $L^p$ norms of $A''$, as we will see in the next chapter (see the proof of Lemma \ref{lem:error-EW}).

In order to avoid the apparition of a Dirac mass in $A''$, we add yet another corrector term which has the following form
\be\label{def:psicorr}
\psi^{corr}_\Sigma = -[\d_y^2 \psi^0_{t}]_{\Sigma}\mathbf 1_{y>y_1} (y-y_1)^2 \chi\left(\frac{y-y_1}{\viscosite^{1/4}}\right)\chi\left(\frac{x-x(s_1^W)}{|\log \viscosite|^{-1/5}}\right).
\ee
It is obvious that $\psi^{corr}_\Sigma\in H^2(\Om)$. Moreover, $\psi^{corr}_\Sigma$ satisfies the following estimates:

\begin{Lem}
 For all $k\in \{0,1,2\}$, we have
$$
\|\psi^{corr}_\Sigma\|_{W^{k,\infty}}=O(\delta_y^{-2} \viscosite^{-\frac{k-2}{4}})=o(\viscosite^{-k/4}).
$$
Moreover,
$$
\|\d_x\psi^{corr}_\Sigma\|_{L^2(\Om)}=o(\viscosite^{1/8}),\quad \| \viscosite \Delta^2 \psi^{corr}_\Sigma\|_{H^{-2}(\Om)}=o(\viscosite^{5/8}).
$$

\label{lem:psi-corr}
\end{Lem}
\begin{proof}
The estimates are fairly straightforward, noticing that for $l\in \{0,1,2\}$, $k\in\N$,
$$
\left\| (y-y_1)^l \chi^{(k)}\left(\frac{y-y_1}{\viscosite^{1/4}}\right)\right\|_{L^\infty}= O(\viscosite^{l/4}).
$$
Since $[\d_y^2 \psi^0_{t}]_{\Sigma}$ does not depend on $x$ nor $y$, we have
$$
\d_x^k\psi^{corr}_\Sigma= -[\d_y^2 \psi^0_{t}]_{\Sigma}\mathbf 1_{y>y_1} (y-y_1)^2 \chi\left(\frac{y-y_1}{\viscosite^{1/4}}\right) |\ln \viscosite|^{k/5}\chi^{(k)}\left(\frac{x-x(s_1^W)}{|\log \viscosite|^{-1/5}}\right),
$$
so that
$$
\| \d_x^k  \psi^{corr}_\Sigma\|_\infty=O(\delta_y^{-2} \viscosite^{1/2} |\ln \viscosite|^{k/5})= O(|\ln \viscosite|^\frac{k-2}{5}|\ln |\ln \viscosite|\:|^{2\beta}).
$$
Let us now estimate $\d_y \psi^{corr}_\Sigma$; the rest of the estimates are left to the reader. We have
\begin{multline*}
\d_y \psi^{corr}_\Sigma(x,y)= -[\d_y^2 \psi^0_{t}]_{\Sigma}\chi\left(\frac{x-x(s_1^W)}{|\log \viscosite|^{-1/5}}\right)\mathbf 1_{y>y_1} (y-y_1) \\
\times
\left[2\chi\left(\frac{y-y_1}{\viscosite^{1/4}}\right) + \viscosite^{-1/4}(y-y_1)\chi'\left(\frac{y-y_1}{\viscosite^{1/4}}\right) \right],
\end{multline*}
so that
$$
\| \d_y \psi^{corr}_\Sigma\|_\infty \leq C \left| [\d_y^2 \psi^0_{t}]_{\Sigma}\right| \viscosite^{1/4}\leq C \delta_y^{-2} \viscosite^{1/4}.
$$
More generally, for $k=0,1,2$
$$
\| \d_y^k  \psi^{corr}_\Sigma\|_\infty=O(\delta_y^{-2} \viscosite^\frac{2-k}{4})=O(\viscosite^{-k/4} |\ln \viscosite|^{-2/5} |\ln |\ln \viscosite|\:|^{2\beta}).
$$
The $W^{k,\infty}$ estimates follow easily. Moreover, since the support of $\psi^{corr}_\Sigma$ is included in the rectangle $$[x_1^W-2|\ln \viscosite|^{-1/5}, x_1^W+2|\ln \viscosite|^{-1/5}]\times [y_1-2 \viscosite^{1/4}, y_1+2\viscosite^{1/4}],$$ we infer that
$$\begin{aligned}
\| \d_x \psi^{corr}_\Sigma\|_{L^2(\Om)}=O(|\ln \viscosite|^{-1/5}|\ln |\ln \viscosite|\:|^{2\beta}|\ln \viscosite|^{-1/10}\viscosite^{1/8})= o(\viscosite^{1/8}),\\
\|\Delta\psi^{corr}_\Sigma\|_{L^2(\Om)}=O(\viscosite^{-1/2}  |\ln \viscosite|^{-2/5} |\ln |\ln \viscosite|\:|^{2\beta}|\ln \viscosite|^{-1/10}\viscosite^{1/8})= o(\viscosite^{-3/8}).
\end{aligned}
$$

\end{proof}

\bigskip

In the rest of the paper, we set
$$
\bar \psi_\Sigma:= \psil+ \psi^\Sigma\chi\left(\frac{y-y_1}{\delta_y}\right)+ \psi^{corr}_\Sigma.
$$
\label{psibSigma}
With this definition, $\psi^0_t+\bar \psi_\Sigma \in H^2(\Om)$.

 Moreover,  we have
$$
\d_x (\psi_t^0 +\bar \psi_\Sigma) - \viscosite\Delta^2(\psi_t^0 +\bar \psi_\Sigma)= \tau + \delta\tau+ r_{\mathrm{lift}}^1 + r_{\mathrm{lift}}^2 +\delta\tau_\sigma \quad\text{on }\Om
$$
where  $\delta \tau$ has been estimated in Proposition \ref{lem:trunc}, $ r_{\mathrm{lift}}^1 + r_{\mathrm{lift}}^2$ is controlled by Proposition \ref{lem:eq-psil} and
\be
\label{dtausig}
 \delta \tau_\Sigma = \d_x \psi^{corr}_\Sigma- \viscosite \Delta^2\psi^{corr} +\viscosite (\d_y^4- \Delta^2) \psi^\Sigma ,
\ee

        \begin{Prop}
                Assume that $\Om$ satisfies assumptions (H1)-(H4).
Then $ \delta \tau_\Sigma$ is an admissible remainder in the sense of  Definition \ref{remainder-def}.
\label{lem:sigma}

        \end{Prop}
This proposition will be proved in Chapter \ref{chap:convergence}.

 Note further that, by assumption on $\Sigma$,  the  trace of $\bar \psi_\Sigma+\psi^0_t+\bar \psi_{N,S}$ on the  boundary $\Gamma_W$  vanishes in the vicinity the points $s_i\in \Gamma_W$.

\begin{Rmk}
\label{rem:hyp-h4}
If one looks carefully at the construction above, there are several reasons why the assumption (H4) is necessary: some of them are technical, and others are probably more fundamental. We first notice that if $I_1=[s_0, s_1]$ with $s_0<s_1$, assumption (H4) is always satisfied, as $\d\Om$ is $\mathcal C^4$ near $s_0$. Therefore we focus on the case $I_1=\{s_1\}$.
\begin{itemize}
\item First, assumption (H4) is required so that the coefficients $a$ and $b$ have derivatives of order four; if $\cos \theta$ vanishes only at order $n$ near $s_1$, with $n\leq 3$ in (H2i), then such a regularity is not achieved in general. Therefore, if (H4) is not satisfied, an additional regularization of $a$ and $b$ must be performed.

\item Moreover, the estimates on $\psi^\Sigma$ are more singular if the cancellation of $\cos \theta$ near $s=s_1$ has a lower order. This could have two different consequences: on the one hand, it could be possible that the estimates on $\psi^\Sigma$ and its derivatives in the interior of $\Om$ are not sufficient to prove that the error terms are admissible. This, however,  is unlikely: indeed, the estimates are exactly the same as the ones for North and South boundary layer terms, for which we are able to prove that the error terms are admissible even when $\cos \theta$ vanishes at the lowest possible order.

On the other hand, if $\cos \theta$ vanishes at a low order near $s=s_1$, the estimations on $\|\psi^\Sigma\|_{W^{k,\infty}}$ become much more singular, and therefore the trace of $\psi^\Sigma$ and its derivatives on the West boundary are much larger. In this regard, assumption (H4) does not seem to be merely technical: it is possible that the discontinuity boundary layer may destabilize the West boundary layer if $\psi^\Sigma$ becomes very large. This could in fact be achieved with a low order of cancellation of $\cos \theta$ near $s_1$. However, such considerations go beyond the scope of this paper.

\end{itemize}

\end{Rmk}


\section{Lifting the West boundary conditions}\label{west-connection}
 \label{sec:Wboundary}

We recall that the boundary layer term on  West boundaries takes the form $\psi_W(s,Z)= A(s)^t f(Z)$, where $A(s)$ is  defined by
\be\label{def:AW-fin}
A(s)={1\over \sqrt{3}}
        \begin{pmatrix}
                 e^{i\pi/6} &-i\lambda_W^{-1}\\ e^{-i\pi/6} &i \lambda_{W}^{-1}
        \end{pmatrix}
        \begin{pmatrix}
                - (\psi^0_t + \psib_\Sigma)_{|\d \Om}(s)\\ \d_n(\psi^0_t + \psib_\Sigma)_{|\d \Om}(s)
        \end{pmatrix}
\ee
for $s\in \Gamma_W$ and
$$
f(Z)=\begin{pmatrix}
        \exp\left(-e^{i\pi/3} Z\right)\\  \exp\left(-e^{-i\pi/3} Z\right)
       \end{pmatrix}.
       $$

\begin{Lem} Define the amplitude of the West boundary layer by (\ref{def:AW-fin}).
Let $(s_i, s_{i+1})$ be a portion of $\Gamma_W$.

Then there exists a constant $c_0$ such that the following estimates hold:

\begin{itemize}
\item  Far from the West end of $\Sigma$, i.e. for $|s-s_1^W|\geq c_0 \delta_y$, and for all $s\in \Supp (\varphi_i^+ \varphi_{i-1}^-)$,
\begin{eqnarray}
|A(s)|&\leq & C\nonumber,\\
\label{est:dsA}|A'(s)|&\leq & C \left(1+|\cos \theta|\ccM(y(s)) \right)\\
 \label{est:ds2A}
|A''(s)|&\leq & C \left(\frac{1}{\delta_y}+ \frac{(\cos \theta)^2}{\delta_y^2|\ln \delta_y|}\right).
\end{eqnarray}

\item If $|s-s_1^W|\leq c_0\delta_y$,
\begin{equation}
\label{A-Sigma-est}
|\d_s^k A(s)| \leq C \viscosite^{-\frac{k}{4}- \frac{1+2k}{8\times 19}}\quad\text{for } k\in \{0,1,2\}.
\end{equation}

\end{itemize}
\label{lem:est-A}
\end{Lem}

\begin{proof}
Consider a component $[s_{i-1},s_i]$ of the West boundary. There are two zones on which the $s$ derivatives of $A$ become singular:
\begin{itemize}
 \item near $s_{i-1}$ and  $s_i$, or, more generally, near points of $\d\Om$ such that $|y(s)-y_j|=0$ ($y_j\neq y_1$); on such zones, the derivatives of the function $\psi^0_t$ become unbounded, while $\psib_\Sigma$ is identically zero.

\item near $s_1^W$, the end point of $\Sigma$; we recall that the trace of $\psi^0_t$ is discontinuous at this point, and that this discontinuity is lifted by the term $\psil$. Additionally, the $s$ derivative of $\d_n\psi_t^0$ is also discontinuous, and this discontinuity is corrected by the term $\psi^{corr}_\Sigma$ introduced  in equation \eqref{def:psicorr}. Notice however that on this zone, $\cos \theta$ is bounded away from zero.

\end{itemize}

Away from these zones, $A$  and its derivatives are bounded, and $\cos \theta$ is bounded away from zero.
Whence we now focus on the two pathological zones.
Since the proof of the estimation on each zone is rather different, we separate the two.
\smallskip

$\bullet$ \textbf{Estimate near points such that $|y(s)-y_j|\ll 1$ ($y_j\neq y_1$):}

In the vicinity of such  points, thanks to assumption (H3), $\psi^0_t$ is $\mathcal C^\infty$ (with  derivatives which are not uniformly bounded). We also use the following bound:
$$
\sup_{s\in \supp (\varphi_i^- \varphi_{i-1}^+)} |\lambda_{W}^{-1}(s)|\leq C\left\{ \begin{array}{ll}
                                                                  \viscosite^\frac{n+1}{4n+3}&\text{ in case (H2i)},\\
                                                                \viscosite^{1/4}|\ln \viscosite|^{-1/2}&\text{ in case (H2ii),}
                                                                 \end{array}\right.
$$
so that $\sup_{s\in \supp (\varphi_i^- \varphi_{i-1}^+) } |\lambda_{W}^{-1}(s)|\ll \viscosite^{1/4}\ll \delta_y$. In the rest of the proof, we write $\lambda$ for $\lambda_{W}$.
Using Lemma \ref{lem:est-psi0}, we infer that for all $s\in \supp (\varphi_i^- \varphi_{i-1}^+)$,
\begin{equation}
\label{est:A}
|A(s)|\leq C \left(1+ |\lambda(s)|^{-1}\ccM(y(s))\right)\\
\leq C\,.
\end{equation}
The estimates of the $s$ derivatives are a little more involved. Using \eqref{def:AW-fin}, we have
$$
|A'(s)|\leq C \left( |\d_s \psi^0_{t|\d \Om}(s) |+|\lambda|^{-1} |  \: |\d_s \d_n\psi^0_{t|\d \Om}(s) | + \left|\frac{\lambda'}{\lambda^2}\right|\; |\d_n\psi^0_{t|\d \Om}(s) | \right).
$$
Notice that on the support of $\varphi_i^- \varphi_{i-1}^+$,
$$
\left|\frac{\lambda'}{\lambda^2}\right|\leq \frac{\viscosite^{1/3} |\theta'|}{|\cos \theta|^{4/3}}\leq |\cos \theta|
$$
thanks to the validity condition \eqref{hyp:valid-E/W-2}.
Lemma \ref{lem:est-psi0} then implies, for $s\in \linebreak\supp (\varphi_i^- \varphi_{i-1}^+)$,
$$|A'(s)|\leq  C \left(1+ \frac{|\cos \theta|}{\delta_y|\ln \delta_y|}\sum_{j\in I_+} \mathbf{1}_{|y(s)-y_j|\leq \delta_y}+|\cos \theta| \ccM(y(s))\right)\,.$$

In a similar fashion,
\begin{eqnarray*}
|A''(s)|&\leq & C \left(|(\d_s^2  \psi_t^0)_{|\d \Om}(s) | + \left|\frac{\lambda'}{\lambda^2}\right|\; |(\d_s\d_n \psi_t^0)_{|\d \Om}(s) | \right)\\
&&+C \left( |\lambda|^{-1} |  \: |(\d_s^2 \d_n \psi_t^0)_{|\d \Om}(s) | + \left|\left(\frac{\lambda'}{\lambda^2}\right)'\right|\; |(\d_n \psi_t^0)_{|\d \Om}(s) | \right).
\end{eqnarray*}
We have
$$
 \left|\left(\frac{\lambda'}{\lambda^2}\right)'\right|\leq C \frac{\viscosite^{1/3}}{|\cos \theta|^{7/3}}(|\theta'|^2 + |\theta''| \:|\cos \theta|).
$$
Since $\cos \theta$ vanishes always at a higher order than $\theta'$, we infer that on $\supp(\varphi_i^- \varphi_{i-1}^+)$,
$$
 \left|\left(\frac{\lambda'}{\lambda^2}\right)'\right|\leq C \frac{\viscosite^{1/3}|\theta'|^2}{|\cos \theta|^{7/3}}\leq C.
$$
Finally, we obtain, for $s\in \supp (\varphi_i^- \varphi_{i-1}^+)$,
$$
|A''(s)|\leq  C \left(\frac{1}{\delta_y}+ \frac{(\cos \theta)^2}{\delta_y^2|\ln \delta_y|}\right).
$$

\smallskip

\textbf{$\bullet$ Estimates near the end point of $\Sigma$:}

As explained above, since the interior term is now $\psi_t^0+\bar \psi_\Sigma$, the coefficient $A(s)$ changes abruptly near $s_1^W$, on an interval whose length is of order $\delta_y$. Moreover, since (H3) is satisfied, we have
$$
\cos \theta(s_1^W)\neq 0,
$$
so that $\varphi_i^- \varphi_{i-1}^+$ and its derivatives are $O(1)$ in a neighbourhood $V$ of $s_1^W$ independent of $\viscosite$.
 As before, we use \eqref{def:AW-fin} and we obtain
\begin{eqnarray*}
\|A\|_{L^\infty(V)}&\leq & C (\| \psi_t^0+\bar \psi_\Sigma \|_{L^\infty} + \viscosite^{1/3} \|\psi_t^0+\bar \psi_\Sigma\|_{W^{1,\infty}}),\\
\|A'\|_{L^\infty(V)}&\leq & C (\| \psi_t^0+\bar \psi_\Sigma\|_{W^{1,\infty}} + \viscosite^{1/3} \|\psi_t^0+\bar \psi_\Sigma\|_{W^{2,\infty}}).
  \end{eqnarray*}
The $W^{k,\infty}$  norms of the right-hand side are to be understood in a neighbourhood of size $O(1)$ in $\Om$ of $(x(s_1^W), y(s_1^W))$. We now recall that thanks to the addition of the corrector term $\psi^{corr}_\Sigma$, the amplitude $A'$ is differentiable at $s=s^W_1$, and
$$
\|A''\|_{L^\infty(V)}\leq C( \|\psi_t^0+\bar \psi_\Sigma\|_{W^{2,\infty}} + \viscosite^{1/3} \|\psi_t^0+\bar \psi_\Sigma\|_{W^{3,\infty}}).
$$
We now estimate $\|\psi_t^0+\bar \psi_\Sigma\|_{W^{k,\infty}}$ using Lemmas \ref{lem:est-psi0}, \ref{lem:est-psi-sigma}, \ref{lem:psi-corr}, together with the explicit definition of $\psil$. We infer that
$$
\| \psi_t^0+\bar \psi_\Sigma\|_{W^{k,\infty}}=O(\viscosite^{-\frac{k}{4}- \frac{1+2k}{8\times 19}})\quad\text{for } k\in \{0,\cdots, 3\},
$$
which leads to the desired estimates.

\end{proof}

\section{Approximate solution in the periodic and rectangle case}
\label{sec:rectangle}

\subsection{In the periodic case},  the solution has already been completely defined (see \eqref{psiapp-per} and section \ref{sec:BL-per}). We therefore only give here a few estimates on that solution.

We recall that
$$\psiapp =\psi^{circ} + \psi^0_{per} + \psbl_{N,S} $$
where $\psi^0_{per} $ and  $\psbl_{N,S} $ are defined as before and $\psi^{circ}=\psi^{circ}(y)$  is the circumpolar current, due to the average forcing by the wind, namely
$$
\begin{aligned}
 -\viscosite \d_y^4 \psi^{circ} = \mean{\tau} (y) , \quad y\in (y_-, y_+),\\
  \psi^{circ}_{|y=y_\pm} = 0,\\
  \d_y \psi^{circ}_{|y=y_c^\pm} = 0.
  \end{aligned}
  $$
Provided $\tau $ is sufficiently smooth (i.e. belongs to $H^s$ for $s$ large enough), we have
\be\label{est:psbl-per}
\begin{aligned}
\|\psi_{circ}\|_{H^s}&\leq C\viscosite^{-1} \|\tau\|_{H^{(s-4)_+}},\\
\|\psi^0\|_{H^s}&\leq C \|\tau\|_{H^s},
\end{aligned}
\ee
and consequently by (\ref{est:psblNS-per})
\be
\begin{aligned}
\| \d_x^k \d_Z^l\psbl_{N,S}\|_{L^2(\T\times \R_+)} &\leq C( \|\psi^0_{y=y_\pm}\|_{H^{k + \frac{l}{4}- \frac{1}{8}}(\T)} + \viscosite^{1/4} \|\d_y\psi^0_{y=y_\pm}\|_{H^{k + \frac{l}{4}- \frac{3}{8}}(\T)})\\
&\leq C( \|\tau \|_{H^{k + \frac{l}{4} + \frac{3}{8}}} +\viscosite^{1/4} \|\tau \|_{H^{k + \frac{l}{4} + \frac{9}{8}}} ).
\end{aligned}
\ee

\bigskip
\subsection{ In the rectangle case}, i.e.  when $\Om=(x_-,x_+)\times (y_-, y_+)$, the construction explained in the preceding sections for a smooth domain $\Om\subset \R^2$ still works:
\begin{itemize}[leftmargin=*]
\item Since $\cos \theta$ does not vanish at the north and south ends of the east boundary, there is no need for a truncation. We therefore simply define
$$
\begin{aligned}
\d_x \psi^0=\tau\quad \text{in } \Om,\\
\psi^0_{|x=x_+}=0.
\end{aligned}
$$

\item We also construct an East boundary layer corrector and a macroscopic corrector (see section \ref{sec:macro-corrector}) on the whole East boundary; again, this is possible because there is no singularity in $\psi^0$. We denote by $\bar \psi_E(x,y)$ the sum of these two terms.


\item We then construct north and south boundary layers $\psbl_{N,S}$, lifting the traces of $(\psi^0+ \bar \psi_E)_{|y=y_\pm}, \d_y(\psi^0+ \bar \psi_E)_{|y=y_\pm}$, with zero initial condition at $x=x_+$. When $\tau$ is sufficiently smooth, $\Psi_0 = - (\psi^0+ \bar \psi_E)_{|y=y_\pm}$ and $\Psi_1 = - \viscosite^{1/4}\d_y(\psi^0+ \bar \psi_E)_{|y=y_\pm}$ are smooth functions of $x$ with large derivatives of order 2 or higher (because of the east boundary layer),  which satisfy better estimates than in the case of a regular domain. Typically
$$
\begin{aligned}
\| \Psi_j \|_{L^\infty} +\| \d_x \Psi_j \|_{L^\infty}\leq C,\\
 \d_x^2 \Psi_j = O(1)_{L^\infty} + O( \viscosite^{-1/3} \indc_{|x- x_+| \leq C \viscosite^{1/3}}).
 \end{aligned}
 $$

Furthermore, as $\psbl_{N,S}$ is a smooth function of $x$ and $Z$,
$$
\d_x \psbl_{N,S|x=x_+}=\d_Z^4 \psbl_{N,S|x=x_+}=0.
$$
As a consequence, the trace of $\psbl_{N,S}$ on the east boundary and the trace of its normal derivative (i.e. its $x$ derivative) is identically zero: the {\bf North and South boundary layer correctors do not pollute the East boundary}.

\item The next and last step is the construction of the western boundary layer. The latter lifts the trace of $\psi^0_t + \psbl_{N,S}$ and its normal derivative on $x=x_-$.  Let
$$
\begin{pmatrix}
A_+(y)\\A_-(y)
\end{pmatrix}
=
\frac{1}{\sqrt{3}}\begin{pmatrix}
e^{i\pi/6} & -i\viscosite^{1/3}\\
e^{-i\pi/6} & i\viscosite^{1/3}
\end{pmatrix}
\begin{pmatrix}
f_0(y)\\f_1(y)
\end{pmatrix},
$$
where
$$
\begin{aligned}
f_0(y)= - \left(\psi^0_{t|x=x_-}(y) + \psbl_{N}\left(x_-, \frac{y_+-y}{\viscosite^{1/4}}\right) \chi_0(y)+ \psbl_{S}\left(x_-, \frac{y-y_-}{\viscosite^{1/4}}\right)\chi_0(y) \right),\\
f_1(y)=\d_x\psi^0_{t|x=x_-}(y) + \d_x\left(x_-, \frac{y_+-y}{\viscosite^{1/4}}\right) \chi_0(y)+ \d_x\psbl_{S}\left(x_-, \frac{y-y_-}{\viscosite^{1/4}}\right)\chi_0(y) .
\end{aligned}
$$
as the macroscopic truncation $\chi_0$ can be chosen as a function of $y$ only.
Notice that by definition of $\psbl_{N,S}$ $A_+ (y)$ and $A_-(y)$ both vanish at $y_- $ and $y_+$. Therefore the amplitude of the western boundary layer is zero at $y=y_\pm$: {\bf its trace vanishes on the northern and southern boundaries}. Thus there is no need for additional south and north boundary layers.

Using the energy estimates of Lemma \ref{reg-lem}, we obtain
$$
\|\d_x^k \psbl_{N,S}\|_{L^\infty_x(L^2_Z)} +  \|\d_x^k \d_Z^2 \psbl_{N,S}\|_{L^2_x(L^2_Z)} \leq C\viscosite^{-(k-1)_+/3}
$$
from which we infer, using Lemma \ref{lem:sobolev}, for $k\geq 0$ and $l\in \{0,1,2,3\}$,
\begin{eqnarray*}
\|\d_x^k \d_Z^l \psbl_{N,S}\|_{L^\infty_{x,Z}}&\leq &C \|\d_x^k \psbl_{N,S}\|_{L^\infty_x(L^2_Z)}^{\frac{7-2l}{8}} \|\d_x^k \d_Z^4\psbl_{N,S}\|_{L^\infty_x(L^2_Z)}^{\frac{1+2l}{8}}\\
&+ &C\viscosite^{-(k-1)_+/3} \\
&\leq & C \|\d_x^k \psbl_{N,S}\|_{L^\infty_x(L^2_Z)}^{\frac{7-2l}{8}} \|\d_x^{k+1}\psbl_{N,S}\|_{L^\infty_x(L^2_Z)}^{\frac{1+2l}{8}}\\
&+ &C\viscosite^{-(k-1)_+/3} \\
&\leq & C \viscosite^{-\frac{(7-2l)(k-1)_+}{24}-\frac{(1+2l)k}{24}}.
\end{eqnarray*}
This implies that for $l\in\{0,1,2\}$
\be\label{est:A-rect}
\|\d_y^l A_\pm\|_{L^\infty(y_-,y_+)}\leq C \viscosite^{-l/4}.
\ee

\end{itemize}

\bigskip

The construction of the approximate solution in all the cases considered in the introduction is now complete. There remains to check that the function which we built is indeed an approximate solution in the sense of Definitions \ref{def:app-general} and \ref{def:app-per}. In other words, we still need to prove that the remainder terms created by the various approximations are admissible. This is the goal of the next chapter.

%% file: convergence.tex
\chapter{Proof of  convergence}\label{chap:convergence}

The approximate solution $\psiapp$ we have built in the previous chapter
\begin{itemize}
\item[(i)] satisfies exactly the boundary conditions
$$\psiapp_{|\d \Omega} = (\d_n\psiapp)_{|\d \Omega} = 0$$
by definition of the boundary layer terms;

\item[(ii)] has $H^2$ regularity, but with non uniform bounds on the derivatives (see for instance the definition of the interior term near East corners and surface discontinuities);

\item[(iii)] but does not satisfy the Munk equation, only the approximate version
$$\d_x \psiapp -\viscosite \Delta^2 \psiapp = \tau +\delta \tau,$$
where $\delta \tau$ comes from both the interior and boundary terms.

\end{itemize}

The goal of this final chapter is to estimate the different contributions to this error term $\delta \tau$, and to check that it is admissible in the sense of Definition \ref{remainder-def} (or satisfies \eqref{app-def-circ} in the periodic case). At this stage, there are no more conceptual difficulties: all the proofs are based on technical computations.  Exactly as in the previous chapters, we focus mainly on the case when  $\Omega$ is a smooth domain in $\R^2$, and we have gathered in a paragraph at the end of the chapter the estimates corresponding to the periodic and rectangular cases.

\section{Remainders stemming from the interior term $\psi^{int}= \psi^0_t + \psil$}$ $

Going back to the construction of the interior term, we see that its contribution $\delta \tau_{int}$ to the remainder is made of two kinds of terms.
\begin{itemize}
\item Because of the truncation, the solution $\psi_t^0$ to the transport equation satisfies the Munk equation with  source terms $\tau (1-\chi_\viscosite)$ and $-\viscosite \Delta^2  \psi_t^0$, which are admissible remainders according to Proposition \ref{lem:trunc}.


\item The lifting term which restores the $H^2$ regularity near discontinuity surfaces also introduces  correctors.
Part of them, namely $\dtl$, are considered as source terms for the singular surface boundary layers.
The other ones will be proved to be admissible remainders as stated in Proposition \ref{lem:eq-psil}.

\end{itemize}

\subsection{Error terms due to the truncation $\chi_\viscosite$.}$ $

Let us first prove  Proposition \ref{lem:trunc}, i.e. the following estimates~:
$$
\begin{aligned}
\| \tau (\chi_\viscosite-1)\|_{H^{-2}(\Om)}= o (\viscosite^{5/8}),\\
\| \Delta \psi^0_t\|_{L^2(\Om^\pm)}= o(\viscosite^{-3/8}).
\end{aligned}
$$

\medskip
\noindent
                $\bullet$ We start by localizing the problem around $(x_i,y_i)$ thanks to a partition of unity:
 for $i\in I_+$,  let $f_i \in\mathcal C^\infty_0(\R^2)$ such that
                $$
                \sum_{i\in I_+} f_i=1\quad \text{on }\bar \Om,
                $$
                and
                \begin{itemize}
                        \item for all $i\in I_+$, $f_i\equiv 1$ in a neighbourhood of $(x_i,y_i)$,
                        \item for all $i\neq j$, $d((x_i,y_i), \Supp f_j)>0$.
                \end{itemize}
Then, according to the definition of $\chi_\viscosite$, for $\viscosite$ small enough,
\begin{eqnarray*}
        \tau \chi_\viscosite&= & \sum_{i\in I_+} f_i \tau \chi_\viscosite\\
        &= &  \sum_{i\in I_+}  f_i \tau \chi_i\left({x-x_i},{y-y_i}\right).
\end{eqnarray*}
Thus, defining $\psi^0_i$ by
$$\begin{aligned}
        \d_x \psi^0_i=  f_i \tau \chi_i\left({x-x_i},{y-y_i}\right),\\
        \psi^0_{i|\Gamma_E}=0,
  \end{aligned}
$$
we have $\psi^0_t=\sum \psi^0_i$, and it suffices to
 prove that under the assumptions of  Proposition \ref{lem:trunc},
 $$\begin{aligned}
  \Big\|f_i \tau \left[\chi_i\left(x-x_i,y-y_i\right)-1\right]\Big\|_{H^{-2}(\Om)}=o(\viscosite^{5/8}),\\
          \|  \viscosite\Delta^2 \psi^0_i\|_{H^{-2}(\Om^\pm)}=\viscosite\|  \Delta \psi^0_i\|_{L^2(\Om^\pm)}=o(\viscosite^{5/8}).
   \end{aligned}
$$

In order not to burden the notation, we drop the indices $i$ and we shift the origin of the arc-length parametrization so that $s_i=0$.  We also shift the origin of the axes so that $(x_i,y_i)=(0,0)$, and we rename respectively $\tau$ and $\chi_\viscosite$  the functions $f_i\tau$ and $\chi_i(x-x_i, y-y_i)$  (which have the same regularity as the original functions).

Note that, by definition of $\chi_\viscosite$, we have the obvious estimates
\begin{equation}
\label{chi-est}
\begin{aligned}
\| \d_x^{s_1} \d_y^{s_2} \chi_\viscosite \|_{L^\infty} \leq C \max_\pm \left({|\ln \delta_y|\over \delta_x^\pm}\right)^{s_1}\left({1\over \delta_y}\right)^{s_2},\\
\| \d_x^{s_1} \d_y^{s_2} (1-\chi_\viscosite) \|_{L^2} \leq C\max_\pm\left(  \left( {|\ln \delta_y|\over \delta_x^\pm}\right)^{s_1}\left({1\over \delta_y}\right)^{s_2}(\delta_x^\pm \delta_y)^{1/2}\right) .
\end{aligned}
\end{equation}

\bigskip
\noindent
$\bullet$ We now give an upper-bound for the source term $\tau (1-\chi_\viscosite)$ in terms of $\delta_x^\pm$, $\delta_y$.
Without loss of generality, we can assume that  $s=0$ is a North boundary point of $\Gamma_E$.
Define, for $(x,y)\in \Om$,
$$
\phi_\viscosite(x,y)=\int_{-\infty}^y dy'\int_{-\infty}^{y'} dz (\chi_\viscosite-1)(z).
$$
so  that
$$
\d_y^2 \phi_\viscosite=\chi_\viscosite-1,\quad \phi_\viscosite(x,y)=0\text{ for }(x,y)\in \Om, \ d((x,y), \d\Om)\gg \delta_x,\delta_y.
$$
Since the support  of $1-\chi_\viscosite$ is included in the rectangle $$ [-\delta_x^- (1+|\ln \delta_y|^{-1}),  \delta_x^+ (1+|\ln \delta_y|^{-1})] \times [-\delta_y,\delta_y],$$
 we have
$$
\begin{aligned}
\phi_\viscosite(x,y)=0\text{ if } x\geq \delta_x^+(1+|\ln \delta_y|^{-1})\text{ or }  x\leq -\delta_x^-(1+|\ln \delta_y|^{-1}),\\
\text{and }  \phi_\viscosite(x,y)=0\text{ if } y<-\delta_y.
\end{aligned}
$$
Notice also that
$$
|\phi_\viscosite(x,y)|\leq C(\delta_y^2 + |y|\delta_y)\quad \forall (x,y)\in \Om.
$$
Therefore, parametrizing locally $\d\Om$ by a graph $(x,y_N(x))$, we obtain
$$
\|\phi_\viscosite\|_{L^2(\Om)}^2\leq C\int_{-\delta_x^- (1+|\ln \delta_y|^{-1})\leq x \leq  \delta_x^+ (1+|\ln \delta_y|^{-1})} \int_{-\delta_y}^{y_N(x)} (\delta_y^4+ \delta_y^2 y^2)\:dy\:dx.
$$
In order to bound the right-hand side of the above inequality, we have to distinguish between two cases, depending on whether $\cos \theta$ has a local extremum at $s=0$.
\begin{itemize}
\item
{\bf No local extremum at $s=0$:} either  $\cos \theta$ changes sign across $s=0$, or  $\cos \theta$ is identically zero in a neighbourhood on the left or on the right of $s=0$. Notice that in this case, $x_E(y)$ is only defined on a neighbourhood of the left of $y=0$, and thus we have $\delta_x^-=\delta_x^+=:\delta_x$.
Moreover  the function $y_N(x)$ remains non-positive for $|x|\leq   \delta_x(1+|\ln \delta_y|^{-1})$. In this case, we simply obtain
$$
\|\phi_\viscosite\|_{L^2(\Om)}^2\leq C \delta_y^5\delta_x.
$$

\begin{figure}
 \begin{center}
 	 \includegraphics[height = 6cm]{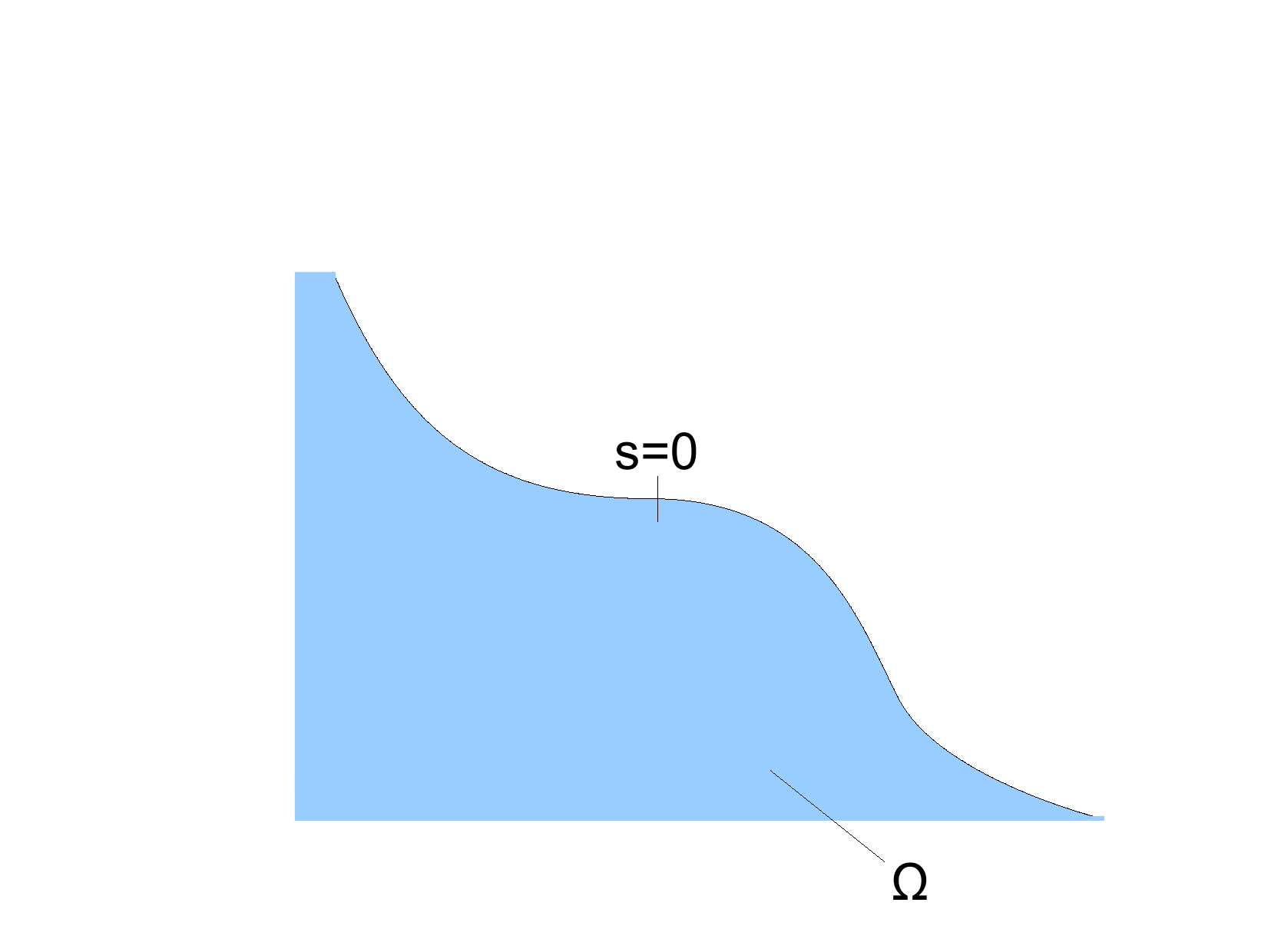}
\end{center}
\caption{$\cos \theta$ has a local minimum at $s=0$}\label{fig:inflexion}
\end{figure}

\item
{\bf Local extremum at $s=0$} (see Figure \ref{fig:inflexion}) : we have
$$
\|\phi_\viscosite\|_{L^2(\Om) }^2\leq C \delta_y^5 \delta_x^+ + C\int^0_{-\delta_x^-(1+|\ln\delta_y|^{-1})}(\delta_y^4 |y_N(x)| + \delta_y^2 |y_N(x)|^3)\:dx.
$$
We now estimate the right-hand side  using the formulas in Appendix B:
$$
\begin{aligned}
y_N(x)\sim C x^{n+1}\text{ in case (i)},\\
y_N(x)\sim C x^2 \exp\left(-\frac{\alpha}{|x|}\right)\text{ in case (ii)}.
\end{aligned}
$$
As a consequence,
$$
 \int^0_{-\delta_x^-(1+|\ln\delta_y|^{-1})}|y_N(x)|\:dx\leq C\left\{ \begin{array}{ll}
                        (\delta_x^-)^{n+2}&\text{ in case (i)},\\
                        (\delta_x^-)^4 \exp \left( -{\alpha \over \delta_x^- (1+|\ln \delta_y|^{-1})}\right)  &\text{ in case (ii)}\end{array}
\right.
$$
and
$$
 \int^0_{-\delta_x^-(1+|\ln\delta_y|^{-1})}|y_N(x)|^3\:dx\leq C\left\{ \begin{array}{ll}
                        (\delta_x^-)^{3n+4}&\text{ in case (i)},\\
                        (\delta_x^-)^8 \exp \left( -{3\alpha \over \delta_x^- (1+|\ln \delta_y|^{-1})}\right) &\text{ in case (ii)}\end{array}
\right.$$
\end{itemize}

We infer that in all cases,
$$
\|\phi_\viscosite\|_{L^2(\Om) }^2 \leq C  \delta_y^5( \delta_x^+ + \delta_x^-).
$$
An integration by part shows that
$$\left\|\tau \left(\chi_\viscosite-1\right) \right\|_{H^{-2}(\Om)}\leq \| \phi_\viscosite\tau\|_{L^2(\Om)} + 2 \|  \phi_\viscosite\d_y \tau\|_{L^2(\Om)} + \|  \phi_\viscosite \d_y^2\tau\|_{L^2(\Om)}$$
Then, by definition of $\delta_x^\pm$ and $\delta_y$, we get
\begin{equation}
\label{est-reste}
\left \|\tau \left(\chi_\viscosite-1\right) \right\|_{H^{-2}(\Om)}\leq C\|\tau\|_{W^{2,\infty }}\viscosite^{5/8} |\ln \viscosite|^{1/2} (\ln|\ln \viscosite|)^{-5\beta/2} \delta_x^{1/2} = o (\viscosite^{5/8}) ,
\end{equation}
since in the worse case
$$\delta_x = \max(\delta_x^+, \delta_x^-)\sim {\alpha \over |\ln \delta_y|} \sim {4\alpha \over |\ln \viscosite|}.$$
Note that the choice of $\delta_y$ is dictated by this estimate (which fixes both the power of $\viscosite$ and the power of $|\ln \viscosite|$).

\bigskip\noindent
$\bullet$
It remains then to estimate the viscous term $\viscosite \Delta^2 \psi_t^0$ in $H^{-2}(\Omega^\pm)$,  that is to estimate $\viscosite \Delta \psi_t^0$ in $L^2(\Omega^\pm)$. We indeed recall that, by Proposition \ref{H2-prop},
$$\Delta \psi^{int}= \chi^+\Delta\psi^{int} + \chi^-\Delta \psi^{int}.$$
 We then start from the identity
\begin{eqnarray}
        \Delta \psi^0_t &=&\d_x (\tau \chi_\viscosite) - \d_y^2 \int_x^{x_E(y)} \tau \chi_\viscosite\nonumber\\
        &=&\chi_\viscosite \d_x\tau + \tau \d_x\chi_\viscosite -\int_x^{x_E(y)} (\chi_\viscosite\d_y^2\tau + 2\d_y\tau\d_y\chi_\viscosite + \tau \d_y^2\chi_\viscosite)\label{delta_psi_1}\\
        &&-2x_E'(y) (\chi_\viscosite\d_y\tau + \tau\d_y\chi_\viscosite)_{|x=x_E(y)} - x_E''(y) (\tau\chi_\viscosite)(x_E(y),y)\label{delta_psi_2}\\
        &&- |x_E'(y)|^2\left(\chi_\viscosite \d_x\tau + \tau \d_x\chi_\viscosite\right)_{|x=x_E(y)}.      \label{delta_psi_3}
\end{eqnarray}

Using (\ref{chi-est}), we can easily check that
$$
\| (\ref{delta_psi_1})\|_{L^2} \leq \max_\pm\left(1+ \delta_y^{1/2}(\delta_x^\pm)^{-1/2}|\ln \delta_y| + \delta_x^\pm\delta_y^{-3/2}\right).
$$
Note that we use here in a crucial way the fact that there is no discontinuity line in the domains $\Omega^\pm$~: the singularity  is localized in a rectangle of width $\delta_x$.

As for \eqref{delta_psi_2}, \eqref{delta_psi_3}, we use \eqref{est:xE-suppchi} together with the formulas of Appendix B:
\begin{itemize}

\item {\bf Case (i)}: we deduce that
\begin{eqnarray*}
        \|\eqref{delta_psi_2}\|_2&\leq & \frac{C}{\delta_y}\left(\int_{|y|\geq \delta_y/2} |x_E'(y)|^2\:dy\right)^{1/2}+ C \left(\int_{|y|\geq \delta_y/2} |x_E''(y)|^2\:dy\right)^{1/2}\\
        &\leq & C \delta_y^{-\frac{n}{n+1}- \frac{1}{2}}=  C \delta_y^{-\frac{3n+1}{2(n+1)}},
\\      \|\eqref{delta_psi_3}\|_2&\leq &\max_\pm\frac{C}{\delta_x^\pm}\left(\int_{|y|\geq \delta_y/2} |x_E'(y)|^4\:dy\right)^{1/2}\\
        &\leq & C \delta_y^{-\frac{1}{n+1}} \delta_y^{\frac{1}{2}- \frac{2n}{n+1}}=  C \delta_y^{-\frac{3n+1}{2(n+1)}}.
\end{eqnarray*}
Hence the choice $\delta_y= \viscosite^{1/4} |\ln(\viscosite)|^{1/5}(\ln|\ln(\viscosite)|)^{-\beta}$, with $\beta>0$ arbitrary, ensures that $\delta \tau$ is an admissible remainder. Notice that any choice of the type $\delta_y= \viscosite^{\gamma}$, with $\gamma \in \left(\frac{5(n+1)}{4(5n+6)}, \frac{3(n+1)}{4(3n+1)}\right)$, also works in this case.

\item {\bf Case (ii)}: similarly, we have
\begin{eqnarray*}
\|\eqref{delta_psi_2}\|_2&\leq  & \frac{C}{\delta_y}\left(\int_{|y|\geq c\delta_y}\frac{1}{y^2(\ln|y|)^4}dy\right)^{1/2} + C \left(\int_{|y|\geq c \delta_y}\frac{1}{y^4(\ln|y|)^4}dy\right)^{1/2} \\
&\leq & \delta_y^{-3/2} |\ln|\delta_y||^{-2}\\
\|\eqref{delta_psi_3}\|_2&\leq  &       \max_\pm\frac{C}{\delta_x^\pm}\left(\int_{|y|\geq c\delta_y}\frac{1}{y^4(\ln|y|)^8}dy\right)^{1/2}\\
&\leq & C \delta_y^{-3/2} |\ln|\delta_y||^{-3}.
\end{eqnarray*}
Thus, with the choice  $\delta_y= \viscosite^{1/4} |\ln(\viscosite)|^{1/5}(\ln|\ln(\viscosite)|)^{-\beta}$, we infer
that  $\delta \tau$ is an admissible remainder.

\end{itemize}

\bigskip
This completes the proof of Proposition \ref{lem:trunc}.

\subsection{Error terms due to the lifting term $\psil$}$ $
\label{ssec:proof-eq-psil}
We give the rest of the proof of Proposition \ref{lem:eq-psil}. Since we have already checked that all the terms are well-defined, there only remains  to estimate the size of the remainders coming from $\psil$. More precisely, we prove that
\be\label{est:rlift}
\| r_{\mathrm{lift}}^1\|_{L^2(\Om)} =o(\viscosite^{1/8}), \quad \| r_{\mathrm{lift}}^2\|_{H^{-2}(\Om)} =o(\viscosite^{5/8}),
\ee
where
\begin{eqnarray*}
	r_{\mathrm{lift}}^1&=  &4 \viscosite^{3/4}\chi^+ \chi'\left(\frac{y-y_1}{\viscosite^{1/4}}\right)b''(x)+ 2 \viscosite^{1/2}\chi^+\chi''\left(\frac{y-y_1}{\viscosite^{1/4}}\right)(a'' (x)+ b''(x)(y-y_1))\\
	r_{\mathrm{lift}}^2&=&\tau (\chi_\viscosite - 1)+\viscosite\chi^+ \chi\left(\frac{y-y_1}{\viscosite^{1/4}}\right)(a^{(4)}(x)+ b^{(4)}(x)(y-y_1)) \\&&- \viscosite \Delta (\chi^+ \Delta \psi^0_t + \chi^- \Delta \psi^0_t ).
\end{eqnarray*}

We rely on Lemma \ref{lem:est-A-B} in Appendix C, which yields
\begin{eqnarray*}
 &&\left\|2 \viscosite^{1/2}\chi''\left(\frac{y-y_1}{\viscosite^{1/4}}\right)(a'' + b''(y-y_1)) + 4 \viscosite^{3/4} \chi'\left(\frac{y-y_1}{\viscosite^{1/4}}\right)b''\right\|_{L^2(\Om)}\\
&=& O(\viscosite^{5/8}(\|a''\|_{L^2}+ \viscosite^{1/4} \|b''\|_{L^2}))=o(\viscosite^{1/8}).
\end{eqnarray*}
As for $r_{\mathrm{lift}}^2$, we have already proved (see Proposition \ref{lem:trunc}) that $\tau(\chi_\viscosite-1)=o(\viscosite^{5/8})$ in $H^{-2}(\Om)$ and that $\|\viscosite \Delta \psi^0_t\|_{L^2(\Om^\pm)}= o(\viscosite^{5/8})$. And we have clearly
\begin{eqnarray*}
&&\viscosite\chi^+ \chi\left(\frac{y-y_1}{\viscosite^{1/4}}\right)(a^{(4)}(x)+ b^{(4)}(x)(y-y_1))\\&=&\d_x^2\left(\viscosite\chi^+ \chi\left(\frac{y-y_1}{\viscosite^{1/4}}\right)(a''(x)+ b''(x)(y-y_1))\right),
\end{eqnarray*}
so that
$$
\|r_{\mathrm{lift}}^2\|_{H^{-2}(\Om)}= o(\viscosite^{5/8}) + O\left(\viscosite^{9/8}(\|a''\|_{L^2}+ \viscosite^{1/4} \|b''\|_{L^2})\right)=o(\viscosite^{5/8}).
$$
Gathering all the terms, we obtain \eqref{est:rlift}.


\section{Remainders coming from the boundary terms}$ $

Following the  construction of the approximate solution, we now turn to the remainders coming from boundary layer terms. This step is more technical since boundary terms are defined through rescaled curvilinear coordinates, so that estimating the bilaplacian requires tedious computations.

Actually most of the terms will be estimated in $H^{-2}$ norm, so that we will need to compute only one iteration of the laplacian.

\subsection{Laplacian in curvilinear coordinates}$ $
\label{ssec:BL-est}
In order to get the boundary layer equation (\ref{BL}), we have kept only the leading order term of the bilaplacian, i.e. the fourth derivative with respect to the scaled normal variable $Z$. Rewriting this boundary layer equation as a Munk equation with remainder terms, we have to estimate the lower order terms of the bilaplacian, which can be done thanks to the following

\begin{Prop}
Denote as previously by $(z,s)$ the local coordinates on a tubular neighbourhood of the boundary.
Let $f\equiv F(s,\lambda(s) z)$ be a smooth function with $\Supp f \subset \Supp \chi_0$.

Then, we have the following error estimate
$$
\begin{aligned}
\left| \Delta f - \d_{zz} f\right| &\leq C\left( |\d_{ss} F| + (\lambda'/\lambda)^2 | \d_{ZZ} F| + |\lambda'/\lambda|\;|Z\d_{sZ} F|\right)\\
& + C \left( \left| {\theta''\over \lambda}\right| | Z \d_s F| +   \left| {\lambda' \theta''\over \lambda^2}\right| |Z^2 \d_Z F| +\left| {\lambda''\over \lambda}\right| |Z \d_Z F|+|\lambda \theta'| | \d_Z F|\right).
\end{aligned}
 $$

In particular,

\begin{itemize}
\item if  $\lambda$ does not depend on $s$
\begin{equation}
\label{NS-laplacian}
\left| \Delta f - \d_{zz} f\right| \leq C\left( |\d_{ss} F| + \left| {\theta''\over \lambda}\right| | Z \d_s F| +  |\lambda \theta'| | \d_Z F|\right);
\end{equation}

\item if $F$ is an exponential profile with respect to $Z$, namely if $F(s,Z) = A(s) \exp (-\mu Z)$, with $\mu\in \bC$ a constant of order one with positive real part,
\begin{eqnarray}
\label{EW-laplacian}
\left| \Delta f - \d_{zz} f\right| &\leq& C\left( |\d_{ss} A| +\left( \left| {\theta''\over \lambda}\right| + \left| \frac{\lambda'}{\lambda}\right|\right)|  \d_s A| \right)e^{-cZ}\\
&&+C  \left( \left| {\lambda' \theta''\over \lambda^2}\right|  +\left| {\lambda''\over \lambda}\right| +|\lambda \theta'|+\left( {\lambda'\over \lambda}\right)^2 \right) |A|e^{-cZ} \nonumber
\end{eqnarray}
where $c=\Re(\mu)/2$.
\end{itemize}
\label{prop:dev-lapla}
\end{Prop}

\begin{proof}

We start from the expression of the Laplacian in curvilinear coordinates
$$\Delta  = {1\over 1+z\theta'} \left[ \d_z ((1+z\theta') \d_z) +\d_s \left( {1\over 1+z\theta'} \d_s\right) \right]\,.$$
Given the specific form of the profile $f\equiv F(s,\lambda(s) z)$, we thus have
\begin{equation}
\label{delta}
\begin{aligned}
 \Delta f &= \lambda^2 \d_{ZZ}F + {\lambda \theta' \over 1+z\theta'} \d_Z F  - {z\theta'' \over (1+z\theta')^3} (\d_s F + \lambda' z \d_Z F)\\
& +{1\over (1+z\theta')^2} \left[ \d_{ss} F +2\lambda'z \d_{sZ} F +\lambda'' z\d_Z F+ (\lambda 'z)^2 \d_{ZZ} F\right].
\end{aligned}
\end{equation}
Then, using the fact that on the support of $\chi_0$ the jacobian of the change of variables is uniformly bounded from above and from below
$$ 0<\frac1C \leq 1+z\theta' \leq C,$$
we immediately get the first, general, estimate.

\medskip
For North, South   and discontinuity layers, we always choose
$$\lambda =\viscosite^{-1/4}$$
so that  $\lambda' =\lambda''=0$. The remainder consists then only on three terms.

\medskip
For East and West boundary layers, $f$ is a combination of exponential profiles.
Therefore, multiplying by $Z$ or differentiating with respect to $Z$ does not change the size of the error term.

\end{proof}

\bigskip

In order to get the expected estimates on the remainders, we first note that
$$\begin{aligned}
\Delta^2 &= (\lambda^2 \d_{ZZ} )^2 +  \Delta (\Delta- \lambda^2 \d_{ZZ}) +  (\Delta- \lambda^2 \d_{ZZ}) (  \lambda^2 \d_{ZZ})\\
& = (\lambda^2 \d_{ZZ} )^2 +  (\Delta+\lambda^2 \d_{ZZ})  (\Delta- \lambda^2 \d_{ZZ}) +[\Delta- \lambda^2 \d_{ZZ},\lambda^2 \d_{ZZ}]
\end{aligned}
$$
\begin{itemize}
\item
The first term of the right-hand side has been considered in the equation for the boundary layers, so it does not introduce any remainder.
\item
For the second term, we expect that it can be estimated in $H^{-2}$ using the previous proposition together with the controls on the boundary layer profiles and their derivatives up to second order.
\item
The difficulty comes then from the third term, which is a commutator between two second-order differential operators. Roughly speaking, it should be a third-order differential operator, so that we could estimate the corresponding remainder in $L^2$ in terms of the third-order derivative of the boundary layer profiles.
\end{itemize}

However we will not proceed exactly in that way, first of all because the commutator $[\Delta- \lambda^2 \d^2_{ZZ},\lambda^2 \d^2_{ZZ}]$ involves too many terms so that the computations would be very tedious, and overall because some terms are not small in $L^2$ and need to be estimated in $H^{-2}$. The proof will actually be  adapted to the case to be considered: recall indeed that
 for North, South   and discontinuity layers, $\lambda$ does not depend on $s$, whereas for East and West boundary layers, $F$ is an exponential profile with respect to $Z$.

A useful tool in the case when $\lambda$ does not depend on $s$  is the following commutator estimate~:

\begin{Lem}\label{commutator-lem}
Let $\rho \in H^2(\Om)$, with support in $\{d(x, \d\Om)\leq \delta\}$, and let $f\in L^2(\d \Om\times\R_+)$. We set $\lambda:=\viscosite^{-1/4}$.
Then
$$
\begin{aligned}
\left\| \rho(s,z) (\lambda^2\d_Z^2) f(s,\lambda z)\right\|_{H^{-2}(\Om)}\leq C|\lambda|^{-1/2}\|f\|_{L^2(\d\Om\times \R_+)} \sum_{k=0}^2\|\d_z^k \rho \|_{L^\infty(\Om)},\\
\left\| \rho(s,z) \Delta (f(s,\lambda z))\right\|_{H^{-2}(\Om)}
\leq  C|\lambda|^{-1/2}\|f\|_{L^2(\d\Om\times \R_+)} \|\rho\|_{W^{2,\infty}(\Om)}.
\end{aligned}
$$

\end{Lem}

\begin{proof}

$\bullet$ Let us first recall that $$\lambda \d_Z:f(s,\lambda Z) \in H^s(\R_+)\to \lambda  \d_Zf(s,\lambda Z)  \in  H^{s-1}(\R_+)$$ is a bounded operator (uniformly with respect to $\viscosite$) for all $s$.

Let $\varphi \in H^2(\Om )$ be an arbitrary test function such that $\varphi_{|\d \Om}= \d_n \varphi_{|\d \Om}=0$. Then, using the Cauchy-Schwarz inequality,
\begin{eqnarray*}
&&\int_{\Om}  \rho(s,z) (\lambda^2\d_Z^2) f(s,\lambda z) \varphi(s,z)(1+z\theta') \:ds\:dz\\
&=& \int_\Om f(s,\lambda z)\d_z^2 ( \rho(s,z) \varphi(s,z)(1+z\theta') )\:ds\:dz\\
&\leq & C |\lambda|^{-1/2}\|f\|_{L^2(\d\Om\times \R_+)} \sum_{k=0}^2\|\d_z^k (\rho\varphi) \|_{L^2(\Om)}\\
&\leq &  C|\lambda|^{-1/2}\|f\|_{L^2(\d\Om\times \R_+)} \|\varphi\|_{H^2(\Om)}\sum_{k=0}^2\|\d_z^k \rho \|_{L^\infty(\Om)}.
\end{eqnarray*}
By definition of $H^{-2}(\Omega)$,
\begin{eqnarray*}
&& \left\| \rho(s,z) (\lambda^2\d_Z^2) f(s,\lambda z)\right\|_{H^{-2}(\Om)}\\&=&\sup_{\substack{\varphi\in H^2_0(\Om),\\ \|\varphi\|_{H^2(\Om)}\leq 1}}\int_{\Om}  \rho(s,z) (\lambda^2\d_Z^2) f(s,\lambda z) \varphi(s,z)(1+z\theta') \:ds\:dz,
\end{eqnarray*}
from which we deduce  the first commutator estimate.

$\bullet$ Concerning the second inequality, we merely develop the commutator $[\Delta, \rho]$. We have
\begin{eqnarray*}
\rho(s,z) \Delta f(s,\lambda z)- \Delta(\rho(s,z)  f(s,\lambda z))&=& - 2 \nabla \rho \nabla(f(s,\lambda z) )- \Delta \rho f(s,\lambda z)\\
&=&-2 \Div (f(s,\lambda z)\nabla \rho ) + \Delta \rho f(s,\lambda z).\end{eqnarray*}
The estimate follows.

\end{proof}

\bigskip
When the profile is exponential, we will rather use the following Lemma~:

\begin{Lem}

Let $I\subset \Gamma_E\cup \Gamma_W$ be a closed interval, and let $\lambda\in L^\infty(I, \bC)$ such that $\inf_I \Re(\lambda)>0$.

 Let $\phi\in L^\infty(\Om)$ such that $\supp \phi\subset I$. Then
$$ \left \| \phi(s) \chi_0(z) \exp(- \lambda z)\right\|_{L^2(\Om)}\leq C \left(\int_I\phi(s)^2 (\Re\lambda(s))^{-1}\:ds\right)^{1/2}
$$
and
\begin{eqnarray*}
&& \left \| \phi(s) \chi_0(z) \exp(- \lambda z)\right\|_{H^{-2}(\Om)}\\
&\leq & C \left(\int_I\phi(s)^2 (\Re\lambda(s))^{-5}\:ds\right)^{1/2}+C \|\phi\|_{\infty} \exp(-C\inf_I \Re\lambda).
\end{eqnarray*}

\label{lem:est-BL-L2H2}

\end{Lem}

\begin{proof}

We begin with the $L^2$ estimate. We recall that the jacobian of the change of variables $(x,y)\to (s,z)$ is equal to $(1+z\theta')$, and is thus bounded  (from above and below) in $L^{\infty}$. As a consequence,
\begin{eqnarray*}
 &&\left \| \phi(s) \chi_0(z) \exp(- \lambda z)\right\|_{L^2(\Om)}^2\\&=& \int_{\d \Om} \int_0^\delta \phi(s)^2 \chi_0(z)^2 \exp(-2 \Re(\lambda) z)(1+z\theta')\:ds\: dz\\
&\leq&C \int_I\phi^2\left(\int_0^\infty \exp(-2 \Re(\lambda)z)\right)dsdz\\
&\leq & C  \int_I\phi^2 (\Re\lambda)^{-1}.
\end{eqnarray*}

As for the $H^{-2}$ estimate, we have
\begin{eqnarray*}
\phi(s)\chi_0(z)\exp(-  \lambda z)&= & \frac{\d^2}{\d z^2} \left( \phi(s) \chi_0(z) (  \lambda(s))^{-2}\exp(-  \lambda z)\right)\\
&&+ 2 \phi(s) \chi_0'(z)  (  \lambda)^{-1}\exp(-  \lambda z)\\
&&- \phi(s) \chi_0''(z) \exp(-  \lambda z).
\end{eqnarray*}
The last two terms are supported in $\{\delta/2\leq z \leq \delta\}$, hence they are exponentially small in $L^2$. More precisely, they are bounded in $L^\infty(\Om)$ by
$$
C\|\phi\|_\infty \exp(-\delta\inf_I \Re(\lambda)/4).
$$
As for the first term, if $\zeta\in H^2_0(\Om)$ is an arbitrary test function, then
\begin{eqnarray*}
&& \int_{\Om}  \frac{\d^2}{\d z^2} \left( \phi(s) \chi_0(z) (  \lambda(s))^{-2}\exp(-  \lambda z)\right) \zeta(s,z)(1+z\theta')ds\:dz \\
&=& \int_{\d \Om} \int_0^\delta  \frac{\d^2}{\d z^2} \left( \phi(s) \chi_0(z) (  \lambda(s))^{-2}\exp(-  \lambda z)\right) \zeta(s,z)(1+z\theta')ds\:dz\\
&=& \int_{\d \Om} \int_0^\delta \phi(s) \chi_0(z) (  \lambda)^{-2}\exp(-  \lambda z) \frac{\d^2}{\d z^2}\left(\zeta(s,z)(1+z\theta')\right)ds\:dz\\
&\leq&C\|\zeta\|_{H^2(\Om)}\left(\int_{\d \Om}\int_0^\delta \phi(s)^2|\lambda|^{-4}\exp(-2  \Re(\lambda) z)ds\:dz\right)^{1/2}\\
&\leq & C\|\zeta\|_{H^2(\Om)} \left(\int_{\d \Om}\phi^2(\Re\lambda)^{-5}\right)^{1/2}.
\end{eqnarray*}

By definition of the $H^{-2}$ norm, we infer the estimate of Lemma \ref{lem:est-BL-L2H2}.

\end{proof}

\subsection{Error terms associated with North and South    layers}$ $

In order to simplify the presentation, we focus on a South boundary layer, as we did in the preceding chapter. Of course, the case of North boundary layers is strictly identical.
We denote by $s_{i-1}, s_i$  the end points of the South boundary under consideration, and we introduce a macroscopic truncation function $\tilde \gamma_{i-1,i}$ (with bounded derivatives) as in \eqref{def:tilde-gamma-i}.

\begin{Lem} Let $\psi_{S}$ be the solution to \eqref{eq:South}.

Then
$$
\left[\d_x - \viscosite \Delta^2\right] \left(\tilde \gamma_{i-1,i}(s) \psi_S(s,\viscosite^{-1/4} z) \chi_0(z) \right)= r_i^1 + r_i^2,
$$
with
$$
\| r_i^1\|_{L^2(\Om)}=o(\viscosite^{1/8}), \quad \| r_i^2\|_{H^{-2}(\Om)}=o(\viscosite^{5/8}).
$$

\label{lem:approx-BL-NS}
\end{Lem}

\begin{proof}[Proof of Lemma \ref{lem:approx-BL-NS}]

Since $\psi_S$ satisfies \eqref{eq:South}, there are several kinds of terms  in
$$
\left[\d_x - \viscosite \Delta^2\right] \left(\psi_S(s,\viscosite^{-1/4}  z)\tilde \gamma_{i-1,i}(s) \chi_0(z) \right):
$$

\medskip
$\bullet$ All terms in which at least one derivative of $\chi_0$ or $\tilde \gamma_{i-1,i}$ occurs are  small: indeed, we have one derivative with respect to $z$  less acting on $\psi_S$,
so that we gain a power $\viscosite^{1/4}$. By Lemma \ref{commutator-lem}, we indeed have
\begin{eqnarray}
 \nonumber \viscosite \d_z^4(\psi_S \chi_0) -\viscosite \chi_0 \d_z^4 \psi_S&=&\viscosite  \sum  _{j=0}^3 C^j_4 \lambda^j \d_Z^j \psi_{S}(s,\lambda z) \d_z^{4-j}\chi_0\\
\label{reste}&=&O(\viscosite^{5/8}\| \psi_S\|_{L^2_sH^2_Z})\text{ in }L^2(\Om) \\&&+ O(\viscosite^{7/8} \|\d_Z \psi_S\|_{L^2_{s,Z}}) \text{ in }H^{-2}(\Om).\nonumber
\end{eqnarray}
In the same way,
\begin{eqnarray*}
 &&\psi_S(s,\viscosite^{-1/4} z) \d_x (\chi_0(z) \tilde \gamma_{i-1,i}(s))\\&=& \psi_S(s,\viscosite^{-1/4} z) \left[\frac{\sin \theta(s)}{1+ z \theta'} \d_s\tilde \gamma_{i-1,i}(s) \chi_0(z) + \cos \theta(s) \tilde \gamma_{i-1,i}(s)\d_z \chi_0(z) \right].
\end{eqnarray*}
Using the estimates of paragraph \ref{ssec:extinction}, we infer that
$$
\left\| \psi_S(s,\viscosite^{-1/4} z)\cos \theta(s) \tilde \gamma_{i-1,i}(s)\d_z \chi_0(z) \right\|_{L^2(\Om)}=O(\viscosite^{1/8}|\ln \viscosite|^{-1}).
$$
Moreover, on the support  of $\d_s \tilde \gamma_{i-1,i}$,   either $\psi_S$ is zero  or $\|\psi_S\|_{L^\infty_s(L^2_Z)}=O(|\ln \viscosite|^{-1})$. We deduce that
$$
\left\| \psi_S(s,\viscosite^{-1/4} z)\sin \theta(s) \d_s\tilde \gamma_{i-1,i}(s) \chi_0(z)  \right\|_{L^2(\Om)}=O(\viscosite^{1/8}|\ln \viscosite|^{-1}),
$$
so that eventually
\be
\| \psi_S(s,\viscosite^{-1/4} z) \d_x (\chi_0(z) \tilde \gamma_{i-1,i}(s))\|_{L^2(\Om)}= O(\viscosite^{1/8}|\ln \viscosite|^{-1}).\label{reste00}
\ee

\medskip
$\bullet$  The terms stemming from the $\d_x$ derivative applied to $\psi_S$ can be split into a remainder term coming from the jacobian
$$
\frac{z\theta' \sin \theta}{1+z \theta'} \chi_0(z) \tilde \gamma_{i-1,i}(s)\d_s  \psi_S(s,\viscosite^{-1/4} z),
$$
and
$$
\chi_0(z) \tilde \gamma_{i-1,i}(s) \left( \viscosite^{-1/4} \cos \theta \d_Z  \psi_S(s,\viscosite^{-1/4} z)+\sin \theta \d_s  \psi_S(s,\viscosite^{-1/4} z)\right) ,
$$
which will simplify with the term $\viscosite \lambda^4  \chi_0(z)  \tilde \gamma_{i-1,i}(s) (\d_Z^4)\psi_S(s,\viscosite^{-1/4} z)$ coming from the bilaplacian.

From the moment estimate in (\ref{moment}), we deduce that
\begin{equation}
\label{reste0}
\left\|\frac{z\theta' \sin \theta}{1+z \theta'} \tilde \gamma_{i-1,i}(s) \chi_0(z)\d_s  \psi_S(s,\viscosite^{-1/4} z) \right\|_{L^2}= O (\viscosite^{1/8}\viscosite^{1/4} \delta_y^{-1} )=o(\viscosite^{1/8}).
\end{equation}

\medskip
$\bullet$ The most technical part is the control of the remainder terms coming from the bilaplacian. For the sake of simplicity, instead of commuting the whole laplacian with $\lambda^2 \d_Z^2$, we will first get rid of the corrections coming from the jacobian, then commute the remaining part of the laplacian. More precisely, we start from (\ref{delta}) with a constant $\lambda$
$$\begin{aligned}
&\Delta^2 - (\lambda^2 \d^2_{Z} )^2 \\
&=  \Delta (\Delta- \lambda^2 \d^2_{Z}) +(\Delta- \lambda^2 \d^2_{Z})\lambda^2 \d^2_{Z}\\
&= \Delta (\Delta- \lambda^2 \d_{ZZ})  +\left(  {\lambda \theta' \over 1+z\theta'} \d_Z   - {z\theta'' \over (1+z\theta')^3} \d_s  +{1\over (1+z\theta')^2}  \d_{ss}  \right) \lambda^2 \d_{ZZ}\\
&= \Delta (\Delta- \lambda^2 \d_{ZZ}) +\lambda^2 \d_{ZZ}  \Big( \lambda \theta' \d_Z - z\theta'' \d_s +\d_{ss}\Big) + 2\lambda \theta'' \d_{Zs} \\&  -\left(  {\lambda z( \theta')^2 \over 1+z\theta'} \d_Z   - {z\theta'' ((1+z\theta')^3-1)) \over (1+z\theta')^3} \d_s  +{(1+z\theta')^2-1\over (1+z\theta')^2}  \d_{ss}  \right) \lambda^2 \d_{ZZ}.\end{aligned}
$$

We have, using the regularity and moment estimates (\ref{moment})
$$
\begin{aligned}
\viscosite &(\Delta- \lambda^2 \d_{ZZ})\psi_S\\
 &=   \viscosite {\lambda \theta' \over 1+z\theta'} \d_Z \psi_S  - \viscosite {\lambda^{-1}Z \theta'' \over (1+z\theta')^3} \d_s  \psi_S+{\viscosite \over (1+z\theta')^2}  \d_{ss} \psi_S\\
& = O( \viscosite \viscosite^{-1/4} \viscosite^{1/8})_{ L^2} + O( \viscosite \viscosite^{1/4}\delta_y^{-1} \viscosite^{1/8})_{L^2}  + O( \viscosite\viscosite^{-1/2}|\ln \viscosite|^{-1} \viscosite^{1/8} )_{L^2}.
\end{aligned}
$$
We now use the second commutator estimates stated in Lemma \ref{commutator-lem} to get
\be \label{reste1}
\|\chi_0(z) \tilde \gamma_{i-1,i}(s) \viscosite \Delta (\Delta- \lambda^2 \d_{ZZ})\psi_S (s,\lambda z)\|_{H^{-2}(\Om)}=o(\viscosite^{5/8}).
\ee

In the same way, we have
$$\begin{aligned}
 \viscosite& \lambda \theta' \d_Z \psi_S - \viscosite z\theta'' \d_s\psi_S +\viscosite \d_{ss}\psi_S \\
 &= O( \viscosite \viscosite^{-1/4} \viscosite^{1/8})_{ L^2} + O( \viscosite \viscosite^{1/4}\delta_y^{-1} \viscosite^{1/8})_{L^2}  + O( \viscosite\viscosite^{-1/2}|\ln \viscosite|^{-1} \viscosite^{1/8}  )_{L^2}
 \end{aligned}$$
so that, using some integration by parts together with the first estimate of Lemma \ref{commutator-lem},
\begin{equation}
\label{reste2}
\viscosite\chi_0(z)\tilde \gamma_{i-1,i}(s)  \lambda^2 \d_{ZZ}  ( \lambda \theta' \d_Z - z\theta'' \d_s +\d_{ss})\psi_S = o(\viscosite^{5/8})_{H^{-2}}\,.
\end{equation}

For the corrections coming from the jacobian, we use the first estimate of the commutator lemma \ref{commutator-lem}: {\small$$
\begin{aligned}
&\left\| \chi_0 \tilde \gamma_{i-1,i}\left(  {\lambda z( \theta')^2 \over 1+z\theta'} \d_Z   - {z\theta'' ((1+z\theta')^3-1)) \over (1+z\theta')^3} \d_s  +{(1+z\theta')^2-1\over (1+z\theta')^2}  \d_{ss}  \right) \lambda^2 \d_{ZZ}\psi_S\right\|_{H^{-2}} \\
&= O( \viscosite^{1/8}) \| Z \d_Z \psi_S\|_{L^2_{s,Z}} \left\| {( \theta')^2 \over 1+z\theta'}\chi_0 \tilde \gamma_{i-1,i}\right\|_{W^{2,\infty}(\Om)} \\
&+ O( \viscosite^{1/8}\viscosite^{1/2} ) \| Z^2 \d_s \psi_S\|_{L^2_{s,Z}} \left\| {\theta''( ( 1+z \theta')^3 - 1) \over z( 1+z\theta')^3}\chi_0 \tilde \gamma_{i-1,i}\right\|_{W^{2,\infty}(\Om)} \\
& + O( \viscosite^{1/8} \viscosite^{1/4} ) \| Z \d_s^2  \psi_S\|_{L^2_{s,Z}} \left\| {\theta''( ( 1+z \theta')^2 - 1) \over z( 1+z\theta')^2}\chi_0 \tilde \gamma_{i-1,i}\right\|_{W^{2,\infty}(\Om)}
\end{aligned}
$$}
which implies that
{\small \begin{multline}
\label{reste3}
\Big\| \viscosite \chi_0 \tilde \gamma_{i-1,i}\left(  {\lambda z( \theta')^2 \over 1+z\theta'} \d_Z   - {z\theta'' ((1+z\theta')^3-1)) \over (1+z\theta')^3} \d_s  +{(1+z\theta')^2-1\over (1+z\theta')^2}  \d_{ss}  \right) \lambda^2 \d_{ZZ}\psi_S\Big\| _{H^{-2}}  \\= o(\viscosite^{5/8})\,.
\end{multline}}

The only remaining term is
\begin{equation}
\label{reste4}
2\viscosite \lambda \theta'' \chi_0(z) \tilde \gamma_{i-1,i}(s)\d_{Zs} \psi_S =  O(\viscosite \viscosite^{-1/4} \delta_y^{-1} \viscosite^{1/8})_{L^2}.
\end{equation}

\bigskip
Combining estimates (\ref{reste})-(\ref{reste4}) shows that the remainders coming from the South and North boundary layers are admissible in the sense of Definition \ref{def:admissible}.
\end{proof}

\subsection{Error terms associated with East and West boundary layers}$ $

We recall that the boundary layer term on the East or West boundary $[s_i,s_{i+1}]$  takes the form
$$
\psi_{E,W}(s,Z)=(\varphi_i^+\varphi_{i+1}^-) (s) A(s)^t f(Z),
$$
where
\begin{itemize}
\item
 $A(s)=-\lambda_E^{-1}(s) \d_n \psi^0_{t|\Gamma_E}(s)$ and $f(Z) =\exp(-Z)$ on East coasts,
 \item  $A$ is the matrix defined by (\ref{def:AW-fin}) and
 $$f(Z)=\begin{pmatrix}
        \exp\left(-e^{i\pi/3} Z\right)\\  \exp\left(-e^{-i\pi/3} Z\right)
       \end{pmatrix}$$
on West coasts,
\item $\varphi_i^\pm$ are defined by (\ref{Evarphii}) on the East coast, and by (\ref{Wsigmaipm})-(\ref{Wvarphii}) on the West coast (see also Lemma \ref{sigma-s-est-bis}).
\end{itemize}


Combining the definitions of $s_i^\pm$, $\sigma_i^\pm$ and the equivalents given in Lemma \ref{sigma-s-est} and Appendix B, we can check that  the eigenvalues $\lambda_{E,W}$ satisfy the following estimates on the domain of validity of the East and West boundary layers:

\begin{Lem}\label{lEW-est}
Under  the compatibility condition \eqref{hyp:valid-E/W-2}, for $s\in \supp \varphi_i,$
$$
\left|\frac{\lambda'}{\lambda^3}\right|\ll 1,\ \left|\frac{\lambda''}{\lambda^2}\right|\ll 1,\ \left|\frac{(\lambda')^2}{\lambda^3}\right|
\ll 1,
$$
and
$$
\left|\viscosite \lambda'\lambda^{1/2}\right|\leq C |\theta'| \; \left|\frac{\viscosite}{\cos \theta}\right|^{1/2}.
$$
\end{Lem}

We can then prove that the error terms associated with East and West boundary layers are admissible.

\begin{Prop} Let $A,f$ be defined by the expressions above and \eqref{def:AW-fin}.

Then
$$
\left[\d_x - \viscosite \Delta^2\right] \Big(\psi_{E,W}(s,\lambda_{E,W} z)\chi_0(z) \Big)=r_{E,W}^1 + r_{E,W}^2,
$$
with
$$
\|r^1_{E,W}\|_{L^2(\Om)}=o(\viscosite^{1/8}),\\
\|r^2_{E,W}\|_{H^{-2}(\Om)}=o(\viscosite^{5/8}).
$$
\label{lem:error-EW}

\end{Prop}

\begin{proof}
Throughout the proof, we write $\lambda$ instead of $\lambda_{E,W}$, and we set $\varphi_i= \varphi_i^+\varphi_{i+1}^-$.
Since $f$ is an exponential profile, we will use the second part of Proposition \ref{prop:dev-lapla} to estimate the remainder terms in the bilaplacian.
In a first step, we bound the error terms $r_{E,W}^1$ and $r_{E,W}^2$ by expressions involving $\viscosite, \theta$ and $A$, and in a second step, we prove that these bounds yield the desired estimates.

\noindent
\textbf{ First step: estimates in terms of the amplitude $A$}

We claim that for all $i$ such that $(s_i, s_{i+1})\subset \Gamma_W \cap \Gamma_E$,
$$
(\d_x- \viscosite \Delta^2) (\varphi_i(s) A(s) f(\lambda z))= r_i^1 + r_i^2,
$$
where

\begin{eqnarray}
\nonumber\| r_i^1\|_{L^2(\Om)}&\leq & C\viscosite^{5/6} \| (\varphi_i A)''\|_{L^2(\d \Om)}\\
\label{borne-ri1}&&+C\viscosite^{1/2} |\ln \viscosite|^{1/2}\| (\varphi_i A)'\|_{L^\infty(\d \Om)}\\
\nonumber&&+ C \viscosite^{1/6} \|A\|_{L^\infty(\d \Om)} +C \viscosite^{-K} \exp(-C/\viscosite^{1/4}),
\end{eqnarray}
and
\begin{eqnarray}
\nonumber\| r_i^2\|_{H^{-2}(\Om)}&\leq & C\viscosite \left(\int(\d^2_s(\varphi_i A))^2 \left| \frac{\viscosite}{\cos (\theta(s))}\right|^{1/3}\:ds\right)^{1/2}\\
\label{borne-ri2}&&+C \left(\int(\d_s(\varphi_i A))^2 \left| \frac{\viscosite}{\cos (\theta(s))}\right|^{5/3}\:ds\right)^{1/2}\\
\nonumber&&+ C  \left(\int_{\supp \varphi_i} A(s)^2 ( \theta')^2 \viscosite^{5/3} |\cos \theta|^{-11/3} \:ds\right)^{1/2}\\
\nonumber&&+C\|A\|_\infty \viscosite^{5/6} + C \viscosite^{-K} \exp(-C/\viscosite^{1/4}).
\end{eqnarray}

We recall the expression of $\d_x$ in curvilinear coordinates:
$$
 \frac{\d}{\d x}= - \cos \theta  \frac{\d}{\d z} + \frac{\sin \theta }{1+z \theta'} \frac{\d}{\d s}
$$

We now estimate
$$
[\d_x-\viscosite \Delta^2]\left(\varphi_i(s) A(s) \exp(-\lambda(s) z)\chi_0(z)\right).
$$

$\bullet$ We start by commuting the differential operator $\d_x-\viscosite \Delta^2$ with the multiplication by $\chi_0$.
Since $\chi_0(z)\equiv 1$ for $z$ in a neighbourhood of zero, all terms in which at least one derivative of $\chi_0$ appears are exponentially small: indeed,  for all $k\geq 1$, and for $\mu\in \{1, e^{\pm i \pi/3}\}$,
$$
\left| \d_z^k\chi_0 (z) \exp\left(- \mu \lambda(s) z\right)\right|\leq C_k \exp(-C\viscosite^{-1/4})\quad \forall (s,z)\in \supp \varphi_i\times \supp \chi_0.
$$
We also have (see Lemmas \ref{lem:est-psi0},  \ref{lem:est-psi-sigma}, \ref{lem:psi-corr})
$$
\left|\d_s^k, \d_z^k (\varphi_i(s) A(s) \exp(-\mu \lambda(s) z))\right|\leq C \viscosite^{-K} \exp(- c\lambda(s) z)),
$$
for some positive constants $c,C, K$. Hence we infer that
\begin{eqnarray*}
&&[\d_x-\viscosite \Delta^2]\left(\varphi_i(s) A(s) \exp(-\lambda(s) z)\chi_0(z)\right)\\&=& \chi_0(z)[\d_x-\viscosite \Delta^2]\left(\varphi_i(s) A(s) \exp(-\lambda(s) z)\right)\\
&&+ O(\viscosite^{-K}  \exp(-C\viscosite^{-1/4}))\text{ in }L^2(\Om).
\end{eqnarray*}

$\bullet$  The terms stemming from $\chi_0\d_x(\varphi_i(s) A(s) \exp(-\lambda(s) z))$ can be split into
$$
 \lambda\cos \theta \varphi_i(s) A(s) (\d_Z f)(\lambda z)\chi_0,
$$
which will simplify with the term $\viscosite \lambda^4  \varphi_i(s) A(s) (\d_Z^4)f(\lambda z)\chi_0$ coming from the bilaplacian, and a remainder term
\be\label{resteEW}
\frac{\sin \theta}{1+z \theta'} \left[\d_s(\varphi_i A) f(\lambda z) + A(s) \varphi_i(s)z \lambda' (\d_Zf)(\lambda z)\right]\chi_0.
\ee

By Lemma \ref{lem:est-BL-L2H2}, replacing $(1+z\theta')^{-1}$ by $1$ in \eqref{resteEW} yields an error term whose square $L^2$ norm is bounded by
\begin{eqnarray*}
&&C \int_{\supp \varphi_i}  \theta'^2 \frac{\viscosite}{|\cos\theta|}  \left( |\d_s(A \varphi_i)|^2 +\frac{\theta'^2}{\cos^2\theta}  |A(s)|^2 \right) \: ds\\
&\leq &C\viscosite |\ln \viscosite| \|\d_s(\varphi_i A)\|_{L^\infty(\d \Om)}^2+ C \viscosite^{5/7} \|A\|_{L^\infty(\d \Om)}^2.
\end{eqnarray*}
Indeed, by definition of $\varphi_i$,
$$
\supp \varphi_i \subset  (s_{i}^+, s_{i+1}^-) ,
$$
and using assumption (H2), it can be easily proved that for all $i$ such that $(s_i, s_{i+1})$ is a western or eastern boundary,
\be
\label{in:theta'}
\begin{aligned}
   \int_{s_{i}^{+}}^{s_{i+1}^{-}}  \theta'^2 |\cos \theta|^{-1}= O(|\ln \viscosite|),\\
  \int_{s_{i}^{+}}^{s_{i+1}^{-}}  \theta'^4 |\cos \theta|^{-3}= O (\viscosite^{-2/7}).
  \end{aligned}
\ee

Then, using Lemma \ref{lem:est-BL-L2H2}, we infer that
\begin{eqnarray*}
&&\left\| \sin\theta\chi_0(z)\d_s(\varphi_i A) f(\lambda z) \right\|_{H^{-2}(\Om)}^2\\
&\leq &C \int_{\d\Om}(\d_s(\varphi_i A))^2 \left| \frac{\viscosite}{\cos (\theta(s))}\right|^{5/3}\:ds+ O(\viscosite^{-K}  \exp(-C\viscosite^{-1/4}))
\end{eqnarray*}
and
\begin{eqnarray*}
&& \left\| \sin\theta\chi_0(z)A(s) \varphi_i(s)z\lambda' f(\lambda z) \right\|_{H^{-2}(\Om)}^2\\
&\leq & C  \int_{\supp \varphi_i}A(s)^2 \theta'^2 \viscosite^{5/3} |\cos \theta|^{-11/3} \:ds+ O(\viscosite^{-K}  \exp(-C\viscosite^{-1/4})).
\end{eqnarray*}

$\bullet$ We now develop the bilaplacian of $(\varphi_i A)(s)f(\lambda z)$, using Proposition \ref{prop:dev-lapla}. We bound all the derivatives of $\theta$ by a positive constant. Moreover, by Lemma \ref{lEW-est} and Proposition \ref{prop:dev-lapla},
\begin{eqnarray*}
&&\left| (\Delta-\lambda^2\d_{ZZ})(\varphi_i A)(s)f(\lambda z) \right|\\
&\leq & C \left(|(A\varphi_i)''| + \left|\frac{\lambda'}{\lambda}  \right| \;| (A\varphi_i)'| + |\lambda A \varphi_i|\right)\exp(-c|\lambda| z).
\end{eqnarray*}
We now evaluate the right-hand side of the above inequality in $L^2(\Om)$ (after multiplication by $\chi_0$). By Lemma \ref{lem:est-BL-L2H2}, we have
$$
\left\|  (\lambda A \varphi_i)(s)\exp(-c|\lambda| z)\chi_0(z)\right\|_{L^2(\Om)}\leq C \|A\varphi_i \lambda^{1/2}\|_{L^\infty}\leq C \|A\|_\infty \viscosite^{-1/6}.
$$
Therefore
$$
\begin{aligned}
\Delta  &\left( (A\varphi_i)(s)f(\lambda z)\right)
= (A\varphi_i)(s)\lambda^2 f''(\lambda z) \\
&+O\left( \|A\|_{L^\infty}\viscosite^{-1/6} +\|  (A\varphi_i)'\lambda'\lambda^{-3/2}\|_{L^2} + \|(A\varphi_i)''\lambda^{-1/2}\|_{L^2}\right)_{L^2}
\end{aligned}
$$
Commuting once again $\chi_0$ and $\Delta$, we infer that
$$
\begin{aligned}
\viscosite \chi_0(z)\Delta^2 & \left( (A\varphi_i)(s)f(\lambda z)\right)\\
&= \viscosite \chi_0(z) \Delta\left((A\varphi_i)(s)\lambda^2 f''(\lambda z) \right)+ O( \exp (-C\viscosite^{1/4}))_{L^2}\\
&+O(\|A\|_{\infty}\viscosite^{5/6}+\viscosite\|(A\varphi_i)'\lambda'\lambda^{-3/2}\|_{L^2}+\viscosite\|(A\varphi_i)''\lambda^{-1/2}\|_{L^2})_{H^{-2}}
\end{aligned}
$$
 Lemma \ref{lEW-est} then implies
$$
\viscosite\|(A\varphi_i)'\lambda'\lambda^{-3/2}\|_{L^2}=o\left( \left(\int_{ \d \Om}(\d_s(\varphi_i A))^2 \left| \frac{\viscosite}{\cos (\theta(s))}\right|^{5/3}\:ds\right)^{1/2}\right).
$$

It remains then to deal with the first term $ \viscosite \chi_0(z) \Delta\left((A\varphi_i)(s)\lambda^2 f''(\lambda z) \right)$.
We  develop the laplacian one more time. Since $f$ is an exponential, we can use the computations above and simply replace $A$ by $A \lambda^2$.
Then, by Lemma \ref{lEW-est},
\begin{eqnarray*}
&&\viscosite\chi_0(z) \Delta \left( (A\varphi_i)(s)\lambda^2 f''(\lambda z) \right)\\
&=&\viscosite\chi_0(z)  (A\varphi_i)(s)\lambda^4 f^{(4)}(\lambda z) +  \frac{\viscosite\chi_0(z)}{(1+z \theta')^2}(A\varphi_i)'' \lambda^2 f''(\lambda z) \\
&&+ \viscosite O(\|A\varphi_i\lambda^{5/2}\|_{L^2}+ \| (A\varphi_i)'\lambda'\lambda^{1/2}\|_{L^2})\text{ in }L^2(\Om).
\end{eqnarray*}
We estimate the second term of the right-hand side in $H^{-2}(\Om)$. This is the only time where we really need to commute the jacobian terms with the $z$ derivatives, as explained at the end of paragraph \ref{ssec:BL-est}. Using Lemma \ref{lem:est-BL-L2H2}, we have
\begin{eqnarray*}
&&\frac{\chi_0(z)}{(1+z \theta')^2}(A\varphi_i)'' \lambda^2 f''(\lambda z)\\
&=&\frac{\d^2}{\d z^2} \left(\frac{\chi_0(z)}{(1+z \theta')^2}(A\varphi_i)'' f(\lambda z)\right)\\
&&+4 \frac{\chi_0(z)\theta'}{(1+z\theta')^3}\lambda f'(\lambda z) (A\varphi_i)'' -6 \frac{\chi_0(z){\theta'}^2}{(1+z\theta')^4} f(\lambda z) (A\varphi_i)''\\
&&+O\Big( \|(A\varphi_i)'' \lambda^{1/2}\|_{L^2}\exp(-C/\viscosite^{1/4})\Big)_{L^2(\Om)}\\
&=&O\Big(\|(A\varphi_i)'' \lambda^{-1/2}\|_{L^2}\Big ) _{H^{-2}(\Om)}+ O\Big( \|(A\varphi_i)'' \lambda^{1/2}\|_{L^2}\Big)_{L^{2}(\Om)}.
\end{eqnarray*}
Inequality \eqref{in:theta'}  leads  then to
$$
 \| (A\varphi_i)'\lambda'\lambda^{1/2}\|_{L^2}\leq C \| (A\varphi_i)'\|_\infty \viscosite^{1/2} |\ln \viscosite|^{1/2}.
$$
Gathering all the terms, we obtain \eqref{borne-ri1}, \eqref{borne-ri2}.

\bigskip

\noindent
\textbf{ Second step: quantitative bounds for $r_i^1$, $r_i^2$:}

We now use the definition of $A$ together with Lemma \ref{lem:est-psi0} and Lemma \ref{lem:est-A} in order to estimate the right-hand sides of \eqref{borne-ri1}, \eqref{borne-ri2}. Notice that the estimates on West coasts are always more singular than the ones on East coasts: indeed, $\varphi_i$ has bounded derivatives on East coasts, and unbounded on West coasts. Moreover, the estimates of Lemma \ref{lem:est-psi0} show that $A$ is always larger on West coasts, together with its derivatives. Therefore, we focus on a West part of the boundary in the following, i.e. we assume that $(s_i, s_{i+1}) \subset\Gamma_W$.

Following Lemma \ref{lem:est-A}, we treat separately the estimates in the vicinity of $s_1^W$ and in the vicinity of points $s$ such that $\cos \theta=0$ or $y(s)=y_j$ for some $j\in I_+$. In the rest of the domain, the estimates on $r_i^1$ and $r_i^2$ are easily proved.

\medskip
$\bullet$ {\bf Estimate near points such that $|y(s)-y_j|\ll 1$ ($y_j\neq y_1$)}:

Let us first recall (see Lemma \ref{sigma-s-est-bis})
$$
\|\d_s^k\varphi_i\|_\infty\leq C \viscosite^{-k/7}\quad k\in \{0,1,2\}.
$$
Therefore
\begin{multline*}
\|r_i^1\|_{L^2}\\ \leq C \left(\viscosite^{1/6}\|A\|_{\infty} + \viscosite^{1/2}|\ln \viscosite|^{1/2}\|A'\|_\infty+ \viscosite^{5/6}\|A''\|_2+  \viscosite^{-K} \exp(-C/\viscosite^{1/4})\right).
\end{multline*}
Inequalities \eqref{est:A}, \eqref{est:dsA}, \eqref{est:ds2A} yield respectively
\be\label{est:A-infty}\|A\|_{\infty}=O(1),\quad\|A'\|_\infty= (\delta_y^{-1}|\ln \delta_y|^{-1}),\quad\|A''\|_\infty=O(\delta_y^{-2}|\ln \delta_y|^{-1}),
\ee
so that $\|r_i^1\|_{L^2}=o(\viscosite^{1/8})$. We now address the bound of $r_i^2$. Using once again \eqref{est:A-infty}, we have
\begin{eqnarray*}
 &&\viscosite \left(\int_{\d \Om}(\d^2_s(\varphi_i A))^2 \left| \frac{\viscosite}{\cos (\theta(s))}\right|^{1/3}\:ds\right)^{1/2}\\
&\leq & C\viscosite^{7/6}\left(\|\varphi_i''\|_\infty + \|\varphi_i'\|_\infty \|A'\|_\infty + \|A''\|_\infty \right)\left(\int \mathbf{1}_{\supp \varphi_i} |\cos \theta|^{-1/3}\right)^{1/2}\\
&\leq & C \viscosite^{2/3}|\ln \viscosite|^{-1}\left(\int \mathbf{1}_{\supp \varphi_i} |\cos \theta|^{-1/3}\right)^{1/2}.
\end{eqnarray*}
Using the formulas in Appendix B together with the definition of $s_i^+, s_{i+1}^-$, it can be proved (treating cases (i) and (ii) in assumption (H2) separately) that
\be\label{est:sing-cos-1}
\int_{\supp \varphi_i}|\cos \theta(s)|^{-\gamma}\:ds=O(\viscosite^{-\gamma/4} |\ln \viscosite|^{-(2+3\gamma/2)}).
\ee
We infer that
$$
\viscosite \left(\int_{\d \Om}(\d^2_s(\varphi_i A))^2 \left| \frac{\viscosite}{\cos (\theta(s))}\right|^{1/3}\:ds\right)^{1/2}=o(\viscosite^{5/8}).
$$
In a similar way, we have
$$
\int_{\supp \varphi_i} (\theta')^2 |\cos \theta|^{-11/3}= O\left( \viscosite^{-5/12} |\ln \viscosite|^{-1/2}\right),
$$
so that
$$
\left(\int_{\supp \varphi_i}A(s)^2 ( \theta')^2 \viscosite^{5/3} |\cos \theta|^{-11/3} \:ds\right)^{1/2}=O(\viscosite^{\frac{5}{6}- \frac{5}{24}}|\ln \viscosite|^{-1/4})=o(\viscosite^{5/8}).
$$
We now tackle the term
\be\label{ri-embetant-1}\int_{\d\Om}A^2(s) (\d_s \varphi_i(s))^2 \left| \frac{\viscosite}{\cos (\theta(s))}\right|^{5/3}\:ds.
\ee
By definition,
the support of $\d_s\varphi_i$ is located in  neighbourhoods of $s_i$ and $s_{i+1}$, so that we can use assumption (H2).
We distinguish between  the cases (i) and (ii). By definition of $\varphi_i^+$,
$$\Supp \d_s\varphi_i^+ \subset \left(s_i^+, s_i ^++ \frac{|s_i-s_i^+|}{2}\right).$$
 In a neighbourhood of $s=s_i$, in case (i), we have
$$
|\d_s \varphi_i^+ (s)|\leq C \viscosite^{-\frac{1}{4n+3}}\mathbf 1_{C_0\viscosite^{\frac{1}{4n+3}}\leq s-s_i\leq 2C_0\viscosite^{\frac{1}{4n+3}}},
$$
so that \eqref{ri-embetant-1} is bounded by
$$
C\viscosite^{\frac{5}{3}- \frac{1}{2(4n+3)}+ \frac{1}{4n+3}\left(1- \frac{5n}{3}\right)}= O\left(\viscosite^\frac{10n+11}{8n+6}\right)=o(\viscosite^{5/4}).
$$
In case (ii), using Lemmas \ref{sigma-s-est} and \ref{sigma-s-est-bis}, we infer that
$$
|s_i^\pm- \sigma_i^\pm|\sim C \frac{\ln |\ln \viscosite|}{(\ln \viscosite)^2},
$$
while
$$
\int_{s_i^\pm}^{\sigma_i^\pm}|\cos \theta|^{-5/3}\leq C \viscosite^{-5/12} |\ln \viscosite|^{-9/2}.
$$
We deduce that the contribution of \eqref{ri-embetant-1} is bounded by
$$
|s_i^\pm- \sigma_i^\pm|^{-2} \viscosite^{5/3} \int_{s_i^\pm}^{\sigma_i^\pm}|\cos \theta|^{-5/3}= O(\viscosite^{5/4} |\ln \viscosite|^{-1/2}(\ln |\ln \viscosite|)^{-2})=o(\viscosite^{5/4}).
$$
There only remains to check that
$$
\int_{\supp \varphi_i}(A'(s))^2 \left| \frac{\viscosite}{\cos (\theta(s))}\right|^{5/3}\:ds= o(\viscosite^{5/4}).
$$
We use inequality \eqref{est:dsA}. We have
\begin{eqnarray*}
&&\int_{\supp \varphi_i}(A'(s))^2 |\cos \theta|^{-5/3}\:ds \\
&\leq & C\int_{\supp \varphi_i} |\cos \theta|^{-5/3}\:ds+C \int_{\supp \varphi_i}|\cos \theta|^{1 /3}(\ccM(y(s)))^2 \,.
\end{eqnarray*}
Using \eqref{est:sing-cos-1}, it can be checked that the first term is $o(\viscosite^{-5/12})$. There remains to estimate the second term of the right-hand side. We have to distinguish between two cases:
\begin{itemize}
 \item If $|\cos \theta|\ll 1$ in a neighbourhood of the point where $|y(s)-y_j|\ll 1$, then we can use assumption (H2) and the estimates of Appendix B. We infer that
$$
 \int_{\supp \varphi_i}|\cos \theta|^{1 /3}(\ccM(y(s)))^2  \leq C \delta_y^{-5/3}|\ln \delta_y|^{-10/3}.
$$

\item If $\cos \theta$ is bounded away from zero at the point where $y(s)=y_j$ for some $s\in \supp \varphi_i$, $j\in I_+$, then we can use the change of variables $s\to y$ in the integrals. Since $dy(s)/ds=-\cos \theta(s)$, the jacobian is not singular, and we have
$$
\begin{aligned}
 \int |\cos \theta|^{1/3} (\ccM(y(s)))^2
&\leq  C \int_{\delta_y/2\leq |y|\leq 1} |y|^{-2} (\ln |y|)^{-2}\\
&\leq  C \delta_y^{-1} |\ln \delta_y|^{-2}.
\end{aligned}
$$

\end{itemize}
In both cases, we infer that
$$\viscosite^{5/3}\int_{\supp \varphi_i}(A'(s))^2 |\cos \theta|^{-5/3}\:ds =o(\viscosite^{5/4}).
$$
This concludes the estimates in the neighbourhood of a point $s\in \d\Om$ such that $\cos \theta$ vanishes or $y(s)=y_j$ for some $j\in I_+$.

\smallskip

$\bullet$ {\bf Estimates near the end point of $\Sigma$:}

Using once again Lemma \ref{lem:est-A}, we infer that the contribution of this zone to $\|r_i^1\|_{L^2}$ is bounded by
$$
C(\viscosite^{\frac{1}{6}- \frac{1}{8\times 19}} + \viscosite^{1/2}|\ln \viscosite|^{1/2}\viscosite^{-\frac{1}{4}-\frac{3}{8\times19}}+ \viscosite^{\frac{5}{6}- \frac{1}{2}- \frac{5}{8\times19}})=o(\viscosite^{1/8}),
$$
while its contribution to $\|r_i^2\|_{H^{-2}}$ is bounded by
$$
C(\viscosite^{\frac{5}{6}- \frac{1}{8\times19}}+ \viscosite^{\frac{5}{6}- \frac{1}{4}- \frac{3}{8\times19}}\delta_y^{1/2}+\viscosite^{\frac{7}{6}- \frac{1}{2}- \frac{5}{8\times19}}\delta_y^{1/2})=o(\viscosite^{5/8}).
$$
This concludes the proof.
\end{proof}

\subsection{Error terms associated with discontinuity    layers}$ $

Since we have already proved that the error terms associated with $\psi^{corr}_\Sigma $ are admissible in the sense of Definition \ref{def:admissible}, in order to establish Proposition \ref{lem:sigma}, there only remains to prove that
$$
(\d_x-\viscosite \Delta^2) \left(\psi^\Sigma\left(x, \frac{y-y_1}{\viscosite^{1/4}}\right) \chi \left(\frac{y-y_1}{\delta_y}\right)\right)= \dtl\left(x, \frac{y-y_1}{\viscosite^{1/4}}\right) + r^1_\Sigma+ r_\Sigma^2,
$$
where
$$
\|r^1_\Sigma \|_{L^2(\Om)}=o(\viscosite^{1/8}),\quad \| r^2_\Sigma\|_{H^{-2}(\Om)}=o(\viscosite^{5/8}).
$$
Note that the computations here are much simpler since all the formulas are given in cartesian coordinates.

By definition of $\psi^\Sigma$, we have
\begin{eqnarray}
\nonumber& &(\d_x-\viscosite \Delta^2) \left(\psi^\Sigma\left(x, \frac{y-y_1}{\viscosite^{1/4}}\right) \chi \left(\frac{y-y_1}{\delta_y}\right)\right)-  \dtl\left(x, \frac{y-y_1}{\viscosite^{1/4}}\right)\\
&=& (\chi-1) \left(\frac{y-y_1}{\delta_y}\right) \dtl\left(x, \frac{y-y_1}{\viscosite^{1/4}}\right)\label{psid1}\\
&&- \viscosite ( 2\d_y^2 +\d_x^2) \d_x^2 \left(\psi^\Sigma\left(x, \frac{y-y_1}{\viscosite^{1/4}}\right) \chi \left(\frac{y-y_1}{\delta_y}\right)\right)\label{psid2}\\
&&-\viscosite \sum_{j=0}^3C_4^j \viscosite^{-j/4}\delta_y^{j-4}\d_Y^j\psi^\Sigma\left(x,\frac{y-y_1}{\viscosite^{1/4}}\right)\chi^{(4-j)}\left(\frac{y-y_1}{\delta_y}\right).\label{psid3}
\end{eqnarray}
The term \eqref{psid1} is supported in
$$
\{|y-y_1|\leq \viscosite^{1/4}\}\cap \{|y-y_1|\geq \delta_y\}.
$$
Since $\delta_y\gg \viscosite^{1/4}$, \eqref{psid1} is zero for $\viscosite$ small enough. Moreover, using Lemma \ref{lem:est-psi-sigma}, we infer that \eqref{psid2} is $O(\viscosite \viscosite^{-2/19} \viscosite^{1/8})$ in $H^{-2}(\Om)$. Eventually, using the same type of commutations as in the proof of Lemma \ref{lem:approx-BL-NS}, we prove that \eqref{psid3} is $o(\viscosite^{1/8})_{L^2} + o(\viscosite^{5/8})_{H^{-2}}$.

\section{Remainders in the periodic and rectangular cases}

$\bullet$ We begin with the rectangular case, which is closer to the case of a smooth domain in $\R^2$. The remainder term coming from the interior is $-\viscosite \Delta^2 \psi^0$, which is $O(\viscosite)$ in $L^2$ if $\tau $ is smooth (say, $\tau \in H^4$), and is therefore admissible.
The north and south boundary layers give rise to the same type of error terms as in the case of a smooth domains; hence these error terms are admissible as well. Eventually, the error terms coming from the western boundary layer are
$$
(\viscosite \d_y^4 + 2\viscosite \d_y^2 \d_x^2)\left(\sum_\pm A_\pm(y) \exp\left(-e^{\pm i \pi/3} \frac{x-x_-}{\viscosite^{1/3}}\right) \chi_0(x)\right)  ,
$$
which are bounded in $H^{-2}$ by
$$
\viscosite \left\| \d_y^2  \sum_\pm A_\pm(y) \exp\left(-e^{\pm i \pi/3} \frac{x-x_-}{\viscosite^{1/3}}\right) \chi_0(x) \right\|_{L^2}.
$$
Using estimates \eqref{est:A-rect}, we infer that the above quantity is bounded by
\begin{eqnarray*}
&&C \viscosite^{7/6} \sum_\pm\|\d_y^2 A_\pm \|_{L^2(y_-,y_+)}\\
&\leq & C \viscosite^{7/6} \viscosite^{-1/2}=o(\viscosite^{5/8}),
\end{eqnarray*}
hence the remainder term is admissible.

$\bullet$ In the periodic case, there are very few remainder terms. Concerning the interior part, we have
$$
(\d_x - \viscosite \Delta^2)(\psi^0+ \psi_{circ})= \tau - \viscosite \Delta^2 \psi^0,
$$
and
$$
\|\viscosite \Delta^2 \psi^0\|_{H^{-2}(\T\times (y_-,y_+))}\leq C \viscosite \|\tau\|_{H^2(\T\times (y_-,y_+))}.
$$
As for the boundary layer terms, we have
\begin{eqnarray*}
&&(\d_x - \viscosite \Delta^2)(\psbl_N(x,(y_+-y)/\viscosite^{1/4})\chi_0(y))\\
&=&-\viscosite (\d_x^4 + 2\d_x^2\d_y^2)(\psbl_N(x,(y_+-y)/\viscosite^{1/4})\chi_0(y))\\
&+&\text{exponentially small terms in }L^2
\end{eqnarray*}
and using \eqref{est:psbl-per}, the first term in the right-hand side is bounded in $H^{-2}$ by
\begin{eqnarray*}
&&C\viscosite \left\|\d_x^2 \psbl_N(x,(y_+-y)/\viscosite^{1/4})\chi_0(y)\right\|_{L^2(\Om)}\\
&\leq & C \viscosite^{7/8}( \|\psi^0_{y=y_\pm}\|_{H^{2 - \frac{1}{8}}(\T)} + \viscosite^{1/4} \|\d_y\psi^0_{y=y_\pm}\|_{H^{2 - \frac{3}{8}}(\T)})\\
&\leq &  C\viscosite^{7/8}( \|\tau \|_{H^{2+ \frac{3}{8}} }+\viscosite^{1/4} \|\tau \|_{H^{2 + \frac{9}{8}}} ).
\end{eqnarray*}
Hence for $\tau\in H^{25/8}$, the remainder terms are all $O(\viscosite)$ in $H^{-2}$, which completes the proof of Corollary \ref{cor:per-rect}.

%% file: conclusion.tex
\chapter{Discussion: Physical relevance of the model}

Of course it is not completely clear whether our simplified model
is really relevant from the physical point of view. Let us therefore explain briefly the derivation of this model (following Desjardins and Grenier \cite{DG}), and track the different simplifications that should be compared with experimental data. For further discussions regarding the physical relevance of such models, we refer to \cite{pedlosky},\cite{pedlosky2} or \cite{gill}.

$\bullet$ As usual in large-scale oceanography, we start with the 3D incompressible Navier-Stokes equations in a rotating frame: the velocity field $u$ is assumed to be divergence-free
$$\nabla \cdot u =0\,,$$
 and to satisfy the dynamical equation
$$
\d_t u + (u\cdot \nabla) u+2\omega \wedge u ={\nabla p\over \rho} +g +{1\over \rho} \cF u$$
expressing the fact that the fluid evolves under the combined effects of the Coriolis force, the pressure, the gravity and some turbulent dissipation mechanism. The precise formulation of this last contribution (involving a turbulent viscosity)
$${1\over \rho}  \cF u =\mu_h \Delta_h u +\mu_3 \d_{33} u $$
 even commonly used by physicists, has no real justification, which is a first limitation of our study.

These equations are set in a bounded domain
$$ \cD = \{ x\in \Omega \times \R\,/\, h_B(x_h) \leq x_3\leq 0\}$$
and supplemented by boundary conditions. On the bottom (which is described by the topography $h_B$) and on the lateral boundaries, we assume that the fluid-solid  interaction can be catched through a no slip condition (Dirichlet condition)~:
$$ u_{|x_3=h_B}=0,\quad u_{|\d \Omega_h} =0\,.$$
For the sake of simplicity, the free surface is replaced by a prescribed spherical boundary corresponding to the depth $x_3=0$ (rigid lid approximation), and the effect of the wind is modeled by some non homogeneous Navier condition
$$ u_{n|x_3=0} =0,\quad \mu_h (\nabla u+(\nabla u)^T)_{t|x_3=0} =\cT.$$
This drastic but standard simplification is investigated for instance in \cite{LOT,Du}, but has no rigorous justification, which is the second weakness of the model.

\bigskip
$\bullet$ Far from the poles and the equator, i.e. around a latitude $\vartheta_0\in (0,\pi/2)$, these equations can be rewritten in scaled cartesian-like coordinates. This involves many non dimensional parameters characterizing the physical properties of the  flow, especially
\begin{itemize}[label=$-$]
\item the Rossby number which measures the  size of the Coriolis force
$$\eps = {U\over 2\omega \sin \vartheta_0 L}$$
where $U$ and $L$ are  the typical  velocity  of the fluid and (horizontal) length of observation;
\item the aspect ratio and the curvature parameter which characterize the geometry of the domain, defined respectively by
$$\rho = {\bar h \over L},\quad r_*=\frac{R_0}{L}$$
 $\bar h$ being the typical height of the ocean and  $R_0$  the radius of the Earth;
 {
 \item
the horizontal and vertical Ekman numbers which account for the viscous effects
$$\viscosite_3 ={\mu_3 \over 2 \omega \sin \vartheta_0 \bar h ^2}  ,\quad \viscosite_h ={\mu_h \over 2 \omega \sin \vartheta_0 L^2}. $$
 }
\end{itemize}
The Munk equation is obtained in the fast rotating limit with thin layer approximation, that is  when $\eps \to 0$, with the following choice of scaling
$$
\begin{aligned}
 \rho, r_*^{-1}=O( \eps) , \quad \viscosite_3=O(\eps)\\
h_B = -\bar h (1+\eps \eta_B),\\
  \hbox{ and } \sin \vartheta= \sin \vartheta_0(1+ \eps \beta y)
\,.
  \end{aligned}$$
where $$\beta=\frac{r_* \cos \vartheta_0}{\eps\sin \vartheta_0}.$$

\bigskip
$\bullet$ A formal asymptotic analysis (based for instance on asymptotic expansions)   shows that the limiting flow is purely two-dimensional: we indeed obtain at leading order
\begin{equation}
\label{2D}
 u_3=0 ,\quad \d_3 u_h =0, \quad \nabla_h \cdot u_h =0\,,
 \end{equation}
 {
  or equivalently
 $$ u =\nabla_h^\perp \psi$$
 for some scalar stream function $\psi$ depending only on the horizontal variables.}
Since these conditions are not compatible with the bottom and surface boundary conditions, one has to introduce boundary layer corrections, referred to as Ekman layers, which contribute to the global energy balance via Ekman pumping.

The dynamical equation is then obtained at next order (more or less as a solvability condition): since $u_h$ is divergence-free, it is indeed completely determined by its vorticity:
\begin{equation}
\label{vorticity}
\nabla^\perp\cdot  ( \d_t u +u\cdot \nabla_h u -\beta y u^\perp) =\nabla^\perp\cdot ( \viscosite_h \Delta_h  u -\beta\cT+\eta_B u^\perp -r u)\end{equation}
where $-r u$ is the Ekman pumping associated to the energy dissipation by friction on the bottom, $\eta_B u^\perp$ accounts  for   the effects of the topography, and  $-\beta \cT $ is the source term resulting from the wind forcing. When the domain $\Om$ is not simply connected, \eqref{vorticity} must be supplemented with further compatibility conditions (see paragraph \ref{ssec:islands}).

Note that this formal derivation can be justified by classical energy methods \cite{DG}, at least for well-prepared initial data (i.e. for initial data satisfying the constraint equations (\ref{2D})). Starting from the limiting system (\ref{2D})(\ref{vorticity}), one can indeed build a smooth approximate solution, then control its $L^2$ distance to the solution $u_\eps$ of the Navier-Stokes equations using some strong-weak stability principle. Note that one even obtains a rate of convergence.

\bigskip
$\bullet$ The apparition of Munk boundary layers and the intensification of oceanic currents on western coasts we would like to describe with such a model  are physical phenomena which are typically linear and which result from the fact that the no-slip Dirichlet  condition on the coasts of the basin becomes a non admissible boundary condition if the viscous dissipation is not a leading order term in the equation (\ref{vorticity}).

More precisely, we expect to exhibit such a behaviour when the remaining rotating term $\beta \nabla^\perp \cdot (y u^\perp)= \beta u_y $ (due to the inhomogeneities of the local rotation vector) is large compared to $\viscosite_h \Delta \nabla^\perp \cdot u$.
In particular, we expect that boundary layers
\begin{itemize}
\item should not depend on the topographical effects and Ekman pumping (if $\beta$ is large enough)
\item decouple from  the convection (at least if boundary layers are  stable in the sense that their kinetic energy remain small)
\item are quasi-stationary, meaning that the equation governing boundary layers does not involve time.
\end{itemize}
We emphasize that our study remains entirely valid if  in equation \eqref{M}, $\beta$ is not a parameter, but a smooth function which remains bounded away from zero.

%% file: app-islands.tex
\section*[Appendix A]{Appendix A: The case of islands: derivation of the compatibility condition \eqref{compatibility} and proof of Lemma \ref{lem:def-M,D}}
 \renewcommand{\theequation}{A.\arabic{equation}}
  \setcounter{equation}{0}  

$\bullet$ The compatibility condition \eqref{compatibility} is inherited from the Navier-Stokes system satisfied by $u=\nabla^\bot \psi$ (see \cite{kikuchi} for a similar argument in the inviscid case). Indeed, we start from the stationary Stokes-Coriolis system in dimension 2, with $\beta$-plane approximation, namely
\be\label{SC}
\begin{aligned}
\frac{1}{\eps}(1+ \beta \eps  y) u^\bot - \viscosite \Delta u+ \nabla p=\mathcal T,\text{ in } \Om,\\
\Div u=0\text{ in } \Om,\\
u_{|\d\Om}=0.
\end{aligned}
\ee
The idea is to take the curl of the first equation in order to get rid of the pressure term. However, for $\Phi\in L^2(\Om)^2$, the identity $\curl \Phi=0$ does not imply the existence of $q\in H^1(\Om)$ satisfying $\Phi=\nabla q$.
Indeed, $\Phi$ is a gradient if and only if its circulation around any closed contour $\mathcal C$ in $\Om$ vanishes:
$$
\oint_\mathcal C \Phi\cdot t=0,
$$
where $t$ is the tangent vector to the curve $\mathcal C$. If $\curl \Phi=0$, this condition becomes $\int_{C_i}  \Phi\cdot t=0$ for all $i\geq 2$.
Hence we obtain
$$
\exists q\in H^1(\Om),\quad \Phi=\nabla q \ \iff \ \left\{\begin{array}{l}
 \curl \Phi=0\\
\text{and }\int_{C_i}  \Phi\cdot t=0\text{ for all }i\geq 2.
\end{array}\right.
$$
Therefore \eqref{SC} is equivalent to
$$
\begin{aligned}
\beta u_2 - \viscosite \Delta \curl u = \curl \cT\text{ in }\Om,\\
\mathrm{div} u=0,\\
u_{|\d\Om}=0,\\
\viscosite \int_{C_i} \Delta u \cdot t + \int_{C_i} \cT \cdot t=0\quad \forall i \geq 2.
\end{aligned}
$$
This amounts to \eqref{Munk-islands}-\eqref{compatibility} with $u=\nabla^\bot \psi$ (notice that, as above, the existence of $\psi$ is ensured by the divergence free condition on $u$ and by the Dirichlet boundary conditions).

$\bullet$ Proof of Lemma \ref{lem:def-M,D}: the existence and uniqueness of the functions $\psi_i\in H^2(\Om)$ follow for instance from the Lax-Milgram Lemma. Identity \eqref{dec:psi} is a consequence of the linearity of the equation \eqref{Munk-islands} and of uniqueness, and the equality $M^\viscosite c=D^\viscosite$ follows from the compatibility condition \eqref{compatibility} and the decomposition \eqref{dec:psi}.

Concerning the invertibility of $M^\viscosite$, we prove in fact that for all $(a_2, \cdots, a_K)\in \R^{K-1}$
\be\label{in:psia}
\sum_{i,j} a_i a_j\viscosite \int_{C_j}\d_n \Delta \psi_i= -\viscosite\int_{\Om} |\Delta \psi_a|^2\leq 0,
\ee
where $\psi_a=\sum_{i=2}^K a_i \psi_i$.

We have
$$
\begin{aligned}
\d_x \psi_a-\viscosite \Delta^2 \psi_a=0,\\
\d_n \psi_{a|\d\Om}=0,\\
\psi_{a|C_1}=0,\ \psi_{a|C_i}=a_i, \ i\geq 2.
\end{aligned}
$$
Therefore, since $C_i$ is a closed contour,
$$
\int_{\Om} \d_x \psi_a\psi_a= \frac{1}{2}\sum_{i=2}^K a_i^2 \int_{C_i} e_x\cdot n=0,
$$
and
\begin{eqnarray*}
\viscosite\int_{ \Om} \Delta^2 \psi_a \psi_a&=&\viscosite \int_{\Om} |\Delta \psi_a|^2+ \sum_{i=2}^K a_i \viscosite \int_{C_i} \d_n \Delta \psi_a\\
&=&\viscosite \int_{\Om} |\Delta \psi_a|^2+ \viscosite \sum_{i,j}a_ia_j  \int_{C_i} \d_n \Delta \psi_j.
\end{eqnarray*}
The identity \eqref{in:psia} follows.

%% file: appendixB.tex
\section*[Appendix B]{Appendix B: Equivalents for the coordinates of boundary points near horizontal parts }
  \renewcommand{\theequation}{B.\arabic{equation}}
  \setcounter{equation}{0}  

In the vicinity of a point $s_0\in  \d \Gamma_E$ such that $\cos (\theta(s_0))=0$, the East part of the boundary can be described as the graph of some function $y\mapsto x_E(y)$.
Assumption (H2) provides then asymptotic expansions of the function $x_E$ and its derivatives.

 Let $(x(s), y(s))$ be the coordinates of the point with arc-length $s$ on $\d\Om$. By definition, if $\cos(\theta(s))>0$,
$$
x(s)=x_E(y(s)).
$$
Moreover,
\be\label{ds-xy}
\frac{d}{ds}\begin{pmatrix}
                x(s)\\y(s)
            \end{pmatrix}
= \begin{pmatrix}
        \sin (\theta(s))\\
        -\cos(\theta(s))
  \end{pmatrix}
\ee
and setting the origin of the axes so that $(x(s_0), y(s_0))= (0,0)$, we get
\be\label{x_E'}\begin{aligned}
x_E'(y(s))=-\tan(\theta(s)),\\
x_E''(y(s))= \frac{\theta'(s)}{\cos^3(\theta(s))}.
\end{aligned}
\ee
Similar formulas hold for $x_E^{(3)}$ and $x_E^{(4)}$.

\textbf{In case (H2i),} we have
$$
y(s)\sim C(s-s_0)^{n+1}\sim C x(s)^{n+1} ,
$$
and thus
\be
\label{bound-est-i}
\begin{aligned}
x_E(y)\sim C|y|^\frac{1}{n+1},\\
        x_E^{(k)}(y)\sim C_k |y|^{\frac{1}{n+1}-k}\quad \text{for } k\in\{1,2,3,4\}.
\end{aligned}
\ee

\textbf{In case (H2ii), }using the same kind of calculations as above, we infer that as $y$ vanishes
\be
\label{bound-est-ii}
\begin{aligned}
        x_E(y)\sim \frac{C}{\ln |y|},\\
        x_E^{(k)}(y)\sim \frac{C_k}{y^k(\ln |y|)^2} \quad \text{for } k\in\{1,2,3,4\}.
\end{aligned}
\ee

We have used the following fact:
$$
\int_{s_0}^s  \exp\left(-\frac{\alpha}{s'-s_0}\right)ds'\sim \frac{(s-s_0)^2 }{\alpha} \exp\left(-\frac{\alpha}{s-s_0}\right)\quad \text{as }s\to s_0.
$$

%% file: appendixC.tex
\section*[Appendix C]{Appendix C: Estimates on the coefficients $a$ and $b$.}
  \renewcommand{\theequation}{C.\arabic{equation}}
  \setcounter{equation}{0}  

We recall that (see \eqref{def:psil}, \eqref{def:ab1}, \eqref{def:ab2})
$$
\psil(x,y)=\chi^+ \chi\left(\frac{y-y_1}{\viscosite^{1/4}}\right)\left[a(x) + b(x) (y-y_1)\right],
$$
where $a$ and $b$ are defined by
\begin{eqnarray*}
b(x)&=&(1-\varphi(s))\left[-\sin \theta(s) \d_n \psi^0_{t|\d\Om}(s) +\cos \theta(s)\d_s\psi^0_{t|\d\Om}(s)\right]\\&& - \cos \theta(s) \varphi' \psi^0_{t|\d\Om}(s),\\
 a(x)&=&-(1-\varphi(s) )\psi^0_{t|\d\Om}(s) - b(x(s))(y(s)-y_1),
\end{eqnarray*}
for $x_1<x=x(s)<x(\sigma_{in})$,
and by
$$
\begin{array}{l}
a(x)= -[\psi^0_t]_\Sigma,\\
b(x)=- [\d_y\psi^0_t]_\Sigma
\end{array}
\quad\text{for }x\leq x_1.
$$
For $x\geq x(\sigma_{in})$, we merely take $a(x)=b(x)=0$.

\begin{Lemma}
Let
$$\begin{aligned}
   x_{min}:=\inf \{x\in \bR,\ \exists y\in \bR,\ (x,y)\in \Om\},\\
   x_{max}:=\sup \{x\in \bR,\ \exists y\in \bR,\ (x,y)\in \Om\}
  \end{aligned}
$$
Then for all $k\in \{1,\cdots, 4\},$
$$
A_k:=\|a^{(k)}\|_{L^1(x_{min}, x_{max})} + \viscosite^{1/4} \|b^{(k)}\|_{L^1(x_{min}, x_{max})}=O(\viscosite^\frac{1-k}{19}).
$$
Additionally,
$$
\begin{aligned}
\|a\|_{L^\infty(x_{min}, x_{max})} + \viscosite^{1/4} \|b\|_{L^\infty(x_{min}, x_{max})} =O(1),\\
\|a''\|_{L^2(x_{min}, x_{max})} + \viscosite^{1/4} \|b''\|_{L^2(x_{min}, x_{max})} = o(\viscosite^{-1/2}).\\
\end{aligned}
$$

\label{lem:est-A-B}
\end{Lemma}

\begin{proof}
We recall that $a$ and $b$ were defined in paragraph \ref{ssec:psil}.
Since $a$ and $b$ are constant for $x<x_1$ and for $x>x(\sigma_{in})$, it suffices to prove the estimates on $(x_1, x(\sigma_{in}))$. On this interval,
 $a$ and $b$ are defined through formulas of the type
$$
a(x(s))=\ga (s),\quad b(x(s))=\gb(s).
$$
Differentiating these inequalities with respect to $s$, we obtain
$$
x'(s)a'(x(s))=\ga'(s),
$$
and thus, using the formula \eqref{ds-xy},
$$
a'(x(s))=\frac{\ga'(s)}{\sin \theta(s)}.
$$
Iterating this process and using the fact that $\sin \theta$ remains bounded away from zero in the interval under consideration, we infer that for all $k\in\{0,\cdots 4\},$
$$
\left|a^{(k)}(x(s))\right|\leq C\sum_{l=0}^k|\ga^{(l)}(s)|,
$$
and eventually, for all $p\in [1,\infty]$,
$$
\|a^{(k)}\|_{L^p(x_1 ,x(\sigma_{in})}\leq C\sum_{l=0}^k\| \ga^{(l)}\|_{L^p(\sigma_{in}, s_1)}.
$$
Of course the same type of inequality holds for $b$ as well. Therefore we now compute the derivatives of $\ga$ and $\gb$ with respect to $s$ up to order 4. Since
$$
\ga(s)=-(1-\varphi(s))\psi^0_{t|\d\Om}(s) - \gb(s)(y(s)-y_1).
$$
we start with the derivatives of $\gb$. Recall that
\begin{eqnarray*}
\gb(s)&=&(1-\varphi(s))\left[-\sin \theta(s) \d_n \psi^0_{t|\d\Om}(s) +\cos \theta(s)\d_s\psi^0_{t|\d\Om}(s)\right]\\&& - \cos \theta(s) \varphi' (s)\psi^0_{t|\d\Om}(s).
\end{eqnarray*}
 The computations are lengthy but do not raise any difficulty. It can be easily checked that the most singular estimates correspond to the case when $\cos \theta$ vanishes algebraically near $s=\inf I_1$ (assumption (H2i)), with the lowest possible exponent $n$.
Indeed, if $\cos\theta$ vanishes exponentially near $s=\inf I_1$, then for all $k\geq 0$
$$
\|\d_s^k\varphi\|_{L^\infty}=O(|\ln \viscosite|^{2k}),
$$
while in case (i)
$$
\|\d_s^k\varphi\|_{L^\infty}=O(\viscosite^{-\frac{k}{4n+3}}).
$$
Similar formulas hold for $\d_s^k \psi^0_{t|\d\Om}$, $\d_s^k \d_n\psi^0_{t|\d\Om}$.

Therefore in the rest of the proof, we only treat the case (H2i) with $n=4$ (Remember that because of assumption (H4), $\cos \theta$ vanishes at least like $(s-\inf I_1)^4)$ near $\inf I_1$). We explain with some detail the estimates of $\gb$ and $\gb'$, and we leave the rest of the derivatives to the reader.

First, we recall that by definition of $\sigma_{in}$ (see \eqref{def:sigma-}),
$$
|y(s)-y_1|\leq \viscosite^{1/4}\ll \delta_y
$$
for all $s\in (\sigma_{in}, s_1)$. Therefore, using the estimates of Lemma \ref{lem:est-psi0}, we have
$$
|\gb(s)|\leq C \left(\frac{1}{\delta_y|\ln \delta_y| } + \viscosite^{-1/19}\cos \theta(s)\right),
$$
so that
$$
\viscosite^{1/4}\|\gb\|_{\infty}=O(|\ln \viscosite|^{-1})=o(1).
$$
Differentiating  $\gb$ with respect to $s$, we obtain, for $s\in (\sigma_{in}, s_1)$,
\begin{eqnarray*}
\gb'(s)&=& (1-\varphi)\left[-\theta' \cos \theta \d_n \psi^0_{t|\d\Om} - \sin \theta \d_s  \d_n \psi^0_{t|\d\Om}\right.\\&&\qquad\qquad\left.- \theta' \sin \theta  \d_s \psi^0_{t|\d\Om}+ \cos \theta  \d_s^2 \psi^0_{t|\d\Om}\right]\\
&&+ \varphi'\left[\sin \theta  \d_n \psi^0_{t|\d\Om}- 2 \cos \theta  \d_s \psi^0_{t|\d\Om}+ \theta' \sin \theta   \psi^0_{t|\d\Om}\right]\\
&&-\varphi''\cos \theta \psi^0_{t|\d\Om},
\end{eqnarray*}
so that
\begin{eqnarray*}
|\gb'(s)|&\leq & C(1-\varphi)\left[1 +\frac{|\cos \theta|}{\delta_y^2 |\ln \delta_y|}\right]\\
&&+ C|\varphi'|\frac{1}{\delta_y |\ln \delta_y|}+ C |\varphi''|\; |\cos \theta|.
\end{eqnarray*}
Eventually, we retrieve
$$
\|\gb'\|_{L^1}=O(\viscosite^{-1/4}|\ln \viscosite|^{-1}).
$$
The same kind of estimate yields, for $k=2,3,4$,
$$
\|\gb^{(k)}\|_{L^1}=O(\viscosite^{-\frac{1}{4}+ \frac{1-k}{19}}|\ln \viscosite|^{-1}).
$$
Now, since
$$
\ga^{(k)}= \sum_{l=0}^kC_k^l \left(\varphi^{(l)}\d_s^{k-l}\psi^0_{t|\d\Om} - \gb^{(l)}\frac{d^{k-l}(y(s)-y_1)}{d s^{k-l}}\right),
$$
we infer, for $k\geq 1$,
\begin{eqnarray*}
\|\ga^{(k)}\|_{L^1}&\leq &C \| \varphi^{(k)}\|_{L^1}\\
&&+C\sum_{l=0^{k-1}} \viscosite^{-l/19}\| \d_s^{k-l}\psi^0_{t|\d\Om}\|_{L^1(\sigma_{in}, s_1)} \\
&&+ C\sum_{l=0^{k}} \viscosite^{-\frac{1}{4}+ \frac{1-l}{19}}|\ln \viscosite|^{-1}\left\| \frac{d^{k-l}(y(s)-y_1)}{d s^{k-l}}\right\|_{L^\infty(\sigma_{in}, s_1)} .
\end{eqnarray*}

It can be checked that for $l\in \{1,\cdots, 4\}$,
$$
\begin{aligned}
 \| \d_s^{l}\psi^0_{t|\d\Om}\|_{L^1(s_-^N, s_+)}\leq C |\ln \viscosite|^{-1} \viscosite^{\frac{1-l}{24}},\\
\left\| \frac{d^{l}(y(s)-y_1)}{d s^{l}}\right\|_{L^\infty(s_-^N, s_+)}\leq C \viscosite^{\frac{6-l}{24}}.
\end{aligned}
$$
Gathering all the terms, we obtain, for $k\geq 1$,
$$
\|\ga^{(k)}\|_{L^1}=O(\viscosite^\frac{1-k}{19}).
$$

The estimates of $\ga''$ and $\gb''$ in $L^2$ go along the same lines.

\end{proof}

%% file: appendixD.tex
\section*[Appendix D]{Appendix D: Proof of Lemma \ref{lem:sobolev}}
   \renewcommand{\theequation}{D.\arabic{equation}}
   \setcounter{equation}{0}  

The estimates on $f$ and $f'$ are derived in a classical fashion: for instance, write
$$
f(Z)^2=-2\int_Z^\infty f(Z')\d_Z f(Z')\:dZ',
$$
so that
$$
\|f\|_\infty^2 \leq 2 \|f\|_2\|\d_Zf\|_2.
$$
Integrating by parts and using the fact that $f_{|Z=0}=\d_Zf_{|Z=0}=0$,
we infer that
$$
\|\d_Zf\|_2^2= -\int_{\R_+} f \d_Z^2 f
$$
and similarly
\be\label{est:dZ2f}
\|\d_Z^2f\|_2^2=\int_{\R_+}f \d_Z^4 f.
\ee
Using several times the Cauchy-Schwarz inequality leads to
$$
\|f\|_\infty \leq \sqrt{2}\|f\|_2^{7/8}\|\d_Z^4 f \|_2^{1/8}.
$$
The estimate on $\d_Z f$ goes along the same lines. The estimate on $\d_Z^2 f$ is a little more tricky: we write
\begin{eqnarray}
\nonumber\|\d_Z^2 f\|_\infty^2&\leq & 2\|\d_Z^2 f\|_2\|\d_Z^3 f\|_2\\
&\leq &2\|f \|_2^{1/2}  \|\d_Z^3 f\|_2\|\d_Z^4 f \|_2^{1/2}.\label{in:dz2}
\end{eqnarray}
On the other hand, an integration by parts yields
$$
\|\d_Z^3 f\|_2^2=-\int_{\R_+}\d_Z^2f \d_Z^4 f - \d_Z^2 f_{|Z=0}\d_Z^3 f_{|Z=0}
$$
and
$$\d_Z^3 f_{|Z=0}^2=-2\int_{\R_+}\d_Z^3f \d_Z^4 f.$$
Therefore, using \eqref{est:dZ2f}
$$
\|\d_Z^3 f\|_2^2\leq  \|f\|_2^{1/2}\|\d_Z^4 f \|_2^{3/2}+ \sqrt{2} \|\d_Z^2 f\|_\infty \|\d_Z^4 f \|_2^{1/2}\|\d_Z^3 f\|_2^{1/2}.
$$
We deduce that
$$
\|\d_Z^3 f\|_2^2\leq C \|f\|_2^{1/2}\|\d_Z^4 f \|_2^{3/2}+ C\|\d_Z^2 f\|_\infty^{4/3} \|\d_Z^4 f \|_2^{2/3}.
$$
Inserting the above inequality into \eqref{in:dz2} and using Young's inequality  leads eventually to
$$
\|\d_Z^2 f\|_\infty^2 \leq C\|f\|_2^{3/4}\|\d_Z^4 f \|_2^{5/4},
$$
so that
$$
\|\d_Z^3 f\|_2\leq C\|f\|_2^{1/4}\|\d_Z^4 f \|_2^{3/4}.
$$
The last estimate on $\|\d_Z^3f\|_\infty$ follows easily.

\qed

%% file: notations.tex
\section*{Notations}

\bigskip
\noindent
$A_i$, $B_i$ subdomains of $\Omega$ where the solution to the Sverdrup equation is continuous,  p. \pageref{A-def}, \pageref{B-def}

\bigskip
\noindent
$\Gamma_E$, $\Gamma_W$ East and West parts of the boundary, p. \pageref{GammaEW}

\bigskip
\noindent
$\Gamma_N$, $\Gamma_S$ North and South parts of the boundary, p. \pageref{GammaNS}

\bigskip
\noindent $\gamma_{i,i-1}$ truncation of the boundary conditions to be lifted by the North/South boundary layers,  p. \pageref{gamma-i}

\bigskip
\noindent $\tilde \gamma_{i,i-1}\geq \gamma_{i,i-1}$ truncation of  the North/South boundary layers, p. \pageref{def:tilde-gamma-i}

\bigskip
\noindent
$\delta_y = \viscosite^{1/4} |\ln \viscosite|^{1/5} (\ln |\ln \viscosite|)^{-\beta}$ truncation parameter with respect to $y$, p.  \pageref{deltay}

\bigskip
\noindent
$\delta_x^{i,\pm} = x_E( y_i\pm \delta_y)- x_E(y_i)$ local truncation parameter with respect to $x$, p.  \pageref{deltax}

\bigskip
\noindent
$\theta(s)$ angle between the horizontal vector $e_x$ and the exterior normal to the boundary, p.  \pageref{theta}

\bigskip
\noindent
$I_+$ set of indices corresponding to East corners ($s_i \in \d \Gamma_E$), p.  \pageref{I+}

\bigskip
\noindent%
$\lambda_E$, $\lambda_W$ inverse sizes of the East and West boundary layers, p.  \pageref{def:lambda_EW}

\bigskip
\noindent%
$\lambda_N=\lambda_S=\lambda_\Sigma=\viscosite^{-1/4}$ inverse size of the horizontal  boundary layers, p.  \pageref{def:lambdaNS}

\bigskip
\noindent%
$\ccM$ majorizing function, p. \pageref{cM-def}

\bigskip
\noindent
$\rho_i$ partition of unity, p.  \pageref{rho-def}

\bigskip
\noindent
$s$ curvilinear abcissa, p.  \pageref{localcoord}

\bigskip
\noindent
$s_i$ curvilinear abcissa of the singular points ($s_i \in \d \Gamma_E\cup \d \Gamma_W $), p.  \pageref{si}

\bigskip
\noindent%
$s_i^\pm $  extremal points of the domain of validity of East and West boundary layers (if relevant), p.  \pageref{sipm}

\bigskip
\noindent%
$t_1^+ $ curvilinear abcissa of the East end of $\Sigma$, p.  \pageref{t1+}

\bigskip
\noindent
$\sigma_i^\pm $  extremal points where energy is injected in  North and South boundary layers (if relevant), p.  \pageref{Esigmaipm}, \pageref{Wsigmaipm}

\bigskip
\noindent
$\sigma_{in}$ initial point for the discontinuity boundary layer $\psil, \psisig$, p.  \pageref{def:sigma-}

\bigskip
\noindent
$\Sigma$ surfaces of discontinuity (ordinates of critical points), p.  \pageref{Sigma}

\bigskip
\noindent
$\tau$ forcing term (coming from Ekman pumping), p.  \pageref{tau}

%

\bigskip
\noindent
$\varphi_i^\pm$ truncation dealing with the transition between East/West and North/South boundary layer terms, p.  \pageref{Evarphii}, \pageref{Wvarphii}

\bigskip
\noindent
$\chi_0$ localization  near the boundary, p.  \pageref{chib}

\bigskip
\noindent
$\chi_\pm$ localization on $\Omega^\pm$, p. \pageref{chipm}

\bigskip
\noindent
$\chi_\viscosite$ truncation of $\tau$ near the singular point $s_i$, p.  \pageref{chii}

\bigskip
\noindent
$x_E^\pm\equiv x_E^\pm(y)$ graphs of the East boundaries of $ \Omega^\pm$, p. \pageref{x-E-pm}

\bigskip
\noindent%
$x_1^W $  abcissa of the West end of $\Sigma$, p.  \pageref{projW-def}

\bigskip
\noindent
$(x_i,y_i)$ coordinates of the singular points ($(x_i,y_i) = (x(s_i), y(s_i))  $), p.  \pageref{si}

\bigskip
\noindent
$\psi^0$ solution to the transport equation, p.  \pageref{psi0}

\bigskip
\noindent
$\psi^0_t$  solution to the transport equation  with truncated source term, p.  \pageref{psi0t}

\bigskip
\noindent
$\psi_{N,S}$ North and South boundary layer terms, p.  \pageref{psiNS}

\bigskip
\noindent%
$\psi_{E,W}$ East and West boundary layer terms, p.  \pageref{psiEW}

\bigskip
\noindent
$\psil$ interior singular terms lifting the discontinuities of $\psi^0$, p.  \pageref{psil}

\bigskip
\noindent
$\psi^\Sigma$  interior singular layer terms, p.   \pageref{eq:psisig}

\bigskip
\noindent
$\psi_{int} = \psi^0_t + \psil$ regularization of the solution to the Sverdrup equation, p. \pageref{psiint-def}

\bigskip
\noindent
$\psib=\psi^0_t + \psi_E+ \psi^{cor}_E$ sum of all terms whose trace is lifted by the North/South boundary layers, p.  \pageref{def:psib}

\bigskip
\noindent
$\psib_\Sigma=\psil+ \psi^\Sigma\chi \left(\frac{y-y_1}{\delta_y}\right) + \psi^{corr}_W$ sum of all singular layer terms whose trace is lifted on the West boundary, p.  \pageref{psibSigma}

\bigskip
\noindent
$Y^-$ parametrization of the interior singular layer, p.  \pageref{Y-}

\bigskip
\noindent
$\Omega = \Omega_1\setminus \cup_{i=2}^K \Omega_i$ domain with islands $\Omega_i$, p.\pageref{omega-i}

\bigskip
\noindent
$\Omega^\pm$ subdomains of $\Omega$ where the transport equation has a continuous solution, p. \pageref{Sigma}

\bigskip
\noindent
$z$ distance to the boundary, p.  \pageref{localcoord}

\bigskip
\noindent%
$Z= \lambda(s) z$ scaled distance to the boundary, p.  \pageref{capZ}

\section*{Sizes of parameters and terms}

$$
\begin{array}{|c|c|c|}
\hline
\textbf{Parameter}&\multicolumn{2}{|c|}{\textbf{Size/Definition}}\\
\cline{2-3}
& \text{\it Case (H2i)} & \text{\it Case (H2ii)}\\
\hline\hline
s_i^{\pm} - s_i & C_i \viscosite^\frac{1}{4n+3} & \ds\frac{\mp 4 \alpha}{\ln \viscosite} \left( 1 - 6 \frac{\ln |\ln \viscosite|}{\ln \viscosite} + O\left(\frac{1}{\ln \viscosite}\right)\right)\\
\hline
\sigma_i^{\pm} - s_i \text{ \it (East coast)} & \multicolumn{2}{|c|}{C_i''} \\
\hline
\sigma_i^{\pm} - s_i \text{ \it (West coast)} &C_i' \viscosite^\frac{1}{4n+4} & \ds\frac{\mp 4 \alpha}{\ln \viscosite} \left( 1 - 8 \frac{\ln |\ln \viscosite|}{\ln \viscosite} + O\left(\frac{1}{\ln \viscosite}\right)\right)\\
\hline
\sigma_{in}&\multicolumn{2}{|c|}{\ds \left\{
\begin{array}{ll}
\sigma_1^- & \text{ if } I_1 =\{s_1\},\\
\sigma_0^- &\text{ if } I_1 =[s_0, s_1]\text{ and }(\sigma_0^-, s_0)\subset \Gamma_W,\\
s_0&\text{ if } I_1 =[s_0, s_1]\text{ and }(\sigma_0^-, s_0)\subset \Gamma_\viscosite.
\end{array}
\right.}\\
\hline
\lambda_{N,S}&\multicolumn{2}{|c|}{\viscosite^{-1/4}}\\
\hline
\lambda_{E,W}(s)&\multicolumn{2}{|c|}{\ds \left(\frac{|\cos \theta(s)|}{\viscosite}\right)^{1/3}}\\
\hline
\lambda_{E,W}'(s)&\multicolumn{2}{|c|}{\ds -\frac{\viscosite^{1/3} \sgn(\cos \theta(s)) \theta'(s) \sin \theta(s)}{3 |\cos \theta(s)|^{4/3}}}\\
\hline
\lambda_{E,W}(s_i^\pm)&\viscosite^{- \frac{n+1}{4n+3}}&\viscosite^{-1/4}|\ln \viscosite|^{1/2}\\
\hline
\lambda_{E,W}(\sigma_i^\pm)&\viscosite^{- \frac{3n+4}{12(n+1)}}&\viscosite^{-1/4}|\ln \viscosite|^{2/3}\\
\hline
\end{array}
$$

\vskip5mm

\label{tableaux-tailles}
$$
\begin{array}{|c|c|c|}
\hline
\textbf{Term}& \textbf{Typical size in }L^2(\Om)& \textbf{Typical size in }H^2(\Om^\pm)\\
\hline \hline
\psi^0_t&1&o(\viscosite^{-3/8})\\
\hline
\psi_{N,S}& \viscosite^{1/8} & \viscosite^{-3/8}\\
\hline
\psi_W & \viscosite^{1/6} & \viscosite^{-1/2} \\
\hline
\psi_E & \viscosite^{1/2} & o(\viscosite^{-3/8})\\
\hline
\psil& \viscosite^{1/8} & \viscosite^{-3/8}\\
\hline
\psisig & \viscosite^{1/8} & \viscosite^{-3/8}\\
\hline
\psi^{corr}_\Sigma & o(\viscosite^{1/8} )& o(\viscosite^{-3/8})\\
\cline{2-3}
&\multicolumn{2}{|c|}{\|\psi^{corr}_\Sigma\|_{W^{k,\infty}}=O(\delta_y^{-2} \viscosite^{-\frac{k-2}{4}})}\\
\hline
\end{array}
$$

%% file: munk-def-HAL.bbl
\begin{thebibliography}{99}
	
	
	
	
	
	
	
	
	\bibitem{BC}
	D. Bresch and T. Colin,  Some remarks on the derivation of the Sverdrup relation. {\sl J. Math. Fluid Mech.} \textbf{4} (2002), no. 2, 95--108.
	
	\bibitem{BGGRB}
	D. Bresch, F. Guillon-Gonzalez, and M. A. Rodriguez-Bellido, A corrector for the Sverdrup solution for a domain with islands. {\sl Appl. Anal.} \textbf{83} (2004), no. 3, 217--230.
	
	
	\bibitem{Brezis}
	H. Brezis,
	{\sl Functional analysis, Sobolev spaces and partial differential equations,}
	Universitext. Springer, New York, 2011. 
	
	\bibitem{CDGG} J.-Y. Chemin, B. Desjardins, I. Gallagher and
	E. Grenier, {\sl Basics of Mathematical Geophysics}, Oxford Lecture
	Series in Mathematics and its Applications, {\bf 32},  {\sl Oxford
		University Press}, 2006.
	
	
	\bibitem{DSR}  A.-L. Dalibard and L. Saint-Raymond,  Mathematical study of the $\beta$-plane model for rotating fluids in a thin layer. {\it J. Math. Pures Appl. } {\bf 94} (2010), no. 2, 131--169.
	
	
	\bibitem{DG} B. Desjardins and
	E. Grenier,  On the Homogeneous Model of Wind-Driven Ocean
	Circulation, {\sl SIAM Journal on Applied Mathematics}
	{\bf 60} (1999),     43--60.
	
	\bibitem{Du} V. Duch\^ene, On the rigid-lid approximation for two shallow layers of immiscible fluids with small density contrast, to appear in {\sl Journal of Nonlinear Science} (2014).
	
	
	\bibitem{eckhaus-jager} W. Eckhaus, E. M.  de Jager, 
	Asymptotic solutions of singular perturbation problems for linear differential equations of elliptic type. 
	{\sl Arch. Rational Mech. Anal.} {\bf  23} (1966),  26--86. 
	
	%
	%
	%
	%
	%
	%
	
	\bibitem{GVP} D. G\'erard-Varet and T. Paul, Remarks on boundary layer expansions. {\it Comm. Partial Differential Equations } {\bf 33} (2008), no. 1-3, 97--130.
	%
	
	
	\bibitem{gill} A. E. Gill, {\it Atmosphere-Ocean Dynamics},
	{ International
		Geophysics Series}, {\bf Vol. 30}, 1982.
	
	\bibitem{grasman} J. Grasman, {\it  On the birth of boundary layers}. Mathematical Centre Tracts, {\bf No. 36}. Mathematisch Centrum, Amsterdam, 1971.  
	
	
	\bibitem{grenier} E. Grenier, Oscillatory perturbations of the Navier--Stokes  equations. {\sl Journal de
		Math{\'e}matiques Pures et Appliqu{\'e}es}, {\bf 76}  (1997), pages
	477--498.
	%
	
	\bibitem{jung-temam} C. Jung, R. Temam, Convection-diffusion equations in a circle~: the compatible case, {\it J. Math. Pures Appl. } {\bf  96} (2011),  88--107.
	
	\bibitem{kikuchi} K. Kikuchi, Exterior problem for the two-dimensional Euler equation, {\it J. Fac. Sci. Univ. Tokyo
		Sect. IA Math.}, {\bf 30} (1983), pp. 63--92.
	
	%
	%
	\bibitem{LOT}  D. C. Levermore, M. Oliver and E. S.Titi, Global well-posedness for models of shallow water in a basin with a
	varying bottom, {\it Indiana Univ. Math. J.} \textbf{45} (1996), 479--510.
	%
	%
	%
	%
	%
	%
	%
	
	\bibitem{pedlosky} J. Pedlosky, {\it Geophysical fluid dynamics},
	{ Springer},
	$1979$.
	
	\bibitem{pedlosky2} J. Pedlosky, {\it Ocean Circulation Theory},
	{ Springer},
	$1996$.
	
	\bibitem{pedlosky3} J. Pedlosky, L.J. Pratt, M.A. Spall and K.R. Helrich,  Circulation around islands and ridges.
	{\it J. Marine Research}, {\bf 55} (1997), 1199--1251.
	
	%
	%
	%
	%
	
	\bibitem{Ro2} F. Rousset, Asymptotic behavior of geophysical fluids in highly rotating balls,
	{\sl Z. Angew. Math. Phys.} \textbf{58} (2007), no. 1, 53--67. 
	
	\bibitem{DR} W. P. M. de Ruijter, \textit{On the Asymptotic Analysis of Large-Scale Ocean Circulation}, Mathematical Centre Tracts \textbf{120}, Amsterdam, 1980.
	
	\bibitem{weak-bl} L. Saint-Raymond, Weak compactness methods for singular penalization problems with boundary layers. {\it SIAM J. Math. Anal. } {\bf 41} (2009), no. 1, 153--177.
	
	
	\bibitem{schochet} S.
	Schochet, Fast singular
	limits of hyperbolic
	PDEs. {\sl Journal of Differential Equations} {\bf
		114}  (1994), $476-512$.
	
	%
	
\end{thebibliography}
